\documentclass[11pt]{amsart}
\usepackage{amssymb,amscd,amsmath,amsthm}
\usepackage[all]{xy}
\usepackage{enumerate}

\newcommand\cA{{\mathcal{A}}}
\newcommand\cB{{\mathcal{B}}}
\newcommand\cC{{\mathcal{C}}}
\newcommand\cD{{\mathcal{D}}}
\newcommand\cE{{\mathcal{E}}}
\newcommand\cF{{\mathcal{F}}}
\newcommand\cG{{\mathcal{G}}}
\newcommand\cH{{\mathcal{H}}}

\newcommand\cK{{\mathcal{K}}}
\newcommand\cL{{\mathcal{L}}}
\newcommand\cM{{\mathcal{M}}}
\newcommand\cN{{\mathcal{N}}}

\newcommand\cR{{\mathcal{R}}}
\newcommand\cS{{\mathcal{S}}}

\newcommand\cU{{\mathcal{U}}}
\newcommand\cV{{\mathcal{V}}}

\newcommand\cX{{\mathcal{X}}}

\newcommand\RR{{\mathbb{R}}}
\newcommand\CC{{\mathbb{C}}}
\newcommand\ZZ{{\mathbb{Z}}}
\newcommand\NN{{\mathbb{N}}}
\newcommand\QQ{{\mathbb{Q}}}

\newcommand{\Hom}{\operatorname{Hom}}
\newcommand{\ch}{\operatorname{ch}}
\newcommand{\End}{\operatorname{End}}

\newcommand{\codim}{\operatorname{codim}}
\newcommand{\im}{\operatorname{im}}
\newcommand{\dom}{\operatorname{dom}}

\newcommand{\id}{\operatorname{id}}

\newcommand{\Pf}{\mbox{\rm Pf}}

\newcommand{\Ind}{\operatorname{Ind}}
\newcommand{\Id}{\operatorname{Id}}
\newcommand{\Td}{\operatorname{Td}}

\newcommand{\ev}{\text{ev}}
\newcommand{\odd}{\text{odd}}

\newcommand{\tr}{{\rm tr}\;}
\newcommand{\trLambda}{{\rm tr}_\Lambda}

\newcommand{\Tr}{{\rm Tr}\;}

\newcommand{\res}{{\rm res}\;}

\newcommand\Ker{\mathop{\rm Ker}\nolimits}
\newcommand\Dom{\mathop{\rm Dom}\nolimits}
\newcommand{\bint}{\ensuremath{-\hspace{-2,4ex}\int}}

\newcommand{\A}{\mathbb{A}}
\newcommand{\fg}{{\mathfrak{g}}}

\theoremstyle{plain}

\newtheorem{thm}{Theorem}[section]

\newtheorem{theorem}[thm]{Theorem}

\theoremstyle{definition}

\newtheorem{defn}[thm]{Definition}

\newtheorem{ex}[thm]{Example}

\theoremstyle{remark}

\numberwithin{equation}{section}
\begin{document}
\bibliographystyle{plain}

\title[Index Theory And Geometry
On Manifolds]{INDEX THEORY AND NON-COMMUTATIVE GEOMETRY\\ ON
FOLIATED MANIFOLDS}
\author{Yuri A. Kordyukov}
\address{Institute of Mathematics, Russian Academy of Sciences, 112
Chernyshevsky street, 450077 Ufa, Russia}
\email{yurikor@matem.anrb.ru}\thanks{Supported by the Russian
Foundation of Basic Research (grant no. 06-01-00208)}

\subjclass[2000]{Primary: 58B34. Secondary: 19K56, 46L87, 58J42}
\keywords{non-commutative geometry, manifolds, foliations,
transversally elliptic operators, tangentially elliptic operators,
index, $K$-theory, operator algebras}

\begin{abstract}
This paper gives a survey of the index theory of tangentially
elliptic and transversally elliptic operators on foliated manifolds
as well as of related notions and results in non-commutative
geometry.
\end{abstract}
\maketitle

\tableofcontents

\section*{Introduction}
The index theory of elliptic operators is one of the most
significant achievements in the mathematics of the twentieth
century. It began with a question raised by Gelfand in 1959. An
arbitrary elliptic operator $A$ on a compact manifold without
boundary determines a Fredholm problem in Sobolev spaces, and
therefore its integer-valued index
\[
\Ind A = \dim \Ker A-\dim {\rm Coker}\, A.
\]
is defined. Gel'fand observed that this index depends only on
topological properties of the operator and raised the question of
finding an explicit formula for the index of an elliptic operator in
topological terms. An answer to this question was given by Atiyah
and Singer in 1963. We refer the reader to the original papers
\cite{ASI,ASII,ASIII,ASIV}, and the books
\cite{Soloviev-Troitsky,BGV,Gilkey95,Melrose93,Palais,Roe98}, and
the bibliography cited therein for the classical Atiyah-Singer
theorem and its proofs.

Later on, the development of the index theory went in several
directions. The present survey is devoted to two of these
directions. One of them originated in the papers of Atiyah
\cite{At:teo} and Singer \cite{Singer:recent} and concerns a class
of non-Fredholm operators. This class consists of differential
operators on a compact manifold $M$ which are invariant under an
action of some compact Lie group $G$ on the manifold and elliptic in
directions conormal to the orbits of the action. Such operators are
said to be transversally elliptic. The index of a transversally
elliptic operator is no longer an integer, but a distribution on the
Lie group $G$. Transversally elliptic operators are naturally
regarded as an analogue of elliptic operators on the space $M/G$ of
orbits of the group action. In the papers \cite{Co:survey,Co:nc},
Connes considered transversally elliptic operators on compact
foliated manifolds.

Another direction in the index theory comes from the Atiyah-Singer
index theorem for families of elliptic operators \cite{ASIV}. The
index of a family of elliptic operators parameterized by points of
some topological space $X$ is defined as an element of the group
$K(X)$ of the topological $K$-theory. A generalization of this
theory is the index theory of tangentially elliptic operators on a
foliated manifold. A differential operator on a foliated manifold
$(M,\cF)$ is said to be tangentially elliptic if it contains
differentiations only along the leaves of the foliation (and
therefore can be restricted to any leaf of the foliation), and its
restriction to each leaf is an elliptic operator. Any tangentially
elliptic operator can be regarded as a family of elliptic operators
on leaves of the foliation, parameterized by the points of the leaf
space $M/\cF$.

Both transversally elliptic and tangentially elliptic operators on a
foliated manifold have a natural interpretation in terms of its leaf
space. The leaf space, generally speaking, is a very singular
object, and it is poorly described by means of classical tools of
geometry, topology, and analysis. Here we get help from
non-commutative geometry, one of the main goals of which is the
development of methods for the study of geometry, topology and
analysis on singular spaces such as the leaf space of a foliation.

There are several fundamental ideas which lie in the basis of
non-com\-mu\-ta\-ti\-ve geometry. The first of them is to pass from
geometric spaces to algebras of functions on these spaces and
translate basic geometric and analytic notions and constructions
into algebraic language. Such a procedure is well known and has been
used for a long time, for instance, in algebraic geometry. The next
idea is that in many cases, especially when the classical algebra of
functions is small or has a bad structure, it is useful to consider
some non-commutative algebra as an analogue of it. This necessitates
extending the basic geometric and analytic definitions to the case
of a general non-commutative algebra. Such a point of view has been
well known since the time of the Gel'fand-Naimark theory of
commutative $C^*$-algebras. For instance, the theory of
$C^*$-algebras is a far reaching generalization of the theory of
topological spaces and is often called non-commutative topology, and
the theory of von Neumann algebras is a generalization of the
classical theory of measure and integration. These ideas have turned
out to be very fruitful in the index theory as well.

In \cite{Co79} Connes defined the $C^*$-algebra $C^*(M,\cF)$ of a
foliation $(M,\cF)$, which it is natural to regard as an analogue of
the algebra of continuous functions on the leaf space $M/\cF$. For
instance, the index of a tangentially elliptic operator on $(M,\cF)$
is well defined as an element of the group $K(C^*(M,\cF))$ of the
operator topological $K$-theory.

The main purpose of non-commutative differential geometry, which was
initiated by Connes \cite{Co:nc} and is actively developing at
present time (cf. the recent surveys
\cite{Connes2000,ConnesLNM1831,survey} and the books
\cite{Co,Gracia:book,Landi:book,invitation,Naz-Savin-Sternin,Varilly,Singha-Goswami}
in regard to various aspects of non-commutative geometry), is to
extend the methods described above to analytic objects on geometric
spaces and to non-commutative algebras. Here the main attention is
focused on the facts that, first, a correct non-commutative
generalization applied in the classical setting, that is, to an
algebra of functions on a compact manifold must agree with its
classical analogue, and second it must inherit basic algebraic and
analytic properties of its classical analogue. Nevertheless, it
should be said that, as a rule, such non-commutative generalizations
are quite non-trivial and have a richer structure and essentially
new features in comparison with their commutative analogues.

It should be noted that the emergence of non-commutative
differential geometry itself in the paper \cite{Co:nc} is directly
connected with the index theory, because the notion of cyclic
cohomology introduced there was invented by Connes in attempts to
define the index of transversally elliptic operators on foliated
manifolds. Namely, the index of a transversally elliptic operator on
a foliated manifolds $(M,\cF)$ is defined as a cyclic cohomology
class on some subalgebra of the $C^*$-algebra $C^*(M,\cF)$. On the
other hand, cyclic cohomology and corresponding non-commutative
differential calculus play an important role in the index theory of
tangentially elliptic operators. First of all, we observe that in
the classical case the Chern character $\ch(\Ind P)$ of the index of
a family $P$ of elliptic operators parameterized by points of a
smooth manifold $B$ can be considered as a de Rham cohomology class
of the base $B$. Moreover, as stated by a local index theorem proved
by Bismut \cite{BismutInv85}, for families of Dirac operators it is
not only the cohomology class $\ch(\Ind P)\in H^*(B)$ that has a
geometrical meaning, but also a certain differential form on $B$
representing this class. In order to get numerical invariants from
$\ch(\Ind P)$, one can consider its pairings with arbitrary de Rham
currents on $B$. A non-commutative generalization of such a
construction is given by the higher index theory for tangentially
elliptic operators on foliated manifolds, which studies higher
indices of a tangentially elliptic operator defined as pairings of
its index with cyclic cohomology classes of some smooth subalgebras
of the algebra $C^*(M,\cF)$.

This paper is devoted to an exposition of the aspects of index
theory and non-commutative geometry mentioned above. We begin in
Section 1 with a survey of necessary notions of the classical index
theory of elliptic operators. In Section 2 we give some information
from non-commutative topology --- the notion of $C^*$-algebra as a
non-commutative topological space, the non-commutative analogues of
a vector bundle and of a field of Hilbert spaces, the simplest (and
the basic) invariant --- the $K$-theory and the $K$-homology. At the
end of the Section we give some necessary information from the
non-commutative theory of measure and integration. Section 3 is
devoted to basic notions of non-commutative differential geometry
--- cyclic (co)homology and non-commutative differential calculus,
spectral triples as an analogue of a non-commutative Riemannian
structure and the non-commutative local index theorem. In Section 4
we turn to foliations, starting with a brief summary of necessary
information from foliation theory. In Section 5 we describe the
construction of the operator algebras associated with a foliation,
which makes use of the notion of holonomy groupoid of a foliation.
Then basic objects of non-commutative topology of foliations are
discussed --- holonomy equivariant bundles, fields of Hilbert
spaces, $K$-theory and the Baum-Connes conjecture, the
non-commutative integration theory. Section 6 is devoted to
non-commutative differential calculus on the leaf space of a
foliation and to constructions of cyclic cocycles on the operator
algebras associated with a foliation. Finally, Sections 7 and 8 are
devoted, respectively, to the index theory of transversally elliptic
and tangentially elliptic operators.

This paper has some overlap with the author's previous survey
\cite{survey}, which is devoted to various aspects of
non-commutative geometry of foliations. In the present survey the
main emphasis is placed on the index theory of differential
operators associated with a foliation, in particular, on
applications of methods of non-commutative geometry to the index
theory.

The author is grateful to A.M. Vershik, who prompted the writing of
this paper, as well as to N.I. Zhukova and M.A. Shubin for useful
remarks.

\section{Survey of the classical index theory}

\subsection{Elliptic operators and the index}
Let $M$ be an $n$-dimensional smooth compact manifold without
boundary, and $E$ and $F$ smooth complex vector bundles on $M$ of
rank $N_E$ and $N_F$, respectively. (Here and subsequently,
``smooth'' means of class $C^\infty$. We will always assume that all
objects under consideration are of class $C^\infty$.) A linear
operator $D : C^\infty(M,E)\to C^\infty(M,F)$ is called a
differential operator of order $m$ if in any local chart $\phi :
U\subset M\to \RR^n$ and any trivializations $E\left|_U\right. \cong
U\times \CC^{N_E}$ and $F\left|_U\right. \cong U\times \CC^{N_F}$ of
the bundles $E$ and $F$ over it the operator $D$ has the form
\begin{equation}\label{e:diff}
D =\sum _{|\alpha|\leq m}a_{\alpha}(x)
\frac{\partial^{|\alpha|}}{\partial x_1^{\alpha_1}\ldots \partial
x_n^{\alpha_n}}, \quad x\in \phi(U)\subset \RR^{n},
\end{equation}
where $\alpha=(\alpha_1,\ldots,\alpha_n)\in \ZZ_+^n$ is a
multi-index, $|\alpha|=\alpha_1+\ldots+\alpha_n$, and the
$a_{\alpha}$ are smooth functions on $\RR^n$ with values in the
space $\cL(\CC^{N_E},\CC^{N_F})$ of complex $N_E\times N_F$
matrices.

For a differential operator $D$ given by the formula (\ref{e:diff})
in some local chart $\phi : U\subset M\to \RR^n$ and for
trivializations of the bundles $E$ and $F$ over it, define its
(complete) symbol
$$ \sigma (x,\xi) = \sum_{| \alpha | \leq m} a_{\alpha}(x)
(i\xi)^{\alpha},\quad x\in \phi(U)\subset \RR^{n},\quad \xi\in\RR^n,
$$ and its principal symbol $$ \sigma _{m}(x,\xi) = \sum _{ | \alpha |
=m} a_{\alpha}(x)(i\xi)^{\alpha},\quad x\in \phi(U)\subset \RR^{n},
\quad \xi\in\RR^n.
$$ The principal symbol is invariantly defined as a smooth section of the bundle
$\Hom(\pi^*E,\pi^*F)$ on the cotangent bundle $T^*M$, where
$\pi:T^*M\to M$ denotes the canonical projection, and $\pi^*E$ and
$\pi^*F$ are the pull-backs of the bundles $E$ and $F$,
respectively, on $T^*M$ by the map $\pi$. A differential operator $D
: C^\infty(M,E)\to C^\infty(M,F)$ is said to be elliptic if
$N_E=N_F$ and the principal symbol $\sigma_{m}(\xi)$ is invertible
for $\xi \in T^*M\setminus 0$.

In the index theory it is necessary to consider a wider class of
operators, the class of pseudodifferential operators. We recall that
a linear operator $P: C^{\infty}(M,E)\to {\cD}'(M,F)$ belongs to the
class $\Psi^m(M,E,F)$ of pseudodifferential operators of order $m$
if in a coordinate domain $X\subset \RR^n$ it can be represented as
\[
Pu(x)=\int e^{(x-y)\xi}p(x,\xi)\,u(y)\,dy\,d\xi, \quad x\in X,
\]
where $u\in C^\infty_c(X,\CC^{N_E})$, the function $p$ called the
complete symbol of $P$ belongs to the symbol class
$S^m(X\times\RR^n, \cL(\CC^{N_E},\CC^{N_F}))$.

The principal symbol $\sigma(P)$ of a $P\in \Psi^m(M,E,F)$ is an
element of the symbol space $S^m(T^*M,\Hom(\pi^*E,\pi^*F))$,
uniquely determined up to elements of the space
$S^{m-1}(T^*M,\Hom(\pi^*E,\pi^*F))$. In local coordinates
$\sigma(P)$ is given by the complete symbol $p$ of $P$. We say that
$P$ is elliptic if its principal symbol $\sigma(P)$ has a
representative $p\in S^m(T^*M,\Hom(\pi^*E,\pi^*F))$, which is
pointwise invertible outside some compact set in $T^*M$ and
satisfies the estimate
$$
\|p(x,\xi)^{-1}\|\le C(1+|\xi|)^{-m}
$$
for some constant $C$ and some Riemannian metric on $M$.

For an elliptic operator $P$ there exists an operator
$Q\in\Psi^{-m}(M;F,E)$ such that $QP-\Id$ and $PQ-\Id$ are smoothing
operators. The operator $Q$ is called a parametrix for $P$. The
existence of a parametrix implies that $P$ defines a Fredholm
operator acting in Sobolev spaces
$$
P_{(s)}: H^{s+m}(M,E)\to H^s(M,F)
$$
for any $s\in\RR$. The kernel of $P_{(s)}$ is finite-dimensional and
lies in $C^\infty(M,E)$, and its image is closed in $H^s(M,F)$ and
coincides with the orthogonal complement in $H^s(M,F)$ of the kernel
of $P^*$. Therefore, the index
\[
\Ind P_{(s)} =\dim \Ker P_{(s)}-\dim {\rm Coker}\, P_{(s)} =\dim
\Ker P_{(s)}-\dim \Ker P^*_{(s)}
\]
of $P_{(s)}$ is well-defined and independent of $s$; it is called
the index of $P$.

Examples of elliptic operators are given by Dirac operators. Let us
recall their definition (for an exposition of basic facts of spin
geometry and the theory of Dirac operators, see, for instance,
\cite{BGV,Lawson-M,Roe98}).

Let $M$ be a compact manifold of even dimension $n$, $g_M$ a
Riemannian metric on $M$, and $\nabla$ the Levi-Civita connection on
$TM$. For any $x\in M$ denote by ${\rm Cl}(T_xM)$ the complex
Clifford algebra of the Euclidean space $T_xM$. If one chooses an
orthonormal base $\{e_1,e_2,\ldots,e_n\}$ in $T_xM$, then ${\rm
Cl}(T_xM)$ is defined as an algebra over $\CC$ generated by the
elements $1$ and $e_1,e_2,\ldots,e_n$ satisfying the relations
\[
e_\alpha e_\beta+e_\beta e_\alpha=-2\delta_{\alpha\beta}, \quad
\alpha, \beta=1,2,\ldots,n.
\]

The Clifford algebra ${\rm Cl}(T_xM)$ has a natural $\ZZ_2$-grading.
Recall that a vector space $V$ is said to be $\ZZ_2$-graded if it
has a decomposition $V=V_0\oplus V_1$ into a direct sum of
subspaces. Equivalently, a $\ZZ_2$-grading on $V$ is determined by
an operator $\gamma\in \cL(V)$ such that $\gamma^2=1$. With respect
to the decomposition $V=V_0\oplus V_1$ the operator $\gamma$ has a
block form
$\begin{pmatrix} 1 & 0\\
0 & -1
\end{pmatrix}.$
If $V$ has a Euclidean structure, then it is natural to assume that
the subspaces $V_0$ and $V_1$ are orthogonal, which is equivalent to
the self-adjointness of $\gamma$. For any linear operator $T$ on a
$\ZZ_2$-graded space $V$ its supertrace is defined by
\[
{\rm Tr}_s(T)=\Tr \gamma T = \Tr T_{11}-\Tr T_{22},
\]
where $T$ is written in the block form $T=\begin{pmatrix} T_{11} & T_{12}\\
T_{21} & T_{22} \end{pmatrix}$ determined by the decomposition
$V=V_0\oplus V_1$.

We consider the vector bundle ${\rm Cl}(TM)$ on $M$ whose fibre at
$x\in M$ coincides with ${\rm Cl}(T_xM)$. This bundle is associated
with the principal $SO(n)$-bundle $O(TM)$ of oriented orthonormal
frames in $TM$: ${\rm Cl}(TM)=O(TM)\times_{O(n)} {\rm Cl}(\RR^n)$.
Therefore, the Levi-Civita connection $\nabla$ induces a natural
connection $\nabla^{{\rm Cl}(TM)}$ in ${\rm Cl}(TM)$ which is
compatible with Clifford multiplication and preserves the
$\ZZ_2$-grading on ${\rm Cl}(TM)$. If $\{e_1, e_2, \ldots , e_n\}$
is a local orthonormal frame in $TM$ and
$\omega^\gamma_{\alpha\beta}$ are the coefficients of $\nabla$:
$\nabla_{e_\alpha}e_\beta
=\sum_\gamma\omega^\gamma_{\alpha\beta}e_\gamma$, then
\[
\nabla^{{\rm Cl}(TM)}_{e_\alpha} =
e_\alpha+\frac{1}{4}\sum_{\gamma=1}^n
\omega^\gamma_{\alpha\beta}c(e_\beta)c(e_\gamma),
\]
where for any $a\in C^\infty(M,{\rm Cl}(TM))$, $c(a)$ denotes the
operator of pointwise left multiplication by $a$ in $C^\infty(M,{\rm
Cl}(TM))$.

A complex vector bundle $\cE$ on $M$ is called a Clifford module if
for any $x\in M$ there is a representation of the algebra ${\rm
Cl}(T_xM)$ in $\cE_x$ depending smoothly on $x$. The action of an
$a\in {\rm Cl}(T_xM)$ on an $s\in \cE_x$ will be denoted by
$c(a)s\in \cE_x$. A Clifford module $\cE$ is said to be self-adjoint
if it is endowed with a Hermitian metric such that the operator
$c(f) : \cE_x\to \cE_x$ is skew-symmetric for any $x\in M$ and $f\in
T_xM$. An arbitrary Clifford module $\cE$ has a natural
$\ZZ_2$-grading $\cE=\cE_+\oplus\cE_-$.

A connection $\nabla^\cE$ on a Clifford module $\cE$ is called a
Clifford connection, if for any $f\in C^\infty(M,TM)$ and $a\in
C^\infty(M,{\rm Cl}(TM))$ the following relation holds:
\[
[\nabla^{\cE}_f, c(a)]=c(\nabla^{{\rm Cl}(TM)}_fa).
\]

A self-adjoint Clifford module $\cE$ endowed with a Hermitian
Clifford connection $\nabla^\cE$ is called a Clifford bundle. Let
$e_1,\ldots,e_n$ be a local orthonormal basis in $TM$. The Dirac
operator $D_\cE$ associated with a Clifford bundle $\cE$ is defined
by
\[
D_\cE= \sum_{i=1}^n c(e_i)\nabla^{\cE}_{e_i}.
\]
The Dirac operator $D_\cE$ is formally self-adjoint in $L^2(M,
\cE)$. It is an odd operator with respect to the $Z_2$-grading
$C^\infty(M,\cE)=C^\infty(M,\cE_+)\oplus C^\infty(M,\cE_-)$. Thus,
it can be written in the form
\[
D_\cE=\begin{bmatrix}
      0 & D^{-}_\cE \\
      D^{+}_\cE & 0 \\
    \end{bmatrix}.
\]

\begin{ex}
An example of a Clifford bundle is the complexified exterior bundle
$\Lambda T^{*}M\otimes \CC$ on a Riemannian manifold $M$. An action
of ${\rm Cl}(TM)$ on it is given by the formula
\[
c(v)=\varepsilon_{v^*}-i_v, \quad v\in
T_xM,
\]
where $v^*\in T^*_xM$ is the covector dual to $v$,
$\varepsilon_{v^*}$ is the exterior multiplication by $v^*$, and
$i_v$ is the inner multiplication by $v$. A Clifford connection is
the Hermitian connection determined by the Riemannian metric. The
corresponding Dirac operator is the de Rham operator acting in
$C^{\infty}(M,\Lambda TM^{*})$ by the formula
\[
D_{\Lambda TM^{*}}=d + d^*,
\]
where $d$ is the de Rham differential.
\end{ex}

\begin{ex}
Recall that the group ${\rm Spin}(n)$ is the non-trivial double
covering of the group $SO(n)$. A spin structure on a Riemannian
manifold $M$ is defined to be a principal ${\rm Spin}(n)$-bundle
$O'(TM)$ on $M$ which is a double covering of the principal
$SO(n)$-bundle $O(TM)$ of oriented orthonormal frames in $TM$ such
that the map $O'(TM)\to O(TM)$ induces the double covering ${\rm
Spin}(n)\to SO(n)$ in each fibre. A manifold $M$ is called a spin
manifold if it admits a spin structure.

There is a unique (up to an isomorphism) non-trivial irreducible
unitary representation $S$ of the group ${\rm Spin}(n)$, called the
spin representation. The space of spinors has a natural
$\ZZ_2$-grading. For a Riemannian spin manifold $M$, denote by
$F(TM)$ the associated Hermitian vector bundle of spinors on $M$:
$F(TM)=O'(TM)\times_{{\rm Spin}(n)}S$. This bundle is a self-adjoint
Clifford module. The Levi-Civita connection $\nabla$ has a lift to a
Clifford connection $\nabla^{F(TM)}$ on $F(TM)$. The corresponding
Dirac operator is called the spin Dirac operator or simply the Dirac
operator.

More generally, one can take a Hermitian vector bundle $E$ endowed
with a Hermitian connection $\nabla^E$. Then $F(TM)\otimes E$ is a
Clifford module: an action of $a\in {\rm Cl}(TM)$ on $F(TM)\otimes
E$ is defined by the operator $c(a)\otimes 1$ ($c(a)$ denotes the
action of $a$ on $F(TM)$). The connection $\nabla^{F(TM)\otimes
E}=\nabla^{F(TM)}\otimes 1 + 1\otimes \nabla^{E}$ on $F(TM)\otimes
E$ is a Clifford connection. The corresponding Dirac operator
$D_{F(TM)\otimes E}=D_E$ is called the twisted spin Dirac operator
(or the spin Dirac operator with coefficients in $E$).
\end{ex}

\begin{ex}
By definition, the group ${\rm Spin}^c(n)$ is the subgroup of the
complex Clifford algebra ${\rm Cl}(\RR^n)$ generated by the group
${\rm Spin}(n)$ and the group $S^1=\{z\in\CC : |z|=1\}$.

Let $L$ be a principal $S^1$-bundle on $M$. The natural
representation $S^1$ in $\CC$ allows one to regard $L$ as a complex
line bundle on $M$. A ${\rm Spin}^c$-structure on $M$ is defined to
be a principal ${\rm Spin}^c(n)$-bundle $O'(TM)$ on $M$ which is a
double covering of the principal $SO(n)\times S^1$-bundle
$O(TM)\times L$ such that the map $O'(TM)\to O(TM)$ induces the
double covering ${\rm Spin}^c(n)\to SO(n)\times S^1$ in each fibre.
The bundle $L$ is called the fundamental line bundle associated with
the ${\rm Spin}^c$-structure. The corresponding Dirac operator is
called the ${\rm Spin}^c$ Dirac operator.
\end{ex}

\subsection{$K$-theory}\label{s:Ktheory}
The topological $K$-theory plays a very important role in the index
theory. Let $X$ be a compact topological space. The set of
isomorphism classes of finite-dimensional complex vector bundles on
$X$ endowed the operation of direct sum is an Abelian semigroup.
This semigroup generates an Abelian group $K(X)$, the Grothendieck
group, which consists of formal differences of vector bundles
(virtual vector bundles). A continuous map $f:X\to Y$ induces a
natural homomorphism $f^* : K(Y)\to K(X)$ which depends only on the
homotopy class of $f$.

In the case when $X$ is a locally compact space we will use the
group $K(X)$ of $K$-theory with compact supports, which is
conveniently described as follows. The group $K(X)$ is generated by
complexes of the form
\begin{equation}\label{e:Ed}
0 \longrightarrow E_0 \stackrel{d_0}{\longrightarrow} E_1
\stackrel{d_1}{\longrightarrow}\ldots
\stackrel{d_{N-1}}{\longrightarrow} E_N
\stackrel{d_N}{\longrightarrow}0,
\end{equation}
where $E_0,\ldots, E_N$ are vector bundles over $X$ and $d_0,\ldots,
d_N$ are morphisms of bundles. The support of such a complex is the
closure of the set of all $x\in X$, for which the sequence
(\ref{e:Ed}) is not exact. Denote by $S(X)$ the set of homotopy
classes of complexes of the form (\ref{e:Ed}) with compact support
and by $S_\varnothing(X)$ the subset in $S(X)$ which consists of
complexes of the form (\ref{e:Ed}) with empty support. The direct
sum operation defines a semigroup structure on $S(X)$. The group
$K(X)$ is defined as the quotient $S(X)/S_\varnothing(X)$, which is
in fact an Abelian group. Equivalently, one can consider the set
$S_N(X)$ of homotopy classes of complexes with compact support of
fixed length $N$ instead of $S(X)$. In particular, for $N=1$ we
obtain a description of $K(X)$ in terms of triples $(E_0,E_1,d_0)$,
where the morphism $d_0:E_0\to E_1$ is an isomorphism outside of
some compact.

For any integer $n\geq 0$ put $K^n(X)=K(X\times \RR^n)$. There is a
fundamental fact (Bott periodicity) stating that
\[
K^2(X)=K(X\times \RR^2)\cong K^0(X)=K(X).
\]
Thus, we have only two essentially different groups $K^0(X)=K(X)$
and $K^1(X)=K(X\times \RR)$ in topological $K$-theory.

One can give a definition of $K^1(X)$ in terms of the algebra
$C(X)$. Consider the group $GL(N,\CC)$ of invertible complex
$N\times N$ matrices, and assume that $GL(N,\CC)$ is embedded into
$GL(N+1,\CC)$ by means of the map
\[
X\mapsto
\begin{pmatrix} X & 0\\ 0 & 1 \end{pmatrix}.
\]
Let $GL_\infty(\CC)=\underset{\rightarrow}{\lim}\, GL(N,\CC)$.
Then
\[
K^1(X)=C(X,GL_\infty(\CC))/C(X,GL_\infty(\CC))_0,
\]
where $C(X,GL_\infty(\CC))$ denotes the group of continuous
functions on $X$ with values in $GL_\infty(\CC)$ and
$C(X,GL_\infty(\CC))_0$ denotes the identity component in
$C(X,GL_\infty(\CC))$. Equivalently, one can consider the group
$U(N)$ of unitary $N\times N$ matrices instead of $GL(N,\CC)$.

The role of orientation in $K$-theory is played by a complex spin
(${\rm Spin}^c$) structure. A real vector bundle $V$ of rank $m$
over $X$ is said to be $K$-orientable if its structure group reduces
to the group ${\rm Spin}^c(m)$. Any complex vector bundle is
$K$-orientable. There is the following description of ${\rm
Spin}^c(m)$-structures on a bundle $V$. Choose an arbitrary
Riemannian metric in the fibres of $V$. Let ${\rm Cl}(V)$ be the
corresponding bundle of complex`Clifford algebras. If $m$ is even,
then a ${\rm Spin}^c(m)$-structure on $V$ is given by a choice of an
orientation on $V$ along with a bundle $S$ of irreducible Clifford
modules. For odd $m$ one should replace ${\rm Cl}(V)$ by its even
part. For any $K$-orientable bundle $V$ over a compact manifold $X$
one has a Thom isomorphism $K(V)\cong K(X)$.

The $K$-theory and the usual cohomology groups (say, singular or
\v{C}ech) of a compact topological space $X$ are connected by the
Chern character
\begin{gather}
\ch : K^0(X)\to H^{{\rm ev}}(X,\QQ)=\bigoplus_{k\ \text{even}}
H^k(X,\QQ), \label{e:ch0}\\
\ch : K^1(X)\to H^{{\rm odd}}(X,\QQ)=\bigoplus_{k\ \text{odd}}
H^k(X,\QQ), \label{e:ch1}
\end{gather}
which becomes an isomorphism after tensoring by $\QQ$.

If $X$ is a smooth manifold, then one has an explicit
differential-geometric construction of the Chern character. The
construction of the even Chern character (\ref{e:ch0}) is a
particular case of the Chern-Weil construction of characteristic
classes of vector bundles. If $E$ is a smooth vector bundle on $X$,
then the Chern character $\ch(E)\in H^\ev(X,\RR)$ of the
corresponding class $[E]$ in $K^0(X)$ is represented by the de Rham
cohomology class of the closed differential form
\[
\ch(E,\nabla)=\Tr \exp \left(\frac{F}{2\pi i}\right)\in \Omega^{\rm
ev}(X,\CC)
\]
for any connection $\nabla : C^\infty(X,E)\to C^\infty(X,E\otimes
T^*X)$ in $E$, where $F=\nabla^2$ is the curvature of $\nabla$.

The odd Chern character (\ref{e:ch1}) is obtained from the even
Chern character (\ref{e:ch0}) by transgression. This construction is
a particular case of the construction of Chern-Simons classes (see
\cite{Baum-D82,Getzler93}). If $U \in C^{\infty}(X,U(N))$, then the
Chern character $\ch(U)\in H^\odd(X,\RR)$ of the corresponding class
$[U]$ in $K^1(X)$ is given by the de Rham cohomology class of the
closed differential form
\[
\ch(U)=\sum_{k=0}^{+\infty} (-1)^{k} \frac{k!}{(2k+1)!} \Tr
(U^{-1}dU)^{2k+1}\in \Omega^{\rm odd}(X,\CC).
\]

We note two more important particular cases of characteristic
classes given by the Chern-Weil construction, which we will need
below:
\begin{enumerate}[(i)]
  \item $\operatorname{Td}(E)\in H^{\rm ev}(X,\CC)$ is the Todd class
  of a complex vector bundle $E$, which is represented by the de Rham
cohomology class of the closed differential form
\[
\operatorname{Td}(E,\nabla)=\det\left(\frac{F}{e^F-1}\right)=\exp\Tr\left(\log
\left(\frac{F}{e^F-1}\right) \right)\in \Omega^{\rm ev}(X,\CC);
\]
  \item $\hat{A}(E)\in H^{4\ast}(X,\RR)$ is the $\hat{A}$-genus (the reduced
  Atiyah-Hirzebruch class) of a real vector bundle $E$, which is represented by the de
Rham cohomology class of the closed differential form
\[
\hat{A}(E,\nabla)={\det}^{1/2}\left(\frac{F/2}{\sinh(F/2)}\right)
=\exp\Tr\left(\frac12\log\frac{F/2}{\sinh(F/2)}\right)\in
\Omega^{4\ast}(X,\RR).
\]
\end{enumerate}

The definition of the Todd class $\operatorname{Td}(E)\in H^{\rm
ev}(X,\CC)$ can be extended to an arbitrary $K$-oriented vector
bundle $E$ by the formula
\[
\operatorname{Td}(E)=e^{c_1(L)/2}\hat{A}(E),
\]
where $c_1(L)$ is the first Chern class of the fundamental line
bundle $L$ associated with the ${\rm Spin}^c$-structure on $E$.

\subsection{The Atiyah-Singer theorem}\label{s:AS}
Let $M$ be a closed oriented manifold of dimension $n$, and consider
an elliptic pseudodifferential operator $P: C^{\infty}(M,E)\to
C^\infty(M,F)$, where $E$, $F$ are complex vector bundles over $M$.
By definition, the principal symbol $\sigma(P)$ of $P$ is an
isomorphism of the bundles $\pi^*E$ and $\pi^*F$ outside some
compact neighbourhood of the zero section $M\subset TM$. Therefore,
it gives rise to a well-defined element $[\pi^* E, \pi^*
F,\sigma(P)]\in K(T M)$ of the $K$-theory with compact supports of
the tangent bundle $\pi:TM\to M$ of $M$ (which can be identified
with the cotangent bundle $T^*M$ by means of a Riemannian metric on
$M$). One can prove that the index of the elliptic operator $P$
depends only on the class in $K(T M)$ defined by its principal
symbol. Moreover, any element of $K(T M)$ can be obtained by means
of this construction from the principal symbol of some elliptic
operator. Thus, a homomorphism (the analytic index) $\Ind_a: K(T
M)\to\ZZ$ is well defined.

On the other hand, by using topological constructions one can define
a homomorphism (the topological index) $\Ind_t: K(TM)\to \ZZ$. Let
us briefly describe its definition. Choose an embedding $i:M\to
\RR^n$ (such an $i$ exists for sufficiently large $n$). Denote by
$di:TM \to T\RR^n$ its differential, which in this case is a proper
embedding. Its normal bundle coincides with the lift of the bundle
$N\oplus N$ by the map $di$, where $N$ is the normal bundle of $i$.
Choose a diffeomorphism $N \oplus N \to W$, where $W$ is a tubular
neighborhood of $TM$ in $T\RR^n$. One has the Thom isomorphism
$\varphi^K : K(TM)\to K(N\oplus N)\cong K(W)$ for the Hermitian
complex vector bundle $N\oplus N\to TM$. For an open subset $W$ in
$T\RR^n$ there is a natural map $K(W)\to K(T\RR^n)$, called the
Gysin homomorphism. The composition of the Thom homomorphism and the
Gysin homomorphism is a map
$$
i_!:K(TM)\to K(T\RR^n)=K(\RR^{2n}).
$$
This construction holds for any smooth proper embedding $M\to V$ of
manifolds, and the resulting map $i_! : K(M)\to K(V)$, also called
the Gysin homomorphism, is independent on the choice of $W$ and
other auxiliary elements of the construction. We regard
$\RR^{2n}=\RR^{n}\oplus \RR^{n}=\CC^n \to pt$ as a complex vector
bundle. Then the Thom isomorphism $K(pt)\cong K(\RR^{2n})$ is
defined. Its inverse is the Bott periodicity isomorphism
$\beta:K(\RR^{2n})\cong K(pt)=\ZZ$. The topological index is defined
by the formula
$$
\Ind_t=\beta\circ i_!:K(TM)\to \ZZ.
$$

The Atiyah-Singer index theorem in the $K$-theoretic form \cite{ASI}
states the following.

\begin{thm}\label{t:ASKth}
One has the identity
$$
\Ind_a=\Ind_t :K(TM)\to \ZZ.
$$
\end{thm}

The following cohomological formula holds for $\Ind_t$ \cite{ASIII}.
Let $\pi_!: H^{*}(TM)\to H^{*}(M)$ be the map given by the
integration along the fibres of the bundle $\pi:TM\to M$. This map
is inverse to the Thom isomorphism $\Phi: H^{*}(M)\to H^{*}(TM)$.

\begin{thm}\label{t:AScoh}
For any $x\in K(TM)$ one has the identity
$$
\Ind_t(x)=(-1)^{n(n+1)/2} \int_M (\pi_! \ch(x))
\operatorname{Td}(TM\otimes \CC).
$$
\end{thm}

An immediate consequence of this formula is that the index of any
elliptic operator on a compact oriented odd-dimensional manifold
equals zero.

Let $M$ be an even-dimensional oriented Riemannian spin compact
manifold and $E$ a complex vector bundle on $M$ endowed with a
Hermitian structure $g^E$ and a unitary connection $\nabla^E$. The
spin Dirac operator $D_{E,+} : C^\infty(M, F_+(TM)\otimes E)\to
C^\infty(M, F_-(TM)\otimes E)$ with coefficients in $E$ has index
\[
\Ind(D_{E,+})=\int_M\hat{A}(TM)\ch(E).
\]

In this case one has a stronger form of the index theorem, which we
will state now. To start with, we recall the McKean-Singer formula
\cite{McKeanSinger}
\[
\Ind(D_{E,+})={\rm Tr}_s(\exp(-t(D_E)^2)), \quad t>0.
\]
Here we regard the space $L^2(M,F(TM)\otimes E)$ as a $\ZZ_2$-graded
Hilbert space and denote by ${\rm Tr}_s$ the supertrace of an
operator.

Let $P_t(x,y)$ be the smooth kernel of the operator
$\exp(-t(D_E)^2)$ with respect to the Riemannian volume form $dy$.
Then
\[
\Ind(D_{E,+})=\int_M{\rm Tr}_s(P_t(x,x))\,dx, \quad t>0.
\]
In \cite{AtiyahBottPatodi,Gilkey73,Patodi71}, the following theorem
is proved (the local index theorem):

\begin{thm}
The following pointwise limit relation holds as $t\to 0$
\[
{\rm Tr}_s(P_t(x,x))\to
\{\hat{A}(TM,\nabla^{TM})\ch(E,\nabla^E)\}^{\rm max}.
\]
\end{thm}

Another proof of the local index theorem, which significantly
improved its geometric understanding, was given by Getzler
\cite{Getzler86}. The Atiyah-Singer index theorem in the
cohomological form is an immediate consequence of the local index
theorem and the McKean-Singer formula.

A generalization of the index theorem to manifolds with boundary was
obtained by Atiyah, Patodi and Singer in \cite{APS:SARG1,APS:SARG3}.
In \cite{ASI} an equivariant index theorem for elliptic operators
invariant under an action of a compact Lie group was also proved. In
\cite{ASII} a formula was proved which provides an expression of the
$G$-index of an elliptic $G$-complex in terms of the fixed points of
the action (an analogue of the Atiyah-Bott-Lefschetz formula
\cite{Atiyah-Bott,Atiyah-BottII}).

\subsection{The index theory for self-adjoint operators}
Let us consider a first-order self-adjoint elliptic operator $D :
C^\infty(M,E)\to C^\infty(M,E)$ on a closed manifold $M$. Then its
index equals zero. Nevertheless, as discovered in
\cite{Baum-D82,Kasparov80}, the operator $D$ gives rise to some
analytic index type topological invariants. For this, one makes use
of Toeplitz operators $T_u$ associated with unitary multipliers $u
\in C^\infty(M, U(N))$. The Toeplitz operator $T_u$ is the bounded
operator on $L^2(M,E\otimes \CC^N)$ given by
\[
T_u=P_+M_uP_+,
\]
where $P_+$ is the spectral projection corresponding to the positive
semi-axis for the operator $D\otimes I_N$ acting in $L^2(M,E\otimes
\CC^N)$, and $M_u$ is the operator of multiplication by $1\otimes
u$. The operator $T_u$ is a Fredholm operator. For its index there
is a topological formula \cite{Baum-D82} which is derived from the
Atiyah-Singer index theorem:
\[
\Ind T_u = \int_{ST^*M} \pi^*\ch(u) \ch(E_+) \pi^*\Td(TM\otimes
\CC),
\]
where $\ch(u)\in H^{\rm odd}(M)$ is the Chern character of the class
$[u]\in K^1(M)$ defined by (\ref{e:ch1}), $\pi^*\ch(u)\in H^{\rm
odd}(ST^*M)$ is its lift to the cosphere bundle $ST^*M$ by the
natural projection $\pi :ST^*M\to M$, $E_+$ is the subbundle of the
bundle $\pi^*E$ on $ST^*M$ generated by the positive eigenvectors of
the principal symbol $\sigma(D)$, and $\ch(E_+)\in K(ST^*M)$ is its
Chern character.

One can show that the index of the Toeplitz operator $T_u$ depends
only on the class $[u]\in K^1(M)$ of the unitary multiplier $u \in
C^\infty(M, U(N))$, thus yielding a map
\[
K^1(M)\to \ZZ, \quad K^1(M) \ni [u] \mapsto \Ind T_u.
\]
As a result, the index theory of self-adjoint operators is often
called the odd index theory.

\subsection{The families index theory}\label{s:fam}
The index theory for families of elliptic operators was developed in
\cite{ASIV}. Let $Z$ be a fibration over a Hausdorff topological
space $Y$ with fibre $X$ and structure group ${\rm Diff}(X)$ (a
manifold over $Y$), and let $E$ and $F$ be fibrations over $Y$.
Denote by $Z_y$ the fibre of $Z$ over $y$ and by $E_y$ and $F_y$ the
restrictions of the bundles $E$ and $F$ to $Z_y$. We consider a
continuous family $\{ P_y\in \Psi^d(Z_y, E_y, F_y) : y\in Y\}$ of
elliptic pseudodifferential operators. If $\dim \Ker P_y$ is
independent of $y$, the family $\{\Ker P_y : y\in Y\}$ of vector
bundles defines a vector bundle $\Ker P$ over $Y$. The same holds
for ${\rm Coker}\,P$. In this case the index of the family $P$ is
defined by
\[
\Ind_a(P)=[\Ker P]-[{\rm Coker}\,P]\in K(Y).
\]
In the general case when $\dim \Ker P_y$ varies, the definition of
the index of $P$ as an element of $K(Y)$ is given by a slight
modification of this definition.

Denote by $TZ/Y$ the vertical tangent space (the tangent space along
the fibres of the fibration). The symbol of the family $P$ defines
an element $[\sigma(P)]\in K(TZ/Y)$. The analytic index of $P$
depends only on $[\sigma(P)]\in K(TZ/Y)$, thus defining an analytic
index map
\[
\Ind_a : K(TZ/Y)\to K(Y).
\]
The construction of a topological index $\Ind_t : K(TZ/Y)\to K(Y)$
for families of elliptic operators is a direct generalization of the
construction of the topological index for elliptic operators. One
simply needs to perform all constructions used there fibrewise over
$Y$. The family index theorem proved in \cite{ASIV} establishes the
coincidence of the analytic and topological indices:

\begin{thm}\label{t:ASfamily}
One has the equality
\[
\Ind_a=\Ind_t : K(TZ/Y)\to K(Y).
\]
\end{thm}

We also formulate the index theorem for families of elliptic
operators in the cohomological form.

\begin{thm}\label{t:ASfamily2}
Let $P$ be a family of elliptic operators parameterized by a
manifold $Y$ and $u\in K(TZ/Y)$ the class of the symbol of $P$. Then
\[
\ch \Ind P= (-1)^n\pi_* (\ch u \Td(TZ\otimes \CC))\in H^*(Y,\QQ),
\]
where $n$ is the dimension of fibres, and $\pi_* : H^*(TZ/Y)\to
H^*(Y)$ is integration along the fibres of $TZ/Y\to Y$.
\end{thm}

An important example is the index theorem for families of Dirac
operators. Let $\pi: M\to B$ be a fibration of compact manifolds
with compact fibres $Z_y, y\in B,$ of even dimension $n=2l$. Suppose
that the vertical tangent bundle $TM/B$ admits a spin structure. Let
$g_{M/B}$ be a smooth metric on $TM/B$. Denote by $F=F_+(TM/B)\oplus
F_-(TM/B)$ the spinor bundle for $TM/B$. Let $E$ be a Hermitian
bundle on $M$ endowed with a unitary connection $\nabla^E$. For any
$y\in B$ one has a well defined Dirac operator
\[
D_{E,y}=\begin{bmatrix}
      0 & D_{E,y,-} \\
      D_{E,y,+} & 0 \\
    \end{bmatrix},
\]
acting in $C^\infty(Z_y,F(TZ_y)\otimes E)$. Thus, the family of
elliptic operators $\{D_{E,y,+} : y\in B\}$ is well defined. For the
index of this family one has the formula
\[
\ch(\Ind D_{E,+}) =\pi_*[\hat{A}(TM/B)\ch E]\in H^{\rm ev}(B,\QQ).
\]

Bismut \cite{BismutInv85} proved a stronger version of this result,
a so-called local index theorem for families.

Let $T^HM$ be a subbundle of $TM$ such that $TM=T^HM\oplus TM/B$.
The bundle $T^HM$ can be identified with $\pi^*TB$. For any $U\in
TB$, denote by $U^H$ its lift to $T^HM$. Any metric $g_B$ on $TB$
lifts to $T^HM$.

For any $U,V\in TB$ put
\[
T(U,V)=-P[U^H,V^H] \in TM/B,
\]
where $P: TM\to TM/B$ is the orthogonal projection.

For $U\in C^\infty(B,TB)$ let $\operatorname{div}_Z(U^H)$ denote the
divergence of the vector field $U^H$ with respect to the vertical
Riemannian volume form $dv_Z$.

Let $\nabla^L$ be the Levi-Civita connection on $TM$ associated with
$g_B\oplus g_{M/B}$. We introduce a Euclidean connection
$\nabla^{TM/B}$ on $TM/B$ by the formula
\[
\nabla^{TM/B}=P \nabla^L.
\]
The connection $\nabla^{TM/B}$ is independent of the choice of $g_B$
and is canonically determined by $T^HM$ and $g_{M/B}$, and
$\nabla^{TM/B}$ and $\nabla^E$ determine a connection
$\nabla^{F\otimes E} : C^\infty(M,F\otimes E)\to C^\infty(M,T^*M
\otimes F\otimes E)$ on $F\otimes E$. We define a connection
\[
\nabla^H : C^\infty(M,F\otimes E) \to C^\infty(M,(T^HM)^* \otimes
F\otimes E)
\]
as follows: for $U\in C^\infty(B,TB)$ and $s\in C^\infty(M,F\otimes
E)$
\[
\nabla^H_Us=\nabla^{F\otimes E}_{U^H}s+\frac12
\operatorname{div}_Z(U^H) s.
\]

Let $f_1,f_2,\ldots,f_m$ be a local orthonormal base in $TB$ and
$f^1,f^2,\ldots,f^m$ its dual base in $T^*B$. For any $t>0$ consider
the operator
\[
A_t : C^\infty(M,F\otimes E)\to C^\infty(M,(T^HM)^* \otimes F\otimes
E)
\]
given by
\[
A_t=\nabla^H +\sqrt{t} D_E -\frac{1}{8\sqrt{t}}f^\alpha f^\beta
c(T(f_\alpha,f_\beta)).
\]
It defines a superconnection on $C^\infty(M,F\otimes E)$, regarded
as an infinite-dimensional vector bundle over $B$ with fibre
$C^\infty(Z_y,F(TZ_y)\otimes E)$ at $y\in B$, in the sense of the
following definition.

Denote by $\Omega(M)$ the space of smooth differential forms on $M$
and, for any vector bundle $\cE$ on $M$, by $\Omega(M,\cE)$ the
space of smooth differential forms on $M$ with coefficients in
$\cE$.

\begin{defn}
Let $\cE$ be a $\ZZ_2$-graded vector bundle on a manifold $M$. A
superconnection on $\cE$ is defined to be an odd first-order
differential operator $\A : \Omega^{\pm}(M,\cE) \to
\Omega^{\mp}(M,\cE)$ satisfying the $\ZZ_2$-graded Leibniz rule: if
$\alpha\in \Omega(M)$ and $\theta\in \Omega(M,\cE)$, then
\[
\A(\alpha \wedge \theta)=d\alpha\wedge
\theta+(-1)^{|\alpha|}\alpha\wedge\A(\theta).
\]
\end{defn}

A superconnection $\A$ on $\cE$ gives rise to an action on the space
$\Omega(M,\End\ \cE)$, which is compatible with the Leibniz rule:
\[
\A\alpha=[\A,\alpha], \quad \alpha \in \Omega(M,\End\ \cE).
\]
The curvature of the superconnection $\A$ is defined as the operator
$\A^2$ acting in the space $\Omega(M,\End \cE)$.

In our case the curvature $A_t^2$ of the superconnection $A_t$ is
given by a family of second-order elliptic operators with
coefficients in $\Lambda^2(T^HM)^*$ acting along the fibres of $\pi:
M\to B$:
\[
A_t^2 : C^\infty(M,F\otimes E)\to C^\infty(M,\Lambda^2(T^HM)^*
\otimes F\otimes E).
\]
The operator
\begin{align*}
\exp(-A_t^2)\colon C^{\infty}(M,F\otimes E) &\to
C^{\infty}\bigl(M,\Lambda^*(T^HM)^*\otimes F\otimes E\bigr)
\\
& \qquad\cong C^{\infty}\bigl(M,\Lambda^*(T^*B)\otimes F\otimes
E\bigr),
\end{align*}
is given by a family of smoothing operators acting along the fibres
of the fibration, with coefficients in differential forms on the
base $B$.

For any $t>0$ consider the even form on $B$
\[
\alpha_t=\phi{\rm Tr}_s[\exp(-A_t^2)],
\]
where the linear endomorphism $\phi : \Lambda(T^*B)\to
\Lambda(T^*B)$ is given by the formula $\omega \mapsto
(2i\pi)^{-deg\omega/2}\omega$, and ${\rm Tr}_s$ denotes the
fibrewise supertrace.

The following facts hold \cite{BismutInv85}:
\begin{enumerate}
  \item $\alpha_t$ is a real even closed differential form on $B$;
  \item the de Rham cohomology class $[\alpha_t]\in
  H^{ev}(B,\QQ)$ is equal to $\ch(\Ind D_{E,+})$;
  \item (the family local index theorem) as $t\to 0$,
\[
\alpha_t=\pi^*[\hat{A}(TM/B,\nabla^{TM/B})\ch(E,\nabla^E)] +O(t).
\]
\end{enumerate}

Moreover, it is proved in \cite{BerlineVergne87}, \cite{BGV} that if
$\Ker D_E$ is a vector bundle and $\nabla^{\Ker D_E}$ is the
orthogonal projection of the connection $\nabla^H$ on $\Ker D_E$,
then as $t\to +\infty$,
\[
\alpha_t=\ch(\Ker D_E,\nabla^{\Ker
D_E})+O\left(\frac{1}{\sqrt{t}}\right).
\]

\subsection{The higher index theory}\label{s:high-ind0}
For an even-dimensional closed connected Riemannian spin manifold
$M$ and a Hermitian vector bundle $E$ on it, the Atiyah-Singer index
theorem establishes a connection between the index of the spin Dirac
operator $D_E$ with coefficients in $E$ and the topological
expression $\int_M \hat{A}(TM) \ch(E)$. If $M$ is not simply
connected, one can modify the index theorem, taking into account the
fundamental group $\Gamma$ of $M$. Denote by $p : M\to B\Gamma$ the
classifying map for the universal covering $\widetilde{M}\to M$. The
higher index theory attempts to give an analytic interpretation of
the expression
\[
\int_M \hat{A}(TM) \ch(E) p^*[\eta],
\]
where $[\eta]\in H^*(B\Gamma,\CC)$. As examples of applications of
the higher index theory, we can mention the Novikov conjecture on
homotopy invariance of non-simply connected manifolds and questions
on the existence of metrics of positive scalar curvature (see, for
instance, the paper
\cite{Ferry:RR,Gromov95,Kas93,Leichtnam-Piazza04} and the references
cited therein).

The higher index theory is directly connected with the $L^2$-index
theory on coverings of compact manifolds, which originated in
\cite{At:discrete}. There Atiyah proved a $\Gamma$-index theorem for
$\Gamma$-invariant elliptic operators on a covering of a compact
manifold. In \cite{Co-M90} Connes and Moscovici proposed an approach
to the higher index theory based on the use of non-commutative
differential geometry, in particular, of cyclic cohomology. Their
approach, like some others, is based on the idea or regarding a
$\Gamma$-invariant elliptic operator on a covering of a compact
manifold as a family of elliptic operators parameterized by points
of a non-commutative space $B$ whose algebra of continuous functions
is the reduced group $C^*$-algebra $C^*_r\Gamma$. In
\cite{Lott:gafa92} there is a proof of a higher index theorem, using
methods of the local index theory for families of elliptic
operators. We also mention the paper \cite{Lott:eta}, where a higher
analogue of the eta-invariant for $\Gamma$-invariant operators on a
$\Gamma$-covering of a compact manifold was introduced.

\section{Basic notions of non-commutative topology}
\subsection{Non-commutative spaces and bundles}
Recall that a $C^*$-algebra is an involutive Banach algebra $A$ with
\[
\|a^*a\|=\|a\|^2, \quad a\in A.
\]
A simplest example of a $C^*$-algebra is the algebra $C_0(X)$ of
continuous functions on a locally compact Hausdorff topological
space $X$ vanishing at infinity, endowed with the operations of the
pointwise addition and multiplication, with the standard involution
given by the complex conjugation, and with the uniform norm
\[
\|f\|=\sup_{x\in X}|f(x)|, \quad f\in C_0(X).
\]

A Gel'fand-Naimark theorem enables us to reconstruct from a
commutative $C^*$-algebra $A$ a unique locally compact Hausdorff
topological space $X$ such that $A\cong C_0(X)$. More precisely, $X$
coincides with the set $\widehat{A}$ of all characters of the
algebra $A$, that is, of all continuous homomorphisms $A\to \CC$,
endowed with the topology of pointwise convergence. This fact allows
one to regard an arbitrary $C^*$-algebra as the algebra of
continuous functions on some virtual space. For this reason, the
theory of $C^*$-algebras is often called non-commutative topology.

The algebra $\cL(H)$ of bounded operators on a Hilbert space $H$
equipped with the involution given by taking adjoints and with the
operator norm is a $C^*$-algebra. By the second Gel'fand-Naimark
theorem, any $C^*$-algebra is isometrically $\ast$-isomorphic to
some norm-closed $\ast$-subalgebra of the algebra $\cL(H)$ for some
Hilbert space $H$.

We recall that a right module $\cE$ over a unital algebra $A$ is
said to be finitely generated if it is generated by a finite family
$\{x_i\in\cE : i=1,2,\ldots,k\}$ of elements, that is, the submodule
of finite $A$-linear combinations of the form $\sum_{i=1}^k x_ia_i$,
where $a_i\in A$, is dense in $A$. A right $A$-module $\cE$ is said
to be projective if there exists a right $A$-module $\cE^\prime$
such that the direct sum $\cE \oplus \cE^\prime$ is isomorphic to
the free module $A^N$ for some $N$.

If $E$ is a continuous complex vector bundle on a compact
topological space $X$, then the space $C(X,E)$ of its continuous
sections is a finitely generated projective module over the algebra
$C(X)$ of continuous functions on $X$. The action of $C(X)$ on
$C(X,E)$ is given by the formula
\[
(a\cdot s)(x)=a(x)s(x),\quad x\in X.
\]
The Serre-Swan theorem states that any finitely generated projective
$C(X)$-module is isomorphic to the $C(X)$-module $C(X,E)$ for some
finite-di\-men\-si\-o\-nal complex vector bundle $E$. Thus, an
arbitrary finitely generated projective module over a $C^*$-algebra
can be regarded as an analogue of a finite-dimensional complex
vector bundle over the corresponding non-commutative space, or in
other words, as a non-commutative vector bundle.

In many problems of index theory and non-commutative geometry, it is
useful to consider more general objects, namely, continuous fields
of Hilbert spaces. A natural example of a continuous field of
Hilbert spaces arises in the index theory for families of elliptic
operators. Let $Z$ be a fibration over a Hausdorff topological space
$Y$, with fibre $X$ and structure group ${\rm Diff}(X)$, and let $E$
be a bundle over $Y$. As above, denote by $Z_y$ the fibre of $Z$
over $y$ and by $E_y$ the restriction of $E$ to $Z_y$. Then the
family $\{ L^2(Z_y,E_y), y\in Y\}$ of Hilbert spaces is a continuous
field of Hilbert spaces on $Y$.

A non-commutative analogue of a continuous field of Hilbert spaces
is the notion of Hilbert $C^*$-module. Hilbert $C^*$-modules can be
regarded as a natural generalization of Hilbert spaces, which arises
if one replaces the field of scalars $\CC$ by an arbitrary
$C^*$-algebra.

\begin{defn}
Let $B$ be a $C^*$-algebra. A pre-Hilbert $B$-module is defined to
be a right $B$-module $\cE$ equipped with a sesquilinear map
$\langle \cdot, \cdot\rangle_B : \cE\times \cE\to B$ (linear in the
second argument), satisfying the following conditions:
\begin{enumerate}
  \item $\langle x,x\rangle_B \geq 0$ for any $x\in \cE$;
  \item $\langle x,x\rangle_B = 0$ if and only if $x=0$;
  \item $\langle y,x\rangle_B = \langle x,y\rangle^*_B$ for any $x,y\in \cE$;
  \item $\langle x,yb\rangle_B=\langle x,y\rangle_Bb$ for any $x,y\in
\cE, b\in B$.
\end{enumerate}
The map $\langle \cdot, \cdot\rangle_B$ is called a $B$-valued inner
product.
\end{defn}

For a pre-Hilbert $B$-module $\cE$ the formula $\|x\|_B=\|\langle
x,x\rangle_B\|^{1/2}$ defines a norm on $\cE$. If $\cE$ is complete
in the norm $\| \cdot \|_B$, then $\cE$ is called a Hilbert
$C^*$-module. In the general case the action of $B$ and the inner
product on $\cE$ are extended to its completion $\tilde{\cE}$,
making $\tilde{\cE}$ into a Hilbert $C^*$-module.

\begin{ex}
If $\cH=\{H_x: x\in X\}$ is a continuous field of Hilbert spaces
over a compact topological space $X$, then the space $C(X,\cH)$ of
its continuous sections is a Hilbert module over $C(X)$. The action
of $C(X)$ on $C(X,\cH)$ is given by the formula
\[
(a\cdot s)(x)=a(x)s(x),\quad x\in X,
\]
and the inner product by
\[
\langle s_1, s_2\rangle_{C(X)}(x)= \langle s_1 (x),
s_2(x)\rangle_{H_x}, \quad x\in X.
\]
\end{ex}

Let $\cE$ be a Hilbert $B$-module. Denote by $\cB(\cE)$ the set of
all endomorphisms $T$ of $\cE$ such that there exists an adjoint
endomorphism, that is, an endomorphism $T^*:\cE\to\cE$ such that
$\langle Tx,y\rangle_B=\langle x,T^*y\rangle_B$ for any $x,y\in
\cE$. Any operator in $\cB(\cE)$ is a bounded operator in $\cE$, and
the algebra $\cB(\cE)$ is a $C^*$-algebra with respect to the
uniform norm.

For any $x,y\in \cE$ denote by $\theta_{x,y}$ the operator defined
in $\cE$ by $\theta_{x,y}(z)=x\langle y,z\rangle_B, z\in \cE$. It is
easy to see that $\theta_{x,y}\in \cB(\cE)$. The closure $\cK(\cE)$
of the linear span of $\{ \theta_{x,y}\in \cB(\cE) : x,y\in \cE\}$
is a closed ideal in $\cB(\cE)$. Its elements are called compact
endomorphisms of $\cE$.

\begin{ex}\label{ex:field1}
Let $Z$ be a fibration over a compact topological space $Y$ with
fibre $X$ and structure group ${\rm Diff}(X)$, and let $E$ be a
bundle over $Y$. Consider the Hilbert $C(Y)$-module $\cE$ of
continuous sections of the continuous field $\{L^2(Z_y,E_y), y\in
Y\}$ of Hilbert spaces on $Y$. Then an arbitrary continuous family
$\{ P_y\in \Psi^d(Z_y, E_y) : y\in Y\}$ of pseudodifferential
operators with $d\leq 0$ defines a bounded endomorphism of the
Hilbert $C(Y)$-module $\cE$. If $d<0$, this endomorphism is a
compact endomorphism.
\end{ex}

The notion of isometric $\ast$-iso\-morp\-hism between
$C^*$-algebras is a natural analogue of the notion of homeomorphism
of topological spaces. There is a broader equivalence relation for
$C^*$-algebras called strong Morita equivalence. It preserves many
invariants of $C^*$-algebras, for instance, the K-theory, the space
of irreducible representations, the cyclic cohomology, and it
coincides with the relation of isometric $\ast$-isomorphism on the
class of commutative $C^*$-algebras. We briefly recall some
information about the strong Morita equivalence (a more detailed
exposition can be found in \cite{Rieffel82}).

\begin{defn}
Let $A$ and $B$ be $C^*$-algebras.

An $A$-$B$-bimodule is a vector space $X$ endowed with structures of
a left $A$-module and a right $B$-module, which are compatible in
the sense that $(ax)b=a(xb)$ for any $x\in X$, $a\in A$, and $b\in
B$.

An $A$-$B$-equivalence bimodule is defined to be an $A$-$B$-bimodule
$X$ endowed with $A$-valued and $B$-valued inner products $\langle
\cdot, \cdot\rangle_A$ and $\langle \cdot, \cdot\rangle_B$,
respectively, such that $X$ is a right Hilbert $B$-module and a left
Hilbert $A$-module with respect to these inner products and,
moreover,
\begin{enumerate}
  \item $\langle x,y\rangle_A z=x\langle y, z\rangle_B$  for any
  $x, y, z\in X$,
  \item the set $\langle X,X\rangle_A$ spans a dense subset in
  $A$, and the set $\langle X,X\rangle_B$ spans a dense subset in
  $B$.
\end{enumerate}

Algebras $A$ and $B$ for which there is an $A$-$B$-equivalence
bimodule are said to be strongly Morita equivalent.
\end{defn}

It is not difficult to show that the strong Morita equivalence is an
equivalence relation.

For any linear space $L$, denote by $\tilde{L}$ the complex
conjugate linear space, which coincides with $L$ as a set and has
the same addition operation, but has multiplication by scalars given
by the formula $\lambda\tilde{x}=(\bar{\lambda} x)^\sim $. If $X$ is
an $A$-$B$-equivalence bimodule, then $\tilde{X}$ is endowed with
the structure of a $B$-$A$-equivalence bimodule. For instance,
$b\tilde{x}a=(a^*xb^*)^\sim .$

\begin{thm} \label{t:morita}
Let $X$ be an $A$-$B$-equivalence bimodule. Then the map $E \to
X\otimes_BE$ defines an equivalence of the category of Hermitian
$B$-modules and the category of Hermitian $A$-modules, with inverse
given by the map $F \to F\otimes_B\tilde{X}$.
\end{thm}

In particular, Theorem~\ref{t:morita} implies that two commutative
$C^*$-algebras are strongly Morita equivalent if and only if they
are isomorphic.

The following theorem relates the notion of strong Morita
equivalence with the notion of stable equivalence.

\begin{thm}
Let $A$ and $B$ be $C^*$-algebras with countable approximate units.
Then these algebras are strongly Morita equivalent if and only if
they are stably equivalent, i.e. the algebras $A\otimes \cK$ and
$B\otimes \cK$ are isomorphic. (Here $\cK$ denotes the algebra of
compact operators in a separable Hilbert space.)
\end{thm}

\subsection{The operator $K$-theory}
We begin this subsection with some facts from $K$-theory for
$C^*$-algebras, the non-commutative analogue of topological
$K$-theory.

Let $A$ be a unital $C^*$-algebra. The group $K_0(A)$ is defined as
the Gro\-then\-dieck group of the semigroup of isomorphism classes
of finitely generated projective modules over $A$, with the direct
sum operation. Thus, elements of $K_0(A)$ can be regarded as formal
differences of isomorphism classes of finitely generated projective
modules over $A$. Equivalently, one can consider isomorphism classes
of orthogonal projections in the algebra of matrices over $A$.

Another definition of the group $K_0(A)$ is given as follows. Denote
by $M_n(A)$ the algebra of $n\times n$ matrices with entries in $A$.
Let us assume that $M_n(A)$ is embedded in $M_{n+1}(A)$ by means of
the
map $X\to \begin{pmatrix} X & 0\\
0 & 0 \end{pmatrix}$. Let
$M_\infty(A)=\underset{\rightarrow}{\lim}\, M_n(A)$. The group
$K_0(A)$ is defined as the set of homotopy equivalence classes of
projections ($p^2=p=p^*$) in $M_\infty(A)$ equipped with the direct
sum operation
\[
p_1\oplus p_2=\begin{pmatrix} p_1 & 0\\ 0 & p_2 \end{pmatrix}.
\]

Denote by $GL_n(A)$ the group of invertible $n\times n$ matrices
with entries in $A$, and assume that $GL_n(A)$ is embedded in
$GL_{n+1}(A)$ by means of the map $X\to
\begin{pmatrix} X & 0\\ 0 & 1 \end{pmatrix}$. Let
$GL_\infty(A)=\underset{\rightarrow}{\lim}\, GL_n(A)$. The group
$K_1(A)$ is defined as the set of homotopy equivalence classes of
unitary matrices ($u^*u=uu^*=1$) in $GL_\infty(A)$ equipped with the
direct sum operation.

If $A$ has no unit and $A^+$ is the algebra obtained by adjoining
the unit to $A$, then one has the homomorphism $i:\CC\to A^+ :
\lambda\mapsto \lambda\cdot 1$, which induces a homomorphism $i_*:
K_0(\CC)\to K_0(A^+)$, and $K_0(A)$ is defined as the kernel of this
homomorphism. Moreover, $K_1(A)=K_1(A^+)$, by definition.

These definitions agree with those given in the commutative case:
the isomorphisms $K_i(C_0(X))\cong K^i(X), i=0,1$, hold for any
locally compact topological space $X$.

For an arbitrary algebra $\cA$ over $\CC$ we will consider the
groups $K_0(\cA)$ and $K_1(\cA)$ of algebraic $K$-theory. The group
$K_0(\cA)$ is defined like the group $K_0(A)$ of topological
$K$-theory with the use of idempotents ($e^2=e$) in $M_\infty(\cA)$
instead of projections. The group $K_1(\cA)$ is defined as the
quotient of $GL_\infty(\cA)$ by the commutator subgroup
$[GL_\infty(\cA), GL_\infty(\cA)]$.

The homological $K$-functor --- the object dual to the topological
$K$-theory --- was introduced by purely homotopic methods by
Whitehead in 1962. Atiyah \cite{Atiyah-global} observed that an
elliptic operator on a smooth manifold can be regarded in some sense
as an element of a $K$-homology group. He formulated
functional-analytic axioms for two basic properties of an elliptic
pseudodifferential operator on a compact manifold. Using them, he
defined, for any compact topological space $X$, a class of objects
${\rm Ell}(X)$ and, in the case when $X$ is a CW-complex, an
epimorphism ${\rm Ell}(X)\to K_0(X)$. He proposed regarding elements
of ${\rm Ell}(X)$ as representing cycles for $K_0(X)$. Atiyah's
ideas were completely realized by Kasparov in \cite{Kasparov75} (see
also \cite{BDF}). The analytic construction in \cite{Kasparov75} of
the $K$-homology groups is applicable to an arbitrary
non-commutative $C^*$-algebra and is based on the notion of Fredholm
module.

\begin{defn}
A Fredholm module (or a $K$-cycle) over an algebra $A$ is a pair
$(H,F)$, where
\begin{enumerate}
  \item $H$ is a Hilbert space equipped with a
$\ast$-representation $\rho$ of the algebra $A$;
  \item $F$ is a bounded operator in $H$ such that for any
  $a\in A$ the operators $(F^2-1)\rho(a)$, $(F-F^*)\rho(a)$, and
  $[F,\rho(a)]$ are compact on $H$.
\end{enumerate}
A Fredholm module $(H,F)$ is said to be even if the Hilbert space
$H$ is endowed with a $\ZZ_2$-grading $\gamma$, the operators
$\rho(a)$ are even, $\gamma\rho(a)=\rho(a)\gamma$, and the operator
$F$ is odd, $\gamma F=-F\gamma$. Otherwise it is said to be odd.
\end{defn}

The homology groups $K^0(A)$ (respectively, $K^1(A)$) of a
$C^*$-algebra $A$ are defined as the sets of homotopy equivalence
classes of even (respectively, odd) Fredholm modules over $A$. The
direct sum operation defines an Abelian group structure on $K^0(A)$
and $K^1(A)$.

\begin{ex}\label{ex:fred0}
Let $M$ be a compact manifold, $E_0$ and $E_1$ Hermitian vector
bundles on $M$, $P\in \Psi^0(M,E_0,E_1)$ a zero order elliptic
operator whose principal symbol $\sigma_P$ satisfies the condition
$\sigma_P\sigma_P^*=1$ (for instance, if $D\in \Psi^d(M,E_0,E_1)$ is
an elliptic operator of order $d>0$, then one can take for $P$ the
operator $D(1+D^*D)^{-1/2}$), and $Q\in \Psi^0(M,E_1,E_0)$ a
parametrix for $P$. Define the $\ZZ_2$-graded Hilbert space
$H=L^2(M,E_0)\oplus L^2(M,E_1)$ endowed with a natural action of the
algebra $C(M)$, and the bounded operator $F$ in $H$ given by the
matrix $\begin{pmatrix} 0 & Q\\ P & 0
\end{pmatrix}.$
Then the pair $(H,F)$ is an even Fredholm module over $C(M)$.

For $D$ one can take the Dirac operator $D^{+}_\cE :
C^\infty(M,\cE_+)\to C^\infty(M,\cE_-)$ associated with an arbitrary
Clifford bundle $\cE$ over an even-di\-men\-si\-o\-nal compact
Riemannian manifold $M$.
\end{ex}

\begin{ex}\label{ex:fred1}
Let $M$ be a compact manifold, $E$ a Hermitian vector bundle on $M$,
and $F\in \Psi^0(M,E)$ a zero order elliptic operator whose
principal symbol $\sigma_F$ satisfies the condition
$\sigma^*_F=\sigma_F$, $\sigma^2_F=1$ (for example, if $D\in
\Psi^d(M,E)$ is a self-adjoint elliptic operator of order $d>0$,
then one can take for $F$ the operator $D(1+D^2)^{-1/2}$). Define
the $\ZZ_2$-graded Hilbert space $H=L^2(M,E)\oplus 0$ endowed with
the natural action of the algebra $C(M)$. Then the pair $(H,F)$ is
an odd Fredholm module over $C(M)$.

For an operator $D$ one can take the Dirac operator $D^{+}_\cE :
C^\infty(M,\cE_+)\to C^\infty(M,\cE_-)$ associated with an arbitrary
Clifford bundle $\cE$ over an odd-di\-men\-si\-o\-nal compact
Riemannian manifold $M$.
\end{ex}

Corresponding to a Fredholm module $(H,F)$ over an algebra $A$ is an
index map $\operatorname{ind}:K_{*}(A) \rightarrow \ZZ$. In the even
case the operator $F$ takes the form
\[
     F= \left(
        \begin{array}{cc}
            0& F^-  \\
            F^+& 0
        \end{array}
        \right),    \quad F_{\pm} : H^{\pm} \rightarrow H^{\mp},
\]
with respect to the decomposition $H=H^{+}\oplus H^-$ given by the
$\ZZ_{2}$-grading of $H$. For an idempotent $e \in M_{q}(A)$ the
operator $e(F^{+}\otimes 1)e$ acting from $e(H^{+}\otimes \CC^q)$ to
$e(H^{-}\otimes \CC^q)$ is Fredholm, and its index depends only on
the class of $e$ in $K_{0}(A)$. Therefore, a map
$\operatorname{Ind}:K_{0}(A) \rightarrow \ZZ$ is well defined by the
formula
    \begin{equation}
        \operatorname{Ind}[e]= \operatorname{Ind} e(F^{+}\otimes 1)e.
       \label{eq:LIFNCG.index-even}
    \end{equation}

In the odd case, for a unitary matrix $U\in GL_{q}(A)$ the operator
$(P\otimes 1) U (P\otimes 1)$ with $P=(1+F)/2$ is a Fredholm
operator. Moreover, the index of $(P\otimes 1) U (P\otimes 1)$
depends only on the class of $U$ in $K_{1}(A)$. Thus, one obtains
the map $\operatorname{Ind}:K_{1}(A) \rightarrow \ZZ$ given by the
formula
\begin{equation}
     \operatorname{Ind}[U]=  \operatorname{Ind} (P\otimes 1) U (P\otimes 1).
        \label{eq:LIFNCG.index-odd}
\end{equation}

If $A$ is a $C^*$-algebra, then in both the even and odd cases the
map $\operatorname{Ind}$ determines a map of the group $K_i(A)$ of
topological $K$-theory to $\CC$, and it depends only on the class
determined by the Fredholm module $(H,F)$ in the $K$-homology group
$K^i(A)$.

In \cite{Baum-D82} a geometric definition of the $K$-homology groups
is given. We briefly recall this definition. A $K$-cycle on a
topological space $X$ is defined to be a triple $(M,E,\phi)$, where
$M$ is a compact ${\rm Spin}^c$ manifold without boundary, $E$ is a
complex vector bundle on $M$, and $\phi$ is a continuous map from
$M$ to $X$. We remark that $M$ is not necessarily connected, and the
fibres of $E$ may have different dimensions on different components.
One can define a natural notion of isomorphism of $K$-cycles. On the
set of isomorphism classes of $K$-cycles on $X$, there is an
equivalence relation generated by three elementary relations, called
bordism, direct sum, and vector bundle modification.

The geometric $K$-homology group $K^{\rm geom}(X)$ is defined as the
set of equivalence classes of $K$-cycles on $X$. An Abelian group
structure on $K^{\rm geom}(X)$ is given by the obvious operation of
disjoint union of $K$-cycles.

Denote by $K^{\rm geom}_0(X)$ (respectively, $K^{\rm geom}_1(X)$)
the subgroup of $K^{\rm geom}(X)$, which consists of equivalence
classes of $K$-cycles $(M,E,\phi)$ such that each component of $M$
is even-dimensional (respectively, odd-dimensional).

An isomorphism of geometric and analytic $K$-homology groups is
constructed as follows. Let $(M,E,\phi)$ be a $K$-cycle on $X$. Let
$D_E$ be the ${\rm Spin}^c$ Dirac operator on $M$ with coefficients
in $E$, and let $[D_E]\in K_i(M)$, $i=\dim M \mod 2$ be the
corresponding class in the $K$-homology of $M$ (see
Examples~\ref{ex:fred0} and~\ref{ex:fred1}). The map $\phi:M\to X$
induces a map $\phi_*: K_i(M)\to K_i(X)$ in the $K$-homology. The
correspondence $(M,E,\phi)\to \phi_*[D_E]\in K_i(X)$ defines a map
$K^{\rm geom}_i(X)\to K_i(X), i=0,1,$ which is an isomorphism
\cite{baum-higson-schick}.

There is a natural transformation of homology theories, the homology
Chern character $\ch : K_*(X)\to H_*(X,\QQ)$, defined as follows
(see \cite{Baum-D82}). Consider a $K$-cycle $(M,E,\phi)$ on $X$. The
map $\phi:M\to X$ induces a map $\phi_* : H_*(M,\QQ)\to H_*(X,\QQ)$
of rational homology groups. Put
\begin{equation}\label{e:chern}
\ch (M,E,\phi)=\phi_*(\ch(E)\cup \Td(M)\cap [M])\in H_*(X,\QQ),
\end{equation}
where the cap product $\ch(E)\cup \Td(M)\cap [M]\in H_*(M,\QQ)$
coincides with the Poincar\'e dual to $\ch(E)\cup \Td(M)\in
H^*(M,\QQ)$. In particular, we observe that, for the Dirac operator
$D_E$ on a ${\rm Spin}^c$-manifold $M$ with coefficients in $E$, we
have
\[
\ch [D_E]=\ch(E)\cup \Td(M)\cap [M]\in H_*(M,\QQ).
\]
If $X$ is a finite CW-complex, then the Chern character $\ch$
determines an isomorphism $K_*(X)\otimes \QQ \to H_*(X,\QQ)$. For an
arbitrary CW-complex $X$, we will consider singular homology
$H_*(X,\QQ)$ and $K$-homology $K_*(X)$ with compact supports.
Therefore, the Chern character $\ch : K_*(X)\otimes \QQ \to
H_*(X,\QQ)$ is well defined and is also a rational isomorphism.

We will also need $K$-homology groups $K_{*,F}(X)$ of $X$ twisted by
a real vector bundle $F$ on $X$. They are defined by the formula
\[
K_{j,F}(X)=K_j(F,F\setminus\{0\}), \quad j=0,1.
\]
One can give an equivalent definition, choosing a Euclidean
structure in the fibres of $F$ and introducing the unit ball and
unit sphere bundles $BF$ and $SF$ of $F$. Then
\[
K_{j,F}(X)=K_j(BF,SF), \quad j=0,1.
\]

The $K$-homology fundamental class of a compact ${\rm
Spin}^c$-manifold $M$ is the class $[D]\in K_i(M)$, where $i=\dim M
\mod 2$, defined by an arbitrary ${\rm Spin}^c$ Dirac operator $D$
on $M$. The Poincar\'e duality in $K$-theory states that the cap
product with $[D]$ gives an isomorphism
\[
K^*(M)\cong K_*^a(M).
\]
An analogous statement also holds for an arbitrary compact ${\rm
Spin}^c$-manifold $M$ with boundary $\partial M$:
\[
K^*(M,\partial M)\cong K_*^a(M).
\]
If $X$ is an arbitrary closed Riemannian manifold, then an
application of these statements to the ${\rm Spin}^c$-manifold
$T^*X$ leads to the Poincar\'e duality isomorphisms
\[
K^*(T^*X)=K^*(BT^*X,ST^*X)\cong K_*^a(X).
\]
This isomorphism has a natural interpretation in terms of index
theory. Namely, consider an arbitrary element $x\in
K^0(BT^*X,ST^*X)$. It is given by a triple $(E_0,E_1,\sigma)$, where
$E_0$ and $E_1$ are smooth vector bundles on $X$ and $\sigma$ is an
isomorphism $\sigma : \pi^*E_0\to \pi^*E_1$ of vector bundles.
Without loss of generality, one can assume that $\sigma$ is the
principal symbol of some pseudodifferential operator $D_\sigma$.
Then the class $[D_\sigma]$ in $K_*^a(X)$, determined by $D_\sigma$,
coincides with the Poincar\'e dual of $x\in K^0(BT^*X,ST^*X)$. A
similar construction holds for elements $x\in K^1(BT^*X,ST^*X)$.

For an arbitrary compact ${\rm Spin}^c$-manifold $M$ with boundary
$\partial M$ there is another Poincar\'e duality isomorphism
\[
K^*(M)\cong K_*(M,\partial M).
\]
If $X$ is an arbitrary closed Riemannian manifold, then an
application of this isomorphism to the ${\rm Spin}^c$-manifold
$T^*X$ leads to the Poincar\'e duality isomorphisms
\[
K_{*,T^*X}(X)=K_*(BT^*X,ST^*X)\cong K^*(X).
\]

At the end of this subsection we present some information about the
bivariant Kasparov $KK$-functor introduced in \cite{Kasparov80}.

Let $A$ and $B$ be $C^*$-algebras. In this subsection we take an
$(A,B)$-bimodule to be a $\mathbb Z_2$-graded countably generated
Hilbert module $H$ equipped with an action of $A$ given by a
representation $\rho: A\to\mathcal B(H)$ of it by even operators on
$H$. Denote by $Q_A(H)$ the subalgebra in $\mathcal B(H)$ formed by
operators $T\in\mathcal B(H)$ such that $[T,\rho(a)]\in\mathcal
K(H)$ for any $a\in A$, and by $I_A(H)$ the ideal in $Q_A(H)$ which
consists of operators $T\in\mathcal B(H)$ such that
$T\rho(a)\in\mathcal K(H)$ and $\rho(a)T\in\mathcal K(H)$ for any
$a\in A$.

We consider pairs $(H,F)$, where $H$ is an $(A,B)$-bimodule and
$F\in Q_A(H)$ is an odd operator such that the operators $F-F^*$ and
$F^2-I$ belong to $I_A(H)$. We say that pairs $(H_1,F_1)$ and
$(H_2,F_2)$ are isomorphic if there exists a grading-preserving,
isometric isomorphism $u: H_1\to H_2$ of $(A,B)$-bimodules such that
$F_2=uF_1u^{-1}$. The set of isomorphism classes of pairs $(H,F)$ is
denoted by $E(A,B)$.

A homotopy between elements $(H_0,F_0)\in E(A,B)$ and $(H_1,F_1)\in
E(A,B)$ is defined to be an element $(H,F)\in E(A,B[0,1])$ whose
restrictions to the extreme points $t=0$ and $t=1$ of the segment
coincide with $(H_0,F_0)$ and $(H_1,F_1)$, respectively:
$H_t=H\otimes_{B[0,1]}B, T_t=T\otimes 1$. We define $KK(A,B)$ as the
set of homotopy equivalence classes in $E(A,B)$.

The set $KK(A,B)$ is an Abelian group with respect to the direct sum
operation: $(H_1,F_1)\oplus (H_2,F_2)=(H_1\oplus H_2, F_1\oplus
F_2)$.

We`denote by $KK_0(A,B)$ (respectively, $KK_1(A,B)$) the subset in
$KK(A,B)$ given by elements $(H,F)\in E(A,B)$ such that the operator
$F$ is odd (respectively, even).

The groups $KK_i(\CC,B)$ are naturally isomorphic to the groups
$K_i(B)$ of topological $K$-theory. In some sense, this isomorphism
is an algebraic reformulation of the definition of the index of a
family of elliptic operators. In the case when $B=C(X)$ is the
algebra of continuous functions on a compact topological space $X$,
a family $P$ of elliptic operators parameterized by the points of
$X$ defines naturally an element of the group $KK(\CC,B)$. The
element of the group $K_0(B)=K(X)$ corresponding to this element
under the isomorphism $KK(\CC,B)\cong K(B)$ coincides with the index
of the family $P$. For an arbitrary $C^*$-algebra $B$, elements of
the group $KK_i(\CC,B)$ can be constructed using elliptic
pseudodifferential operators invariant under the action of $B$. An
index theorem for elliptic operators over a $C^*$-algebra was proved
by Mishchenko and Fomenko in \cite{Mishch-Fomenko} (see also
\cite{Soloviev-Troitsky}). It is also clear that in the case of
$B=\CC$ the definition of $KK_i(A,\CC)$ coincides with the
definition of the $K$-homology groups $K^i(A)$.

The main technical tool is a bilinear pairing (the cup-cap product)
\[
KK_i(A_1,B_1\otimes D)\otimes KK_j(D\otimes A_2,B_2)\to
KK_{i+j}(A_1\otimes A_2,B_1\otimes B_2),
\]
which generalizes all known products in $K$-theory and $K$-homology.
We will denote the cup-cap product of elements $x\in
KK_i(A_1,B_1\otimes D)$ and $y\in KK_j(D\otimes A_2,B_2)$ by
$x\otimes_Dy$. An idea of the construction of the cup-cap product is
borrowed from the theory of elliptic operators. It has many natural
properties of a product: it is contravariant in $A_1$ and $A_2$,
covariant in $B_1$ and $B_2$, functorial in $D$ and associative. The
latter means that
\[
(x_1\otimes_{D_1} x_2)\otimes_{D_2}x_3=x_1\otimes_{D_1}
(x_2\otimes_{D_2}x_3)
\]
for any $x_1\in KK_i(A_1,B_1\otimes D_1)$, $x_2\in KK_j(D_1\otimes
A_2,B_2\otimes D_2)$, $x_3\in KK_{\ell}(D_2\otimes A_3,B_3)$.

\subsection{Non-commutative measure theory} The initial data of
non-com\-mu\-ta\-ti\-ve measure theory are a pair $(\cM,\phi)$,
consisting of a von Neumann algebra $\cM$ and a weight $\phi$ on
$\cM$.

\begin{defn}
A von Neumann algebra is an involutive subalgebra of the algebra
$\cL(H)$ of bounded operators on a Hilbert space $H$, which is
closed in the weak operator topology.
\end{defn}

\begin{defn}
A weight on a von Neumann algebra $\cM$ is a function $\phi$ defined
on the set $\cM_+$ of positive elements of $\cM$ with values in
$\bar{\RR}_+=[0.+\infty]$, which satisfies the conditions
\begin{gather*}
\phi(a+b)=\phi(a)+\phi(b), \quad a,b \in \cM_+,\\ \phi(\alpha
a)=\alpha \phi(a), \quad \alpha \in \RR_+,\quad a\in \cM_+.
\end{gather*}

A weight on a von Neumann algebra $\cM$ is called a trace, if
\[
\phi(a^*a)=\phi(aa^*), \quad a\in \cM_+.
\]
\end{defn}

\begin{defn}
A weight $\phi$ on a von Neumann algebra $\cM$ is
\begin{enumerate}
  \item faithful if for any $a\in \cM_+$ the equality $\phi(a)=0$
  implies that $a=0$;
  \item normal if for any bounded increasing net
  $\{a_\alpha\}$ of elements of $\cM_+$ with the least upper bound
  $a$, one has the equality
$\phi(a)=\sup_\alpha\phi(a_\alpha)$.
  \item semifinite if the linear span of the set $\{x\in
  \cM_+:\phi(x)<\infty\}$ is $\sigma$-weakly dense in $\cM$.
\end{enumerate}
\end{defn}

Any von Neumann algebra has a faithful, normal, semifinite weight. A
von Neumann algebra is said to be semifinite if it has a faithful
normal semifinite trace.

\begin{ex}
The usual trace $\tr $ on the von Neumann algebra $\cL(\cH)$ of
bounded linear operators on a Hilbert space $\cH$ is a faithful
normal semifinite trace. Moreover, for any bounded positive operator
$T$ in $\cH$ the functional $\phi_T(A)=\tr AT, A\in \cL(\cH)$, is a
faithful normal semifinite weight on $\cL(\cH)$.
\end{ex}

\begin{ex}
If there is given a $\sigma$-finite measure $\mu$ on a measurable
space $X$, then the elements of the space $L^\infty(X,\mu)$ regarded
as multiplication operators acting in the Hilbert space
$L^2(X,\mu)$, form a von Neumann algebra. Moreover, the equality
\[
\phi(f)=\int_X f(x)\,d\mu(x), \quad f\in L^\infty(X,\mu),
\]
defines a faithful normal semifinite trace on $L^\infty(X,\mu)$.
\end{ex}

We recall that for any subset $S\subset \cL(\cH)$ its commutant is
the set of bounded operators on $\cH$ which commute with all
operators from $S$. An involutive subalgebra $\cM$ in $\cL(\cH)$ is
a von Neumann algebra if and only if $\cM^{\prime\prime}=\cM$. An
unbounded linear operator $T$ acting in a Hilbert space $\cH$ is
said to be affiliated with a von Neumann algebra $\cM$ acting in
$\cH$ if $Tu=uT$ for any unitary operator on $\cH$ belonging to the
commutant $\cM^\prime$ of $\cM$.

Let $\cN$ be a semifinite von Neumann algebra, and $\tau$ a faithful
normal semifinite trace on $\cN$. The norm-closed two-sided ideal in
$\cN$ generated by elements $E\in \cN$ with $\tau(E)<\infty$ will be
denoted by $\cK_\cN$. Its elements are called $\tau$-compact
operators. A Breuer-Fredholm operator is defined to be any operator
$P\in\cN$ whose image under the canonical map $\pi : \cN\to
\cN/\cK_\cN$ is invertible in the algebra $\cN/\cK_\cN$. The
Breuer-Fredholm index of a Breuer-Fredholm operator $P\in\cN$ is
defined by the formula
\[
\Ind (F)=\tau(P_{\Ker F})-\tau(P_{{\rm Coker}\ F}),
\]
where $P_{\Ker F}$ and $P_{{\rm Coker}\ F}$ are the projections to
the kernel and cokernel of the operator $F$ respectively. The theory
of Breuer-Fredholm operators was developed in the papers
\cite{Breuer1,Breuer2} in the case when $\cN$ is a factor, and was
extended to the case when $\cN$ is not a factor in
\cite{Phillips-Raeburn}.

\section{Non-commutative differential geometry}
\subsection{Cyclic cohomology and homology}
In this subsection we give a definition of cyclic cohomology and
homology, which play the role of non-commutative analogue of the
homology and cohomology of topological spaces (for cyclic
cohomology, see the books \cite{Co,Karoubi87,Loday} and the
references cited therein). It is important to remark that the
definition of the de Rham homology and cohomology of a topological
space requires the introduction of an additional structure on this
space, for instance, the structure of a smooth manifold. In the
non-commutative case this results in the fact that cyclic cocycles
are usually defined not on a $C^*$-algebra, the analogue of the
algebra of continuous functions, but on a subalgebra of it which
consists of ``smooth'' functions. For the time being, we postpone
(see Subsection~\ref{s:smooth}) a discussion of the non-commutative
analogue of the algebra of smooth functions on a smooth manifold,
the notion of smooth algebra, and turn to the definition of the
cyclic cohomology for an arbitrary algebra.

Let $\cA$ be an algebra over $\CC$. Consider the complex
$(C^*(\cA,\cA^*), b)$, where:
\begin{itemize}
  \item[(i)] $C^k(\cA,\cA^*), k\in \NN$, is the space of $(k+1)$-linear forms on $\cA$;
  \item[(ii)] the coboundary $b\psi\in
  C^{k+1}(\cA,\cA^*)$ of an element $\psi\in C^k(\cA,\cA^*)$ is
  given by
\begin{align*}
b\psi(a^0, \cdots, a^{k+1}) = &  \sum_{j=0}^k (-1)^j \psi(a^0,
\cdots, a^j a^{j+1}, \cdots, a^{k+1})  \\ &  +
(-1)^{k+1}\psi(a^{k+1}a^0, \cdots, a^k), \quad  a^0,a^1,\ldots,
a^{k+1} \in \cA.
  \end{align*}
\end{itemize}
The cohomology of this complex is called the Hochschild cohomology
of the algebra $\cA$ with coefficients in the bimodule $\cA^*$ and
is denoted by $HH(\cA)$.

Let $C^k_\lambda(\cA)$ be the subspace of $C^k(\cA,\cA^*)$ which
consists of all $\psi\in C^k(\cA,\cA^*)$, satisfying the cyclicity
condition
\begin{equation}
    \psi(a^1, \cdots, a^k,a^0) =(-1)^k \psi(a^0, a^{1},\cdots, a^k),
    \quad a^0,a^1,\ldots, a^{k} \in
    \cA.
\end{equation}
The differential $b$ maps the subspace $C^k_\lambda(\cA)$ to
$C^{k+1}_\lambda(\cA)$, and the cyclic cohomology $HC^{*}(\cA)$ of
$\cA$ is defined as the cohomology of the complex
$(C^*_\lambda(\cA),b)$.

\begin{ex}
For $k=0$ the space $HC^{0}(\cA)$ coincides with the space of all
trace functionals on $\cA$. For this reason cyclic $k$-cocycles on
$\cA$ are called $k$-traces on $\cA$ (for $k>0$, higher traces).
\end{ex}

\begin{ex}
If $\cA=\CC$, then $HC^{n}(\CC)=0$ if $n$ is odd and
$HC^{n}(\CC)=\CC$ if $n$ is even. For $n$ odd a non-trivial cocycle
$\phi\in C^n_\lambda(\CC)$ is given by
\[
\phi(a^0, a^{1},\cdots, a^n)=a^0a^{1}\cdots a^n, \quad a^0,
a^{1},\cdots, a^n \in \CC.
\]
\end{ex}

Equivalently, the cyclic cohomology can be described using a
$(b,B)$-bi\-comp\-lex. We define an operator $B : C^k(\cA,\cA^*)
\rightarrow C^{k-1}(\cA,\cA^*)$ by
\begin{equation}\label{e:defB}
B=AB_{0},
\end{equation}
where, for any $a^0,a^1,\ldots, a^{k-1}\in \cA$,
\begin{gather*}
     A\psi (a^0, \cdots, a^{k-1})=\sum_{j=0}^{k-1} (-1)^{(k-1)j}
    \psi(a^j,a^{j+1}\ldots, a^{k-1},a^0, a^1, \ldots,
    a^{j-1}),  \\
    B_{0}\psi (a^0, \ldots, a^{k-1}) =\psi(1, a^0, \ldots, a^{k-1}) -
    (-1)^k\psi(a^0, \ldots, a^{k-1},1).
\end{gather*}
One has $B^2=0$ and $b B=-B b$.

Consider the following double complex:
\[
C^{n,m}=C^{n-m}(\cA,\cA^*),\quad n,m\in\ZZ,
\]
with differentials $d_1:C^{n,m}\to C^{n+1,m}$ and $d_2:C^{n,m}\to
C^{n,m+1}$ given by
\[
d_1\psi=(n-m+1)b\psi,\quad d_2\psi=\frac{1}{n-m}B\psi,\quad \psi \in
C^{n,m}.
\]
For any $q\in \NN$ consider the complex $(F^qC,d)$, where
\[
(F^qC)^p=\bigoplus_{\substack{m\geq q,\\ n+m=p}}C^{n,m}, \quad
p\in\NN,\quad d=d_1+d_2.
\]
Then one has the isomorphism
\[
HC^n(\cA)\cong H^p(F^qC), \quad n=p-2q.
\]
This isomorphism associates with any $\psi \in HC^n(\cA)$ a cocycle
$\phi\in H^p(F^qC)$ with arbitrary fixed $p$ and $q$ satisfying
$n=p-2q$ which has a single non-vanishing component
$\phi_{p,q}=(-1)^{[n/2]}\psi$. In particular, any cocycle of the
complex $(F^qC,d)$ is cohomologous to a cocycle of the above form.

The periodic cyclic cohomology groups $HP^\ev (\cA)$ and $HP^\odd
(\cA)$ can also be defined by taking the inductive limit of the
groups $HC^k(\cA)$, $k\geq 0$, with respect to the periodicity
operator $S : HC^k(\cA)\to HC^{k+2}(\cA)$. For any $\psi\in
C^k_\lambda(\cA)$ one has
\begin{multline*}
S\psi(a^0,\ldots,a^{k+2})=-\frac{1}{(k+1)(k+2)}\sum_{j=1}^{k+1}\psi(a^0,\ldots,
a^{j-1}a^ja^{j+1},\ldots, a^{k+2})\\
-\frac{1}{(k+1)(k+2)}\sum_{0\leq i< j\leq k+1}\psi(a^0,\ldots,
a^{i-1}a^i,\ldots,a^ja^{j+1},\ldots, a^{k+2}).
\end{multline*}

In terms of the $(b,B)$-bicomplex the periodic cyclic cohomology is
described as the cohomology of the complex
\[
C^{\ev} (\cA) \stackrel{b+B}{\longrightarrow} C^{\odd}(\cA)
\stackrel{b+B}{\longrightarrow} C^{\ev}(\cA),
\]
where
\[
C^{\ev/\odd} (\cA) = \bigoplus_{k \ \text{even}/\text{odd}}
C^{k}(\cA).
\]

\begin{ex}\label{ex:C}
Let $\cA$ be the locally convex topological algebra $C^\infty(M)$ of
smooth functions on an $n$-dimensional compact manifold $M$, and let
$\cD_{k}(M)$ denote the space of $k$-dimensional de Rham currents on
$M$. The Hochschild cohomology and the cyclic cohomology of $\cA$
are computed in \cite{Co:nc}. The following map $\varphi\mapsto
C_\varphi$ defines an isomorphism of the continuous Hochschild
cohomology group $HH^k(\cA)$ with the space $\cD_{k}(M)$:
\begin{gather*}
\langle C_{\varphi},f^{0}df^{1}\wedge\cdots\wedge df^k\rangle=
\frac{1}{k!}\sum_{\sigma\in S_k}\varepsilon(\sigma)
\varphi\bigl(f^{0},f^{\sigma(1)},\dots,f^{\sigma(k)}\bigr),
\\
f^{0},f^{1},\dots,f^k\in C^{\infty}(M).
\end{gather*}
One has $C_{B\varphi}=k d^t C_\varphi$, where $d^t$ is the de Rham
boundary for currents. Therefore, for any $k$ the continuous cyclic
cohomology group $HC^k(\cA)$ is canonically isomorphic to the direct
sum
\[
\Ker d^t\oplus H_{k-2}(M,\CC)\oplus H_{k-4}(M,\CC)\oplus\cdots,
\]
where $H_{k}(M,\CC)$ denotes the usual de Rham homology. The
continuous periodic cyclic cohomology $HP^{\ev/\odd}(\cA)$ is
canonically isomorphic to the de Rham homology
\[
H_{\ev/\odd}(M,\CC)=\bigoplus_{k\ \text{even/odd}}H_{k}(M,\CC).
\]
\end{ex}

\begin{ex}\label{ex:phi}
We recall that a (homogeneous) $k$-cocycle on a discrete group
$\Gamma$ is a map $h:\Gamma^{k+1} \to \CC$ satisfying
\begin{gather*}
h(\gamma \gamma_0, \ldots ,\gamma \gamma_k) = h(\gamma_0, \ldots
,\gamma_k),\quad \gamma, \gamma_0, \ldots , \gamma_k \in \Gamma;\\
\sum_{i=0}^{k+1} (-1)^i h(\gamma_0, \ldots,
\gamma_{i-1},\gamma_{i+1},\ldots, \gamma_{k+1}) = 0, \quad \gamma_0,
\ldots, \gamma_{k+1}\in \Gamma.
\end{gather*}
With any homogeneous $k$-cocycle $h$ we can associate a
(non-homogeneous) $k$-cocycle $c\in Z^k(\Gamma, \mathbb C)$ by
$$ c(\gamma_1, \ldots, \gamma_k) = h (e, \gamma_1, \gamma_1
\gamma_2, \ldots, \gamma_1\ldots \gamma_k). $$ It is easy to check
that $c$ satisfies the following condition:
\begin{multline*}
c(\gamma_1, \gamma_2,\ldots, \gamma_k)
+\sum_{i=0}^{k-1}(-1)^{i+1}c(\gamma_0,\ldots, \gamma_{i-1},
\gamma_i\gamma_{i+1}, \gamma_{i+2},\ldots,\gamma_k)\\ + (-1)^{k+1}
c(\gamma_0,\gamma_1,\ldots,\gamma_{k-1})=0.
\end{multline*}

We say that a cocycle $c\in Z^k(\Gamma, \mathbb C)$ is normalized
(in the sense of Connes) if $c(\gamma_1,\gamma_2,\ldots,\gamma_k)$
equals zero in the case when either $\gamma_i = e$ for some $i$ or
$\gamma_1\ldots \gamma_k = e$. Every cohomology class in
$H^k(\Gamma, \mathbb C)$ can be represented by a normalized cocycle.

The group ring $\CC\Gamma $ consists of all functions
$f:\Gamma\to\CC$ with finite support. Multiplication in $\CC\Gamma $
is given by the convolution
\[
f_1\ast f_2(\gamma)=\sum_{\gamma_1\gamma_2=\gamma}f_1(\gamma_1)
f_2(\gamma_2), \quad \gamma \in \Gamma.
\]

A normalized cocycle $c\in Z^k(\Gamma, \mathbb C)$ determines a
cyclic $k$-cocycle $\tau_c$ on $\CC\Gamma $ by the formula
\[
\tau_c(f_0, \ldots, f_k) = \sum_{\gamma_0 \ldots \gamma_k
=e}f_0(\gamma_0)\ldots f_k(\gamma_k) c(\gamma_1,\ldots, \gamma_k),
\quad f_0,f_1,\ldots,f_k\in \CC\Gamma.
\]
\end{ex}

Let us now recall the definition of cyclic homology for an algebra
$\cA$ over $\CC$. Denote by $A^{\otimes,k+1}$ the tensor product of
$k+1$ copies of $A$ and consider the endomorphism $t$ of
$A^{\otimes,k+1}$ given by
\[
    t(a^0\otimes a^{1}\otimes\cdots \otimes a^k)=(-1)^k a^1\otimes \cdots \otimes a^k \otimes a^0,
    \quad a^0,a^1,\ldots, a^{k} \in
    A.
\]
Let us consider also the map $b$ from $A^{\otimes,k+1}$ to
$A^{\otimes,k}$ defined by the formula
\begin{align*}
b(a^0\otimes \cdots\otimes a^{k}) = &  \sum_{j=0}^{k-1} (-1)^j
a^0\otimes \cdots\otimes a^j a^{j+1}\otimes \cdots\otimes a^{k}  \\
&  + (-1)^{k} a^{k}a^0\otimes \cdots\otimes a^{k-1}, \quad
a^0,a^1,\ldots, a^{k} \in \cA.
  \end{align*}
Put
\[
C_k^\lambda(\cA)=\frac{A^{\otimes,k+1}}{\im (\Id-t)}
\]
The differential $b$ defines a map $b$ from $C_{k+1}^\lambda(\cA)$
to $C_{k}^\lambda(\cA)$, and cyclic homology $HC_{*}(\cA)$ of $\cA$
is defined as the homology of the complex $(C_*^\lambda(\cA),b)$.

If $A$ is a unital locally $m$-convex Frechet algebra, that is, a
unital algebra, which is a locally convex topological Frechet vector
space such that the product is continuous, then the topological
cyclic cohomology groups $HC^k(A)$ are defined in the same way as
above, using continuous $(k+1)$-linear functionals. Similarly, the
topological cyclic homology groups $HC_k(A)$ are defined in the same
way as above, using completed projective tensor products.

\subsection{Non-commutative differential forms}
Let $A$ be a unital algebra. A differential graded algebra is a
graded algebra
\[
\Omega_*(A)=\Omega_0(A)\oplus \Omega_1(A)\oplus \Omega_2(A)\oplus
\cdots
\]
endowed with a linear differentiation $d$ of degree $1$. Thus, for
any $j$ and $k$, one has $\Omega_j(A) \Omega_k(A)\subset
\Omega_{j+k}(A)$. The operator $d$ defines a map
\[
d : \Omega_j(A)\to \Omega_{j+1}(A),\quad j\geq 0
\]
and satisfies the conditions $d^2=0$ and
\[
d(\omega_j\cdot\omega_k)=d\omega_j\cdot\omega_k+(-1)^j\omega_j\cdot
d\omega_k, \quad \omega_j\in \Omega_j(A), \quad \omega_k\in
\Omega_k(A).
\]
We denote by $[\Omega_*(A),\Omega_*(A)]_l$ the linear subspace
spanned by the graded commutators
$[\omega_j,\omega_k]=\omega_j\cdot\omega_k-(-1)^{jk}\omega_k\cdot\omega_j$,
where $j+k=l$ and $\omega_j\in \Omega_j(A)$, $\omega_k\in
\Omega_k(A)$. Let
\[
\overline{\Omega_l(A)}=\frac{\Omega_l(A)}{[\Omega_*(A),\Omega_*(A)]_l}.
\]
The differential $d$ induces a linear differential, also denoted by
$d$, on the graded vector space
$\overline{\Omega_*(A)}=\bigoplus_l\overline{\Omega_l(A)}$. Denote
by $\overline{H}_*(A)$ the homology of this complex and call it the
non-commutative de Rham homology of the algebra $\Omega_*(A)$.

For a unital locally $m$-convex Frechet algebra $A$ one can define
the completion $\hat{\Omega}_*(A)$ of $\Omega_*(A)$, which is a
differential graded Frechet algebra. Then one can define the
non-commutative topological de Rham cohomology $\widehat{H}_*(A)$ as
the homology of the complex $
(\hat{\Omega}_*(A)/\overline{[\hat{\Omega}_*(A),\hat{\Omega}_*(A)]},d)$.

An example of a differential graded algebra is the universal
differential graded algebra $\Omega A$ of a unital algebra $A$. We
recall its construction. The algebra $\Omega A$ is a graded algebra:
$\Omega A=\bigoplus_{p=0}^\infty \Omega^p A$. In degree $0$ the
space $\Omega^0 A$ coincides with $A$. In degree $1$ the space
$\Omega^1 A$ is generated as a left $A$-module by symbols of degree
$1$ of the form $\delta a$, where $a\in A$, satisfying the relations
\begin{gather*}
    \delta(ab)=\delta(a)b+a\delta(b), \quad a,b\in A,\\
    \delta(\alpha a+\beta b)=\alpha \delta(a)+\beta \delta(b), \quad
    a,b\in A, \quad \alpha,\beta\in\CC.
\end{gather*}
In particular, one has $\delta(1)=0$. Thus, a general element of the
space $\Omega^1 A$ has the form $a=\sum_{i}a_i\delta b_i$, where
$a_i,b_i\in A$.

The algebra $\Omega A$ is generated as an algebra by the elements of
$\Omega^1 A$. In particular,
\[
(a_0\delta a_1)(b_0\delta b_1)=a_0\delta (a_1b_0)\delta b_1 - a_0
a_1 \delta b_0\delta b_1.
\]
Therefore, an arbitrary element of $\Omega A$ is represented as a
finite linear combination of the form $a_0\delta a_1\delta
a_2\ldots\delta a_p$ with some $a_0,a_1,a_2,\ldots,a_p\in A$. Left
and right multiplications by elements of $\Omega^0 A=A$ yield a an
$A$-$A$-bimodule structure on $\Omega A$.

A differential $\delta:\Omega^pA\to \Omega^{p+1}A$ is well defined
by
\[
\delta(a_0\delta a_1\delta a_2\ldots\delta a_p)=\delta a_0\delta
a_1\delta a_2\ldots\delta a_p, \quad a_0,a_1,a_2,\ldots,a_p\in A.
\]

If $A$ has an involution $a\mapsto a^*$, then the algebra $\Omega A$
has also a natural involutive algebra structure.

The non-commutative de Rham homology of the universal enveloping
algebra $\Omega A$ is called the non-commutative de Rham homology of
$A$ and denoted by $HDR_*(A)$. It is closely related with the cyclic
homology $HC_*(A)$ of $A$. Namely, for $n>0$ the non-commutative de
Rham homology group $HDR_n(A)$ coincides with the kernel of $B:
\overline{HC}_n(A)\to \overline{HH}_{n+1}(A)$ (see (\ref{e:defB})).
Here $\overline{HC}_n(A)$ and $\overline{HH}_{n+1}(A)$ denote the
reduced cyclic homology and reduced Hochschild homology respectively
(for more details, see \cite{Karoubi87,Loday}).

The use of differential graded algebras enables one to give a
general construction of cyclic cocycles on an arbitrary algebra
\cite{Co:nc}. We define a cycle of dimension $n$ to be a triple
$(\Omega, d, \int)$, where $(\Omega=\oplus_{j=0}^n\Omega^j, d)$ is a
differential graded algebra and $\int$ is a closed graded trace on
$(\Omega, d)$ of degree $n$. Here by a closed graded trace of degree
$n$ we mean a linear functional $\int:\Omega^n\to \CC$ satisfying
the following conditions:
\medskip\par
  (1) $\int \omega_2\omega_1=(-1)^{jk}\int
  \omega_1\omega_2$ for $\omega_1\in\Omega^{j}$,
  $\omega_2\in\Omega^{k}$;

  (2) $\int d\omega =0$ for $\omega\in\Omega^{n-1}$.
\medskip\par
Let $A$ be an algebra over $\CC$. A cycle over $A$ is defined to be
a cycle $(\Omega, d, \int)$ along with a homomorphism $\rho: A
\to\Omega^0$. For any cycle $(A\stackrel{\rho}{\to}\Omega, d, \int)$
over $A$, we can define its character by the formula
\[
\tau(a_0, a_1, \ldots, a_n)=\int \rho(a_0)d(\rho(a_1))\ldots
d(\rho(a_n)), \quad a_0,a_1,a_2,\ldots,a_n\in A.
\]
It is easy to check that $\tau $ is a cyclic cocycle on $A$.
Moreover, one can show that any cyclic cocycle on $A$ is the
character of some cycle over $A$.

\begin{ex}
Let $M$ be a smooth manifold without boundary, and consider the
graded algebra $\Omega$ of smooth differential forms on $M$:
$\Omega=\bigoplus_{p=0}^n\Omega^p$, where
$\Omega^p=C^\infty_c(M,\Lambda^pT^*M\otimes \CC)$. The de Rham
differential $d$ makes $\Omega$ into a differential graded algebra.
Finally, the linear functional $I(\omega)=\int_M\omega$ is a closed
graded trace of degree $n$ on $\Omega$. The corresponding cyclic
cocycle on $C^\infty_c(M)$ is given by
\[
\tau(f^0, f^1, \ldots, f^n)=\int_M f^{0}df^{1}\wedge \ldots
  \wedge df^n, \quad f^{0},f^{1}, \ldots, f^n \in C^\infty_c(M).
\]
Moreover, any closed $k$-dimensional de Rham current
$C\in\cD_{k}(M)$ on $M$ defines a closed graded trace of degree $k$
on $\Omega$:
\[
I_C(\omega)=\langle C, \omega\rangle,\quad  \omega \in \Omega^k.
\]
The character of the cycle $(\omega,d,I_C)$ over $C^\infty_c(M)$
coincides with the cyclic cocycle on $C^\infty_c(M)$ given by
    \begin{equation}\label{e:C}
  \psi_{C}(f^{0},f^{1}, \ldots, f^k)= \langle C, f^{0}df^{1}\wedge \ldots
  \wedge df^k\rangle,
        \quad f^{0},f^{1}, \ldots, f^k \in C^\infty_c(M).
    \end{equation}.
\end{ex}

\begin{ex}
Let $\Gamma$ be a discrete group. The universal differential graded
algebra $\Omega^*(\Gamma)=\Omega\CC\Gamma$ of the group algebra
$\CC\Gamma$ consists of finite linear combinations of symbols of the
form $\gamma_0d\gamma_1\ldots d\gamma_n$ with some $\gamma_0,
\gamma_1\ldots, \gamma_n\in\Gamma$. Any normalized cocycle $c\in
Z^k(\Gamma, \mathbb C)$ (see Example~\ref{ex:phi}) defines a closed
graded trace on $\Omega^*(\Gamma)$ by the formula
\[
\int \gamma_0d\gamma_1\ldots d\gamma_n =\begin{cases}
c(\gamma_1,\gamma_2,\ldots,\gamma_k), & \text{if}\ n=k\ \text{and}\
\gamma_0\gamma_1\ldots \gamma_n=e,\\
0 & \text{otherwise}
\end{cases}
\]
The character of this cycle coincides with the cyclic cocycle
$\tau_c$.
\end{ex}

\subsection{Non-commutative Chern-Weil construction}
For an algebra $\cA$ over $\CC$ the Chern character in $K$-homology
was constructed in \cite{Karoubi87} as a map $\ch : K_0(A)\to
HDR_{ev}(\cA)$. This construction is a straightforward
generalization of the classical Chern-Weil construction and makes
use of the notions of connection and curvature for finitely
generated projective $\cA$-modules.

In the cohomological setting a non-commutative analogue of the
Chern-Weil construction \cite{Co:nc} is a construction of a pairing
between $HC^*(\cA)$ and $K_*(\cA)$ for an arbitrary algebra $\cA$.

The pairing between $HC^\ev(\cA)$ and $K_{0}(\cA)$ is defined as
follows. For any cocycle $\varphi =(\varphi_{2k})$ in $C^\ev(\cA)$
and for any idempotent $e$ in $M_{q}(\cA)$ put
\begin{equation}
\langle[\varphi],[e]\rangle = \sum_{k \geq 0}(-1)^k \frac{(2k)!}{k!}
\varphi_{2k}\# \Tr \left(e-\frac{1}{2},e,\cdots,e\right),
\end{equation}
where $\varphi_{2k}\# \Tr$ is the $(2k+1)$-linear map on
$M_{q}(\cA)=M_{q}(\CC)\otimes \cA$ given by
\begin{equation}
  \varphi_{2k}\# \Tr(\mu^{0}\otimes a^0, \cdots, \mu^{2k}\otimes a^{2k}) =
    \Tr(\mu^{0}\ldots \mu^{2k}) \varphi_{2k}(a^0, \cdots, a^{2k}),
\end{equation}
for any $\mu^j \in M_{q}(\CC)$ and $a^j \in \cA$.

The pairing between $HC^\odd(\cA)$ and $K_{1}(\cA)$ is given by
 \begin{equation}
    \langle [\varphi], [U]\rangle = \frac1{\sqrt{2i\pi}} \sum_{k \geq 0}
    (-1)^k k! \varphi_{2k+1}\# \Tr
    (U^{-1}-1, U-1,\cdots, U^{-1}-1, U-1),
\end{equation}
where $\varphi =(\varphi_{2k+1})\in C^\odd(\cA)$ and $U\in
U_{q}(\cA)$.

It is important in index theory that the index
maps~(\ref{eq:LIFNCG.index-even}) and~(\ref{eq:LIFNCG.index-odd})
associated with a Fredholm module $(H,F)$ satisfying an additional
$p$-summability condition can be computed in terms of the pairing of
elements of $K^{*}(\cA)$ with certain cyclic cohomology class
$\tau=\ch_\ast (H,F) \in HC^n(\cA)$ called the Chern character of
the Fredholm module $(H,F)$ \cite{Co:nc}.

We recall that for any $p\geq 1$ the Schatten class $\cL^p(H)$
consists of all compact operators $T$ on a Hilbert space $H$ such
that $|T|^p$ is a trace class operator. Let $\mu_1(T)\geq
\mu_2(T)\geq \ldots$ be the singular numbers (the $s$-numbers) of a
compact operator $T$ on $H$, that is, the eigenvalues of
$|T|=\sqrt{T^*T}$ taken with multiplicities. Then
\[
T\in \cL^p(H) \Leftrightarrow \tr |T|^p = \sum_{n=1}^\infty
|\mu_n(T)|^p <\infty.
\]

\begin{defn}
A Fredholm module $(H,F)$ over an algebra $\cA$ is $p$-summable if
$(F^2-1)\rho(a)$, $(F-F^*)\rho(a)$ and $[F,\rho(a)]$ belong to
$\cL^p(H)$ for any $a\in \cA$.
\end{defn}

Let $(H,F)$ be a Fredholm module over an algebra $\cA$, and assume
that the module is $(p+1)$-summable and even if $p$ is even. It
determines a cycle over the algebra $\cA$ as follows. First of all,
one can assume without loss of generality that $F^2=1$. WE construct
a graded algebra $\Omega=\bigoplus_{j=0}^n\Omega^j$. For $k=0$ put
$\Omega^0=\cA$. For $k>0$ the space $\Omega^k$ is the linear span of
the bounded operators on $H$ of the form $\omega=a^0[F,a^1]\ldots
[F,a^k]$, where $a^0, a^1, \ldots, a^k\in \cA$. The product in
$\Omega$ is given by the product of operators. The differential
$d:\Omega\to \Omega$ is defined by
\[
d\omega=F\omega-(-1)^k\omega F, \quad \omega \in \Omega^k.
\]
Let us define a closed graded trace ${\rm Tr}_s : \Omega^n\to \CC$,
$n>p$. For any operator $T$ on $H$ such that $FT+TF\in \cL^1(H)$ we
put
\[
{\rm Tr}^\prime(T)=\frac12 \Tr (F(FT+TF)).
\]
Note that ${\rm Tr}^\prime(T)=\Tr(T)$ if $T\in \cL^1(H)$. For any
$\omega\in\Omega^n$ put ${\rm Tr}_s\omega = {\rm Tr}^\prime(\omega)$
if $n$ is odd and ${\rm Tr}_s\omega = {\rm Tr}^\prime(\gamma\omega)$
if $n$ is even, where $\gamma$ denotes the grading operator on $H$.
The character $\tau_n$ of the cycle described above is called the
Chern character of the Fredholm module $(H,F)$. Under the condition
$F^2=1$, $\tau_n$ is given for odd $n>p$ by
\begin{equation}\label{e:odd1}
\tau_n(a^0,a^1,\ldots,a^n)=\lambda_n\tr (a^0[F,a^1]\ldots [F,a^n]),
\quad a^0, a^1, \ldots, a^n\in\cA,
\end{equation}
and for even $n>p$ by
\begin{equation}\label{e:even1}
\tau_n(a^0,a^1,\ldots,a^n)=\lambda_n\tr (\gamma a^0[F,a^1]\ldots
[F,a^n]), \quad a^0, a^1, \ldots, a^n\in\cA,
\end{equation}
where the $\lambda_n$ are some constants depending only on $n$.

For different $n>p$ the characters $\tau_n$ agree in the sense that
one has the relation $S\tau_n=\tau_{n+2}$. Therefore, the class
$\tau=\ch_\ast (H,F) \in HP^*(\cA)$ in the periodic cyclic
cohomology of the algebra $\cA$ is well defined.

\subsection{Smooth algebras}\label{s:smooth}
In this subsection we give some general facts about smooth
subalgebras of $C^*$-algebras that amount to a non-commutative
analogue of the algebra of smooth functions on a smooth manifold.
Let $A$ be a $C^*$-algebra and $A^+$ the algebra obtained by
adjoining the unit to $A$. Let $\cA$ be a $*$-subalgebra of the
algebra $A$ and $\cA^+$ the algebra obtained by adjoining the unit
to $\cA$

\begin{defn}
We say that $\cA$ is a smooth subalgebra of a $C^*$-algebra $A$ if
$\cA$ is a dense $*$-subalgebra of $A$ that is stable under the
holomorphic functional calculus, that is, for any $a\in \cA^+ $ and
for any function $f$ holomorphic in a neighborhood of the spectrum
of $a$ (regarded as an element of the algebra $A^+$), we have $f(a)
\in \cA^+$.
\end{defn}

Suppose that ${\cA}$ is a dense $*$-subalgebra of $A$ endowed with
the structure of a Fr\'echet algebra whose topology is finer than
the topology induced by the topology of ${A}$. A necessary and
sufficient condition for ${\cA}$ to be a smooth subalgebra is given
by the spectral invariance condition (cf. \cite[Lemma 1.2]{Schw}):

\begin{itemize}
\item $\cA^+\cap GL(A^+) = GL(\cA^+)$, where $GL(\cA^+)$ and
$GL(A^+)$ denote the group of invertible elements in $\cA^+$ and
$A^+$, respectively.
\end{itemize}

This fact remains true in the case when ${\cA}$ is a locally
multiplicatively convex Fr\'echet algebra (that is, its topology is
given by a countable family of submultiplicative seminorms) such
that the group $GL(\widetilde {\cA})$ of invertibles is open
\cite[Lemma 1.2]{Schw}.

If $\cA$ is a smooth subalgebra of a $C^*$-algebra $A$, then, for
any $n$, the algebra $M_n(\cA)$ is a smooth subalgebra of the
$C^*$-algebra $M_n(A)$. In particular, $GL_n(\cA^+)$ coincides with
the intersection $M_n(\cA^+)\cap GL_n(A^+)$. Let us regard
$GL_n(\cA^+)$ as a topological group equipped with the induced
topology of the spaces $GL_n(A^+)$. Denote by $GL_\infty(\cA^+)$ the
inductive limit of the topological groups $GL_n(\cA^+)$.

One of the most important properties of smooth subalgebras consists
in the following fact, an analogue of smoothing in the operator
$K$-theory (cf. \cite[Sect. VI.3]{Co81}, \cite{Bost}).

\begin{thm}\label{t:smooth}
If $\cA$ is a smooth subalgebra of a $C^*$-algebra $A$, then the
inclusion $\cA\hookrightarrow A$ induces isomorphisms
\[
K(\cA) \cong K(A)\ \text{and}\ \pi_n(GL_\infty(\cA^+))\to
\pi_n(GL_\infty(A^+))=K_{n+1}(A).
\]
\end{thm}

For any $p$-summable Fredholm module $(H,F)$ over an algebra $\cA$
the algebra
\[
\cC=\{T\in \bar{\cA} : [F,T]\in \cL^p(H)\},
\]
where $\bar{\cA}$ denotes the uniform closure of $\cA$ in $\cL(H)$,
is a smooth subalgebra of the $C^*$-algebra $\bar{\cA}$. Moreover,
one can show that the cocycle $\tau_n$ on $\cA$ defined by
(\ref{e:odd1}) and (\ref{e:even1}) extends by continuity to a cyclic
cocycle on $\cC$. This along with Theorem~\ref{t:smooth} allows us
to assert that the pairing with $\ch_\ast (H,F) \in HC^n(\cA)$
defines a map $K_*(\bar{\cA})\cong K_*(\cC)\to \ZZ$. This fact is
called the topological invariance property of the Chern character
$\ch_\ast (H,F)$ of the Fredholm module $(H,F)$.

\subsection{Non-commutative Riemannian geometry}\label{s:riem}
According to \cite{Co-M,Sp-view}, the initial data of the
non-commutative Riemannian geometry is a spectral triple.
\begin{defn}
A spectral triple is a set $({\mathcal A}, {\mathcal H}, D)$, where:
\begin{enumerate}
\item ${\mathcal A}$ is an involutive algebra;
\item ${\mathcal H}$ is a Hilbert space equipped with a
$\ast$-representation of ${\mathcal A}$;
\item $D$ is an (unbounded) self-adjoint operator acting in ${\mathcal H}$
such that
\begin{enumerate}
\item[(i)] for any $a\in {\mathcal A}$ the operator $a(D-i)^{-1}$
is a compact operator on ${\mathcal H}$;

\item[(ii)]~$D$ almost
commutes with elements of ${\mathcal A}$ in the sense that $[D,a]$
is bounded for any $a\in {\mathcal A}$.
\end{enumerate}
\end{enumerate}

A spectral triple is said to be even if $\cH$ is endowed with a
$\ZZ_{2}$-grading $\gamma \in \cL(\cH)$, $\gamma=\gamma^{*}$,
$\gamma^{2}=1$, and moreover $ \gamma D=-D\gamma$ and $\gamma a =a
\gamma$ for any $a\in \cA$. Otherwise, a spectral triple is said to
be odd.
\end{defn}

Spectral triples were considered for the first time in the paper
\cite{Baaj-Julg}, where they were called unbounded Fredholm modules.
A spectral triple $(\cA, \cH, D)$ determines a Fredholm module
$(\cH,F)$ over $\cA$, where $F=D(I+D^2)^{-1/2}$ \cite{Baaj-Julg}. In
a certain sense, the operator $F$ is connected with measurement of
angles and is responsible for a conformal structure, whereas $|D|$
is connected with measurement of lengths.

\begin{defn}\label{d:dim}
A spectral triple $({\mathcal A}, {\mathcal H}, D)$ is $p$-summable
(or $p$-di\-men\-si\-o\-nal), if for any $a\in {\mathcal A}$ the
operator $a(D-i)^{-1}$ is an element of the Schatten class
$\cL^p({\mathcal H})$.

A spectral triple $({\mathcal A}, {\mathcal H}, D)$ is said to be
finite-dimensional if it is $p$-summable for some $p$.

The greatest lower bound of all $p$ such that a finite-dimensional
spectral triple is $p$-summable is called the dimension of the
spectral triple.
\end{defn}

The dimension of a spectral triple $({\mathcal A}, {\mathcal H}, D)$
coincides with the dimension of the corresponding Fredholm module
$(\cH,F)$, $F=D(I+D^2)^{-1/2}$, over $\cA$ \cite{Baaj-Julg}.

Classical Riemannian geometry is described by the spectral triple
$(\cA,{\mathcal H},D)$ associated with a compact Riemannian spin
manifold $(M,g)$:
\begin{enumerate}
\item the involutive algebra ${\mathcal A}$ is the algebra $C^{\infty}(M)$
of smooth functions on $M$;
\item the Hilbert space ${\mathcal H}$ is the space
$L^2(M,F(TM))$, where the algebra ${\mathcal A}$ acts by
multiplication;
\item the operator $D$ is the spin Dirac operator.
\end{enumerate}

The Weyl asymptotic formula for eigenvalues of self-adjoint elliptic
operators on a compact manifold implies at once that this spectral
triple is finite-dimensional and has dimension equal to the
dimension of $M$.

Let $({\mathcal A}, {\mathcal H}, D)$ be a spectral triple. Assume
for simplicity that the algebra $\cA$ has a unit. We consider the
operator $|D|=(D^2)^{1/2}$. Denote by $\delta$ the (unbounded)
differentiation on ${\mathcal L}({\mathcal H})$ given by
\begin{equation}
\label{derivative} \delta(T)=[|D|,T],\quad T\in \Dom\delta\subset
{\mathcal L}({\mathcal H}).
\end{equation}
For any $T\in {\mathcal L}({\mathcal H})$ denote by $\delta^i(T)$
the $i$th commutator with $|D|$.

\begin{defn}\label{d:op0}
The space ${\rm OP}^0$ consists of all $T\in {\mathcal L}({\mathcal
H})$ such that $\delta^i(T) \in {\mathcal L}({\mathcal H})$ for any
$i\in\NN$:
\[
{\rm OP}^{0}=\bigcap_n {\rm Dom}\;\delta^n.
\]
\end{defn}

The space ${\rm OP}^0$ is a smooth subalgebra of the $C^*$-algebra
${\mathcal L}({\mathcal H})$ (for instance, see \cite[Theorem
1.2]{Ji}). The uniform closure $\bar{\cA}$ of $\cA$ in ${\mathcal
L}({\mathcal H})$ can be regarded as the algebra of continuous
functions on some virtual topological space. In a certain sense the
algebra ${\rm OP}^0\cap \bar{\cA}$ consists of functions on the
given space, which are infinitely differentiable in the quantum
sense. For the spectral triple associated with a compact Riemannian
manifold, the algebra ${\rm OP}^0\cap C(M)$ contains the algebra
$C^\infty(M)$. We refer the reader to \cite{Connes:reconstruction}
for the problem of reconstructing the smooth structure of a manifold
$M$ from the associated spectral triple $(\cA,{\mathcal H},D)$.

We recall that the Riemannian volume form $d\nu$ associated with a
Riemannian metric $g$ is given in local coordinates by the formula
$d\nu=\sqrt{\det g}\,dx$. To define its non-commutative analogue one
uses the trace $\operatorname{Tr}_\omega$ introduced by Dixmier in
\cite{Dixmier66} as an example of a non-standard trace on
$\cL(\cH)$.

Consider the ideal $\cL^{1+}(\cH)$ in the algebra $\cK(\cH)$ of
compact operators, which consists of all $T\in \cK(\cH)$ such that
\[
\sup_{N\in\NN}\frac{1}{\ln N}\sum_{n=1}^N\mu_n(T)<\infty,
\]
where $\mu_1(T)\geq \mu_2(T)\geq \ldots$ are the singular numbers of
$T$. For any invariant mean $\omega$ on the amenable group of upper
triangular $(2\times 2)$-matrices, Dixmier constructed a linear form
$\lim_\omega$ on the space ${\ell}^\infty(\NN)$ of bounded sequences
which has the following properties:
\begin{enumerate}
  \item $\lim_\omega$ coincides with the functional of taking limit
  $\lim$ on the subspace of convergent sequences,
  \item $\lim_\omega\{c_n\}\geq 0$ if $c_n\geq 0$ for any
  $n$,
  \item $\lim_\omega\{c^\prime_{n}\}=\lim_\omega\{c_{n}\}$
  where $\{c^\prime_{n}\}=\{c_1,c_1,c_2,c_2,c_3,c_3,\ldots \}$,
  \item $\lim_\omega\{c_{2n}\}=\lim_\omega\{c_{n}\}$.
\end{enumerate}
For a positive operator $T\in \cL^{1+}(\cH)$ the value of
$\operatorname{Tr}_\omega$ on $T$ is given by
\[
\operatorname{Tr}_\omega(T)={\lim}_\omega \frac{1}{\ln N}
\sum_{n=1}^N\mu_n(T).
\]

Let $M$ be an $n$-dimensional compact manifold, $E$ a vector bundle
on $M$, and $P\in\Psi^m(M,E)$ a classical pseudodifferential
operator. Thus, in any local coordinate system its complete symbol
$p$ can be represented as an asymptotic sum $p\sim
p_m+p_{m-1}+\ldots$, where $p_l(x,\xi)$ is a homogeneous function of
degree $l$ in $\xi$. As shown in \cite{Co-action}, the Dixmier trace
$\operatorname{Tr}_\omega(P)$ does not depend on the choice of
$\omega$ and coincides with the value $\tau(P)$ of a trace $\tau$
introduced by Wodzicki \cite{Wo} and Guillemin \cite{Gu85} on the
algebra $\Psi^{*}(M,E)$ of classical pseudodifferential operators of
arbitrary order. The trace $\tau$ is defined as follows. For
$P\in\Psi^{*}(M,E)$ the density $\rho_P$ is defined in local
coordinates as $$ \rho_P = \left(\int_{|\xi|=1}\Tr
p_{-n}(x,\xi)\,d\xi\right) |dx|. $$ The density $\rho_P$ turns out
to be independent of the choice of a local coordinate system, and
therefore gives a well-defined density on $M$. The integral of
$\rho_P$ over $M$ is the value of the Wodzicki-Guillemin trace:
\begin{equation}\label{e:wodz}
\tau(P)=\frac{1}{(2\pi)^n}\int_M\rho_P
=\frac{1}{(2\pi)^n}\int_{S^*M}\Tr p_{-n}(x,\xi)\, dx d\xi.
\end{equation}
Wodzicki \cite{Wo} showed that $\tau$ is a unique trace on
$\Psi^*(M,E)$.

According to \cite{Co-action} (cf. also \cite{Gracia:book}), any
operator $P\in\Psi^{-n}(M,E)$ belongs to the ideal
$\cL^{1+}(L^2(M,E))$ and, for any invariant mean $\omega$
\[
\operatorname{Tr}_\omega(P)=\tau(P).
\]
The above results imply the formula
\[
\int_Mf\,d\nu=c(n)\operatorname{Tr}_\omega(f|D|^{-n}), \quad f\in
\cA,
\]
where $c(n)=2^{(n-[n/2])} \pi^{n/2} \Gamma(\frac{n}{2}+1)$. Thus,
the Dixmier trace $\operatorname{Tr}_\omega$ can be considered as a
proper non-commutative generalization of the integral.

\subsection{Non-commutative local index theorem}
If one looks at a geometric space as the union of parts of different
dimensions, the notion of dimension introduced in
Definition~\ref{d:dim} gives only the maximum of the dimensions of
the parts of this space. To take into account lower dimensional
parts of the geometric space, Connes and Moscovici \cite{Co-M}
proposed taking as a more correct notion of dimension of a smooth
spectral triple not a single real number $d$, but a subset ${\rm
Sd}\subset {\mathbb C}$, called its dimension spectrum.

\begin{defn}
A spectral triple $({\mathcal A}, {\mathcal H}, D)$ is said to be
smooth if for any $a\in {\mathcal A}$ one has the inclusions $a,
[D,a] \in  {\rm OP}^{0}$.
\end{defn}

The spectral triple associated with a smooth Riemannian manifold $M$
is smooth. The smoothness condition of a spectral triple $({\mathcal
A}, {\mathcal H}, D)$ can be treated in the following sense: the
algebra $\cA$ consists of smooth (in the quantum sense) functions on
the corresponding non-commutative space (see comments after
Definition~\ref{d:op0}).

Let $({\mathcal A}, {\mathcal H}, D)$ be a smooth spectral triple.
Denote by ${\mathcal B}$ the algebra generated by all elements of
the form $\delta^n(a)$, where $a\in{\mathcal A}$ and $n\in {\mathbb
N}$. Thus, $\cB$ is the smallest subalgebra in ${\rm OP}^{0}$
containing $\cA$ and invariant under the action of $\delta$.

\begin{defn}
A spectral triple $({\mathcal A},{\mathcal H},D)$ has discrete
dimension spectrum ${\rm Sd}\subset{\mathbb C}$, if ${\rm Sd}$ is a
discrete subset in ${\mathbb C}$, the triple is smooth, and for any
$b\in {\mathcal B}$ the distributional zeta-function $\zeta_b(z)$ of
$|D|$ given by $$ \zeta_b(z)=\tr b|D|^{-z},
$$ is defined in the half-plane $\{z\in\CC:{\rm Re}\, z>d\}$ and
extends to a holomorphic function on ${\mathbb C}\backslash {\rm
Sd}$.

The dimension spectrum is said to be simple if the singularities of
$\zeta_b(z)$ at $z\in {\rm Sd}$ are at most simple poles.
\end{defn}

\begin{ex}\label{ex:classic}
Let $M$ be a compact manifold of dimension $n$, $E$ a vector bundle
on $M$, and $D\in\Psi^1(M,E)$ a self-adjoint elliptic operator. Then
the triple $(C^{\infty}(M),L^{2}(M,E),D)$ is a smooth
$n$-dimensional spectral triple. The algebra $\mathcal B$ is
contained in the algebra $\Psi^0(M,E)$ of zero order classical
pseudodifferential operators. For any classical pseudodifferential
operator $P\in \Psi^m(M,E)$, $m\in\ZZ$, the function $z \rightarrow
\tr P|D|^{-z}$ has a meromorphic extension to $\CC$ with at most
simple poles at integer points $k$, $k\leq m+n$. The residue of this
function at $z=0$ coincides with the Wodzicki-Guillemin residue
$\tau(P)$ of $P$ (see~(\ref{e:wodz}):
\begin{equation}
\underset{z=0}{\res} \tr P|D|^{-z}=\tau(P).
     \label{eq:LIFNCG-NCR}
\end{equation}

In particular, this spectral triple has simple discrete dimension
spectrum lying in $\{k \in \ZZ : \ k\leq n\}$.
\end{ex}

In \cite{Co-M,Sp-view} a definition is given of the algebra
$\Psi^*({\mathcal A})$ of pseudodifferential operators associated
with a smooth spectral triple $({\mathcal A},{\mathcal H},D)$ in the
case when the algebra $\cA$ is unital. By the spectral theorem, for
any $s\in\RR$ the operator $|D|^s$ is a well-defined positive
self-adjoint operator acting in $\cH$ which is unbounded for $s>0$.
For any $s\geq 0$ we denote by $\cH^s$ the domain of the operator
$\langle D\rangle^s$, and for $s<0$ we put $\cH^s=(\cH^{-s})^*$.
Also, let $\cH^{\infty}=\bigcap_{s\geq 0}\cH^s, \quad
\cH^{-\infty}=(\cH^{\infty})^*$. We say that a bounded operator $P$
on the space $\cH^{\infty}$ belongs to the class ${\rm
OP}^{\alpha}$, if $P\langle D\rangle^{-\alpha}\in {\rm OP}^0$.

\begin{defn}
We say that an operator $P : \cH^{\infty}\to \cH^{-\infty}$ belongs
to the class $\Psi^*({\mathcal A})$ if it admits an asymptotic
expansion:
\[
P\sim \sum_{j=0}^{+\infty}b_{q-j} |D|^{q-j}, \quad b_{q-j}\in
{\mathcal B},
\]
which means that for any $N$
\[
P - \left(b_q|D|^q + b_{q-1}|D|^{q-1}+\ldots+b_{-N}|D|^{-N}
\right)\in {\rm OP}^{-N-1}.
\]
\end{defn}

It was proved in \cite[Appendix B]{Co-M} that $\Psi^*({\mathcal A})$
is an algebra. For the spectral triple
$(C^{\infty}(M),L^{2}(M,E),D)$ described in Example~\ref{ex:classic}
the algebra $\Psi^*({\mathcal A})$ is contained in $\Psi^*(M,E)$.

Recall that a spectral triple $(\cA, \cH, D)$ defines the Fredholm
module $(\cH,F)$ over $\cA$, where $F=D(I+D^2)^{1/2}$, and thereby
the index map $\operatorname{Ind} : K_*(\cA)\to\CC$ (see the
formulae (\ref{eq:LIFNCG.index-even}) and
(\ref{eq:LIFNCG.index-odd})). As shown above, this map can be
expressed in terms of the pairing with the cyclic cohomology class
$\operatorname{ch}_\ast (\cH,F) \in HP^*(\cA)$, the Chern character
of the Fredholm module $(\cH,F)$ (see the formulae (\ref{e:odd1})
and (\ref{e:even1})). Let us call $\operatorname{ch}_\ast (\cH,F)$
by the Chern character of the spectral triple $(\cA, \cH, D)$ and
denote it by $\operatorname{ch}_\ast (\cA, \cH, D)$. The formulae
(\ref{e:odd1}) and (\ref{e:even1}) have a defect, in that they
express the map $\operatorname{ind}$ in terms of the operator
traces. In the classical case these traces are non-local
functionals, and it is impossible to compute them in coordinate
charts. To correct this defect, Connes and Moscovici proved for the
index map another formula, which involves Wodzicki-Guillemin trace
type functionals. These functionals are local in the sense of
non-commutative geometry, because they vanish on any trace class
operator in $\cH$. Therefore, the Connes-Moscovici formula can be
naturally called the non-commutative local index theorem.

Suppose that a spectral triple $(\cA,\cH,D)$ is smooth and, for
simplicity, has simple discrete dimension spectrum. We define the
non-commutative integral determined by this spectral triple by
setting
\begin{equation}
    \bint b = \underset{z=0}{\res} \tr b|D|^{-z}, \quad b\in \cB.
\end{equation}
The functional $\bint$ is a trace on $\cB$ which is local in the
sense of non-commutative geometry.

\begin{thm}[{\cite[Thm.~II.3]{Co-M}}]\label{t:cmeven}
Suppose that $(\cA,\cH,D)$ is an even spectral triple, which is
$p$-summable and has simple discrete dimension spectrum.

\begin{enumerate}
  \item An even cocycle
$\varphi_{\mathop{CM}}^\ev =(\varphi_{2k})$ in the $(b,B)$-bicomplex
of $\cA$ is defined by the following formulae: for $k=0$,
             \begin{equation*}
    \varphi_{0}(a^{0})  = \underset{z=0}{\res}z^{-1}\tr \gamma a^0|D|^{-z},
    \end{equation*}
and for $k\neq0$,
    \begin{equation*}
        \varphi_{2k}(a^{0}, \ldots, a^{2k}) =   \sum_{\alpha\in
        \NN^n}c_{k,\alpha}\bint \gamma a^0 [D,a^1]^{[\alpha_{1}]} \ldots
          [D,a^{2k}]^{[\alpha_{2k}]} |D|^{-2(|\alpha|+k)},
       \end{equation*}
where
\[
c_{k,\alpha}=
\frac{(-1)^{|\alpha|}2\Gamma(|\alpha|+k)}{\alpha!(\alpha_{1}+1)
\cdots (\alpha_{1}+\cdots +\alpha_{2k}+ 2k)},
\]
and the symbol $T^{[j]}$ denotes the $j$th iterated commutator with
$D^2$.
  \item The cohomology class defined by $\varphi_{\mathop{CM}}^\ev$ in $HP^\ev
(\cA)$ coincides with the Chern character $\operatorname{ch}_\ast
(\cA, \cH, D)$.
\end{enumerate}
\end{thm}

\begin{theorem}[{\cite[Thm.~II.2]{Co-M}}]\label{t:cmodd}
Suppose that $(\cA,\cH,D)$ is a spectral triple, which is
$p$-summable and has simple discrete dimension spectrum. Then
\begin{enumerate}
  \item An odd cocycle $\varphi_{\mathop{CM}}^\odd
=(\varphi_{2k+1})$ in the $(b,B)$-bi\-comp\-lex of $\cA$ is defined
by
     \begin{multline*}
             \varphi_{2k+1}(a^{0}, \ldots, a^{2k+1})\\ =
            \sqrt{2i\pi}
               \sum_{\alpha\in \NN^n}c_{k,\alpha}
               \bint a^0 [D,a^1]^{[\alpha_{1}]} \ldots
              [D,a^{2k+1}]^{[\alpha_{2k+1}]} |D|^{-2(|\alpha|+k)-1},
       \end{multline*}
where
\[
c_{k,\alpha}= \frac{(-1)^{|\alpha|}\Gamma(|\alpha|+k
+\frac12)}{\alpha!(\alpha_{1}+1) \cdots (\alpha_{1}+\cdots
+\alpha_{2k+1}+2k +1)}.
\]
  \item the cohomology class defined by $\varphi_{\mathop{CM}}^\odd$
in $HP^\odd (\cA)$ coincides with the Chern character
$\operatorname{ch}_\ast (\cA, \cH, D)$.
\end{enumerate}
 \end{theorem}

\begin{ex}
Let $M$ be a compact manifold of dimension $n$ and $D$ a first
order, self-adjoint, elliptic, pseudodifferential operator on $M$
acting on sections of a vector bundle $E$ on $M$. Then the
non-commutative integral $\bint$ defined by the spectral triple
$(C^{\infty}(M),L^{2}(M,E),D)$ coincides with the Wodzicki-Guillemin
trace $\tau$ (see~(\ref{e:wodz}) and Example~\ref{ex:classic}).

In the case when $D$ is the spin Dirac operator on a compact
Riemannian spin manifold $M$ in Theorem~\ref{t:cmeven}, we have for
any $f^0, f^1,\ldots, f^{m}\in C^{\infty}(M)$ that
\[
\tau (\gamma f^0 [D,f^1]^{[\alpha_{1}]} \ldots
          [D,f^{m}]^{[\alpha_{m}]} |D|^{-(2|\alpha|+m)})=0,
\]
when $|\alpha|\neq 0$, and
\[
\tau (\gamma f^0 [D,f^1] \ldots [D,f^{m}] |D|^{-m}) = c_m\int_M
f^0\,df^1\wedge \ldots df^{m}\wedge \hat{A}(TM,\nabla),
\]
when $\alpha_1=\alpha_2=\ldots=\alpha_{2k}=0$, where $c_m$ is some
constant.

If the dimension of $M$ is even, then the spectral triple is even
and the components of the corresponding even cocycle
$\varphi_{\mathop{CM}}^\ev = (\varphi_{2k})$ are given by
\[
\varphi_{2k}(f^0,\ldots,f^{2k})=\frac{1}{(2k)!}\int_M
f^0\,df^1\wedge \ldots df^{2k}\wedge \hat{A}(TM,\nabla)^{(n-2k)},
\]
where $f^0, f^1,\ldots, f^{2k}\in C^{\infty}(M)$.

If the dimension of $M$ is odd, then the spectral triple is odd and
the components of the corresponding odd cocycle
$\varphi_{\mathop{CM}}^\odd = (\varphi_{2k+1})$ are given by
\begin{multline*}
\varphi_{2k+1}(f^0,\ldots,f^{2k})=\sqrt{2i\pi}
\frac{(2i\pi)^{-[n/2]+1}}{(2k+1)!} \\ \times \int_M f^0\,df^1\wedge
\ldots df^{2k+1}\wedge \hat{A}(TM,\nabla)^{(n-2k-1)},
\end{multline*}
where $f^0, f^1,\ldots, f^{2k+1}\in C^{\infty}(M)$.
\end{ex}

\subsection{Semifinite spectral triples}\label{s:semifinite}
The study of index theory problems in various situations such as
measurable foliations, Galois coverings, almost-periodic operators
has served as a motivation for extending methods of non-commutative
geometry to the case when the algebra of bounded operators on a
Hilbert space is replaced by an arbitrary semifinite von Neumann
algebra. In this subsection we give some information from semifinite
non-commutative geometry (for further information see, for instance,
a survey \cite{Carey-S} and the references therein).

\begin{defn}
A semifinite spectral triple is a set $({\mathcal A}, {\mathcal H},
D)$, where:
\begin{enumerate}
\item ${\mathcal H}$ is a Hilbert space;
\item ${\mathcal A}$ is an involutive subalgebra of a semifinite
von Neumann algebra $\cN$ acting in ${\mathcal H}$
\item $D$ is an (unbounded) self-adjoint operator acting in ${\mathcal
H}$ and affiliated to $\cN$ such that

\begin{enumerate}
\item[(i)] the operator $(D-i)^{-1}$ is a $\tau$-compact operator
in $\cN$ (relative to some faithful, normal, semifinite trace $\tau$
on $\cN$);

\item[(ii)]~the operator $[D,a]$ is a bounded operator,
belonging to $\cN$, for any $a\in {\mathcal A}$.
\end{enumerate}
\end{enumerate}

A spectral triple is said to be even if the space $\cH$ is equipped
with a $\ZZ_{2}$-grading $\gamma \in \cL(\cH)$, $\gamma=\gamma^{*}$,
$\gamma^{2}=1$, and moreover $ \gamma D=-D\gamma$ and $\gamma a =a
\gamma$ for any $a\in \cA$. Otherwise, the spectral triple is said
to be odd.
\end{defn}

For an element $S\in \cN$, its $t$th generalized singular number
($t\in\RR$) is given by
\[
\mu_t(S)=\inf\{\|SE\|\ | E\ \text{is a projection in}\ \cN\
\text{with}\ \tau(1-E)\leq t\}.
\]

The space $\cL^{(1,\infty)}({\mathcal N})$ consists of elements
$T\in\cN$ such that
\[
\|T\|_{\cL^{(1,\infty)}}=\sup_{t>1}\frac{1}{\log(1+t)}\int_0^t\mu_s(T)\,ds<\infty.
\]
For any $p>1$ put
\[
\psi_p(t)=\begin{cases} t & \text{for}\ 0\leq t\leq 1,\\
t^{1-\frac{1}{p}} & \text{for}\ 1\leq t.
\end{cases}
\]
The space $\cL^{(p,\infty)}({\mathcal N})$ consists of $T\in\cN$
such that
\[
\|T\|_{\cL^{(p,\infty)}}=\sup_{t>1}\frac{1}{\psi_p(t)}\int_0^t\mu_s(T)\,ds<\infty.
\]

\begin{defn}
A semifinite spectral triple $({\mathcal A}, {\mathcal H}, D)$ is
$(p,\infty)$-summable, if for any $a\in {\mathcal A}$ the operator
$(1+D^2)^{-1/2}$ is an elements of the class
$\cL^{(p,\infty)}({\mathcal N})$.
\end{defn}

\begin{defn}
A pre-Fredholm module over a unital Banach algebra $\cA$ is a pair
$(\cH,F)$, where
\begin{enumerate}
  \item $\cA$ has a continuous representation in a semifinite von Neumann algebra
  $\cN$ acting in a Hilbert space $\cH$;
  \item $F$ is a self-adjoint Breuer-Fredholm operator acting in $\cH$ such
  that $1-F^2\in \cK_\cN$ and $[F,a]\in \cK_\cN$ for any
  $a\in A$.
\end{enumerate}

If $1-F^2=0$, then $(\cH,F)$ is called a Fredholm module.

A pre-Fredholm module $(\cH,F)$ is said to be even if the Hilbert
space $\cH$ is equipped with a $\ZZ_2$-grading $\gamma$ such that
the operators $\rho(a)$ are even, $\gamma a= a\gamma$, and the
operator $F$ is odd, $\gamma F=-F\gamma$. Otherwise, it is said to
be odd.
\end{defn}

\begin{defn}
A pre-Fredholm module $(\cH,F)$ is $(p,\infty)$-summable, if
$1-F^2\in \cL^{(p/2,\infty)}({\mathcal N})$ and $[F,a]\in
\cL^{(p,\infty)}({\mathcal N})$ for a dense set of elements $a\in
\cA$.
\end{defn}

If $({\mathcal A}, {\mathcal H}, D)$ is a semifinite spectral
triple, then the pair $(\cH,F)$. where $F=D(I+D^2)^{1/2}$, is a
pre-Fredholm module over $\cA$. Using a faithful normal semifinite
trace $\tau$ on the von Neumann algebra $\cN$ instead of the
standard trace on the algebra $\cL(\cH)$, one can define by standard
formulae the Chern character of an arbitrary $(p,\infty)$-summable
pre-Fredholm module $(\cH,F)$, and, hence, the Chern character of
any finite-dimensional semifinite spectral triple $({\mathcal A},
{\mathcal H}, D)$. In the papers \cite{CPRS2,CPRS3} a
non-commutative local index theorem is proved for semifinite
spectral triples. Its formulation is similar to that of the
non-commutative local index theorem (see Theorems~\ref{t:cmeven} and
\ref{t:cmodd}), with the sole difference that the role of the
Dixmier trace $\operatorname{Tr}_\omega$ in these theorems is played
by its generalization to the case of an arbitrary semifinite von
Neumann algebra. For Dixmier traces and general singular traces on
semifinite von Neumann algebras and their applications, see the
survey \cite{Carey-S} and its references.

\subsection{Non-commutative spectral geometry and type III}
The authors of \cite{CoM:twisted} introduced the notion of twisted
(or $\sigma$-spectral) triple, making it possible to apply the
methods of non-commutative geometry for certain type III
non-commutative spaces.

\begin{defn}
Let ${\mathcal A}$ be an algebra equipped with an automorphism
$\sigma$. An ungraded $\sigma$-spectral triple over $\cA$ is defined
to be a set $({\mathcal A}, {\mathcal H}, D)$, where:
\begin{enumerate}
\item ${\mathcal H}$ is a Hilbert space endowed with an action of $\cA$;
\item $D$ is an (unbounded) self-adjoint operator acting in ${\mathcal H}$
such that the operator $(D-i)^{-1}$ is a compact operator in $\cH$
and the operator $Da-\sigma(a)D$ is bounded for any $a\in {\mathcal
A}$.
\end{enumerate}

If the algebra ${\mathcal A}$ is involutive and its representation
in the Hilbert space ${\mathcal H}$ is a $\ast$-representation, then
in addition we impose the unitarity condition
\[
\sigma(a^*)=(\sigma^{-1}(a))^*, \quad a\in \cA.
\]

A graded $\sigma$-spectral triple is defined in a similar way, but
in this case the space $\cH$ is equipped with a $\ZZ_{2}$-grading
$\gamma \in \cL(\cH)$, $\gamma=\gamma^{*}$, $\gamma^{2}=1$, and we
have $\gamma D=-D\gamma$ and $\gamma a =a \gamma$ for any $a\in
\cA$.
\end{defn}

\begin{ex}
Let $({\mathcal A}, {\mathcal H}, D)$ be an arbitrary spectral
triple and $h\in\cA$ a self-adjoint element, $h=h^*$. Consider the
automorphism $\sigma$ of $\cA$ given by
\[
\sigma(a)=e^{2h}ae^{-2h},\quad a\in \cA.
\]
The ``perturbed'' spectral triple $({\mathcal A}, {\mathcal H},
D^\prime)$, where $D^\prime=e^hDe^{-h}$, is a $\sigma$-spectral
triple.
\end{ex}

Another example of $\sigma$-spectral triples arises from transverse
geometry of an arbitrary codimension one foliation (see
Example~\ref{ex:circle}). The most vital open question is to extend
the notions mentioned above to the case of higher codimension
foliations. It is expected that the general case will require the
use of dual actions of Lie groups such as $GL(n)$ and, more
generally, of quantum groups.

\begin{defn}
A $\sigma$-spectral triple $({\mathcal A}, {\mathcal H}, D)$ is
called Lipschitz regular if the operator $|D|a-\sigma(a)|D|$ is
bounded for any $a\in {\mathcal A}$.
\end{defn}

If a $\sigma$-spectral triple $({\mathcal A}, {\mathcal H}, D)$ is
Lipschitz regular and $F=D|D|^{-1}$, then $(\cH,F)$ is a Fredholm
module over $\cA$. Moreover, if $({\mathcal A}, {\mathcal H}, D)$ is
finite-dimensional, then the Fredholm module $(\cH,F)$ is also
finite-dimensional. Thus, the Chern character $\ch (\cH,F)\in
HP^*(\cA)$ is well defined for any finite-dimensional Lipschitz
regular $\sigma$-spectral triple $({\mathcal A}, {\mathcal H}, D)$.
On the other hand, for any $\sigma$-spectral triple $({\mathcal A},
{\mathcal H}, D)$ such that $D^{-1}\in \cL^{(n,\infty)}$ for some
even $n$ one can define a cyclic cocycle $\Psi_{D,\sigma}$ on $\cA$
by
\[
\Psi_{D,\sigma}(a^0,a^1,\ldots,a^n)=\tr(\gamma d_\sigma a^0 D^{-1}
d_\sigma a^1 \ldots D^{-1} d_\sigma a^n), \quad a^0, a^1, \ldots,
a^n\in \cA,
\]
where $d_\sigma a= D a- \sigma(a)D$ for any $a\in \cA$. If a
$\sigma$-spectral triple $({\mathcal A}, {\mathcal H}, D)$ satisfies
the stronger condition
\[
|D|^{-t}(|D|^ta-\sigma^t(a)|D|^t)\in \cL^{(n,\infty)}, \quad a\in
\cA, \quad t\in \RR,
\]
then the cocycle $\Psi_{D,\sigma}$ defines the same cohomology class
as the Chern character $\ch (\cH,F)\in HP^*(\cA)$.

An analogue of the non-commutative local index theorem in this case
remains an open question. For a certain class of twisted spectral
triples of type III, such a theorem was proved very recently in
\cite{Mosc09}.

\section{Some background material from foliation theory}\label{rev}
\subsection{Foliations: definitions and examples}\label{s:fol}
In this subsection we recall the definition of a foliated manifold
and some notions connected with foliations.

\begin{defn}
(1) An atlas $\cA=\{(U_i,\phi_i)\}$, where $\phi_i: U_i\subset
M\rightarrow \RR^n$, of a smooth manifold $M$ of dimension $n$ is
called an atlas of a foliation of dimension $p$ and codimension $q$
($p\leq n, p+q=n$), if, for any $i$ and $j$ such that $U_i\cap
U_j\not=\varnothing$, the coordinate transformations  $
\phi_{ij}=\phi_i\circ\phi_j^{-1}:\phi_j(U_i\cap U_j)\subset
\RR^p\times\RR^q \to \phi_i(U_i\cap U_j)\subset \RR^p\times\RR^q$
have the form
$$ \phi_{ij}(x,y)=(\alpha_{ij}(x,y), \gamma_{ij}(y)), \quad (x,y)\in
\phi_j(U_i\cap U_j)\subset \RR^p\times\RR^q. $$

(2) Two atlases of a foliation of dimension $p$ are equivalent if
their union is again an atlas of a foliation of dimension $p$.
\medskip\par
(3) A manifold $M$ endowed with an equivalence class $\cF$ of
atlases of a foliation of dimension $p$ is called a manifold
equipped with a foliation of dimension $p$.
\end{defn}

An equivalence class $\cF$ of atlases of a foliation is also called
a complete atlas of a foliation. We will also say that $\cF$ is a
foliation on the manifold $M$.

A pair $(U,\phi)$ belonging to some atlas of a foliation $\cF$ and
also the corresponding map $\phi:U\rightarrow \RR^n$ are called a
foliated chart of $\cF$, and $U$ a foliated coordinate neighborhood.

Let $\phi:U\subset M \rightarrow \RR^n$ be a foliated chart. The
components of the set $\phi^{-1}(\RR^p\times \{y\}), y\in \RR^q$,
are called plaques of $\cF$.

The plagues of $\cF$ taken for all possible foliated charts form a
base of a topology on $M$. This topology is called the leaf topology
on $M$. We will also denote by $\cF$ the set $M$ endowed with the
leaf topology. One can introduce a $p$-dimensional smooth manifold
structure on $\cF$.

The connected components of $\cF$ are called leaves of $\cF$. The
leaves are (one-to-one) immersed $p$-dimensional submanifolds in
$M$. For any $x\in M$ there exists a unique leaf passing through
$x$. We will denote this leaf by $L_x$.

One can give an equivalent definition of a foliation by saying that
there is a foliation $\cF$ of dimension $p$ on an $n$-dimensional
manifold $M$ if $M$ is represented as the union of a family
$\{L_\lambda : \lambda\in \cL\}$ of disjoint connected (one-to-one)
immersed submanifolds of dimension $p$, and $M$ has an atlas
$\cA=\{(U_i,\phi_i)\}$ such that, for any coordinate chart
$(U_i,\phi_i)$ with local coordinates $(x_1,x_2,\ldots,x_n)$ and for
any $\lambda\in \cL $, the connected components of the set
$L_\lambda\cap U_i$ are given by equations of the form
$x_{p+1}=c_{p+1}, \ldots, x_n=c_n$ for some constants $c_{p+1},
\ldots, c_n$.

\begin{ex}
Let $M$ be an $n$-di\-men\-si\-o\-nal smooth manifold, $B$ a
$q$-di\-men\-si\-o\-nal smooth manifold, and $\pi: M\to B$ a
submersion (i.e. the differential $d\pi_x:T_xM\to T_{\pi(x)}B$ is
surjective for any $x\in M$). The connected components of the
pre-images of points of $B$ under the map $\pi$ give a codimension
$q$ foliation of $M$ which is called the foliation determined by the
submersion $\pi$. If, in addition, the pre-images $\pi^{-1}(b), b\in
B,$ are connected, then the foliation is called simple.
\end{ex}

\begin{ex}
If $X$ is a non-singular (that is, non-vanishing) smooth vector
field on a manifold $M$, then its phase curves form a dimension one
foliation.

More generally, let a connected Lie group $G$ act smoothly on a
smooth manifold $M$, and the dimension of the stationary subgroup
$G_x=\{x\in G: gx=x\}$ is independent of $x\in M$. In particular,
one can assume that the action is locally free, which means the
discreteness of the isotropy group $G_x$ for any $x\in M$. Then the
orbits of the action of $G$ define a foliation on $M$.
\end{ex}

\begin{ex}[a linear foliation on the torus]
Consider the vector field $\tilde{X}$ on $\RR^2$ given by
\[
\tilde{X}=\alpha\frac{\partial}{\partial
x}+\beta\frac{\partial}{\partial y}
\]
with constant $\alpha$ and $\beta $. Since $\tilde{X}$ is invariant
by all translations, it defines a vector field $X$ on the
two-dimensional torus ${T}^2={\RR}^2/{\ZZ}^2$. The vector field $X$
defines a foliation $\cF$ on ${T}^2$. The leaves of $\cF$ are the
images of the parallel lines $\tilde{L}=\{(x_0+t\alpha, y_0+t\beta):
t\in\RR\}$ with the slope $\theta=\beta/\alpha $ under the
projection $\RR^2\to T^2$.

In the case when $\theta$ is rational all leaves of $\cF$ are closed
and are the circles, and the foliation $\cF$ is defined by the
fibres of a fibration $T^2\to S^1$. In the case when $\theta$ is
irrational, all leaves of $\cF$ are dense in $T^2$.
\end{ex}

\begin{ex}[homogeneous foliations]
Let $G$ be a Lie group and $H\subset G$ a connected Lie subgroup of
it. The family $\{ gH : g\in G\}$ of right cosets of $H$ forms a
foliation $\cH$ on $G$. If $H$ is a closed subgroup, then $G/H$ is a
manifold and $\cH$ is the foliation with the leaves given by the
fibres of the fibration $\pi : G\to G/H$.

Moreover, suppose that $\Gamma\subset G$ is a discrete subgroup $G$.
Then the set $M=\Gamma\backslash G$ of left cosets of $\Gamma $ is a
manifold of the same dimension as $G$. If $\Gamma $ is cocompact in
$G$, then $M$ is compact. In any case, because $\cH$ is invariant
under left translations, and $\Gamma $ acts from the left, the
foliation $\cH$ is mapped by the map $G\to M=\Gamma\backslash G$ to
a well-defined foliation $\cH_\Gamma $ on $M$, which is often
denoted by $\cF(G,H,\Gamma)$ and called a locally homogeneous
foliation. The leaf of $\cH_\Gamma $ through a point $\Gamma g \in
M$ is diffeomorphic to $H/(g\Gamma g^{-1}\cap H)$.
\end{ex}

\begin{ex}[suspension]\label{ex:suspension}
Let $B$ be a connected manifold and $\tilde{B}$ its universal cover
equipped with the action of the fundamental group $\Gamma=\pi_1(B)$
by deck transformations. Suppose that there is given a homomorphism
$\phi : \Gamma \to \operatorname{Diff}(F)$ of $\Gamma$ to the
diffeomorphism group $\operatorname{Diff}(F)$ of a smooth manifold
$F$. We define a manifold $ M=\tilde{B} \times_\Gamma F$ as the
quotient of the manifold $\tilde{B}\times F$ by the action of
$\Gamma$ given, for any $\gamma\in \Gamma$, by
\[
\gamma (b, f)=(\gamma b, \phi(\gamma)f),\quad (b, f) \in
\tilde{B}\times F.
\]
There is a natural foliation $\cF$ on $M$ whose leaves are the
images of the sets $\tilde{B}\times \{f\}, f\in F,$ under the
projection $\tilde{B}\times F\to M$. If for any $\gamma\in \Gamma$
with $\gamma\not= e$ the diffeomorphism $\phi(\gamma)$ has no fixed
points, then all leaves of $\cF$ are diffeomorphic to $\tilde{B}$.

There is defined a bundle $\pi : M\to B : [(b, f)]\mapsto b\mod
\Gamma $ such that the leaves of $\cF$ are transverse to the fibres
of $\pi$. The bundle $\pi :M\to B$ is often said to be a flat
foliated bundle.
\end{ex}

A foliation $\cF$ defines a subbundle $F=T{\mathcal F}$ of the
tangent bundle $TM$, called the tangent bundle of $\cF$. It consists
of all vectors tangent to the leaves of $\cF$. Denote by
$\cX(M)=C^\infty(M,TM)$ the Lie algebra of all smooth vector fields
on $M$ with respect to the Lie bracket and by
$\cX(\cF)=C^\infty(M,F)$ the subspace of vector fields on $M$
tangent to the leaves of $\cF$ at each point. The subspace
$\cX(\cF)$ is a subalgebra of the Lie algebra $\cX(M)$. Moreover, by
the Frobenius theorem, a subbundle $E$ of the bundle $TM$ is the
tangent bundle of some foliation if and only if it is involutive,
that is, the space of sections of this bundle is a Lie subalgebra of
the Lie algebra $\cX(M)$: for any $X,Y\in C^\infty(M,E)$, we have
$[X,Y]\in C^\infty(M,E)$.

We also introduce the following objects: $\tau=TM/T\cF$ is the
normal bundle of ${\mathcal F}$; $P_\tau: TM\to \tau$ is the natural
projection; $N^*{\mathcal F}=\{\nu\in T^*M:\langle\nu,X\rangle=0\
\mbox{\rm for any}\ X\in F\}$ is the conormal bundle of ${\mathcal
F}$. Usually we will denote by $(x,y)\in I^p\times I^q$ ($I=(0,1)$
is an open interval) the local coordinates in a foliated chart
$\phi: U\to I^p \times I^q$ and by $(x,y,\xi,\eta) \in I^p \times
I^q \times {\RR}^p\times {\RR}^q$ the local coordinates in the
corresponding chart on $T^*M$. Then the subset $N^*{\mathcal F}\cap
\pi^{-1}(U)=U_1$ (here $\pi: T^*M\to M$ is the bundle map) is given
by the equation $\xi=0$. Therefore, $\phi$ defines a natural
coordinate chart $\phi_n: U_1\to I^p \times I^q\times {\RR}^q$ on
$N^*{\mathcal F}$ with coordinates $(x,y,\eta)$.

A $q$-dimensional distribution $Q\subset TM$ such that $TM=F\oplus
Q$ is called a distribution transversal to the foliation, or a
connection on the foliated manifold $(M,\cF)$. Any Riemannian metric
$g$ on $M$ defines a transversal distribution $H$ formed by the
orthogonal complement of $F$ with respect to this metric :
\[
H=F^{\bot}=\{X\in TM : g(X,Y)=0\ \mbox{\rm for any}\ Y\in F\}.
\]

\begin{defn}
A vector field $V$ on a foliated manifold $(M,\cF)$ is called an
infinitesimal transformation of $\cF$ if $[V,X]\in \cX(\cF)$ for any
$X\in \cX(\cF)$.
\end{defn}

The set of infinitesimal transformations of $\cF$ is denoted by
$\cX(M/\cF)$. If $V\in \cX(M/\cF)$ and $T_t:M\to M, t\in \RR$, is
the flow of the vector field $V$, then the diffeomorphisms $T_t$ are
automorphisms of the foliated manifold $(M,\cF)$, that is, they take
each leaf of $\cF$ to a (possibly, different) leaf.

\begin{defn}
A vector field $V$ on a foliated manifold $(M,\cF)$ is said to be
projectable, if its normal component $P_\tau(V)$ is locally the lift
of a vector field on the local base.
\end{defn}

In other words, a vector field $V$ on $M$ is projectable, if, in any
foliated chart with local coordinates $(x,y), x\in \RR^p, y\in
\RR^q,$ it has the form $$
V=\sum_{i=1}^pf^i(x,y)\frac{\partial}{\partial x^i} + \sum_{j=1}^q
g^j(y)\frac{\partial}{\partial y^j}. $$

There is a natural action of the Lie algebra $\cX(\cF)$ on the space
$C^\infty(M,\tau)$. The action of a vector field  $X\in \cX(\cF)$ on
$N\in C^\infty(M,\tau)$ is given by
\[
\theta(X) N=P_\tau[X,\widetilde{N}],
\]
where $\widetilde{N}\in \cX(M)$ is any vector field on $M$ such that
$P_\tau(\widetilde{N})= N$. A vector field $N\in \cX(M)$ is
projectable if and only if its transverse component $P_\tau(N)\in
C^\infty(M,\tau)$ is invariant under the $\cX(\cF)$-action $\theta$.
It is easy to see from this that a vector field on a foliated
manifold is projectable if and only if it is an infinitesimal
transformation of the foliation.

\subsection{Holonomy and transverse structures}\label{s:holonomy}
Let $(M, \cF)$ be a foliated manifold. The holonomy map is a
generalization of the first return map (or Poincar\'e map) for flows
to the case of foliations.

\begin{defn}
A smooth transversal is a compact $q$-dimensional manifold $T$
(possibly disconnected and with boundary) and an embedding
$i:T\rightarrow M$ whose image is everywhere transverse to the
leaves of ${\mathcal F}$: $T_{i(t)}i(T)\oplus T_{i(t)}\cF=T_{i(t)}M$
for any $t\in T$.
\end{defn}

We will identify a transversal $T$ with the image $i(T)\subset M$.

\begin{defn}
A transversal is complete, if it meets every leaf of the foliation.
\end{defn}

Take an arbitrary continuous leafwise path $\gamma$ with initial
point $\gamma(0)=x$ and final point $\gamma(1)=y$. (A path $\gamma
:[0,1]\to M$ is said to be leafwise if its image $\gamma([0,1])$ is
entirely contained in one leaf of the foliation). Let $T_0$ and
$T_1$ be smooth transversals such that $x\in T_0$ and $y\in T_1$.

Choose a partition $t_0=0<t_1<\ldots < t_k=1$ of $[0.1]$ such that
for any $i=1,\ldots, k$ the curve $\gamma([t_{i-1},t_i])$ is
contained in some foliated chart $U_i$. Shrinking the neighborhoods
$U_1$ and $U_2$ if necessary, one can assume that for any plaque
$P_1$ of $U_1$ there is a unique plaque $P_2$ of $U_2$ which meets
$P_1$. Shrinking the neighborhoods $U_1$, $U_2$ and $U_3$ if
necessary, one can assume that for any plaque $P_2$ of $U_2$ there
is a unique plaque $P_3$ of $U_3$ which meets $P_2$ and so on. In
the end we get a family $\{U_1, U_2, \ldots, U_k\}$ of foliated
coordinate neighborhoods which covers the curve $\gamma([0,1])$ and
is such that for any $i=1,\ldots, k$ and for any plaque $P_{i-1}$ of
$U_{i-1}$ there is a unique plaque $P_i$ of $U_i$ which meets
$P_{i-1}$. In particular, we get a one-to-one correspondence between
the plaques of $U_1$ and the plaques of $U_k$.

The smooth transversal $T_0$ determines a parametrization of the
plaques of $U_1$ near $x$. Correspondingly, a smooth transversal
$T_1$ determines a parametrization of the plaques of $U_k$ near $y$.
Taking into account the one-to-one correspondence constructed above
between the plaques of $U_1$ and $U_k$, we get a diffeomorphism
$H_{T_0T_1}(\gamma)$ of some neighborhood of $x$ in $T_0$ to some
neighborhood of $y$ in $T_1$, which is called the holonomy map along
the path $\gamma$.

It is easy to see that the germ of $H_{T_0T_1}(\gamma)$ at $x$ does
not depend on the choice of a partition $t_0=0<t_1<\ldots < t_k=1$
of $[0.1]$ and a family $\{U_1, U_2, \ldots, U_k\}$ of foliated
coordinate neighborhoods. Moreover, the germ of $H_{T_0T_1}(\gamma)$
at $x$ is not changed if we replace $\gamma$ by any other continuous
leafwise path $\gamma_1$ from $x$ to $y$ which is homotopic to
$\gamma$ in the class of continuous leafwise paths from $x$ to $y$.

If $\gamma $ is a closed leafwise path starting and ending at $x$,
and $T$ is a smooth transversal such that $x\in T$, then
$H_{TT}(\gamma)$ is a local diffeomorphism of $T$, which leaves $x$
fixed. The correspondence $\gamma \to H_{TT}(\gamma)$ defines a
group homomorphism $H_T$ from the fundamental group $\pi_1(L_x,x)$
of the leaf $L_x$ to the group $\operatorname{Diff}_x(T)$ of germs
at $x$ of local diffeomorphisms of $T$ which leave $x$ fixed. The
image of $H_T$ is called the holonomy group of the leaf $L_x$ at
$x$. The holonomy group of a leaf $L$ at a point $x\in L$ is
independent (up to an isomorphism) of the choice of a transversal
$T$ and of the choice of $x$. A leaf is said to have trivial
holonomy if its holonomy group is trivial.

\begin{ex}
Let $X$ be a complete non-singular vector field on a manifold $M$ of
dimension $n$, $x_0$ a (for simplicity, isolated) periodic point of
the flow $X_t$ of the given vector field, and $C$ the corresponding
closed phase curve. Let $T$ be an $(n-1)$-dimensional submanifold of
$M$ passing through $x_0$ and transverse to the vector $X(x_0)$:
\[
T_xM=T_xT\oplus \RR X(x_0).
\]
For all $x\in T$ close enough to $x_0$ there is a least $t(x)>0$
such that the corresponding positive semi-trajectory
$\{X_t(x):t>0\}$ of the flow meets $T$: $X_{t(x)}(x)\in T$. Thus, we
get a local diffeomorphism $\phi_T : x\mapsto X_{t(x)}(x)$ of $T$
which is defined in a neighborhood of $x_0$ and takes $x_0$ to
itself. This diffeomorphism is called the first return map (or the
Poincar\'e map) along the curve $C$.

If $\cF$ is the foliation on $M$ given by the trajectories of $X$,
then the holonomy group of the leaf $C$ coincides with $\ZZ$, and
the germ of $\phi_T$ at $x_0$ is a generator of this group.
\end{ex}

For any smooth transversal $T$ and for any $x\in T$ there is a
natural isomorphism between the tangent space $T_xT$ and the normal
space $\tau_x$ to $\cF$. Thus, the normal bundle $\tau$ plays the
role of the tangent bundle of the (germs of) transversals to $\cF$.
For any continuous leafwise path $\gamma$ from a point $x$ to a
point $y$ and for any smooth transversals $T_0$ and $T_1$ with $x\in
T_0$ and $y\in T_1$ the differential of the holonomy map
$H_{T_0T_1}(\gamma)$ at $x$ defines a linear map
$dH_{T_0T_1}(\gamma)_x: \tau_x\to\tau_y$. It is easy to check that
this map is independent of the choice of transversals $T_0$ and
$T_1$. It is called the linear holonomy map and denoted by
$dh_{\gamma}:\tau_x\to\tau_y$. Taking the adjoint of $dh_{\gamma}$,
one obtains a linear map $dh_{\gamma}^*:N^*\cF_y\to N^*\cF_x$.

We now turn to another notion connected with holonomy, the notion of
holonomy pseudogroup. First recall the general definition of a
pseudogroup.

\begin{defn}
A family $\Gamma$ of diffeomorphisms between open subsets of a
manifold $X$ (or, in other words, of local diffeomorphisms of $X$)
is called a pseudogroup on $X$ if the following conditions hold:
\begin{enumerate}
  \item if $\Phi\in\Gamma$, then $\Phi^{-1}\in\Gamma$;
  \item if $\Phi_1:U\to U_1$ and $\Phi_2:U_1\to U_2$ belong to $\Gamma$,
  then $\Phi_2 \circ \Phi_1:U\to U_2$ belongs to $\Gamma$;
  \item if $\Phi :U\to U_1$ belongs to $\Gamma$, then its restriction to any open subset
  $V\subset U$ belongs to $\Gamma$;
  \item if a diffeomorphism $\Phi:U\to U_1$ coincides on some neighborhood of each point in $U$
  with an element of $\Gamma$, then $\Phi\in\Gamma$;
  \item the identity diffeomorphism belongs to $\Gamma$.
\end{enumerate}
\end{defn}

\begin{ex}
The set of all local diffeomorphisms of a manifold $X$ forms a
pseudogroup on $X$. One can also consider pseudogroups consisting of
local diffeomorphisms of a manifold $X$ which preserve some
geometric structure, for instance, the pseudogroup of local
isometries of a Riemannian manifold, and so on.
\end{ex}

\begin{defn}
Let $(M,\cF)$ be a smooth foliated manifold and $X$ the disjoint
union of all smooth transversals to $\cF$. The holonomy pseudogroup
of the foliation $\cF$ is defined to be the pseudogroup $\Gamma$ of
all local diffeomorphisms of $X$ whose germ at any point coincides
with the germ of the holonomy map along some leafwise path.
\end{defn}

\begin{defn}
Let $(M,\cF)$ be a smooth foliated manifold and $T$ a smooth
transversal. The holonomy pseudogroup induced by the foliation $\cF$
on $T$ is defined to be the pseudogroup $\Gamma_T$ of all local
diffeomorphisms of $T$ whose germ at each point coincides with the
germ of the holonomy map along some leafwise path.
\end{defn}

There is a special class of smooth transversals given by good covers
of the manifold $M$.

\begin{defn}
A foliated chart $\phi: U\subset M \to \RR^p\times \RR^q$ is said to
be regular if it admits an extension to a foliated chart
$\overline{\phi}: V\to \RR^p\times \RR^q$ such that $\overline{U}
\subset V$.
\end{defn}

\begin{defn}
A cover of a manifold $M$ by foliated neighborhoods $\{U_i\}$ is
said to be good, if:
\begin{enumerate}
\item any chart $(U_i,\phi_i)$ is a regular foliated chart;
\item if $\overline{U}_i\cap \overline{U}_j\not=\varnothing$, then
${U}_i\cap {U}_j\not=\varnothing$ and the set ${U}_i\cap {U}_j$ is
connected, and the same is true for the corresponding foliated
neighborhoods $V_i$;
\item each plaque of $V_i$ meets at most one plaque of $V_j$; a plaque
of $U_i$ meets a plaque of $U_j$ if and only if the intersection of
the corresponding plaques of $V_i$ and $V_j$ is non-empty.
\end{enumerate}
\end{defn}

Good covers always exist.

Let $\cU=\{U_i\}$ be a good cover for the foliation ${\mathcal F}$,
$\phi_i:U_i\overset{\cong}{\to} I^p \times I^q$. For any $i$ put $$
T_i=\phi^{-1}_i(\{0\}\times I^q). $$ Then $T_i$ is a transversal,
and $T=\bigcup T_i$ is a complete transversal. For $y\in T_i$ denote
by $P_i(y)$ the plaque of $U_i$ passing through $y$. For any pair of
indices $i$ and $j$ such that ${U}_i\cap {U}_j\not=\varnothing$ we
define $$ T_{ij}=\{y\in T_i: P_i(y)\cap U_j\not=\varnothing\}. $$
There is defined a transition function $f_{ij}:T_{ij}\rightarrow
T_{ji}$ given for $y\in T_{ij}$ by the formula $f_{ij}(y)=y_1,$
where $y_1\in T_{ji}$ corresponds to the unique plaque $P_j(y_1)$
such that $P_i(y)\cap P_j(y_1)\not=\varnothing.$ The holonomy
pseudogroup $\Gamma_T$ induced by $\cF$ on $T$ coincides with the
pseudogroup generated by the maps $f_{ij}$.

\begin{defn}[cf., for instance, \cite{Haefliger84}]
A transverse structure on a foliation $\cF$ is a structure on a
complete transversal $T$ that is invariant under the action of the
holonomy pseudogroup $\Gamma_T$.
\end{defn}

Using the notion of transverse structure, one can distinguish
classes of foliations with specific transverse properties. For
instance, if a complete transversal $T$ is equipped with a
Riemannian metric, and the holonomy pseudogroup $\Gamma_T$ consists
of local isometries of this Riemannian metric, then we get the class
of Riemannian foliations (see Subsection~\ref{riemann}). Similarly,
if a complete transversal $T$ is equipped with a symplectic
structure, and the holonomy pseudogroup $\Gamma_T$ consists of local
diffeomorphisms preserving this symplectic structure, then we get
the class of symplectic foliations. One can also consider Kaehler
foliations, measurable foliations and so on.

In \cite{Haefliger80} Haefliger defined cohomology groups associated
with the transverse structure of a foliation. As above, let
$\cU=\{U_i\}$ be a good cover for a foliation ${\mathcal F}$,
$T=\cup T_i$ the corresponding complete transversal, and $\Gamma_T$
the holonomy pseudogroup induced by the foliation $\cF$ on $T$.
Denote by $\Omega^k_c(M/\cF)$ the quotient of the space
$\Omega^k_c(T)$ of smooth compactly supported differential $k$-forms
on $T$ by the subspace spanned by forms of the form
$\alpha-h^*\alpha$ with $h\in \Gamma_T$ and the support of
$\alpha\in \Omega^k_c(T)$ contained in the image of $h$. Let us
endow the space $\Omega^k_c(M/\cF)$ with the quotient topology
induced by the usual $C^\infty$-topology on $\Omega^k_c(T)$. In
general, $\Omega^k_c(M/\cF)$ is not a Hausdorff topological space.
The de Rham differential $d_T : \Omega^k_c(T)\to \Omega^{k+1}_c(T)$
induces a continuous differential $d_H : \Omega^k_c(M/\cF)\to
\Omega^{k+1}_c(M/\cF)$. It should be noted that both
$\Omega^k_c(M/\cF)$ and $d_H$ are independent of the choice of a
good cover $\cU=\{U_i\}$. The complex $(\Omega_c(M/\cF),d_H)$ and
its cohomology $H^*_c(M/\cF)$ are called respectively the Haefliger
forms and the Haefliger cohomology of the foliation $\cF$.

If the tangent bundle $T\cF$ is oriented, then there is defined a
continuous open surjective linear map, called integration along the
leaves,
\[
\int_\cF : \Omega^{k+p}(M)\to \Omega^k_c(M/\cF),
\]
which satisfies the condition $d_H \circ \int_\cF = \int_\cF\circ
d$, where $d$ is the de Rham differential on $\Omega(M)$. Thus, one
has the induced map
\[
\int_\cF : H^{k+p}(M)\to H^k_c(M/\cF).
\]

Before turning to a discussion of an analogue of the notion of
measure on the leaf space of a foliation, we recall some basic facts
concerning to densities and integration of densities.

\begin{defn}
Let $L$ be an $n$-dimensional linear space and $\cB(L)$ the set of
bases in $L$. An $\alpha$-density on $L$ ($\alpha\in \RR$) is a
function $\rho : \cB(L)\to \CC$ such that for any $A=(A_{ij})\in
\mathop{GL}(n,\CC)$ and $e=(e_1, e_2,\ldots, e_n)\in \cB(L)$
\[
\rho(e\cdot A)=|\det A|^\alpha \rho(e),
\]
where $(e\cdot A)_i=\sum_{j=1}^ne_jA_{ji}, i=1,2,\ldots,n$.
\end{defn}

We will denote by $|L|^\alpha$ the space of all $\alpha$-densities
on $L$. For any vector bundle $V$ on $M$ denote by $|V|^{\alpha}$
the associated bundle of $\alpha$-densities, $|V|=|V|^1$.

For any smooth, compactly supported density $\rho$ on a smooth
manifold $M$ there is a well-defined integral $\int_M\rho$,
regardless of whether $M$ is orientable or not. This fact enables
one to define a Hilbert space $L^2(M)$ canonically associated with
$M$ and consisting of square integrable half-densities on $M$. The
diffeomorphism group of $M$ acts on the space $L^2(M)$ by unitary
transformations.

\begin{defn}
A (Borel) transversal to a foliation ${\mathcal F}$ is a Borel
subset of $M$ which meets each leaf of the foliation in an at most
countable set.
\end{defn}

\begin{defn}
A transverse measure is a countably additive Radon measure $\Lambda$
defined on the set of all transversals to the foliation.
\end{defn}

\begin{defn}
A transverse measure $\Lambda$ is called holonomy invariant if for
any transversals $B_1$ and $B_2$ and for any bijective Borel map
$\phi :B_1\to B_2$ such that for any $x\in B_1$ the point $\phi(x)$
belongs to the leaf through the point $x$ we have:
$\Lambda(B_1)=\Lambda(B_2)$.
\end{defn}

\begin{ex}
A transverse density is defined to be any section of the bundle
$|\tau|$. Since for any smooth transversal $T$ there is a canonical
isomorphism $T_xT\cong \tau_x$, a continuous positive density
$\rho\in C(M,|\tau |)$ determines a continuous positive density on
$T$, which in its turn determines a transverse measure. This
transverse measure is holonomy invariant if and only if $\rho$ is
invariant under the linear holonomy action.
\end{ex}

\begin{ex}
Any compact leaf $L$ of the foliation $\cF$ determines a holonomy
invariant transverse measure $\Lambda$. For any transversal $T$ and
for any set $A\subset T$ its measure $\Lambda(A)$ equals the number
of elements in $A\cap L$.
\end{ex}

\begin{ex}
Let $\cF$ be the horizontal foliation of a flat foliated bundle
$M=\tilde{B}\times_\Gamma F$ (see Example~\ref{ex:suspension}). Any
measure on $F$ invariant under the action of the group $\Gamma$
determines a holonomy invariant measure for $\cF$.
\end{ex}

Let $\alpha\in C^\infty(M,|T\cF|)$ be a smooth positive leafwise
density on $M$. Starting from a transverse measure $\Lambda$ and the
density $\alpha$ one can construct a Borel measure $\mu$ on $M$ in
the following way. Take a good cover $\{U_i\}$ of $M$ by foliated
coordinate neighborhoods with the corresponding coordinate maps
$\phi_i: U_i\to I^p \times I^q$ and a partition of unity
$\{\psi_i\}$ subordinate to this cover. Consider the corresponding
complete transversal $T=\bigcup_i T_i$, where
$T_i=\phi^{-1}_i(\{0\}\times I^q).$ In any foliated chart
$(U_i,\phi_i)$ the transverse measure $\Lambda$ defines a measure
$\Lambda_i$ on $T_i$, and the smooth positive leafwise density
$\alpha$ defines a family $\{\alpha_{i,y} : y\in T_i\}$, where
$\{\alpha_{i,y}\}$ is a smooth positive density on the plaque
$P_i(y)$. Observe that $\Lambda$ is holonomy invariant if and only
if, for any pair of indices $i$ and $j$ such that ${U}_i\cap
{U}_j\not=\varnothing$, we have the relation
$f_{ij}(\Lambda_i)=\Lambda_j$.

For any $u\in C^\infty_c(M)$, put \begin{equation}\label{def:mu}
 \int_Mu(m)d\mu(m)=\sum_i \int_{T_i}\int_{P_i(y)} \psi_i(x,y)u(x,y)
\alpha_{i,y}(x) d\Lambda_i(y).
\end{equation}
One can show that this formula defines a measure
$\mu$ on $M$, which is independent of the choice of a cover
$\{U_i\}$ and a partition of unity $\{\psi_i\}$.

A measure $\mu$ on $M$ will be said to be holonomy invariant if it
is obtained from a holonomy invariant transverse measure $\Lambda$
by means of the above construction with some choice of a smooth
positive leafwise density $\alpha$.

If in the above construction we take the restrictions to the leaves
of an arbitrary differential $p$-form $\omega$ on $M$ instead of the
leafwise density $u\cdot\alpha$, we obtain a well-defined functional
$C$ on $C^\infty_c(M,\Lambda^pT^*M)$ called the Ruelle-Sullivan
current corresponding to $\Lambda$:
\[
\langle C,\omega\rangle =\sum_i \int_{T_i}\int_{P_i(y)}
\psi_i(x,y)\omega_{i,y}(x) d\Lambda_i(y), \quad \omega\in
C^\infty_c(M,\Lambda^pT^*M),
\]
where $\omega_{i,y}$ is the restriction of $\omega$ to the plaque
$P_i(y), y\in T_i$.

A transverse measure $\Lambda$ is holonomy invariant if and only if
the corresponding Ruelle-Sullivan current $C$ is closed:
\[
\langle C,d\sigma \rangle =0, \quad \sigma \in
C^\infty_c(M,\Lambda^{p-1}T^*M).
\]

\begin{ex}
Suppose that a transverse measure $\Lambda$ is given by a smooth
positive transverse density $\rho\in C^\infty(M,|\tau |)$. Take a
positive leafwise density $\alpha\in C^\infty(M,|T\cF|)$. Then the
corresponding measure $\mu$ on $M$ is given by the smooth positive
density  $\alpha\otimes \rho\in C^\infty(M,|TM|)$, which corresponds
to $\alpha$ and $\rho$ by the canonical isomorphism $|TM|\cong
|T\cF|\otimes |\tau|$ given by the short exact sequence $0\to
T\cF\to TM\to\tau\to 0$.
\end{ex}

\begin{ex}
Suppose that a holonomy invariant transverse measure $\Lambda$ is
given by a compact leaf $L$ of the foliation $\cF$, and $\alpha\in
C^\infty(M,|T\cF|)$ is a smooth positive leafwise density on $M$.
Then the corresponding measure $\mu$ on $M$ is a $\delta$-measure
along $L$:
\[
\int_M f(x)\,d\mu(x)=\int_Lf(x)\,\alpha(x), \quad f\in C_c(M).
\]
\end{ex}

\begin{ex}
Suppose that a foliation $\cF$ is given by the orbits of a locally
free action of a Lie group $H$ on a compact manifold $M$ and a
smooth leafwise density $\alpha $ is given by a fixed Haar measure
$dh$ on $H$. Then the corresponding measure $\mu $ on $M$ is
holonomy invariant if and only if it is invariant under the action
of $H$.
\end{ex}

\subsection{Transverse Riemannian geometry}\label{riemann}
An infinitesimal expression of the holonomy on a foliated manifold
is the canonical flat connection
\[
\stackrel{\circ}{\nabla} : \cX(\cF)\times C^\infty(M,\tau)\to
C^\infty(M,\tau)
\]
defined in the normal bundle $\tau$ along the leaves of $\cF$ (the
Bott connection) \cite{Bott71}. It is given by
\begin{equation}\label{e:Bott}
{\stackrel{\circ}\nabla}_X N=\theta(X)N=
P_\tau[X,\widetilde{N}],\quad X\in \cX(\cF),\quad N\in
C^\infty(M,\tau),
\end{equation}
where $\widetilde{N}\in C^\infty(M, TM)$ is any vector field on $M$
such that $P_\tau(\widetilde{N})= N$. Thus, the restriction of
$\tau$ to any leaf of $\cF$ is a flat vector bundle. The parallel
transport defined by ${\stackrel{\circ}\nabla}$ in $\tau$ along any
leafwise path $\gamma:x\to y$ coincides with the linear holonomy map
$dh_{\gamma}:\tau_x\to\tau_y$.

\begin{defn}
A connection $\nabla : \cX(M)\times C^\infty(M,\tau)\to
C^\infty(M,\tau)$ in the normal bundle $\tau$ is said to be adapted
if its restriction to $\cX(\cF)$ coincides with the Bott connection
${\stackrel{\circ}\nabla}$.
\end{defn}

One can construct an adapted connection starting with an arbitrary
Riemannian metric $g_M$ on $M$. Denote by $\nabla^g$ the Levi-Civita
connection defined by $g_M$. An adapted connection $\nabla$ is given
by
\begin{equation}\label{e:adapt}
\begin{aligned}
\nabla_XN &=P_\tau[X,\widetilde{N}],\quad X\in \cX(\cF),\quad N\in
C^\infty(M,\tau)\\ \nabla_XN&=P_\tau\nabla^g_X\widetilde{N},\quad
X\in C^\infty(M,F^{\bot}),\quad N\in C^\infty(M,\tau),
\end{aligned}
\end{equation}
where $\widetilde{N}\in C^\infty(M,TM)$ is any vector field such
that $P_\tau(\widetilde{N})=N$. One can show that the adapted
connection $\nabla$ described above is torsion-free.

\begin{defn}
An adapted connection $\nabla$ in the normal bundle $\tau$ is said
to be holonomy invariant if for any $X\in \cX(\cF)$, $Y\in \cX(M)$,
and $N\in C^\infty(M,\tau)$ we have
\[
(\theta(X)\nabla)_YN:=\theta(X)[\nabla_YN]-\nabla_{\theta(X)Y}N-
\nabla_Y[\theta_YN]=0.
\]

A holonomy invariant adapted connection in $\tau$ is called a basic
(or projectable) connection.
\end{defn}

A fundamental property of basic connections is the fact that their
curvature $R_\nabla$ is a basic form, that is, $i_XR_\nabla=0,
\theta(X)R_\nabla=0$ for any $X\in \cX(\cF)$. There are topological
obstructions for the existence of basic connections for an arbitrary
foliation.

\begin{defn}
A foliation $\cF$ on a manifold $M$ is said to be Riemannian if it
has a transverse Riemannian structure. In other words, $\cF$ is
Riemannian if there is a cover $\{U_i\}$ of $M$ by foliated
coordinate charts, $\phi_i:U_i\to I^p\times I^q$, and Riemannian
metrics $g^{(i)}(y)= \sum_{\alpha\beta} g^{(i)}_{\alpha\beta}(y)
dy^\alpha dy^\beta$ defined on the local bases $I^q$ of $\cF$ such
that, for any coordinate transformation
$$ \phi_{ij}(x,y)=(\alpha_{ij}(x,y), \gamma_{ij}(y)), \quad
(x,y)\in \phi_j(U_i\cap U_j), $$ the map $\gamma_{ij}$ preserves the
metric on $I^q$, $\gamma_{ij}^*(g^{(j)})=g^{(i)}$.
\end{defn}

\begin{theorem}\label{rev:blike}
A foliation $\mathcal F$ is Riemannian if and only if there exists a
Riemannian metric $g_M$ on $M$ such that the induced metric $g_\tau$
on the normal bundle $\tau$ is holonomy invariant: for any $X\in
\cX(\cF)$ and for any $U, V\in C^\infty(M,\tau)$,
\[
{\stackrel{\circ}\nabla}_X g_\tau(U,V):=X[g_\tau(U,V)]-
g_\tau({\stackrel{\circ}\nabla}_XU, V) -g_\tau(U,
{\stackrel{\circ}\nabla}_XV)=0.
\]
\end{theorem}

\begin{defn}
Any Riemannian metric on $M$ satisfying the conditions of
Theorem~\ref{rev:blike} is called bundle-like.
\end{defn}

For a Riemannian foliation $\cF$ and a bundle-like metric $g_M$ the
adapted connection $\nabla$ on the normal bundle $\tau$ given by
(\ref{e:adapt}) is a Riemannian connection: for any $Y\in \cX(M)$
and $U, V\in C^\infty(M,\tau)$,
\[
Y[g_\tau(U,V)]=g_\tau(\nabla_YU,V)+g_\tau(U,\nabla_YB).
\]
One can prove that a (torsion-free) Riemannian connection on the
normal bundle $\tau$ of a Riemannian foliation $\cF$ is unique. It
is uniquely determined by the transverse metric $g_\tau$ and is
called the transverse Levi-Civita connection for $\cF$. Thus, the
transverse Levi-Civita connection is an adapted connection.
Moreover, it turns out that the transverse Levi-Civita connection is
holonomy invariant, and therefore is a basic connection. In
particular, this proves the existence of a basic connection for any
Riemannian foliation.

The existence of a bundle-like metric on a foliated manifold imposes
strong restrictions on geometry of the foliation. There are
structure theorems for Riemannian foliations obtained by Molino.
Using these structure theorems, one can reduce many questions
concerning Riemannian foliations to the case of Lie foliations, that
is, of foliations with transverse structure modelled by a
finite-dimensional Lie group  (see Example~\ref{ex:lie}).

\begin{ex}
Any foliation defined by a submersion $\pi:M\to B$ is Riemannian.
\end{ex}

\begin{ex}
The orbits of a locally free isometric action of a Lie group on a
Riemannian manifold define a Riemannian foliation. On the other
hand, flows whose orbits form a Riemannian foliation are called
Riemannian flows. There are examples of Riemannian flows which are
not isometric (see, for instance, \cite{Carriere84}).
\end{ex}

\begin{ex}\label{ex:lie}
Let $M$ be a smooth manifold, $\mathfrak g$ a real
finite-dimensional Lie algebra, and $\omega$ a 1-form on $M$ with
values in $\mathfrak g$ satisfying the conditions:
\begin{enumerate}
  \item the map $\omega_x:T_xM\to \mathfrak g$ is surjective for any $x\in M$;
  \item $d\omega+\frac{1}{2}[\omega,\omega]=0$.
\end{enumerate}
The distribution $F_x=\ker \omega_x$ is integrable and, hence,
defines a foliation of codimension $q=\dim \mathfrak g$ on $M$. Such
a foliation is called a Lie $\mathfrak g$-foliation. Any Lie
foliation is Riemannian.
\end{ex}

\begin{ex}
The horizontal foliation $\cF$ on a flat foliated bundle $M$ (see
Example~\ref{ex:suspension}) is Riemannian if and only if for any
$\gamma\in\Gamma$ a diffeomorphism $\phi(\gamma)$ preserves some
Riemannian metric on $F$.
\end{ex}

Let $\cF$ be a transversely oriented Riemannian foliation and $g$ a
bundle-like Riemannian metric. The induced metric on $\tau$ yields
the transverse volume form $v_\tau\in C^\infty(M, \Lambda^q
\tau^*)=C^\infty(M, \Lambda^q N^*\cF)$ which is holonomy invariant
and, therefore, gives rise to a holonomy invariant transverse
measure on ${\mathcal F}$.

\section{Non-commutative topology of foliations}
In this section we will describe the non-commutative algebras
associated with the leaf space of a foliation. First, we will define
an algebra consisting of very nice functions, on which all basic
operations of analysis are defined, then, depending on a problem in
question, we will complete this algebra and obtain an analogue of
the algebra of measurable, continuous, or smooth functions. The role
of a ``nice'' algebra is played by the algebra $C^\infty_c(G)$ of
smooth compactly supported functions on the holonomy groupoid $G$ of
the foliation. Therefore, we start with the notion of the holonomy
groupoid of a foliation.

\subsection{Holonomy groupoid}\label{s:groupoid}
A foliation ${\cF}$ defines an equivalence relation ${\mathcal
R}\subset M\times M$ on $M$: $(x, y)\in {\mathcal R}$ if and only if
$x$ and $y$ lie on the same leaf of the foliation ${\cF}$.
Generally, ${\mathcal R}$ is not a smooth manifold, but one can
resolve its singularity, constructing a smooth manifold $G$, called
the holonomy groupoid or the graph of the foliation, which coincides
``almost everywhere'' with $\cR$ and which can be used in many cases
as a substitute for $\cR$. The idea of the holonomy groupoid
appeared in papers of Ehresmann, Reeb, and Thom and was completely
realized by Winkelnkemper in \cite{Winkeln}. First of all, we give
the general definition of a groupoid

\begin{defn}
We say that a set $G$ carries the structure of groupoid with a set
of units $G^{(0)}$ if there are maps
\begin{enumerate}
\item $\Delta : G^{(0)}\rightarrow G$ (the diagonal map or the unit map);
\item an involution $i:G\rightarrow G$ called the inversion and
written as $i(\gamma)=\gamma^{-1}$;
\item the range map $r:G\rightarrow G^{(0)}$ and the source map $s:G\rightarrow G^{(0)}$;
\item an associative multiplication $m: (\gamma,\gamma')\rightarrow
\gamma\gamma'$ defined on the set
\[
G^{(2)}=\{(\gamma,\gamma')\in G\times G : r(\gamma')=s(\gamma)\},
\]
\end{enumerate}
satisfying the conditions
\begin{itemize}
\item[(i)] $r(\Delta(x))=s(\Delta(x))=x$ and $\gamma\Delta(s(\gamma))=\gamma$,
$\Delta(r(\gamma))\gamma=\gamma$;
\item[(ii)]~$r(\gamma^{-1})=s(\gamma)$ and $\gamma\gamma^{-1}=\Delta(r(\gamma))$.
\end{itemize}
\end{defn}

Alternatively, one can define a groupoid as a small category, in
which each morphism is an isomorphism. In particular, tt is
convenient to represent an element $\gamma\in G$ as an arrow $\gamma
:x \to y$, where $x=s(\gamma)$ and $y=r(\gamma)$. We will use
standard notation (for $x,y\in G^{(0)}$):
\begin{gather*}
G^x=\{\gamma\in G:r(\gamma)=x\} =r^{-1}(x),\quad G_x=\{\gamma\in
G:s(\gamma)=x\} =s^{-1}(x),\\
G^x_y=\{\gamma\in G : s(\gamma)=x, r(\gamma)=y\}.
\end{gather*}

\begin{defn}
A groupoid $G$ is said to be smooth (or a Lie groupoid), if
$G^{(0)}$, $G$ and $G^{(2)}$ are smooth manifolds, $r$, $s$, $i$ and
$m$ are smooth maps, $r$ and $s$ are submersions, and $\Delta$ is an
immersion.
\end{defn}

\begin{ex}[trivial groupoid] \label{ex:trivial}
Let $X$ be an arbitrary set, put $G=X$, $G^{(0)}=X$, and let the
maps $s$ and $r$ be the identity maps (in other words, each element
$x\in G^{(0)}= X$ is identified with the unique element $\gamma :
x\to x$).
\end{ex}

\begin{ex}[equivalence relations] \label{ex:equiv}
 Any equivalence relation $R\subset X\times
X$ defines a groupoid, if one puts $G^{(0)}=X$, $G=R$, and lets the
maps $s:R\to X$ and $r:R\to X$ be given by $s(x,y)=y$, $r(x,y)=x$.
Thus, pairs $(x_1,y_1)$ and $(x_2,y_2)$ can be multiplied if and
only if $y_1=x_2$, and $(x_1,y_1)(x_2,y_2) = (x_1,y_2)$. Moreover,
one has
\begin{gather*}
\Delta(x)=(x,x), \quad x\in X,\\
(x,y)^{-1}=(y,x),\quad (x,y)\in R.
\end{gather*}

In the particular case when $R=X\times X$, we obtain a so called
principal groupoid.
\end{ex}

\begin{ex}[Lie groups] \label{ex:Lie}
 A Lie group $H$ defines a smooth groupoid as follows:
$G=H$, $G^{(0)}$ consists of a single point, and the maps $i$ and
$m$ are defined by the group operations in $H$.
\end{ex}

\begin{ex}[group actions] \label{ex:actions}
 Let a Lie group $H$ act smoothly from the left on a smooth
manifold $X$. The crossed product groupoid $X\rtimes H$ is defined
as follows: $G^{(0)}=X$, $G=X\times H$. The maps $s: X\times H\to X$
and $r:X\times H\to X$ have the form $s(x,h)=h^{-1}x$, $r(x,h)=x$.
Thus, pairs $(x_1,h_1)$ and $(x_2,h_2)$ can be multiplied if and
only if $x_2=h^{-1}_1x_1$, and $(x_1,h_1)(x_2,h_2) = (x_1,h_1h_2)$.
Moreover, one has
\begin{gather*}
\Delta(x)=(x,e), \quad x\in X.\\
(x,h)^{-1}=(h^{-1}x,h^{-1}),\quad x\in X,\quad h\in H.
\end{gather*}
\end{ex}

\begin{ex}[the fundamental groupoid] \label{ex:fund}
 Let $X$ be a topological space,
$G=\Pi(X)$ the set of homotopy classes of paths in $X$ with all
possible endpoints. More precisely, if $\gamma : [0,1] \to X$ is a
path from $x=\gamma(0)$ to $y=\gamma(1)$, then we denote by
$[\gamma]$ the homotopy class of $\gamma$ with fixed $x$ and $y$.
Define the groupoid $\Pi(X)$ as the set of triples $(x, [\gamma],
y)$, where $x,y\in X$, $\gamma $ is a path with the initial point
$x=\gamma(0)$ and the final point $y=\gamma(1)$ and with the
multiplication given by the product of paths. The groupoid $\Pi(X)$
is called the fundamental groupoid of $X$.
\end{ex}

 \begin{ex}[the Haefliger groupoid $\Gamma_n$ \cite{Haefliger72,Haefliger:CIME}] \label{ex:Heafliger}
Let $M$ be a smooth manifold. The groupoid $\Gamma_M$ consists of
germs of local diffeomorphisms of $M$ at various points of $M$.
Thus, $(\Gamma_M)^{(0)}=M$. If $\gamma\in \Gamma_M$ is the germ at
$x\in M$ of a diffeomorphism $f$ from some neighborhood $U$ of $x$
onto the open set $f(U)$, then $s(\gamma)=x$, $r(\gamma)=f(x)$. The
multiplication in $\Gamma_M$ is given by the composition of maps. If
$M=\RR^n$, then the groupoid $\Gamma_M$ is denoted by $\Gamma_n$.
\end{ex}

The holonomy groupoid $G=G(M,{\mathcal F})$ of a foliated manifold
$(M,{\mathcal F})$ is defined in the following way. Let $\sim_h$ be
the equivalence relation on the set of continuous leafwise paths
$\gamma:[0,1]\rightarrow M$, specifying that $\gamma_1\sim_h
\gamma_2$ if $\gamma_1$ and $\gamma_2$ have the same initial and
final points and the same holonomy maps: $h_{\gamma_1} =
h_{\gamma_2}$. The holonomy groupoid $G$ is the set of
$\sim_h$-equivalence classes of leafwise paths. The set of units
$G^{(0)}$ is the manifold $M$. The multiplication in $G$ is given by
the product of paths. The corresponding range and source maps
$s,r:G\rightarrow M$ are given by $s(\gamma)=\gamma(0)$ and
$r(\gamma)=\gamma(1)$. Finally, the diagonal map
$\Delta:M\rightarrow G$ takes any $x\in M$ to the element in $G$
given by the constant path $\gamma(t)=x, t\in [0,1]$. To simplify
the notation we will identify $x\in M$ with $\Delta(x)\in G$.

For any $x\in M$ the map $s$ maps $G^x$ onto the leaf $L_x$ through
$x$. The group $G^x_x$ coincides with the holonomy group of $L_x$.
The map $s:G^x\rightarrow L_x$ is the regular covering with the
covering group $G^x_x$, called the holonomy covering.

The holonomy groupoid $G$ has the structure of a smooth (generally,
non-Hausdorff and non-paracompact) manifold of dimension $2p+q$. We
recall the construction of an atlas on $G$ \cite{Co79}.

Suppose that $\phi: U\to I^p\times I^q$ and $\phi': U'\to I^p\times
I^q$ are two foliated charts, and let $\pi=pr_{nq}\circ\phi: U\to
\RR^q$ and $\pi'=pr_{nq}\circ\phi': U'\to \RR^q$ be the
corresponding distinguished maps. The charts $\phi$ and $\phi'$ are
said to be compatible if for any $m\in U$ and $m'\in U'$ with
$\pi(m)=\pi'(m')$ there is a leafwise path $\gamma$ from $m$ to $m'$
such that the corresponding holonomy map $h_{\gamma}$ takes the germ
$\pi_m$ of $\pi$ at $m$ to the germ $\pi'_{m'}$ of $\pi'$ at $m'$.

For any pair of compatible foliated charts $\phi$ and $\phi'$ denote
by $W(\phi,\phi')$ the subset in $G$ consisting of all $\gamma\in G$
with $s(\gamma)=m=\phi^{-1}(x,y)\in U$ and
$r(\gamma)=m'={\phi'}^{-1}(x',y)\in U'$ such that the corresponding
holonomy map $h_{\gamma}$ takes the germ $\pi_m$ of the map
$\pi=pr_{nq}\circ\phi$ at $m$ to the germ $\pi'_{m'}$ of the map
$\pi'=pr_{nq}\circ\phi'$ at $m'$. There is a coordinate map
\begin{equation}
\label{rev:wchart} \Gamma:W(\phi,\phi')\to I^p\times I^p\times I^q,
\end{equation}
which takes each element $\gamma\in W(\phi,\phi')$ such that
$s(\gamma)=m=\varkappa^{-1}(x,y)$, $r(\gamma)=m'={\phi'}^{-1}(x',y)$
and $h_{\gamma}\pi_m=\pi'_{m'}$ to the triple $(x,x',y)\in I^p\times
I^p\times I^q$. As shown in \cite{Co79}, the coordinate
neighborhoods $W(\phi,\phi')$ form an atlas of a
$(2p+q)$-dimensional manifold (generally, non-Hausdorff and
non-paracompact) on $G$. Moreover, the groupoid $G$ is a smooth
groupoid.

Non-Hausdorffness of the holonomy groupoid is related with the
phenomenon of one-sided holonomy. The simplest example of a
foliation with non-Hausdorff holonomy groupoid is given by the
trajectories of a non-singular vector field on the plane, having a
one-sided limit cycle. As shown in \cite{Winkeln}, the holonomy
groupoid is Hausdorff if and only if the holonomy maps
$H_{T_0T_1}(\gamma_1)$ and $H_{T_0T_1}(\gamma_2)$ along any leafwise
paths $\gamma_1$ and $\gamma_2$ from $x$ to $y$ and given by smooth
transversals $T_0$ and $T_1$ passing through $x$ and $y$,
respectively, coincide, if they coincide on some open subset
$U\subset T_0$ such that $x\in \bar{U}$. In particular, the holonomy
groupoid is Hausdorff, if the holonomy is trivial or real analytic.
Moreover, the holonomy groupoid of a Riemannian foliation is
Hausdorff. In the following we will always assume that $G$ is a
Hausdorff manifold.

\begin{ex}
If $\cF$ is a simple foliation defined by a submersion $\pi :M\to
B$, then its holonomy groupoid $G$ consists of all $(x,y)\in M\times
M$ such that $\pi(x)=\pi(y)$, and, moreover, $G^{(0)}=M$, and the
maps $s:G\to M$ and $r:G\to M$ are given by $s(x,y)=y$, $r(x,y)=x$.
\end{ex}

\begin{ex}
If a foliation $\cF$ is given by the orbits of a free smooth action
of a connected Lie group $H$ on a manifold $M$, then its holonomy
groupoid coincides with the crossed product groupoid $M\rtimes H$.
\end{ex}

\begin{ex}\label{ex:susp}
Consider the horizontal foliation $\cF$ on a flat foliated bundle
$M$ (see Example~\ref{ex:suspension}). Suppose that the following
condition holds: if, for some element $g\in\Gamma$, there exists an
open set $U$ such that $xg=x$ for any $x\in U$, then $g$ is the
identity element of the group $\Gamma$. Under this condition, the
holonomy groupoid $G$ of $\cF$ is isomorphic to the orbit space of
the action of $\Gamma$ given by $(b_1,b_2,f)g=(b_1g,b_2g,\phi(g)f)$
on the manifold $\tilde{B}\times \tilde{B}\times F$, where
$(b_1,b_2,f)\in \tilde{B}\times \tilde{B}\times F$, $g\in\Gamma$:
\[
G\cong (\tilde{B}\times \tilde{B}\times F)/\Gamma.
\]
Denote by $[b_1,b_2,f]$ the equivalence class of an element
$(b_1,b_2,f)\in \tilde{B}\times \tilde{B}\times F$ in
$(\tilde{B}\times \tilde{B}\times F)/\Gamma$. Then the source and
the range maps in the groupoid $G$ are given by
\[
r([b_1,b_2,f])=[b_1,f], \quad s([b_1,b_2,f])=[b_2,f].
\]
Elements $[b_1,b_2,f]$ and $[b^\prime_1,b^\prime_2,f^\prime]$ can be
multiplied if and only if there exists an element $g\in\Gamma$ such
that $b_2=b^\prime_1g$, $f=\phi(g)f^\prime$. In this case,
\[
[b_1,b_2,f][b^\prime_1,b^\prime_2,f^\prime]=[b_1g^{-1},b^\prime_2,f^\prime].
\]
\end{ex}

In addition to the holonomy groupoid, there are another groupoids,
which can be associated with the foliation. First of all, it is the
groupoid given by the equivalence relation on $M$ for which points
$x$ and $y$ are equivalent if they lie on the same leaf of the
foliation (the coarse groupoid). As mentioned above, this groupoid
is not smooth. One can consider the fundamental groupoid of the
foliation $\Pi(M,\cF)$, which also consists of equivalence classes
of leafwise paths, but in this case two leafwise paths are
equivalent if they are homotopic in the class of leafwise paths with
fixed endpoints. The fundamental groupoid of the foliation
$\Pi(M,\cF)$ is a smooth groupoid (cf., for instance,
\cite{phillips87}).

There is a foliation ${\mathcal G}$ of dimension $2p$ on the
holonomy groupoid $G$. In any coordinate chart $W(\phi,\phi')$ given
by a pair of compatible foliated charts $\phi$ and $\phi'$ the
leaves of ${\mathcal G}$ are given by equations of the form $y={\rm
const}$. The leaf of $\cG$ through $\gamma\in G$ consists of all
$\gamma'\in G$ such that $r(\gamma)$ and $r(\gamma')$ lie on the
same leaf of $\cF$, and it coincides with the holonomy groupoid of
this leaf. The holonomy group of a leaf of $\cG$ coincides with the
holonomy group of the corresponding leaf of $\cF$.

The differential of the map $(r,s):G\to M\times M$ maps the tangent
bundle $T{\mathcal G}$ to $\cG$ isomorphically to the bundle
$F\boxtimes F$ on $M\times M$, therefore, there is a canonical
isomorphism $T{\mathcal G}\cong r^*F\oplus s^*F$.

A distribution $H$ on $M$ transverse to $\cF$ determines a
distribution $HG$ on $G$ transverse to $\cG$. For any $X\in H_y$,
there is a unique vector $\widehat{X}\in T_{\gamma}G$ such that $d
s(\widehat{X})= dh^{-1}_{\gamma}(X)$ and $d r(\widehat{X})=X$, where
$dh_{\gamma} : H_x\to H_y$ is the linear holonomy map associated
with $\gamma$. The space $H_\gamma G$ consists of all vectors of the
form $\widehat{X}\in T_{\gamma}G$ for different $X\in H_y$. In any
coordinate chart $W(\phi,\phi')$ on $G$ the tangent space $T_\gamma
{\mathcal G}$ of $\cG$ at some $\gamma $ with the coordinates
$(x,x',y)$ consists of vectors of the form
$X\frac{\partial}{\partial x}+X'\frac{\partial}{\partial x'}$, and
the distribution $H_\gamma G$ consists of vectors
$X\frac{\partial}{\partial x}+X'\frac{\partial}{\partial
x'}+Y\frac{\partial}{\partial y}$ such that
$X\frac{\partial}{\partial x}+Y\frac{\partial}{\partial y}\in
H_{(x,y)}$ and $X'\frac{\partial}{\partial
x'}+Y\frac{\partial}{\partial y}\in H_{(x',y)}$.

Let $g_M$ be a Riemannian metric on $M$ and $H=F^{\bot}$. Then a
Riemannian metric $g_G$ on $G$ is defined as follows. All the
components in $T_{\gamma}G=F_y\oplus F_x\oplus H_\gamma G$ are
mutually orthogonal, and by definition $g_G$ coincides with $g_M$ on
$F_y\oplus H_\gamma G\cong F_y\oplus H_y=T_yM$ and with $g_F$ on
$F_x$.

If $\cF$ is Riemannian and $g_M$ is a bundle-like Riemannian metric,
then $g_G$ is bundle-like, and therefore $\cG$ is Riemannian.
Moreover, in this case the maps $s:G\to M$ and $r:G\to M$ are
Riemannian submersions and locally trivial fibrations. In
particular, all the holonomy coverings $G^x$ of leaves of $\cF$ are
diffeomorphic \cite{Winkeln}.

Let $G$ be an arbitrary Lie groupoid with the set of units
$G^{(0)}=M$. For any (possibly non-Hausdorff) smooth manifold $Z$
and for any smooth map $\rho: Z\to M$, put
\[
Z\times_\rho G=\{(z,\gamma)\in Z\times G : \rho(z)=s(\gamma)\}.
\]
A smooth right action of the groupoid $G$ on $Z$ is defined to be a
map $Z\times_\rho G\to Z, (z,\gamma)\mapsto z\gamma$, satisfying the
conditions
\[
\rho(z\gamma)=r(\gamma), \quad (z\gamma)\gamma^\prime
=z(\gamma\gamma^\prime ), \quad z\cdot x=x.
\]

As an example, one can consider the action of $G$ on $M$ defined by
\[
\rho=\id : M \to M, \quad  M\times_\rho G \ni (y,\gamma) \mapsto
y\gamma=s(\gamma) \in M.
\]

An action of the groupoid $G$ on $Z$ is said to be proper, if:
\begin{enumerate}[(i)]
  \item the map $Z\times_\rho G\to Z\times Z, (z,\gamma)\mapsto
  (z,z\gamma),$ is proper (that is, the pre-image of every compact
  is compact);
  \item the set of equivalence classes $Z/\Gamma$ of the equivalence relation $\Gamma$
  on $Z$ with $z\sim z^\prime$ whenever $z\gamma=z^\prime$ for some $\gamma\in G$
 is Hausdorff.
\end{enumerate}

If $Z$ is a $G$-manifold, then the orbits of the $G$-action define a
foliation on $Z$.

\subsection{The $C^*$-algebra and the von Neumann algebra of a foliation}\label{s:calgebra}
In this subsection we will describe the construction of the
$C^*$-algebra associated with an arbitrary foliation. This algebra
can be regarded as an analogue of the algebra of continuous
functions on the leaf space of the foliation. We will only consider
Hausdorff groupoids. For the definition of the $C^*$-algebra of a
foliation in the case when the holonomy groupoid is not Hausdorff,
see, for instance, \cite{Co:survey}.

There are two ways to define the $C^*$-algebras associated with a
foliation. The first makes use of the auxiliary choice of a smooth
Haar system, the second does not require auxiliary choices and uses
the language of half-densities.

\subsubsection{Definitions using a Haar system}
In this subsection we give the definition of the $C^*$-algebras
associated with an arbitrary smooth groupoid $G$. In fact, the
assumption of smoothness of the groupoid is not essential here, and
all the definitions can be generalized to the case of topological
groupoids.

\begin{defn}
A smooth Haar system on a smooth groupoid $G$ is a family of
positive Radon measures $\{\nu^x:x\in G^{(0)}\}$ on $G$ satisfying
the following conditions:
\begin{enumerate}
  \item the support of $\nu^x$ coincides with $G^x$, and
  $\nu^x$ is a smooth positive measure on $G^x$;
  \item the family $\{\nu^x:x\in G^{(0)}\}$ is left-invariant,
  that is, for any continuous function $f\in
  C_c(G^x), f\geq 0,$ and for any $\gamma\in G$ such that
  $s(\gamma)=x$ and
  $r(\gamma)=y$ we have
\[
\int_{G^y}f(\gamma_1)d\nu^y(\gamma_1)=\int_{G^x}
f(\gamma\gamma_1)\,d\nu^x(\gamma_1);
\]
  \item the family $\{\nu^x:x\in G^{(0)}\}$ is smooth, that is,
  for any $\phi\in C^\infty_c(G)$, the function
\[
G^{(0)}\ni x \mapsto \int_{G^x}\phi(\gamma)d\nu^x(\gamma)
\]
is a smooth function on $G^{(0)}$.
\end{enumerate}
\end{defn}

For a compact foliated manifold $(M,\cF)$, a smooth Haar system
$\{\nu^x:x\in G^{(0)}\}$ on the holonomy groupoid $G$ of $\cF$ is
given by an arbitrary smooth positive leafwise density $\alpha \in
C^\infty(M,|T{\mathcal F}|)$. For any $x\in M$, the positive Radon
measure $\nu ^{x}$ on $G^x$ is defined as the lift of the density
$\alpha $ by the holonomy covering $s:G^x\to M$.

Let $G$ be a smooth groupoid, $G^{(0)}=M$ and $\{\nu^x:x\in M\}$ a
smooth Haar system. Introduce an involutive algebra structure on
$C^{\infty}_c(G)$ by
\begin{align*}
k_1\ast k_2(\gamma)&=\int_{G^x} k_1(\gamma_1)
k_2(\gamma^{-1}_1\gamma)\,d\nu^x(\gamma_1),\quad \gamma\in G^x,\\
k^*(\gamma)&=\overline{k(\gamma^{-1})}, \quad \gamma\in G.
\end{align*}
For any $x\in M$ there is a natural representation of
$C^{\infty}_c(G)$ on $L^2(G^x,\nu^x)$ given for $k\in
C^{\infty}_c(G)$ and $\zeta \in L^2(G^x,\nu^x)$ by
\[
R_x(k)\zeta(\gamma)=\int_{G^x}k(\gamma^{-1}\gamma_1) \zeta(\gamma_1)
d\nu^x(\gamma_1),\quad r(\gamma)=x.
\]
The completion of the involutive algebra $C^{\infty}_c(G)$ in the
norm
\[
\|k\|=\sup_x\|R_x(k)\|
\]
is called the reduced $C^*$-algebra of the groupoid $G$ and denoted
by $C^{\ast}_r(G)$. Also the full $C^{\ast}$-algebra of the groupoid
$C^{\ast}(G)$ is defined as the completion of $C^{\infty}_c(G)$ in
the norm
\[
\|k\|_{\text{max}}=\sup \|\pi(k)\|,
\]
where supremum is taken over the set of all $\ast$-representations
$\pi$ of the algebra $C^{\infty}_c(G)$ on Hilbert spaces.

\begin{ex}
In Example~\ref{ex:trivial} the groupoid $G$ is smooth if $X$ is a
smooth manifold. In this case $G^x=\{x\}$ for any $x\in X$, and a
smooth Haar system is given by an arbitrary smooth function on $X$.
The operator algebras $C^*_r(G)$ and $C^*(G)$ coincide with the
commutative $C^*$-algebra $C_0(X)$.
\end{ex}

\begin{ex}
The principal groupoid introduced in Example~\ref{ex:equiv} is
smooth if $X$ is a smooth manifold. In this case any smooth Haar
system has the following form: $\nu^x$ is a fixed positive smooth
density $\mu$ in each $G^x\cong X$. The operations in
$C^\infty_c(G)$ are given by
\begin{align*} (k_1\ast
k_2)(x,y)&=\int_Xk_1(x,z)k_2(z,y)d\mu(z), \quad (x,y)\in X\times X,
\\ k^*(x,y)&=\overline{k(y,x)}, \quad (x,y)\in X\times X,
\end{align*}
where $k, k_1, k_2\in C^\infty_c(G)$. Thus, elements of
$C^\infty_c(G)$ can be regarded as the kernels of integral operators
in $C^\infty(X)$ (with respect to the density $\mu$). For any $x\in
X$ the representation $R_x$ associates to every $k\in
C^\infty_c(G)\subset C^\infty(X\times X)$ the integral operator
acting in $L^2(G^x,\nu^x)\cong L^2(X,\mu)$ with integral kernel $k$:
\[
R_x(k)u(y)=\int_Xk(y,z)u(z)d\mu(z), \quad u\in L^2(X,\mu).
\]
Finally, $C^*_r(G)$ and $C^*(G)$ coincide with the algebra
$\cK(L^2(X,\mu))$ of compact operators in $L^2(X,\mu)$.
\end{ex}

\begin{ex}
In Example \ref{ex:Lie} a smooth Haar system is given by a
left-invariant Haar measure $dh$ on $H$: $\nu^x=dh$. The product in
$C^\infty_c(G)$ is the classical convolution operation given for any
functions $u, v\in C^\infty_c(H)$ by
\[
(u\ast v)(g)=\int_Hu(h)v(h^{-1}g)\,dh, \quad g\in H,
\]
the involution is given by
\[
u^*(g)=\overline{u(g^{-1})}, \quad u\in C^\infty_c(H),
\]
and the operator algebras $C^*_r(G)$ and $C^*(G)$ are the group
$C^*$-algebras $C^*_r(H)$ and $C^*(H)$.
\end{ex}

\begin{ex}\label{ex:crossed}
In Example \ref{ex:actions} the manifold $G^x=\{(x,h) : h\in H\}$ is
diffeomorphic to $H$ for any $x\in X$ and a smooth Haar system on
$G$ can be defined by using an arbitrary Haar measure on $H$. The
product in $C^\infty_c(G)$ is given for any functions $u, v\in
C^\infty_c(X\times H)$ by
\[
(u\ast v)(x,g)=\int_Hu(x,h)v(h^{-1}x, h^{-1}g)\,dh, \quad (x,g)\in
X\times H,
\]
the involution is given, for a function $u\in C^\infty_c(X\times
H)$, by
\[
u^*(x,g)=\overline{u(g^{-1}x, g^{-1})}, \quad \quad (x,g)\in
X\times H.
\]
The operator algebras $C^*_r(G)$ and $C^*(G)$ corresponding to the
crossed product groupoid $G=X\rtimes H$ coincide with the crossed
products $C_0(X)\rtimes_r H$ and $C_0(X)\rtimes H$ of the algebra
$C_0(X)$ by the group $H$ with respect to the induced action of the
group $H$ on $C_0(X)$.

If the group $H$ is discrete, then elements of the algebra
$C^\infty_c(G)$ are families $\{a_\gamma \in C^\infty_c(X) :
\gamma\in H\}$ such that $a_\gamma\not=0$ for finitely many elements
$\gamma$. It convenient to write them in the form $a=\sum_{\gamma\in
H} a_\gamma U_\gamma $. The product in the algebra $C^\infty_c(G)$
is written as
\[
(a_{\gamma_1} U_{\gamma_1} )(b_{\gamma_2} U_{\gamma_2}) =
(a_{\gamma_1} T_{\gamma_1}(b_{\gamma_2}))U_{\gamma_1\gamma_2},
\]
where $T_\gamma$ denotes the operator in $C_0(X)$ induced by the
action of $\gamma\in H$:
\[
T_\gamma f(x)=f(\gamma^{-1}x),\quad x\in X,\quad f\in C_0(X).
\]
The involution in the algebra $C^\infty_c(G)$ is given by
\[
(a_\gamma U_\gamma )^*=T_{\gamma^{-1}}(\bar{a}_\gamma)
U_{\gamma^{-1}}.
\]
\end{ex}

Let $G$ be the holonomy groupoid of a foliation $\cF$ on a compact
manifold $M$. Elements of the algebra $C^\infty_c(G)$ can be
regarded as families of the kernels of integral operators along the
leaves of the foliation (more precisely, on the holonomy coverings
$G^x$). Namely, each $k\in C^\infty_c(G)$ corresponds to the family
$\{R_x(k):x\in M\}$, where $R_x(k)$ is the integral operator acting
in $L^2(G^x,\nu^x)$ given by integral kernel
\[
K(\gamma_1,\gamma_2)=k(\gamma_1^{-1}\gamma_2), \quad \gamma_1,
\gamma_2\in G^x.
\]
The product of elements $k_1$ and $k_2$ in $C^\infty_c(G)$
corresponds to the composition of integral operators
$\{R_x(k_1)R_x(k_2):x\in M\}$. The $C^*$-algebra $C^*(G)$
(respectively, $C^*_r(G)$) associated with the holonomy groupoid $G$
of the foliation $(M,\cF)$ will be called the $C^*$-algebra of the
foliation $(M,\cF)$ (respectively, the reduced $C^*$-algebra of the
foliation $(M,\cF)$) and denoted by $C^*(M,\cF)$ (respectively,
$C^*_r(M,\cF)$).

\begin{ex}
Let $\cF$ be a simple foliation on a compact manifold $M$ defined by
a submersion $\pi :M\to B$. Fix a smooth Haar system on it. For any
$y\in B$, denote by $\Psi^{-\infty}(Z_y)$ the involutive algebra of
integral operators with smooth kernel acting in the space
$C^\infty(Z_y)$, where $Z_y$ is the fibre of the fibration $\pi$ at
$y$. Let us consider the field $\Psi^{-\infty}(M/B)$ of involutive
algebras on $B$, whose fibre at $y\in B$ is $\Psi^{-\infty}(Z_y)$.
For any section $\sigma$ of the field $\Psi^{-\infty}(M/B)$ the
integral kernels of the operators $\sigma_y$ give rise to a
well-defined function on the holonomy groupoid $G$ of the foliation
$\cF$. We say that a section $\sigma$ is smooth if the corresponding
function on $G$ is smooth. Thus, we obtain a description of the
algebra $C^\infty(G)$ as the algebra of smooth sections of the field
$\Psi^{-\infty}(M/B)$ of fibrewise integral operators with smooth
kernels.
\end{ex}

\begin{ex}
Consider an example of the linear foliation on the torus. Thus, let
$M=T^2=\RR^2/\ZZ^2$ be the two-dimensional torus, and suppose that a
foliation $\cF_\theta$ is given by the trajectories of the vector
field $X=\frac{\partial}{\partial x}+ \theta\frac{\partial}{\partial
y},$ where $\theta\in \RR$ is a fixed irrational number.

Since this foliation is given by the orbits of a free group action
of $\RR$ on $T^2$, its holonomy groupoid coincides with the crossed
product groupoid  $T^2\rtimes \RR$. Thus, $G=T^2\times \RR$,
$G^{(0)}=T^2$, $s(x,y,t)=(x-t,y-\theta t)$, $r(x,y,t)=(x,y)$,
$(x,y)\in T^2, t\in\RR,$ and the multiplication is given by
\[
(x_1,y_1,t_1)(x_2,y_2,t_2)=(x_1,y_1,t_1+t_2),
\]
if $x_2=x_1-t_1, y_2=y_1-\theta t_1$.

The reduced $C^*$-algebra $C^*_r(T^2,\cF_\theta)$ of the linear
foliation $\cF_\theta$ on $T^2$ coincides with the reduced crossed
product $C(T^2)\rtimes_r \RR$. Therefore, the product $k_1\ast k_2$
of $k_1, k_2\in C^\infty_c(T^2\times \RR)\subset C(T^2)\rtimes_r
\RR$ is given by
\begin{multline*}
(k_1\ast k_2)(x,y,t)=\int_{-\infty}^\infty
k_1(x_1,y_1,t_1)k_2(x-t_1, y-\theta t_1, t-t_1)\,dt_1, \\ (x,y)\in
T^2,\quad t\in\RR,
\end{multline*}
and for any $k\in C^\infty_c(T^2\times \RR)$,
\[
k^*(x,y,t)=\overline{k(x-t,y-\theta t,-t)}, \quad (x,y)\in T^2,\quad
t\in\RR.
\]
For any $k\in C^\infty_c(T^2\times \RR)$ and for any $(x,y)\in T^2$
the operator $R_{(x,y)}(k)$ has the following form on
$L^2(G^{(x,y)}, \nu^{(x,y)})\cong L^2(\RR,dt)$: for any $u\in
L^2(\RR, dt)$,
\[
R_{(x,y)}(k)u(t)=\int_{-\infty}^\infty k(x-t_1, y-\theta t_1,
t-t_1)u(t_1)\,dt_1, \quad t\in\RR.
\]

If $\theta$ is rational, then the linear foliation on $T^2$ is given
by the orbits of a free group action of $S^1$ on $T^2$, and its
holonomy groupoid coincides with the crossed product groupoid
$G=T^2\rtimes S^1$.
\end{ex}

We note some facts which connect the structure of the reduced
$C^*$-algebra of a foliation $C^*_r(M,\cF)$ with the topology of
$\cF$ (for more details cf. \cite{F-Sk,Hil-Skan83}).

\begin{thm}[\cite{F-Sk}]
Let $(M,\cF)$ be a foliated manifold.
\begin{enumerate}
  \item The $C^*$-algebra $C^*_r(M,\cF)$ is simple if and only if
  $\cF$ is minimal, that is, every its leaf is dense in $M$.
  \item The $C^*$-algebra $C^*_r(M,\cF)$ is primitive if and only if $\cF$ is
   (topologically) transitive, that is, it has a leaf, which is dense in $M$.
  \item If $\cF$ is amenable in the sense that
  $C^*_r(M,\cF)=C^*(M,\cF)$ the $C^*$-algebra $C^*_r(M,\cF)$ has a
  representation, consisting of compact operators if and only if $\cF$ has
  a compact leaf.
\end{enumerate}
\end{thm}

In the paper \cite{F-Sk} a description is given of the space of
primitive ideals of the $C^*$-algebra $C^*_r(M,\cF)$.

\subsubsection{Definition using half-densities}
In this subsection we will give the definitions of the operator
algebras associated with a foliated manifold without using a choice
of Haar system. For this, we will use the language of
half-densities.

Let $(M,\cF)$ be a compact foliated manifold, and consider the
vector bundle of leafwise half-densities $|T{\mathcal F}|^{1/2}$ on
$M$. Using the source map $s$ and the range map $r$, lift
$|T{\mathcal F}|^{1/2}$ to the vector bundles $s^*(|T{\mathcal
F}|^{1/2})$ and $r^*(|T{\mathcal F}|^{1/2})$ on the holonomy
groupoid $G$. We define a vector bundle $|T\cG|^{1/2}$ on $G$ as $$
|T\cG|^{1/2}=r^*(|T{\mathcal F}|^{1/2})\otimes s^*(|T{\mathcal
F}|^{1/2}). $$ The bundle $|T\cG|^{1/2}$ is naturally identified
with the bundle of leafwise half-densities on the foliated manifold
$(G,{\mathcal G})$.

An involutive algebra structure on $C^{\infty}_c(G,|T\cG|^{1/2})$ is
defined by
\begin{align*}
\sigma_1\ast \sigma_2(\gamma)&=\int_{\gamma_1\gamma_2=\gamma}
\sigma_1(\gamma_1)\sigma_2(\gamma_2),\quad \gamma\in G,\\
\sigma^*(\gamma)&=\overline{\sigma(\gamma^{-1})},\quad \gamma\in G,
\end{align*}
where $\sigma, \sigma_1, \sigma_2\in C^{\infty}_c (G,
|T\cG|^{1/2})$. The formula for $\sigma_1\ast \sigma_2$ should be
interpreted in the following way. If we write $\gamma : x\to y,
\gamma_1: z\to y $ and $\gamma_2: x\to z$, then
\begin{align*}
\sigma_1(\gamma_1)\sigma_2(\gamma_2) & \in |T_y{\mathcal
F}|^{1/2}\otimes |T_z{\mathcal F}|^{1/2}\otimes |T_z{\mathcal
F}|^{1/2}\otimes |T_x{\mathcal F}|^{1/2}\\ & \cong |T_y{\mathcal
F}|^{1/2}\otimes |T_z{\mathcal F}|^{1}\otimes |T_x{\mathcal
F}|^{1/2},
\end{align*}
and integrating the $|T_z{\mathcal F}|^{1}$-component
$\sigma_1(\gamma_1)\sigma_2(\gamma_2)$ with respect to $z\in M$ we
obtain a well-defined section of the bundle $r^*(|T{\mathcal
F}|^{1/2})\otimes s^*(|T{\mathcal F}|^{1/2})=|T\cG|^{1/2}. $

\subsection{Vector bundles and fields of Hilbert spaces}
\label{s:vect} A natural analogue of the notion of vector bundle on
the leaf space of a foliation is the notion of holonomy equivariant
vector bundle, or vector $G$-bundle.

\begin{defn}
A vector bundle $E$ (complex or real) on a foliated manifold
$(M,\cF)$ is said to be holonomy equivariant if there is given a
representation $T$ of the holonomy groupoid $G$ of the foliation
$\cF$ in the fibres of $E$, that is, for any $\gamma\in G,
\gamma:x\rightarrow y$, there is defined a linear operator
$T(\gamma):E_x\rightarrow E_y$ such that
$T(\gamma_1\gamma_2)=T(\gamma_1)T(\gamma_2)$ for any
$\gamma_1,\gamma_2 \in G$ with $r(\gamma_2)=s(\gamma_1)$.

A Hermitian (respectively, Riemannian) vector bundle $E$ on a
foliated manifold $(M,\cF)$ is said to be holonomy equivariant if it
is a holonomy equivariant vector bundle and the representation $T$
is unitary (respectively, orthogonal): $T(\gamma^{-1})=T(\gamma)^*$
for any $\gamma\in G$.
\end{defn}

For any holonomy equivariant vector bundle $E\rightarrow M$ the
action of the groupoid $G$ on $E$ defines a horizontal foliation
$\cF_E$ on $E$ of the same dimension as the foliation $\cF$. The
leaf of $\cF_E$ through a point $v\in E$ consists of all points of
the form $T(\gamma)^{-1}(v)$ with $\gamma\in G$, $r(\gamma)=\pi(v)$.
Thus, any holonomy equivariant vector bundle is foliated in the
sense of \cite{KamberTondLNM}.

\begin{defn}
A vector bundle $p:P\to M$ is said to be foliated if there is a
foliation $\bar{\cF}$ on $P$ of the same dimension as $\cF$ whose
leaves are transversal to the fibres of $p$ and mapped by $p$ to the
leaves of $\cF$.
\end{defn}

Equivalently, one can say that a foliated vector bundle is a vector
bundle $P$ on $M$ such that there exists a flat connection in the
space $C^\infty(M,P)$ defined along the leaves of $\cF$, that is, an
operator
\[
\nabla : \cX(\cF)\times C^\infty(M,P)\to C^\infty(M,P),
\]
satisfying the standard conditions
\[
\nabla_{fX}=f\nabla_{X}, \quad  \nabla_X(fs)=(X f)s +f \nabla_Xs,
\]
for any $f \in C^\infty(M), X\in \cX(\cF), s \in C^\infty(M,P)$, and
also the flatness condition
\[
[\nabla_X, \nabla_Y]=\nabla_{[X,Y]}, \quad X, Y\in \cX(\cF).
\]
The parallel transport along leafwise paths associated with the
connection $\nabla $ defines an action of the fundamental groupoid
of the foliation $\Pi(M,\cF)$ in the fibres of the foliated vector
bundle $P$. In general, the parallel transport may depend on the
holonomy of the corresponding path, therefore, this action does not
necessarily descent to an action of the holonomy groupoid in the
fibres of $P$, and hence a foliated vector bundle is not necessarily
holonomy equivariant.

\begin{ex}
The normal bundle $\tau_x=T_xM/T_x{\mathcal F}, x\in M,$ is a
holonomy equivariant vector bundle, if it is equipped with the
action of the holonomy groupoid $G$ by the linear holonomy map
$dh_{\gamma}:\tau_x\rightarrow \tau_y, \gamma:x\rightarrow y$. The
corresponding partial flat connection defined along the leaves of
$\cF$ is the Bott connection ${\stackrel{\circ}\nabla}$ (see
(\ref{e:Bott})). The normal bundle $\tau$ is a holonomy equivariant
Riemannian vector bundle, if $\cF$ is a Riemannian foliation.

The conormal bundle $N^*{\mathcal F}$ equipped with the action of
the holonomy groupoid $G$ by the linear holonomy map
$(dh^*_{\gamma})^{-1}: N^*_x{\mathcal F}\rightarrow N^*_y{\mathcal
F}$ for $\gamma:x\rightarrow y$, and, more generally, an arbitrary
tensor bundle associated with the normal bundle $\tau$ are also
holonomy equivariant.
\end{ex}

For any holonomy equivariant vector bundle $E$ on a foliated
manifold $(M,\cF)$ there is defined a natural representation $R_E$
of the algebra $C^\infty_c(G)$ in the space $C^{\infty}(M,E)$ of
sections of this bundle. Let $\{\nu^x: x\in G^{(0)}\}$ be a smooth
Haar system on $G$. For any $u\in C^{\infty}(M,E)$ the section
$R_E(k)u\in C^{\infty}(M,E)$ is given by
\begin{equation}
\label{rev:tang1}
R_E(k)u(x)=\int_{G^x}k(\gamma)\,T(\gamma)[u(s(\gamma))]\,d\nu^x(\gamma),
\quad x\in M.
\end{equation}

If $E$ is holonomy equivariant Hermitian vector bundle on $M$, then
the representation $R_E$ is a $\ast$-representation, and therefore
it extends to a $\ast$-representation of the $C^*$-algebra
$C^*(M,\cF)$.

In \cite{Co:survey} the notion of continuous field of Hilbert spaces
on the leaf space of the foliation $\cF$ introduced, and it is shown
that there is a one-to-one correspondence between continuous fields
of Hilbert spaces on the leaf space of the foliation $\cF$ and
Hilbert $C^*$-modules over the $C^*$-algebra $C^*_r(M,\cF)$.

Let $H=\{H_x : x\in M\}$ be a measurable field of Hilbert spaces on
$M$ equipped with a unitary representation of the holonomy groupoid
$G$:
\[
U(\gamma) : H_x\to H_y,\quad \gamma\in G, \quad \gamma :x\to y.
\]

For any measurable sections $\xi,\eta$ of $H$ define a function
$(\xi,\eta)$ on $G$ by
\[
(\xi,\eta)(\gamma)=\langle\xi_y,U(\gamma)\eta_x\rangle, \quad
\gamma\in G, \quad \gamma : x\to y.
\]

We fix an arbitrary smooth Haar system $\{\nu^x : x\in M\}$ on $G$.
For any measurable section $\xi$ of the field $H$ denote by
$\|\xi\|_\infty$ the least number $c\in [0,+\infty]$ such that for
any $y\in M$ and for any $\alpha\in H_y$
\[
\int_{G^y} |\langle \alpha, U(\gamma)\xi_x\rangle
|^2d\nu^y(\gamma)\leq c\|\alpha\|^2.
\]

\begin{defn}
A field $\{H_x : x\in M\}$ is called a continuous field of Hilbert
spaces on the leaf space of the foliation $\cF$ if there is a
distinguished linear space $\Gamma$ of its sections such that
\begin{enumerate}
  \item $\Gamma$ contains a countable total $\|\ \|_\infty$-dense
  subset; in particular, $\|\ \|_\infty$ takes finite values on
  $\Gamma$;
  \item  $(\xi,\eta)\in C^*_r(M,\cF)$ for any $\xi,\eta\in \Gamma$;
  \item $\Gamma$ is closed in the $\|\ \|_\infty$-norm;
  \item for any $f\in C^\infty_c(G)$ and $\xi\in\Gamma$ one has $\xi\ast
  f\in \Gamma$, where (cf. (\ref{rev:tang1}))
\[
 (\xi\ast f)(y)=\int_{G^y}f(\gamma) U(\gamma)\xi(x)d\nu^y(\gamma).
  \]
\end{enumerate}
\end{defn}

We will not present a general construction of the Hilbert
$C^*$-module over the $C^*$-algebra $C^*_r(M,\cF)$ associated with
an arbitrary continuous field of Hilbert spaces on the leaf space of
a foliation $\cF$, but consider only one important particular case.
Let $E$ be a vector bundle on a foliated manifold $(M,\cF)$. It
defines a continuous field of Hilbert spaces on the leaf space of a
foliation $\cF$ and therefore a Hilbert $C^*$-module over the
$C^*$-algebra $C^*_r(M,\cF)$. It can be defined in two ways,
depending on whether we wish to consider left or right modules. If
we wish to work with left modules, then, as in \cite{Co79} we put
$H_x=L^2(G^x,s^*E)$ for any $x\in M$, and define for any $\gamma :
x\to y$ the operator $L(\gamma) : L^2(G^x,s^*E) \to L^2(G^y,s^*E)$
acting on $\xi\in L^2(G^x,s^*E)$ by the formula
\begin{equation}\label{def:L}
L(\gamma)\xi(\gamma_1)=\xi(\gamma^{-1}\gamma_1), \quad \gamma_1\in
G^y.
\end{equation}

To work with right modules, we put $H_x=L^2(G_x,r^*E)$ for any $x\in
M$ (as in \cite{Co:survey,Co-skandal}) and define for any $\gamma :
x\to y$ the operator $R(\gamma) : L^2(G_x,r^*E) \to L^2(G_y,r^*E)$
acting on $\xi\in L^2(G_x,r^*E)$ by
\begin{equation}\label{def:R}
R(\gamma)\xi(\gamma_1)=\xi(\gamma_1\gamma), \quad \gamma_1\in G_x.
\end{equation}

Let us describe the corresponding right Hilbert $C^*$-module over
the $C^*$-algebra $C^*_r(M,\cF)$ \cite{Co:survey,Co-skandal}. We
begin with the definition of a pre-Hilbert $C^\infty_c(G)$-module
$\mathcal E_\infty$. As a linear space $\mathcal E_\infty$ coincides
with $C^\infty_c(G, r^*E)$. The module structure on $\cE_\infty$ is
introduced as follows: the action of $f\in C^\infty_c(G)$ on $s\in
\cE_\infty$ is given by
\[
(s\ast f)(\gamma)=\int_{G^y}
s(\gamma')f({\gamma'}^{-1}\gamma)d\nu^y(\gamma'), \quad \gamma\in
G^y,
\]
and the inner product on $\cE_\infty $ with values in
$C^\infty_c(G)$ is given by
\[
\langle s_1, s_2\rangle (\gamma)=\int_{G^y}\langle
s_1({\gamma'}^{-1}), s_2({\gamma'}^{-1}\gamma )
\rangle_{E_{s(\gamma')}} d\nu^y(\gamma'), \quad s_1, s_2\in
\cE_\infty.
\]
There is also a left action of the algebra $C^\infty(M)$ on
$\cE_\infty$ given by
\[
(a\cdot s)(\gamma)=a(y)s(\gamma), \quad \gamma\in G^y,
\]
where $a\in C^\infty_c(M)$ and $s\in \cE_\infty$.

The completion of the space $\cE_\infty $ in the norm
$\|s\|=\|R(\langle s, s\rangle)\|^{1/2}$ defines a Hilbert
$C^*$-module over the algebra $C^*_r(M,\cF)$, which we denote by
$\cE=\cE_{M,E}$. It has a $C^*_r(M,\cF)$-valued
$C^*_r(M,\cF)$-sesquilinear form $\langle \cdot, \cdot\rangle $
which is the extension by continuity of the sesquilinear form on
$\cE_\infty $. Thus, $\cE_{M,E}$ is a
$C(M)$-$C^*_r(M,\cF)$-bimodule.

If $E$ is holonomy equivariant, there is a left action of the
algebra $C^\infty_c(G)$ on $\cE_\infty$ given by
\[
(f\ast
s)(\gamma)=\int_{G^y}
f(\gamma')T(\gamma')[s({\gamma'}^{-1}\gamma)]d\nu^y(\gamma'), \quad
\gamma\in G^y,
\]
where $f\in C^\infty_c(G)$ and $s\in \cE_\infty$. Unlike the right
action, the left action does not extend to an action of the algebra
$C^*_r(M,\cF)$ by bounded endomorphisms of the $C^*$-Hilbert module
$\cE$ over $C^*_r(M,\cF)$. Nevertheless, using the methods of
reduction to the maximal compact subgroup developed by Kasparov
\cite{Kasparov-conspectus}, one can construct an element $[E]\in
KK(C^*(M,\cF),C^*(M,\cF))$ corresponding to an arbitrary holonomy
equivariant bundle $E$. This correspondence is a generalization of
the map $j^G : KK^G(C_0(X),C_0(Y))\to KK^G(C_0(X)\rtimes G,
C_0(Y)\rtimes G)$ constructed by Kasparov \cite{Kasparov-conspectus}
for any $G$-manifolds $X$ and $Y$.

\subsection{Strong Morita equivalence and transversals}
\label{s:morita1} Consider a compact foliated manifold $(M,\cF)$. A
choice of a complete transversal for the foliation enables us to
reduce the holonomy groupoid $G$ of the foliation $\cF$ to an
equivalent groupoid, which in many cases turns out to be simpler.

For any two subsets $A, B\subset M$ let
\[
G^A_B=\{\gamma\in G : r(\gamma)\in A, s(\gamma)\in B\}.
\]
In particular,
\[
G^M_T=G_T=\{\gamma\in G : s(\gamma)\in T\}.
\]
If $T$ is a complete transversal, then $G^T_T$ is a submanifold and
a subgroupoid of $G$. It is called a reduced holonomy groupoid.

As shown in \cite{Hil-Skan83}, if $T$ is a complete transversal,
then the reduced $C^*$-algebras $C^*_r(M,\cF)$ and $C^*_r(G^T_T)$
are strongly Morita equivalent. In particular, this easily implies
that
\[
C^*_r(M,\cF) \cong \cK \otimes C^*_r(G^T_T).
\]

Following the paper \cite{MRW87}, we describe the construction of a
$C^*_r(M,\cF)$-$C^*_r(G^T_T)$-equivalence bimodule which gives the
strong Morita equivalence of the algebras $C^*_r(M,\cF)$ and
$C^*_r(G^T_T)$, .

Consider the manifold $P=G_T$. There is a natural left action of the
groupoid $G$ on $P$ given by left multiplication in $G$, and a right
action of the groupoid $G^T_T$ given by right multiplication in $G$.
These actions commute. In the language of the paper \cite{MRW87} the
manifold $P$ is a $(G,G^T_T)$-equivalence.

Correspondingly, there is a left action of the algebra
$C^\infty_c(G)$ on $C^\infty_c(G_T)$ given for any $f\in
C^\infty_c(G)$ and $\varphi\in C^\infty_c(G_T)$ by
\[
(f\cdot \varphi)(p)=\int_{G^{r(p)}}
f(\gamma)\varphi(\gamma^{-1}p)d\nu^{r(p)}(\gamma),\quad p\in G_T,
\]
and a right action of the algebra $C^\infty_c(G^T_T)$ on
$C^\infty_c(G_T)$ given for any $g\in C^\infty_c(G^T_T)$ and
$\varphi\in C^\infty_c(G_T)$ by
\[
(\varphi\cdot g)(p)=\sum_{\gamma\in G^T_{s(p)}}
\varphi(p\gamma)g(\gamma^{-1}),\quad p\in G_T.
\]
The inner product on $C^\infty_c(G_T)$ with values in the algebra
$C^\infty_c(G^T_T)$ is given by the following formula: for any
$\varphi, \psi\in C^\infty_c(G_T)$
\[
\langle \varphi,\psi \rangle_{C^\infty_c(G^T_T)}(\gamma)
=\int_{G^{r(p)}} \overline{\varphi(\gamma_1^{-1}p)}
\psi(\gamma_1^{-1}p\gamma) d\nu^{r(p)}(\gamma_1),\quad \gamma\in
G^T_T,
\]
where $p\in G_{r(\gamma)}$ is an arbitrary point (the right hand
side of the formula is independent of the choice of $p$). Similarly,
the inner product on $C^\infty_c(G_T)$ with values in the algebra
$C^\infty_c(G)$ is given by the formula: for any $\varphi, \psi\in
C^\infty_c(G_T)$
\[
\langle \varphi,\psi \rangle_{C^\infty_c(G)}(\gamma) =
\sum_{\gamma_1\in G^T_{s(p)}} \varphi(\gamma^{-1}p\gamma_1)
\overline{\psi(p\gamma_1)},\quad \gamma\in G,
\]
where $p\in G^{r(\gamma)}$ is an arbitrary point (the right-hand
side of the formula is independent of the choice of $p$).

As a result, one has isomorphisms of the Hochschild homology, cyclic
homology and periodic cyclic homology of the algebras
$C^\infty_c(G)$ and $C^\infty_c(G^T_T)$ and of the $K$-theory of the
$C^*$-algebras $C^*_r(M,\cF)$ and $C^*_r(G^T_T)$. It is also proved
in \cite{Haefliger84} that the natural embedding $G^T_T\subset G$
induces a homotopy equivalence $BG^T_T\backsimeq BG$ of the
classifying spaces of $G^T_T$ and $G$ (see Subsection~\ref{s:BC}).

\begin{ex}
If $\cF$ is a simple foliation defined by a fibration $M\to B$, then
the $C^*$-algebra $C^*_r(M,\cF)$ is strongly Morita equivalent to
the $C^*$-algebra $C_0(B)$.
\end{ex}

\begin{ex}
On the manifold $M=\tilde{B}\times_\Gamma F$ we consider a foliation
$\cF$  obtained from the manifold $B$ and a homomorphism $\phi :
\Gamma=\pi_1(B) \to \operatorname{Diff}(F)$ by the suspension
construction (see Example~\ref{ex:suspension}). The image of the set
$\{b_0\}\times F\subset \tilde{B}\times F$ ($b_0\in \tilde{B}$ is an
arbitrary element) under the projection $\tilde{B}\times F \to M$ is
a complete transversal to $\cF$. If the condition given in
Example~\ref{ex:susp} holds, then the algebra $C^*_r(M,\cF)$ is
strongly Morita equivalent to the reduced crossed product
$C(F)\rtimes_r\Gamma$.
\end{ex}

\begin{ex}\label{ex:Atheta}
Consider the linear foliation $\cF_\theta$ on the two-dimensional
torus $T^2$, where $\theta\in \RR$ is a fixed irrational number, and
the transversal $T$ given by the equation $y=0$. The leaf space of
the foliation $\cF_\theta$ is identified with the orbit space  of
the $\ZZ$-action on $S^1=\RR/\ZZ$ generated by the rotation
\[
R_\theta(x)=x-\theta \mod 1,\quad x\in S^1.
\]
The algebra $A_\theta=C^*(G^T_T)$ coincides with the crossed product
$C(S^1)\rtimes \ZZ$ of the algebra $C(S^1)$ by the group $\ZZ$ with
respect to the $\ZZ$-action $R_\theta$ on $C(S^1)$. One can show
that the algebra $A_\theta$ is generated by the elements $U$ and
$V$, satisfying the relation $VU=e^{2\pi i \theta} U V$. It has a
concrete realization as the uniform closure of the subalgebra of
$\cL(L^2(S^1\times \ZZ))$ consisting of finite sums of the form
$\sum_{(n,m)\in\ZZ^2}a_{nm}u^nv^m$, where $a_{nm}\in \CC$ and the
operators $u$ and $v$ have the following form: for $f\in
L^2(S^1\times \ZZ)$
\[
uf(x,n)=f(x,n+1), \quad vf(x,n)=e^{2\pi i (x-n\theta)} f(x,n), \quad
x\in S^1, \quad n\in \ZZ.
\]

In many cases it is convenient to consider as a dense subalgebra of
$A_\theta$ the algebra
\[
\cA_\theta=\left\{a = \sum_{(n,m)\in\ZZ^2}a_{nm}U^nV^m :
\{a_{nm}\}\in \cS(\ZZ^2)\right\},
\]
where $\cS(\ZZ^2)$ is the space of rapidly decreasing sequences
(that is, of sequences such that
$\sup_{(n,m)\in\ZZ^2}(|n|+|m|)^k|a_{nm}|<\infty$ for any natural
$k$).

For $\theta=0$, the algebra $A_\theta$ is commutative and isomorphic
to the commutative $C^*$-algebra of continuous functions on the
two-dimensional torus $T^2$. Therefore. for an arbitrary $\theta$
the algebra $A_\theta$ is often called the algebra of continuous
functions on the non-commutative torus $T^2_\theta$, and the
notation $\cA_\theta= C^\infty(T^2_\theta)$,
$A_\theta=C(T^2_\theta)$ is used. The algebra $A_\theta$ was
introduced in the paper \cite{RieffelRot} (see also \cite{CoCRAS})
and has found many applications in mathematics and physics (see, for
instance, the survey \cite{konechny:schw}).

Thus, the $C^*$-algebra $C^*_r(T^2,\cF_\theta)$ of the linear
foliation on $T^2$ is strongly Morita equivalent to the algebra
$A_\theta$.
\end{ex}

An important property of the groupoid $G^T_T$ associated with a
complete transversal $T$ is (see \cite{CrainicM01}) that it is an
\'etale groupoid, that is, the source map $s: G\to G^{(0)}$ is a
local diffeomorphism. According to \cite{CrainicM01}, a smooth
groupoid is equivalent to an \'etale one if and only if all its
isotropy groups $G^x_x$ are discrete. In the latter case such a
groupoid is called a foliation groupoid.

One can introduce an equivalence relation for groupoids that is
similar to the strong Morita equivalence for $C^*$-algebras
\cite{Haefliger84,MRW87} (see also \cite{CrainicM01}). Roughly
speaking, two groupoids are equivalent if they have the same orbits
spaces, and therefore the same transverse geometries. It was proved
in \cite{MRW87} that if groupoids $G$ and $H$ are equivalent, then
their reduced $C^*$-algebras are strongly Morita equivalent.

\subsection{$K$-theory of foliation $C^*$-algebras and the
Baum-Connes conjecture}\label{s:BC} Computation of the $K$-theory
for foliation $C^*$-algebras is not a simple problem. This problem
has a simple solution for foliations with proper holonomy groupoid.
The holonomy groupoid $G_X$ of a foliation $\cF_X$ on a manifold $X$
is said to be proper if the map $(r,s): G_X\to X\times X$ is proper.
If the holonomy groupoid of a foliation $(X,\cF_X)$ is proper, then
the foliation $\cF_X$ is proper, that is, every its leaf is an
embedded submanifold of $X$. The leaf space $X/\cF_X$ of a foliation
$(X,\cF_X)$ with proper holonomy groupoid is an orbifold, and the
group $K(C^*(X,\cF_X))$ coincides with the $K$-theory defined by
$G_X$-equivariant bundles on $X$ compactly supported in $X/\cF_X$
(see \cite{Hil-Skan}).

In the general case a geometric construction of elements from
$K(C^*(M,\cF))$ was proposed in \cite{Baum-Connes} (see also
\cite{Co-skandal}). It is based on a definition of the geometric
$K$-theory $K^*_{\rm top}(M,\cF)$. Before we give this definition,
we introduce some auxiliary notions.

We will need the notion of smooth map $f$ from the leaf space
$M_1/\cF_1$ of a foliation $\cF_1$ on a manifold $M_1$ to the leaf
space $M_2/\cF_2$ of a foliation $\cF_2$ on a manifold $M_2$. A
smooth map $f: M_1\to M_2/\cF_2$ can be defined as a cocycle  $(U_i,
\gamma_{ij})$ on $M_1$ with values in the holonomy groupoid $G_2$ of
the foliation $\cF_2$. Here $\{U_i\}_{i\in A}$ is an open cover of
the manifold $M_1$, and the $\gamma_{ij} : U_i\cap U_j\to G_2$ are
smooth maps, satisfying the relations
\[
\gamma_{ij}(x)\gamma_{jk}(x)=\gamma_{ik}(x), \quad x\in U_i\cap
U_j\cap U_k.
\]
More precisely, we say that two cocycles $(U_i, \gamma_{ij})$ and
$(U^\prime_i, \gamma^\prime_{ij})$ are equivalent, if they extend to
a cocycle on the disjoint union of the coverings $\{U_i\}$ and
$\{U^\prime_i\}$. A map $f: M_1\to M_2/\cF_2$ is defined as an
equivalence class of cocycles. This definition can be understood as
follows. The maps $\gamma_{ii} : U_i\to (G_2)^{(0)}=M_2$ are the
local lifts of the map $f : U_i\to M_2/\cF_2$. The fact that, for
different $i$ and $j$ and for $x\in U_i\cap U_j$, the points
$f_i(x)$ and $f_j(x)$ lie on the same leaf of the foliation $\cF_2$
is described by the element $\gamma_{ij}(x): f_j(x)\to f_i(x)$ of
the holonomy groupoid $G_2$.

The graph of a map $f: M_1\to M_2/\cF_2$ given by a cocycle $(U_i,
\gamma_{ij})$ is a smooth (not necessarily Hausdorff) manifold
$G_f$, which is defined as the set of equivalence classes on the set
$\{(x,i,\gamma)\in M_1\times A\times G_2 : x\in U_i,
r(\gamma)=f_i(x)\}$ with respect to the equivalence relation,
specifying that $(x,i,\gamma)\sim (x^\prime,j,\gamma^\prime)$, if
$x=x^\prime$ and $\gamma_{ji}(x)=\gamma^\prime\gamma^{-1}$. There is
a smooth map $r_f : G_f\to M_1$, $r_f(x,i,\gamma)=x$, and a right
action of the groupoid $G_2$ on $G_f$ given by the maps $s_f :
G_f\to M_2$, $s_f(x,i,\gamma)=s(\gamma)$ and $G_f\times_{s_f}G_2\to
G_f$, $(x,i,\gamma)\gamma^\prime=(x,i,\gamma\gamma^\prime)$.

The description of a map in terms of its graph can be generalized to
the case of maps from $M_1/\cF_1$ to $M_2/\cF_2$ as follows. Let
$G_1$, $G_2$ be the holonomy groupoids of the foliations $\cF_1$,
$\cF_2$, respectively. The graph $G_f$ of a map $f: M_1/\cF_1 \to
M_2/\cF_2$ is a smooth (not necessarily Hausdorff) manifold equipped
with smooth maps $r_f : G_f\to M_1$ and $s_f: G_f\to M_2$, a left
action $G_1\times_{r_f}G_f\to G_f$ of the groupoid $G_1$ and a right
action $G_f\times_{s_f}G_2\to G_f$ of the groupoid $G_2$. It is
assumed that the actions commute, $r_f : G_f\to M_1$ is a principal
bundle with structure groupoid $G_2$ (that is, $r_f$ is a
submersion, and for any $x$ and $y$ from $G_f$ such that
$r_f(x)=r_f(y)$ there is a unique element $\gamma\in G_1$ with
$x\gamma=y$) and the action of the groupoid $G_2$ is proper.

Equivalently, a map $f : M_1/\cF_1\to M_2/\cF_2$ can be defined
either as a cocycle on $G_1$ with values in $G_2$, or as a
homomorphism $\varphi : (G_1)^{T_1}_{T_1}\to (G_2)^{T_2}_{T_2}$ of
the reduced groupoids, where the $T_j$ are complete transversals for
the foliations $(M_j,\cF_j), j=1,2$ (for more details, see
\cite{Co:survey,Co-skandal,Hil-Skan}).

The composition $f_2\circ f_1: M_1/\cF_1\to M_3/\cF_3$ of maps $f_1
: M_1/\cF_1\to M_2/\cF_2$ and $f_2 : M_2/\cF_2\to M_3/\cF_3$ is
given by the graph $G_{f_2\circ f_1}=G_{f_1}\times_{G_2} G_{f_2}$,
which is defined as the set of equivalence sets on the set $
\{(x_1,x_2)\in G_{f_1}\times G_{f_2} : s_{f_1}(x_1)=r_{f_2}(x_2)\} $
with respect to the equivalence relation, specifying that
$(x_1,x_2)\sim (y_1,y_2)$ if there exists $\gamma\in G_2$ such that
$x_1\gamma=y_1$, $\gamma y_2=x_2$.

\begin{ex}
Let $\cF_i$ be a simple foliation defined by a submersion $\pi_i :
M_i \to B_i$, $i=1,2$. Any map $F: B_1\to B_2$ defines a map $f:
M_1/\cF_1\to M_2/\cF_2$. Its graph $G_f$ consists of all $(x,y)\in
M_1\times M_2$ such that $\pi_2(y)=F(\pi_1(x))$. The maps
$s_f:G_f\to M_2$ and $r_f:G_f\to M_1$ are given by $s_f(x,y)=x$,
$r_f(x,y)=y$.
\end{ex}

\begin{ex}
The graph of the identity map $\id : M/\cF\to M/\cF$ is the holonomy
groupoid $G$.
\end{ex}

\begin{ex}
The graph of the projection $\id : M\to M/\cF$ is also the holonomy
groupoid $G$. More generally, the graphs of a map $f : M_1/\cF_1\to
M_2/\cF_2$ and of the corresponding lift $\bar{f} : M_1\to
M_2/\cF_2$ coincide as sets and differ only by the fact that the
graph of $f$ has the structure of a left $G_1$-manifold, but the
graph of $\bar{f}$ does not.
\end{ex}

\begin{ex}
The graph of the projection $\id : M/\cF\to pt$ is $M$.
\end{ex}

Denote by ${\rm Ml}_n(\RR)$ the metalinear group, that is, the
non-trivial twofold covering of the group ${\rm Gl}^+_n(\RR)$ of
non-singular real $n\times n$ matrices with positive determinant.
Put ${\rm Ml}^c_n={\rm Ml}_n(\RR)\times_{\ZZ_2}S^1$. The maximal
compact subgroup of ${\rm Ml}^c_n$ is ${\rm Spin}^c(n)$. A holonomy
equivariant real vector bundle $E$ on $M$ of rank $r$ is said to be
$K$-orientable if its structure group (as a $G$-bundle) reduces to
${\rm Ml}^c_n$.

Let $f : M_1/\cF_1\to M_2/\cF_2$ and let $E$ be a holonomy
equivariant bundle on $M_2$. Then $s_f^*E$ is a $G_2$-bundle on
$G_f$ such that the action of $G_1$ on it is trivial. Since $G_f$ is
a principal $G_1$-bundle, there exists a unique $G_1$-bundle
$E^\prime$ on $M_1$ such that $r_f^*E^\prime\cong s_f^*E$. Let
$E^\prime=f^*E$.

A smooth map $f : M_1/\cF_1\to M_2/\cF_2$ is said to be
$K$-orientable if the $G_1$-bundle $\tau_1\oplus f^*\tau_2$ on $M_1$
is $K$-orientable.

For any smooth $K$-oriented map $f : M_1/\cF_1\to M_2/\cF_2$ one can
naturally define an element $f!\in
KK(C^*(M_1,\cF_1),C^*(M_2,\cF_2))$, generalizing the Gysin
homomorphism in $K$-theory (see Subsection~\ref{s:AS}). In the case
when the foliation $\cF_1$ is trivial, a construction of the element
$f!\in KK(C(M_1),C^*_r(M_2,\cF_2))$ was given in \cite{Co:survey}
and studied more systematically in \cite{Co-skandal}, where it was
also proved that $(f\circ g)!=g!\otimes f!$. In \cite{Baum-Connes} a
construction of $f!$ was used in the case when the holonomy groupoid
of $\cF_1$ is proper. Finally, in \cite{Hil-Skan} there is a
definition of $f!\in KK(C^*(M_1,\cF_1),C^*(M_2,\cF_2))$ for an
arbitrary smooth $K$-oriented map $f : M_1/\cF_1\to M_2/\cF_2$, and
it is proved that $(f\circ g)!=g!\otimes f!$. For a $K$-oriented
tangent bundle $T\cF$ the natural projection $p : M\to M/\cF$ is a
$K$-oriented map, and the element $p!\in KK(C(M),C^*_r(M,\cF))$
coincides with the element $[D]\in KK(C(M),C^*_r(M,\cF))$, defined
by the corresponding tangential ${\rm Spin}^c$ Dirac operator $D$
(see Subsection~\ref{s:Kth}). For the map $\pi:M/\cF\to pt$, the
construction of $\pi!\in KK(C^*(M,\cF),\CC)$ is directly connected
with a $K$-theoretic analogue of the construction of the transversal
fundamental class of a foliation (see Subsection~\ref{s:tfclass}).
The difficult problem of constructing the element $\pi!$ was solved
in \cite{Hil-Skan}. Its solution makes use of the lift to
para-Riemannian foliations, the Thom isomorphism, and transversally
hypoelliptic operators (see Subsections~\ref{para}
and~\ref{s:high-ind}, in particular, Theorem~\ref{t:HS88}). Finally,
in the case when a map $f : M_1/\cF_1\to M_2/\cF_2$ is an embedding,
the construction of $f!$ is analogous to the classical construction
of the Gysin homomorphism in $K$-theory, in particular, of the
topological index (see Subsection~\ref{s:AS}), and makes use of the
Thom isomorphism and of the construction of the normal groupoid
associated to a foliation.

Generators of the geometric $K$-theory group $K^*_{\rm top}(M,\cF)$
($K$-cocycles) are equivalence classes of quadruples
$(X,\cF_X,x,f)$, where $(X,\cF_X)$ is a foliation with proper
holonomy groupoid, $x\in K_*(C^*(X,\cF_X))$, and $f : X/\cF_X \to
M/\cF$ is a $K$-oriented map. The equivalence relation of
$K$-cocycles is given by
\[
(X,\cF_X,x,f\circ g)\sim (Y,\cF_Y,x\otimes g!,f),
\]
where $g : X/\cF_X \to Y/\cF_Y$ is a $K$-oriented map. The addition
in $K^*_{\rm top}(M,\cF)$ is given by the disjoint sum operation.

In \cite{Baum-Connes} (see also \cite{Co-skandal}), a map $\mu :
K^*_{\rm top}(M,\cF) \to K_*(C^*(M,\cF))$ is defined by
\[
\mu : (X,\cF_X,x,f)\to f_!(x)=x\otimes f!.
\]
The Baum-Connes conjecture asserts that the map $\mu$ is an
isomorphism. Composing $\mu$ with the map $K_*(C^*(M,\cF))\to
K_*(C^*_r(M,\cF))$, induced by the canonical projection
$C^*(M,\cF)\to C^*_r(M,\cF)$, one gets a map $\mu_r : K^*_{\rm
top}(M,\cF) \to K_*(C^*_r(M,\cF))$.

For an arbitrary Lie groupoid $G$ as for any small category the
classifying space $BG$ is defined. It can be constructed by a slight
modification of Milnor's classical construction of the classifying
space of a group (see, for instance, \cite{Haefliger72}). We
describe G. Segal's construction \cite{Segal68}. Consider the
simplicial set $NG$ such that the set of its $n$-simplices $NG_n$ is
defined as
\[
NG_n= G^{(n)}=\{(\gamma_1,\ldots,\gamma_{n})\in G^{n} :
s(\gamma_i)=r(\gamma_{i+1}), \quad i=1,\ldots,n-1\},
\]
for $n>1$ the boundary operators $\delta_j : NG_n\to NG_{n-1},
j=0,1,\ldots,n,$ have the form
\begin{align*}
\delta_0(\gamma_1,\ldots,\gamma_n)&=(\gamma_2,\ldots,\gamma_n),\\
\delta_j(\gamma_1,\ldots,\gamma_n)&=(\gamma_1,\ldots,\gamma_{j}\gamma_{j+1},\ldots,\gamma_n),
\quad 1\leq j\leq n-1,\\
\delta_n(\gamma_1,\ldots,\gamma_n)&=(\gamma_1,\ldots,\gamma_{n-1}),
\end{align*}
and for $n=0$ the boundary operators $\delta_j : NG_1\to NG_{0}=
G^{(0)}, j=0,1,$ have the form
\[
\delta_0(\gamma_1)=r(\gamma_1),\quad \delta_1(\gamma_1)
=s(\gamma_1).
\]
and the degeneration operators $s_j : NG_n\to NG_{n+1}$ have the
form
\[
s_j(\gamma_1,\ldots,\gamma_n)=(\gamma_1,\ldots,\gamma_j,
s(\gamma_j)=r(\gamma_{j+1}), \gamma_{j+1},\ldots, \gamma_n), \quad
0\leq j\leq n.
\]
The simplicial set $NG$ is called the nerve of the groupoid $G$. The
classifying space $BG$ of $G$ is defined as the geometrical
realization of the simplicial set $NG$.

As an important particular case, consider the groupoid $G_{\cU}$
associated with an arbitrary open cover $\cU=\{U_i\}_{i\in A}$ of a
topological space $X$. It is defined by the formulas
\begin{align*}
G_{\cU}^{(0)}=& \bigsqcup_{i\in A} U_i=\{(x,i)\in
M\times A : x\in U_i\},\\
 G_{\cU}=& \bigsqcup_{i,j\in A}
U_i\cap U_j =\{(x,i,j)\in M\times A \times A : x\in U_i\cap U_j\},
\end{align*}
the maps $s,r : G_{\cU}\to G_{\cU}^{(0)}$ have the form
\[
s(x,i,j)=(x,i), \quad r(x,i,j)=(x,j).
\]

The set of $n$-simplices of the nerve $NG_{\cU}$ of this groupoid is
described as
\[
(NG_{\cU})_n =\{(x,i_1,\ldots,i_n) : x\in U_{i_1} \cap\ldots \cap
U_{i_n}\}= \bigsqcup U_{i_1}\cap\ldots\cap U_{i_n},
\]
where the union is taken over all sets $(i_1,\ldots,i_n)$ such that
$U_{i_1}\cap\ldots\cap U_{i_n}$ is non-empty. If a cover $\cU$ is
locally finite and all non-empty intersections
$U_{i_1}\cap\ldots\cap U_{i_n}$ are contractible, then the
classifying space $BG_{\cU}$ is homotopy equivalent to $X$.

A smooth map $f: X\to M/\cF$ given by a cocycle $(U_i, \gamma_{ij})$
defines a morphism from the groupoid $G_{\cU}$ associated with the
cover $\{U_i\}$ to the holonomy groupoid $G$ of the foliation
$(M,\cF)$. This morphism takes any $(x,i)\in G_{\cU}^{(0)}$ to
$\gamma_{ii}(x)\in M$ and any $(x,i,j)\in G_{\cU}$ to
$\gamma_{ij}(x)\in G$. The induced map $BG_{\cU}\to BG$ of the
classifying spaces is well defined. If the cover $\cU$ is locally
finite and all non-empty intersections $U_{i_1}\cap\ldots\cap
U_{i_n}$ are contractible, then this construction gives rise to a
map $\bar{f} : X \to BG$ defined up to homotopy equivalence. The map
$\bar{f}$ is the classifying map for the principal $G$-bundle $r_f :
G_f\to X$ defined by the map $f: X\to M/\cF$.

Thus, the classifying space $BG$ of the holonomy groupoid of the
foliation is an analogue of the leaf space $M/\cF$ in homotopy
theory (the homotopy leaf space). There is also the universal
classifying space for all smooth codimension $q$ foliations. It is
the classifying space $B\Gamma_q$ of the Haefliger groupoid
$\Gamma_q$ (see Example~\ref{ex:Heafliger}) introduced by Haefliger
\cite{Haefliger72,Haefliger84}.

For any codimension $q$ foliation $\cF$ on a manifold $M$, there is
defined (up to homotopy equivalence) the universal map $BG\to
B\Gamma_q$. It is constructed as follows. Let $\cU=\{U_i\}$ be a
good cover of $M$ by foliated charts, $T=\bigcup_iT_i$ the
corresponding complete transversal. Since every local transversal
$T_i$ is diffeomorphic to $I^q$, a natural morphism $G^T_T\to
\Gamma_q$ of groupoids is well defined. The universal map $BG\to
B\Gamma_q$ is obtained as the composition of the induced map $BG^T_T
\to B\Gamma_q$ and the homotopy equivalence $BG\backsimeq BG^T_T$.

We refer the reader to the survey \cite{Hurder08} and its references
for various questions related to the classification problem for
foliations.

Let $(M,\cF)$ be a foliated manifold. The normal bundle $\tau$ on
$M$, being a holonomy equivariant bundle, defines a vector bundle on
$BG$ which will be also denoted by $\tau$. Let us consider the
$K$-homology group $K_{*,\tau}(BG)$ of $BG$ twisted by the bundle
$\tau$. Using Poincar\'e duality, one can show that elements of
$K_{*,\tau}(BG)$ are represented by geometric cycles $(X,E,f)$,
where $X$ is a smooth compact manifold, $E$ is a complex vector
bundles on $X$, and $f : X \to M/\cF$ is a $K$-oriented map.

There is defined a map
\[
K^*_{\rm top}(M,\cF) \to K_{*,\tau}(BG),
\]
which associated to any quadruple $(X,\cF_X,x,f)$ the set
$(X,E,\bar{f})$, where $E$ is a $G_X$-equivariant bundle on $X$
corresponding to the element $x\in K_*(C^*(X,\cF_X))$, and $\bar{f}
: X \to M/\cF$ is a natural lift of $f : X/\cF_X \to M/\cF$. This
map is rationally injective and, in the case when the isotropy
groups $G^x_x$ are torsion free, is an isomorphism.

Denote by $H_{*,\tau}(BG,\QQ)$ the singular homology of the pair
$(B\tau,S\tau)$. One has the Chern character $\ch : K_{*,\tau}(BG)
\to H_{*,\tau}(BG,\QQ)$ and the Thom isomorphism $\Phi :
H_{k+q,\tau}(BG,\QQ)\to H_{k}(BG,\QQ), q=\dim \tau$. Their
composition $\Phi\ch : K_{*,\tau}(BG) \to H_{*}(BG,\QQ)$ has the
following form \cite{Co86} (see also \cite{DGKY} and
(\ref{e:chern})): for any $y=[X,E,f]\in K_{*,\tau}(BG)$
\begin{equation}\label{e:Fch}
\Phi\ch (y)=\bar{f}_*(\ch(E)\cup \Td(TX\oplus \bar{f}^*\tau)\cap
[X]) \in H_{*}(BG,\QQ),
\end{equation}
where $\bar{f} : X\to BG$ is the map defined by the map $f: X\to
M/\cF$.

We refer the reader to the bibliography given in \cite{survey} for
various aspects of the Baum-Connes conjecture for foliations and
related computations of the $K$-theory for foliation $C^*$-algebras.

\subsection{Transverse integration}
The foundations of non-commutative integration theory for foliations
were laid by Connes in \cite{Co79} (see also \cite{M-S,CoLNP80,Co}).
Let $(M,\cF)$ be a compact foliated manifold. Let $\Lambda$ be a
holonomy quasi-invariant transverse measure for $\cF$, $\alpha$ a
strictly positive smooth leafwise density on $M$ and $\nu=s^*\alpha$
the corresponding smooth Haar system. The measure $\Lambda$ and
$\alpha$ enables one to construct a measure $\mu$ on $M$ (see
(\ref{def:mu})). Finally, the measure $\mu$ and the Haar system
$\nu$ define a measure $m$ on $G$:
\[
\int_Gf(\gamma)\,dm(\gamma)=\int_M\left(\int_{G^x}f(\gamma)\,d\nu^x(\gamma)
\right)\,d\mu(x), \quad f\in C_c(G).
\]
In Subsection~\ref{s:calgebra} we defined the representation $R_x$
of the involutive algebra $C^{\infty}_c(G)$ on the Hilbert space
$L^2(G^x,\nu^x)$ for any $x\in M$. Consider the representation $R$
of the algebra $C^{\infty}_c(G)$ in $L^2(G,m)=\int^\oplus_M
L^2(G^x,\nu^x)\,d\mu(x)$ defined as the direct integral of $R_x$:
\[
R=\int^\oplus_M R_x\,d\mu(x).
\]

\begin{defn}
The von Neumann algebra $W^{\ast}_\Lambda(M,\mathcal F)$ of a
foliation $\mathcal F$ is defined as the closure of the image of the
algebra $C^{\infty}_{c}(G)$ by the representation $R$ in the weak
topology of the space ${\mathcal L} (L^2(G,m))$.
\end{defn}

Because the definition of the von Neumann algebra
$W^{\ast}_\Lambda(M,\mathcal F)$ depends only on the measure class
of $m$ (that is, on the family of all sets of $m$-measure zero),
this definition is independent of the choice of $\alpha$.

Elements of the algebra $W^{\ast}_\Lambda(M,\mathcal F)$ can be
described as families $(T_x)_{x\in M}$ of bounded operators in
$H_x=L^2(G^x,\nu^x)$ such that:

(1) for any $\gamma\in G$, $\gamma : x\to y$,
\begin{equation}\label{e:G}
L(\gamma)T_xL(\gamma)^{-1}=T_y,
\end{equation}
where the representation $L$ is given by the formula (\ref{def:L});

(2) the function $M \ni x\mapsto \|T_x\|$ is essentially bounded
with respect to the measure $\mu$;

(3) for any $\xi, \eta \in L^2(G,m)$ the function $M \ni x \mapsto
\langle T_x(\xi_x),\eta_x\rangle_{H_x}$ is measurable.

Operator families satisfying the condition (\ref{e:G}) will be
called (left-in\-va\-ri\-ant) $G$-operators. We will also use
right-invariant $G$-operators, which are families $(T_x)_{x\in M}$
of bounded operators in $H_x=L^2(G_x,\nu_x)$, $\nu_x=r^*\alpha$,
such that $R(\gamma)T_xR(\gamma)^{-1}=T_y$ for any $\gamma\in G$,
$\gamma : x\to y$, where the representation $R$ is given by
(\ref{def:R}).

A holonomy invariant measure $\Lambda$ on $M$ defines a normal
semifinite trace ${\rm tr}_{\Lambda}$ on the von Neumann algebra
$W^{\ast}_\Lambda(M,{\mathcal F})$. For any bounded measurable
function $k$ on $G$ the value ${\rm tr}_{\Lambda}(k)$ is finite and
given by
\begin{equation*}
{\rm tr}_{\Lambda}(k) = \int_{M}k(x)\,d\mu(x).
\end{equation*}
In \cite{Co79} a description of weights on the von Neumann algebra
$W^{\ast}_\Lambda(M,{\mathcal F})$ is given. As explained in
\cite{CoLNP80}, the construction of \cite{Co79} can be interpreted
as a correspondence between weights on $W^{\ast}(M,{\mathcal F})$
and operator-valued densities on the leaf space $M/\cF$.

\begin{thm}
Let $(M,\cF)$ be a foliated manifold. The von Neumann algebra
$W^{\ast}_\Lambda(M,{\mathcal F})$ is a factor if and only if the
foliation is ergodic, i.e., any bounded measurable function constant
along the leaves of $\cF$ is constant on $M$.
\end{thm}

It is known that von Neumann algebras are classified according to
three classes: type I, II, and III. Any von Neumann algebra $\cM$
can be canonically represented as a direct sum $\cM_{\rm I}\oplus
\cM_{\rm II}\oplus \cM_{\rm III}$ of von Neumann algebras $\cM_{\rm
I}$, $\cM_{\rm II}$ and $\cM_{\rm III}$ of type I, II and III,
respectively.

\begin{thm}
Let $(M,\cF)$ be a foliated manifold. The von Neumann algebra
$W^{\ast}_\Lambda(M,{\mathcal F})$ is of:
\begin{enumerate}
  \item type {\rm I} if and only if the leaf space is isomorphic to a
  standard Borel space.
  \item type {\rm II} if and only if there exists a holonomy invariant
  transverse measure and the algebra is not of type {\rm I}.
  \item type {\rm III} if and only if there exists no holonomy invariant
  transverse measure.
\end{enumerate}
\end{thm}

\section{Non-commutative  differential calculus on the leaf space}

\subsection{The transverse de Rham complex of a foliation}\label{s:deRham}
Let $(M,\cF)$ be a foliated manifold. In this subsection we describe
a non-commutative analogue of the de Rham complex on the leaf space
$M/\cF$ \cite{Co,Sau}.

Consider the space $\Omega_\infty=C^\infty_c(G, r^*\Lambda
N^*\cF\otimes |T\cG|^{1/2})$. We define a product on it by the
formula
\[
\omega_1\ast \omega_2(\gamma)=\int_{\gamma_1\gamma_2=\gamma}
\omega_1(\gamma_1)\wedge H(\gamma_1)[\omega_2(\gamma_2)],\quad
\gamma\in G,
\]
where $\omega_1, \omega_2\in C^{\infty}_c (G, |T\cG|^{1/2})$, and
$H(\gamma): \Lambda N^*_x\cF\to \Lambda N^*_y\cF$ is the linear
holonomy map associated with $\gamma\in G$, $\gamma : x\to y$. Thus,
$\Omega_\infty$ is a graded algebra.

We now define the transverse de Rham differential as a linear map
$d_H: \Omega_\infty^0= C^\infty_c(G,|T\cG|^{1/2}) \to
\Omega_\infty^1=C^\infty_c(G, r^*N^*\cF\otimes |T\cG|^{1/2})$
satisfying the condition
\[
d_H(k_1\ast k_2) = d_Hk_1\ast k_2+ k_1\ast d_Hk_2, \quad k_1, k_2
\in C^\infty_c(G,|T\cG|^{1/2}).
\]
The construction of $d_H$ makes essential use of an auxiliary choice
of a distribution $H$ on $M$ transverse to $F=T\cF$.

There is a decomposition of $TM$ into the direct sum $TM=F\oplus H$.
Therefore, one has the corresponding bigrading of the exterior
bundle $\Lambda^*T^*M$: $$
\Lambda^kT^*M=\bigoplus_{i+j=k}\Lambda^{i,j}T^*M, \quad
\Lambda^{i,j}T^*M=\Lambda^iH^*\otimes \Lambda^jF^*, $$  and also
(see, for instance, \cite[Proposition 10.1]{BGV}) the corresponding
decomposition of the de Rham differential $d$ into a sum of bigraded
components of the form
\begin{equation}\label{e:d}
d=d_F+d_H+\theta.
\end{equation}
Here
\begin{enumerate}
\item  $ d_F=d_{0,1}: C^{\infty}(M,\Lambda^{i,j}T^*M)\to
C^{\infty}(M,\Lambda^{i,j+1}T^*M)$ is the tangential de Rham
differential, which is a first order tangentially elliptic operator
independent of the choice of $g$;
\item $d_H=d_{1,0}: C^{\infty}(M,\Lambda^{i,j}T^*M)\to
C^{\infty}(M,\Lambda^{i+1,j}T^*M)$ is the transversal de Rham
differential, which is a first order transversally elliptic
operator;
\item $\theta=d_{2,-1}: C^{\infty}(M,\Lambda^{i,j}T^*M)\to
C^{\infty}(M,\Lambda^{i+2,j-1}T^*M)$ is a zero-order differential
operator which is the operator of inner multiplication by the 2-form
$\theta$ on $M$ with values in $F$, $\theta \in
C^{\infty}(M,F\otimes \Lambda^2\tau^*)$, given by
$$ \theta(X,Y)=p_F([X,Y]),\quad X,Y\in C^{\infty}(M,H), $$ where
$P_F:TM\to F$ is the natural projection (the form $\theta$ is called
the curvature of the transverse distribution $H$; in particular,
$\theta$ vanishes if and only if $H$ is completely integrable).
\end{enumerate}

The transverse distribution $H$ defines naturally a transverse
distribution $HG\cong r^*H$ on the foliated manifold $(G,\cG)$ (see
Subsection~\ref{s:groupoid}) and the corresponding transversal de
Rham differential $d_H: C^\infty_c(G) \to C^\infty_c(G, r^*N^*\cF)$
(see (\ref{e:d})).

For an arbitrary smooth leafwise density $\lambda \in
C^\infty_c(M,|T\cF|)$ we define a 1-form $k(\lambda)\in
C^\infty(M,H^*)\cong C^\infty(M, N^*\cF)$ as follows. Take an
arbitrary point $m \in M$. In a foliated chart $\phi : U \to
I^p\times I^q$ defined in a neighborhood of $m$
($\phi(m)=(x^0,y^0)$) the density $\lambda $ can be written as $
\lambda=f(x,y)|dx|, (x,y)\in I^p\times I^q. $ Then
\[
k(\lambda)=d_Hf+\sum_{i=1}^p{\cL}_{\frac{\partial}{\partial
x_i}}d_Hx_i,
\]
where for any $X=\sum_{i=1}^pX^i\frac{\partial}{\partial x_i} \in
\cX(\cF)$ and for any $\omega =\sum_{j=1}^q \omega_j\, dy_j \in
C^\infty(M,H^*)$ the Lie derivative ${\cL}_X\omega \in
C^\infty(M,H^*)$ is given by
\[
{\cL}_X\omega=\sum_{i=1}^p\sum_{j=1}^q X^i \frac{\partial \omega_j
}{\partial x_i}\,dy_j.
\]
One can give a slightly different description of $k(\lambda)$. For
any $X\in H_{m}$ let $\tilde{X}$ be an arbitrary projectable vector
field, which coincides with $X$ at $m$:
\[
\tilde{X}(x,y)=\sum_{i=1}^pX^i(x,y)\frac{\partial}{\partial
x_i}+\sum_{j=1}^qY^j(y)\frac{\partial}{\partial y_j}.
\]
Put
\begin{multline*}
k(\lambda)(X)=\sum_{i=1}^pX^i(x^0,y^0)\frac{\partial f}{\partial
x_i}(x^0,y^0)\\ +\sum_{j=1}^qY^j(y^0)\frac{\partial f}{\partial
x_j}(x^0,y^0) +\sum_{i=1}^p \frac{\partial X^i}{\partial
x_i}(x^0,y^0) f(x^0,y^0).
\end{multline*}
It can be checked that this definition is independent of the choice
of a foliated chart $\phi : U \to I^p\times I^q$ and an extension
$\tilde{X}$.

If $M$ is Riemannian, $\lambda $ is given by the induced leafwise
Riemannian volume form and $H=F^\bot$, then $k(\lambda)$ coincides
with the mean curvature 1-form of $\cF$.

An arbitrary leafwise half-density $\rho \in
C^\infty(M,|T\cF|^{1/2})$ can be written as $\rho=f|\lambda|^{1/2}$
with $f\in C^\infty(M)$ and $\lambda\in C^\infty(M,|T\cF|)$. Then
$d_H\rho \in C^\infty_c(M, N^*\cF\otimes |T\cF|^{1/2})$ is defined
as
\[
d_H\rho=(d_Hf)|\lambda|^{1/2}+\frac{1}{2} f|\lambda|^{1/2}
k(\lambda).
\]
Any $f\in C^\infty_c(G,|T\cG|^{1/2})$ can be written as
$f=us^*(\rho)r^*(\rho)$, where $u\in C^\infty_c(G)$ and $\rho\in
C^\infty_c(M,|T\cF|^{1/2})$. The element $d_Hf\in C^\infty_c(G,
r^*N^*\cF\otimes |T\cG|^{1/2})$ is defined as
\[
d_Hf = d_Hus^*(\rho)r^*(\rho) + us^*(d_H\rho)r^*(\rho) +
us^*(\rho)r^*(d_H\rho).
\]

The operator $d_H$ has a unique extension to a differentiation of
the differential graded algebra $\Omega_\infty=C^\infty_c(G,
r^*\Lambda N^*\cF\otimes |T\cG|^{1/2})$. By definition, for any
$f\in C^\infty_c(G,|T\cG|^{1/2})$ and $\omega\in
C^\infty_c(M,\Lambda N^*\cF)$, one has
\[
d_H(fr^*\omega )=(d_Hf)r^*\omega + fr^*(d_H\omega ).
\]

We define a closed graded trace $\tau$ on the differential graded
algebra $(\Omega_\infty, d_H)$ by the formula
\[
\tau(\omega)=\int_M\left.\omega\right|_M,
\]
where $\omega\in \Omega_\infty^q=C^\infty_c(G, r^*\Lambda^q
N^*\cF\otimes |T\cG|^{1/2})$. Here $\left.\omega\right|_M$ denotes
the restriction of the form $\omega$ to $M$, which is a section of
the bundle $\Lambda^qN^*\cF\otimes |T\cF|$ on $M$. Since the
foliation is transversally oriented, the integral of
$\left.\omega\right|_M$ over $M$ is well defined,

\subsection{Transverse fundamental class of a foliation}\label{s:tfclass}
Let $(M,\cF)$ be a foliated manifold, and $H$ an arbitrary
distribution on $M$ transverse to $F=T\cF$. In this subsection we
describe, following \cite{Co86} (see also \cite{Co}), the simplest
construction of a cyclic cocycle on the algebra
$C^\infty_c(G,|T\cG|^{1/2})$, namely, the construction of the
transverse fundamental class.

In the previous subsection we constructed a graded algebra
$\Omega_\infty$, a differential $d_H$ and a closed graded trace
$\tau$ on it. The problem is that, since the distribution $H$ is
non-integrable, it is not true in general that $d^2_H=0$. Using the
equality (\ref{e:d}) and computing the type $(0,2)$ component in the
representation of the operator $d^2$ as a sum of bihomogeneous
components, we get that
\[
d^2_H=-(d_F\theta+\theta d_F).
\]
The operator $-(d_F\theta+\theta d_F)$ is a tangential differential
operator, so is given by exterior multiplication by a vector-valued
distribution $\Theta \in {\cD}'(G, r^*\Lambda^2 N^*\cF\otimes
|T\cG|^{1/2})$ supported in $G^{(0)}$. One can show that
\[
d^2_H\omega=\Theta \wedge \omega - \omega \wedge \Theta, \quad
\omega\in \Omega_\infty.
\]
Moreover, $d_H\Theta=0$. Using these facts, one can canonically
construct a differential graded algebra $(\tilde{\Omega}_\infty,
\tilde{d}_H)$ and a closed graded trace $\tilde{\tau}$ on it (see
\cite{Co:nc,Co}). The algebra $\tilde{\Omega}_\infty$ consists of
$2\times 2$ matrices $\omega=\{\omega_{ij}\}$ with entries from
$\Omega_\infty$. An element $\omega\in \tilde{\Omega}_\infty$ has
degree $k$ if $\omega_{11}\in \Omega_\infty^k$, $\omega_{12},
\omega_{21}\in \Omega_\infty^{k-1}$, and $\omega_{22}\in
\Omega_\infty^{k-2}$. The product in $\tilde{\Omega}_\infty$ is
given by
\[
\omega\cdot\omega'=
\begin{bmatrix}
\omega_{11} & \omega_{12} \\ \omega_{21} & \omega_{22}
\end{bmatrix}
\begin{bmatrix}
1 & 0 \\ 0 & \Theta
\end{bmatrix}
\begin{bmatrix}
\omega'_{11} & \omega'_{12} \\ \omega'_{21} & \omega'_{22}
\end{bmatrix},
\]
the differential by
\[
\tilde{d}_H\omega =\begin{bmatrix} d_H\omega_{11} & d_H\omega_{12}
\\ -d_H\omega_{21} & -d_H\omega_{22}
\end{bmatrix}
+
\begin{bmatrix}
0 & -\Theta \\ 1 & 0
\end{bmatrix}
\omega + (-1)^{|\omega|}\omega
\begin{bmatrix}
0 & 1 \\ -\Theta & 0
\end{bmatrix},
\]
and the closed graded trace $\tilde{\tau}:
\tilde{\Omega}_\infty^q\to \CC$ is defined by
\[
\tilde{\tau}\left(
\begin{bmatrix}
\omega_{11} & \omega_{12} \\ \omega_{21} & \omega_{22}
\end{bmatrix}
\right)=\tau(\omega_{11})-(-1)^q\tau(\omega_{22}\Theta).
\]
Finally, the homomorphism $\tilde{\rho}: C^\infty_c(G, |T\cG|^{1/2})
\to \tilde{\Omega}^0_\infty$ is given by
\[
\tilde{\rho}(k)=\begin{bmatrix} k & 0 \\ 0 & 0
\end{bmatrix}.
\]

Thus, the triple $(\tilde{\Omega}_\infty, \tilde{d}_H,
\tilde{\tau})$ is a cycle over the algebra $C^\infty_c(G,
|T\cG|^{1/2})$. It is called the fundamental cycle of the
transversally oriented foliation $(M,\cF)$. The character of this
cycle defines a cyclic cocycle $\phi_H$ on the algebra
$C^\infty_c(G, |T\cG|^{1/2})$. The cocycle $\phi_H$ depends on the
auxiliary choice of a horizontal distribution $H$, but the
corresponding cyclic cohomology class is independent of this choice.
The class $[M/\cF]\in HC^q(C^\infty_c (G, |T\cG|^{1/2}))$ defined by
the cocycle $\phi_H$ is called the transverse fundamental class of
the foliation $(M,\cF)$.

Let $C^\infty_c(G, |T\cG|^{1/2})^+$ be got by adjoining a unit to
the algebra $C^\infty_c(G, |T\cG|^{1/2})$. For an even $q$ let us
extend the cycle $(\tilde{\Omega}_\infty, d_H, \tilde{\tau})$ to a
cycle over $C^\infty_c(G, |T\cG|^{1/2})^+$ by putting
$\tilde{\tau}(\Theta^{q/2})=0$. It is proved in \cite{Gor-CMP} that
the formula
\begin{multline*}
\chi^r(k^0,k^1,\ldots,k^r)\\
=\frac{(-1)^{(q-r)/2}}{\left((q+r)/2\right)!}
\sum_{i_0+\ldots+i_r=(q-r)/2}\int_M
k^0\,\Theta^{i_0}d_Hk^1\,\Theta^{i_1} \ldots d_Hk^r\Theta^{i_r},
\end{multline*}
where $r=q,q-2,\ldots$, $k^0 \in C^\infty_c(G, |T\cG|^{1/2})^+$,
$k^1,\ldots, k^r\in C^\infty_c(G, |T\cG|^{1/2})$ defines a cocycle
in the $(b,B)$-complex of the algebra $C^\infty_c(G,
|T\cG|^{1/2})^+$. The class defined by the cocycle $\chi$ in $HC^q
(C^\infty_c(G, |T\cG|^{1/2}))$ coincides with $[M/\cF]$.

The pairing with the class $[M/\cF]\in HC^q(C^\infty_c(G,
|T\cG|^{1/2}))$ defines an additive map $K(C^\infty_c(G,
|T\cG|^{1/2})) \to \CC$. An important problem is the question of
topological invariance of this map, that is, the question whether it
is possible to extend it to an additive map from $K(C^*_r(G))$ to
$\CC$. This problem was solved in \cite{Co86}.

A standard method of solving the problem of topological invariance
of the map $\varphi$ consists in constructing a smooth subalgebra
$\cB$ in $C^*_r(M,\cF)$, which contains the algebra $C^\infty_c(G,
|T\cG|^{1/2})$ and is such that the cyclic cocycle $\phi_H$ on
$C^\infty_c(G, |T\cG|^{1/2})$, which defines the transverse
fundamental class $[M/\cF]$, extends by continuity to a cyclic
cocycle on $\cB$ and thereby defines a map in the topological
$K$-theory $K(\cB)\cong K(C^*_r(G))\to\CC$. This has been done for
so called para-Riemannian foliations (see Subsection~\ref{para}) by
using some properties of densely defined cyclic cocycles on Banach
algebras.

For an arbitrary foliated manifold $(M,\cF)$ a bundle $P$ over $M$
is constructed in \cite{Co86} whose fibres are connected spin
manifolds of non-positive curvature, and then a natural lift of the
foliation $\cF$ to a para-Riemannian foliation $\cV$ on $P$ (see a
more detailed exposition in Subsection~\ref{para}). Since $\cV$ is
para-Riemannian, its transverse fundamental class defines a map
$K(C^*_r(P,\cV))\to\CC$. On the other hand, since the fibres of $P$
are connected spin manifolds of non-positive curvature, there is an
injective map $K(C^*_r(M,\cF)) \to K(C^*_r(P,\cV))$, making it
possible to construct the desired extension $K(C^*_r(M,\cF))\to\CC$
for the initial foliation $(M,\cF)$. For geometric consequences of
this construction, see Subsection~\ref{s:high-ind}, and also
\cite{Co86}.

Actually, if the tangent bundle $T\cF$ is $K$-oriented, then, for
any element $P$ of the subring of the ring $H^*(M,\RR)$ generated by
Pontryagin classes of the normal bundle $\tau$ and Chern classes of
arbitrary holonomy equivariant bundles on $E$, an additive map
$\varphi_P : K(C^*_r(M,\cF)) \to \CC$ is constructed in \cite{Co86}
such that for any $E\in K^i(M)$, $i=\dim M \mod 2$,
\begin{equation}\label{e:muP}
\varphi_P(E\otimes p!)=\langle \ch (E)\cdot {\rm Td}(TM\otimes
\CC)\cdot P, [M]\rangle,
\end{equation}
where $p!\in KK(C(M),C^*_r(M,\cF))$ is the element associated with
the natural projection $p:M\to M/\cF$ (see a more precise statement
in Theorem~\ref{t:Co86}).

We recall that the $K$-homology fundamental class of a compact ${\rm
Spin}^c$-manifold $M$ is defined to be the class $[D]\in K_i(M)$,
where $i=\dim M\mod 2$, determined by an arbitrary ${\rm Spin}^c$
Dirac operator $D$ on $M$. If the foliation $(M,\cF)$ is Riemannian,
the tangent bundle $T\cF$ is $K$-oriented, and the normal bundle has
a holonomy invariant complex spin structure, then the corresponding
transverse ${\rm Spin}^c$ Dirac operator $D^\pitchfork$ defines a
$K$-cohomology class $[D^\pitchfork]\in K^i(C^*(M,\cF))$, where
$i=\dim (M/\cF) \mod 2$ (see Subsection~\ref{s:steo}), and therefore
a map $K_i(C^*_r(M,\cF))\to \CC$. It is proved in \cite[Theorem
4.2]{DGKY} that for any class $[E]\in K^0(M)$ given by a bundle $E$
on $M$ one has
\[
(E\otimes p!)\otimes [D^\pitchfork]=\Ind (D_E),
\]
where $D_E$ is the ${\rm Spin}^c$ Dirac operator on $M$ with
coefficients in $E$. In particular, by the Atiyah-Singer index
theorem,
\begin{equation}\label{e:Dpitch}
(E\otimes p!)\otimes [D^\pitchfork]=\langle \ch (E)\cdot {\rm
Td}(TM\otimes \CC), [M]\rangle.
\end{equation}
Thus, the right hand side of (\ref{e:Dpitch}) coincides with the
right hand side of~(\ref{e:muP}) with $P={\rm Td}(\tau)$. In the
general case, the class $[D^\pitchfork]\in K^i(C^*(M,\cF))$ is not
well defined, but the construction described above makes it possible
to construct the corresponding map $\varphi_P : K(C^*_r(M,\cF)) \to
\CC$. These arguments may serve as a justification of the name
``transverse fundamental class'' for a map of the form $\varphi_P$.

\begin{defn}
A transverse current is a current $C$ (that is, a continuous linear
functional on the space of smooth, compactly supported differential
forms) defined on the disjoint union of all transversals to the
foliation.
\end{defn}

\begin{defn}
A transverse current $C$ is said to be holonomy invariant, if, for
any transversals $B_1$ and $B_2$ and any map $\phi :B_1\to B_2$,
belonging to the holonomy pseudogroup, $\phi_* (C_{B_1})=C_{B_2}$.
\end{defn}

Let $C$ be a closed, holonomy invariant transverse current of degree
$k$. Define a continuous linear functional $\rho_C$ on the space
$C^\infty_c(M,\Lambda N^*\cF\otimes |T\cF|)$ as follows.

Take a good cover $\{U_i\}$ of the manifold $M$ by foliated
coordinate neighborhoods with the corresponding coordinates maps
$\phi_i: U_i\to I^p \times I^q$ and a partition of unity
$\{\psi_i\}$ subordinate to this covering. We consider the
corresponding complete transversal $T=\bigcup_i T_i$, where
$T_i=\phi^{-1}_i(\{0\}\times I^q).$ In any foliated chart
$(U_i,\phi_i)$ the transverse current $C$ defines a current $C_i$ on
$T_i$. The current $C$ is holonomy invariant if and only if
$f_{ij}(C_i)=C_j$ for any pair of indices $i$ and $j$ such that
${U}_i\cap {U}_j\not=\varnothing$.

For any section $\omega\in C^\infty_c(M,\Lambda^k N^*\cF\otimes
|T\cF|)$ and any $i$ the expression
\[
\int_L \Psi_i\omega=\int_{P_i(y)} \psi_i(x,y)\omega(x,y)
\]
gives a well defined differential $k$-form on $T_i$. Put
\[
\rho_C(\omega)=\sum_i \left\langle C_i, \int_L \Psi_i\omega
\right\rangle.
\]
One can show that the functional $\rho_C$ is well defined by this
formula, that is, the result is independent of the choice of a
covering $\{U_i\}$ and a partition of unity $\{\psi_i\}$. The
formula
\[
\tau_C(\omega)=\rho_C(\left.\omega\right|_M), \quad \omega\in
C^\infty_c(G, r^*\Lambda^k N^*\cF\otimes |T\cG|^{1/2}),
\]
gives a closed graded trace $\tau_C$ on the differential graded
algebra
\[
(\Omega_\infty=C^\infty_c(G, r^*\Lambda N^*\cF\otimes |T\cG|^{1/2}),
d_H).
\]
Thus, any closed holonomy invariant transverse current of degree $k$
defines a cyclic cocycle on the algebra $C^\infty_c(G,
|T\cG|^{1/2})$. Topological invariance of these cocycles follows
from the results of the paper \cite{Co86} (see \cite[Remark
5.12]{Benameur-Heitsch04}).

\subsection{The cyclic cocycle defined by the Godbillon-Vey class}
\label{s:GV} We consider a smooth compact manifold $M$ equipped with
a transversally oriented codimension one foliation $\cF$. The
Godbillon-Vey class is a 3-dimensional cohomology class of $M$. It
is the simplest example of secondary characteristic classes of the
foliation. We recall its definition. Since $\cF$ is transversally
oriented, it is globally defined by a non-vanishing smooth 1-form
$\omega$ (that is, $\ker \omega_x=T_x\cF$ for any $x\in M$). It
follows from the Frobenius theorem that there exists a 1-form
$\alpha$ on $M$ such that $d\omega=\alpha\wedge \omega$. One can
check that the 3-form $\alpha\wedge d\alpha$ is closed, and its
cohomology class does not depend on the choice of $\omega$ and
$\alpha$. This class $GV\in H^3(M,\RR)$ is called the Godbillon-Vey
class of $\cF$.

Let $T$ be a complete smooth transversal given by a good cover of
$M$ by regular foliated charts. Thus, $T$ is an oriented
one-dimensional manifold. In this subsection we construct a cyclic
cocycle on the algebra $C^\infty_c(G^T_T)$, corresponding to the
Godbillon-Vey class. This is done in several steps. To start with,
we describe the construction of the Godbillon-Vey class as a
secondary characteristic class associated with the cohomology
$H^*(W_1,\RR)$ of the Lie algebra $W_1=\RR[[x]]\partial_x$ of formal
vector fields on $\RR$.

The cohomology $H^*(W_1,\RR)$ was computed by Gel'fand and Fuchs
(for instance, see the book \cite{Fuks} and its references). They
are finite-dimensional, and the only non-trivial groups are
$H^0(W_1,\RR)=\RR\cdot 1$ and $H^3(W_1,\RR)=\RR\cdot gv$, where
\[
gv(p_1\partial_x, p_2\partial_x, p_3\partial_x)=
\begin{vmatrix}
p_1(0) & p_2(0) & p_3(0) \\ p'_1(0) & p'_2(0) & p'_3(0)\\ p''_1(0) &
p''_2(0) & p''_3(0)
\end{vmatrix}
, \quad p_1, p_2, p_3 \in \RR[[x]].
\]

Consider the bundle $F^k_+T \rightarrow T$ of positively oriented
frames of order $k$ on $T$ and the bundle
$F^\infty_+T=\underset{\leftarrow}{\lim} F^k_+T$ of positively
oriented frames of infinite order on $T$. By definition, a
positively oriented frame $r$ of order $k$ at $x\in T$ is the
$k$-jet at $0\in \RR$ of an orientation preserving diffeomorphism
$f$ which maps a neighborhood of $0$ in $\RR$ onto some neighborhood
of $x=f(0)$ in $T$. If $y:U\to\RR$ is a local coordinate on $T$
defined in a neighborhood $U$ of $x$, then the numbers
\[
y_0=y(x),\quad y_1=\frac{d (y\circ f)}{dt}(0),\quad \ldots,\quad
y_k=\frac{d^k(y\circ f)}{dt^k}(0),
\]
are coordinates of the frame $r$, and moreover, $y_1>0$.

There is a natural action of the pseudogroup $\Gamma^+(T)$ of
orientation preserving local diffeomorphisms of $T$ on
$F^\infty_+T$. Let $\Omega^*(F^\infty_+T)^{\Gamma^+(T)}$ denote the
space of differential forms on $F^\infty_+T$ invariant under the
action of $\Gamma^+(T)$. There is a natural isomorphism $J:
C^*(W_1)\to \Omega^* (F^\infty_+T)^{\Gamma^+(T)}$ of differential
algebras defined in the following way. First of all, let $v\in W_1$,
and let $h_t$ be any one-parameter group of local diffeomorphisms of
$\RR$ such that $v$ is the $\infty$-jet of the vector field
$\left.\frac{dh_t}{dt}\right|_{t=0}$. Then we define a
$\Gamma^+(T)$-invariant vector field on $F^\infty_+T$ whose value at
$r=j_0f\in J^\infty_+(T)$ is given by
\[
v(r)=j_0 \left(\left.\frac{d(f\circ h_t)}{dt}\right|_{t=0}\right).
\]

For any $c\in C^q(W_1)$ put
\[
J(c)(v_1(r),\ldots,v_q(r))=c(v_1,\ldots,v_q).
\]
One can check that this isomorphism takes the cocycle $gv \in
C^3(W_1,\RR)$ to the three-form
\[
gv =\frac{1}{y_1^3}dy\wedge dy_1\wedge dy_2\in
\Omega^3(F^2_+T)^{\Gamma^+(T)}.
\]
Consider the bundle $F^\infty(M/\cF)$ on $M$, consisting of infinite
order jets of all possible distinguished maps. There is a natural
map $F^\infty(M/\cF) \to F^\infty_+T/\Gamma^+(T)$. Using this map
and the $\Gamma^+(T)$-invariance of $gv \in \Omega^3(F^2_+T)$, one
can lift $gv$ to a closed form $gv(\cF)\in
\Omega^3(F^\infty(M/\cF))$. Since the fibres of the fibration
$F^2(M/\cF)\to M$ are contractible, so
$H^3(F^\infty(M/\cF),\RR)\cong H^3(M,\RR)$, and the cohomology class
in $H^3(M,\RR)$ corresponding to the cohomology class of $gv(\cF)$
in $H^3(F^\infty(M/\cF),\RR)$ under this isomorphism coincides with
the God\-bill\-ion-Vey class of $\cF$.

Let $C^*(G^T_T, \Omega^*(G^T_T))$ denote the space of cochains on
$G^T_T$ with values in the space $\Omega^*(G^T_T)$ of differential
forms on $G^T_T$, By a Van Est type theorem (see
\cite{Haefliger:CIME}) there is an embedding
\begin{equation}\label{e:vanest}
 \Omega^*(F^\infty_+T)^{\Gamma^+(T)}\to C^*(G^T_T,
\Omega^*(G^T_T)).
\end{equation}
This map is a homomorphism of complexes which induces an isomorphism
in cohomology.

Let $\rho$ be an arbitrary smooth positive density on $T$. We define
a homomorphism $\delta : G^T_T\to \RR^*_+$, called the modular
homomorphism by setting $\delta=r^*\rho/s^*\rho$, where $r^*\rho$
(respectively, $s^*\rho$) denotes the lift of the density $\rho$ to
$G^T_T$ by the map $r: G^T_T\to T$ (respectively, $s: G^T_T\to T$),
and also the homomorphism $\ell=\log \delta : G^T_T\to \RR$. The
formula
\[
c(\gamma_1, \gamma_2) =\ell(\gamma_2)\, d\ell (\gamma_1)
-\ell(\gamma_1)\, d\ell(\gamma_2), \quad \gamma_1,\gamma_2\in G^T_T,
\]
defines a 2-cocycle on $G^T_T$ with values in the space of 1-forms
on $G^T_T$. This cocycle, called the Bott-Thurston cocycle,
corresponds to $gv \in \Omega^3(F^2_+T)^{\Gamma^+(T)}$ under the
isomorphism given by the embedding (\ref{e:vanest}).

The Bott-Thurston cocycle $c$ defines a cyclic 2-cocycle $\psi$ on
$C^\infty_c(G^T_T)$ by the formula (cf. Example~\ref{ex:phi}):
\begin{equation}\label{e:GV-cyclic}
\psi(k^0,k^1,k^2)=\int_{\gamma_0\gamma_1\gamma_2\in T}
k^0(\gamma_0)k^1(\gamma_1)k^2(\gamma_2)c(\gamma_1,\gamma_2), \quad
k^0,k^1,k^2 \in C^\infty_c(G^T_T).
\end{equation}
This is the cyclic cocycle corresponding to the Godbillon-Vey class
of $\cF$.

Connes proved the topological invariance of this cocycle, that is,
the fact that the additive map $\varphi : K(C^\infty_c(G^T_T))\to
\CC$ defined by the pairing with the class $[\psi]\in
HC^2(C^\infty_c(G^T_T))$ of $\psi$, $\phi(e)=\langle e, [\psi]
\rangle$ determines a map $\varphi : K(C^*_r(G^T_T)) \cong
K(C^*_r(M,\cF))\to \CC$. Moreover, one has the formula (see also the
formula (\ref{e:mur}))
\begin{equation}\label{e:muGV}
\varphi(\mu_r(x))=\langle\Phi\ch(x), GV \rangle, \quad x\in K^*_{\rm
top}(M,\cF).
\end{equation}

For the cyclic cocycle corresponding to the God\-bill\-i\-on-Vey
class of $\cF$ there is another description which connects it with
invariants of the von Neumann algebra of this foliation \cite{Co86}.
The formula
\[
\tau(k^0,k^1)=\int_{G^T_T}k^0(\gamma^{-1})dk^1(\gamma),\quad k^0,
k^1 \in C^\infty_c(G^T_T)
\]
defines a cyclic $1$-cocycle on $C^\infty_c(G^T_T)$. The class of
$\tau$ in $HC^1(C^\infty_c(G^T_T))$ corresponds to the transverse
fundamental class of $\cF$ in $HC^1(C^\infty_c(G))$ under the
isomorphism $HC^1(C^\infty_c(G^T_T)) \cong HC^1(C^\infty_c(G))$
defined by strong Morita equivalence (see
subsection~\ref{s:morita1}).

The fixed smooth positive density $\rho$ on $T$ defines a faithful
normal semifinite weight $\phi_\rho$ on the von Neumann algebra
$W^*(G^T_T)$ of the groupoid $G^T_T$. For any $k\in
C^\infty_c(G^T_T)$ the value of the weight $\phi_\rho$ is given by
\[
\phi_\rho(k)=\int_Tk(x)\,\rho(x).
\]
Let us consider the one-parameter group $\sigma_t$ of automorphisms
of the von Neumann algebra $W^*(G^T_T)$ given by
\[
\sigma_t(k)(\gamma)=\delta(\gamma)^{it}k(\gamma), \quad k\in
C^\infty_c(G^T_T), \quad t\in \RR.
\]
This group is the group of modular automorphisms associated with the
weight by the Tomita-Takesaki theory.

The significance of the group of modular automorphisms is explained,
in particular, by the following characterization of it: a
one-parameter group of $\ast$-auto\-mor\-phi\-sms $\sigma_t$ of a
von Neumann algebra $M$ is the group of modular automorphisms
associated with a weight $\omega$ if and only if $\omega$ satisfies
the Kubo-Martin-Schwinger conditions with respect to $\sigma_t$,
that is, there is a function $f$ analytic in the strip
$\operatorname{Im} z\in (0,1)$ and continuous in its closure such
that, for any $a,b\in M, t\in\RR$
\[
f(t)=\omega(\sigma_t(a)b), \quad f(t+i)=\omega(b\sigma_t(a)).
\]

Following Connes \cite{Co86}, we define a $1$-trace on a Banach
algebra $B$ to be a bilinear functional $\phi$ defined on a dense
subalgebra $\cA\subset B$ such that
\begin{enumerate}
  \item $\phi$ is a cyclic cocycle on $\cA$;
  \item for any $a^1\in\cA$, there is a constant $C>0$ such that
\[
|\phi(a^0,a^1)|\leq C\|a^0\|, \quad a^0\in \cA;
\]
\end{enumerate}
and we define a $2$-trace on $B$ to be a trilinear functional $\phi$
defined on a dense subalgebra $\cA\subset B$ such that
\begin{enumerate}
  \item $\phi$ is a cyclic cocycle on $\cA$;
  \item for any $a^1,a^2\in\cA$, there is a constant $C>0$ such
  that
\[
|\phi(x^0,a^1x^1,a^2)-\phi(x^0a^1,x^1,a^2)|\leq C\|a^1\|\,\|a^2\|,
\quad x^0, x^1\in \cA.
\]
\end{enumerate}
The formula
\[
\dot{\tau}(k^0,k^1)=\lim_{t\to 0}\frac{1}{t}
\left(\tau(\sigma_t(k^0),\sigma_t(k^1)) -\tau(k^0,k^1) \right),
\quad k^0, k^1\in C^\infty_c(G^T_T)
\]
defines a 1-trace on $C^*_r(G^T_T)$ with domain $C^\infty_c(G^T_T)$
and invariant under the action of the automorphism group $\sigma_t$.

For a $C^*$-algebra $A$ and for any $1$-trace $\phi$ on it that is
invariant under an action of a one-parameter automorphism group
$\alpha_t$ with generator $D$ such that the space $\dom \phi\cap\dom
D$ is dense in $A$ one can define a $2$-trace $\chi=i_D\phi$ on
$C^*_r(G^T_T)$ (an analogue of the contraction) by
\[
\chi(a^0,a^1,a^2)=\phi(D(a^2)a^0,a^1)-\phi(a^0D(a^1),a^2), \quad
a^0,a^1,a^2 \in \dom \phi\cap\dom D.
\]

\begin{thm}\label{t:GV}
Suppose that $(M,\cF)$ is a manifold with a transversally oriented
codimension one foliation, $T$ is a complete smooth transversal, and
$\rho$ is a smooth positive density on $T$. Then the cyclic cocycle
$\psi\in HC^2(C^\infty_c(G^T_T))$ corresponding to the Godbillon-Vey
class of $\cF$ coincides with $i_D\dot{\tau}$.
\end{thm}

One can naturally associate to any von Neumann algebra $M$ an action
(called the flow of weights \cite{Connes-Tak}) of the multiplicative
group $\RR^*_+$ on a certain commutative von Neumann algebra,
namely, the center of the crossed product $M\rtimes \RR$ of $M$ by
$\RR$ relative to the action of $\RR$ on $M$ given by the modular
automorphism group $\sigma_t$. As a consequence of Theorem
\ref{t:GV}, Connes established the following geometric fact.

\begin{thm}\cite{Co86}
Suppose that $(M,\cF)$ is a manifold with a transversally oriented
codimension one foliation. If the Godbillon-Vey class $GV\in
H^3(M,\RR)$ does not vanish, then the flow of weights of the von
Neumann algebra of the foliation $\cF$ has a finite invariant
measure.
\end{thm}

In particular, this implies the following, earlier result.

\begin{thm}\cite{Hurder-KatokBAMS}
Suppose that $(M,\cF)$ is a manifold with a transversally oriented
codimension one foliation. If the Godbillon-Vey class $GV\in
H^3(M,\RR)$ does not vanish, then the von Neumann algebra of the
foliation $\cF$ has a non-trivial type {\rm III} component.
\end{thm}

In \cite{Mor,Mor-Natsume} there is another construction of the
cyclic cocycle associated with the Godbillon-Vey class, as a cyclic
cocycle on the $C^*$-algebra of foliation $C^*_r(M,\cF)$ in the case
when the foliation $\cF$ is the horizontal foliation of a flat
foliated $S^1$-bundle.

\subsection{General constructions of cyclic cocycles}\label{s:second}
Let $G$ be a smooth \'etale groupoid (for instance, the reduced
holonomy groupoid $G^T_T$ of a foliated manifold $(M,\cF)$
associated with a complete transversal $T$). The tangent bundle
$\tau$ on $G^{(0)}$, being a $G$-bundle, defines a vector bundle on
$BG$, which will be also denoted by $\tau$. Consider the cohomology
group $H^{*}_\tau(BG)$ of the space $BG$ twisted by $\tau$:
$H^{*}_\tau(BG)=H^*(B\tau, S\tau)$. Connes \cite[Chapter III,
Section 2 $\delta$, Theorem 14 and Remark b)]{Co} constructed a
natural map
\begin{equation}\label{e:Phi*}
\Phi_*: H^*_\tau(BG)\to
HP^*(C^\infty_c(G)).
\end{equation}
The constructions of the transverse fundamental class of a foliation
and of the cyclic cocycle associated with the Godbillon-Vey class
are particular cases of this general construction.

Let us start a description of the construction of the map $\Phi_*$
with a definition of a bicomplex $(C^{*,*},d_1,d_2)$. For $n
> 0$ and $-q\leq m\leq 0$ ($q=\dim G^{(0)}$) the space $C^{n,m}$
consists of de Rham currents of degree $-m$ on the manifold
\[
G^{(n)}=\{(\gamma_1,\ldots,\gamma_n)\in G^n :
s(\gamma_i)=r(\gamma_{i+1}), \quad i=1,\ldots,n-1\},
\]
which vanish if either some $\gamma_j$ belongs to $G^{(0)}$ or
$\gamma_1\ldots\gamma_n\in G^{(0)}$. The space $C^{0,m}, -q\leq
m\leq 0,$ consists of de Rham currents of degree $-m$ on the
manifold $G^{(0)}$. Otherwise, put $C^{n,m}=\{0\}$. The coboundary
$d_1 : C^{n,m}\to C^{n+1,m}$ is given by
\[
d_1=(-1)^m\sum (-1)^j\delta^*_j,
\]
where the pull-back maps $\delta^*_j : \cD_{-m}(G^{(n-1)}) \to
\cD_{-m}(G^{(n)})$ are induced by the \'etale maps $\delta_j :
G^{(n)}\to G^{(n-1)}$ defined for $n>1$ by
\begin{align*}
\delta_0(\gamma_1,\ldots,\gamma_n)&=(\gamma_2,\ldots,\gamma_n),\\
\delta_j(\gamma_1,\ldots,\gamma_n)&=(\gamma_1,\ldots,\gamma_j\gamma_{j+1},\ldots,\gamma_n),
\quad 1\leq j\leq n-1,\\
\delta_n(\gamma_1,\ldots,\gamma_n)&=(\gamma_1,\ldots,\gamma_{n-1}).
\end{align*}
For $n=1$ the maps $\delta_j : G\to G^{(0)}, j=0,1$, are defined by
\[
\delta_0(\gamma_1)=r(\gamma_1), \quad
\delta_1(\gamma_1)=s(\gamma_1), \quad \gamma_1\in G.
\]
The de Rham boundary $d^t : \cD_{-m}(G^{(n)}) \to
\cD_{-m-1}(G^{(n)})$ gives the coboundary $d_2 : C^{n,m}\to
C^{n,m+1}$. The $k$th cohomology group of the complex
$(C^{*},d=d_1+d_2)$ associated with the bicomplex
$(C^{*,*},d_1,d_2)$ coincides with $H^{k+q}_\tau(BG)$.

We introduce a bicomplex
$(\Omega^{*,*}_c(G),d^\prime,d^{\prime\prime})$. The space
$\Omega^{n,m}_c(G)$ is defined as the quotient of the space of
smooth compactly supported differential forms of degree $m$ on
$G^{(n+1)}$ by the subspace of forms supported in the set of all
$(\gamma_0,\ldots,\gamma_n)$ such that $\gamma_j\in G^{(0)}$ for
some $j>0$. The product in $\Omega^*_c(G)$ is given by
\begin{multline*}
(\omega_1\omega_2)(\gamma_0,\ldots,\gamma_{n_1},\ldots,\gamma_{n_1+n_2})\\
\begin{aligned}
&=\sum_{\gamma\gamma^\prime=\gamma_{n_1}}\omega_1(\gamma_0,\ldots,\gamma_{n_1-1},\gamma)\wedge
\omega_2(\gamma^\prime,\gamma_{n_1+1},\ldots,\gamma_{n_1+n_2})\\
&+\sum_{j=0}^{n_1-2}(-1)^{n_1-j-1}\sum_{\gamma\gamma^\prime=\gamma_{j}}
\omega_1(\gamma_0,\ldots,\gamma_{j-1},\gamma,\gamma^{\prime},\ldots,\gamma_{n_1-1})\wedge
\omega_2(\gamma_{n_1},\ldots,\gamma_{n_1+n_2}),
\end{aligned}
\end{multline*}
where $\omega_1\in \Omega^{n_1,m_1}_c(G)$, $\omega_2\in
\Omega^{n_2,m_2}_c(G)$. In this formula the fact that the maps $r$
and $s$ are \'etale is used to identify the cotangent spaces.

The differential $d^{\prime} : \Omega^{n,m}_c(G)\to
\Omega^{n+1,m}_c(G)$ is given by
\[
d^\prime\omega(\gamma_0,\ldots,\gamma_{n+1})=\chi_{G^{(0)}}(\gamma_0)\omega(\gamma_1,\ldots,\gamma_{n+1}),
\]
where $\chi_{G^{(0)}}\in C^\infty(G)$ is the indicator function of
the set $G^{(0)}$, and the differential $d^{\prime\prime} :
\Omega^{n,m}_c(G)\to \Omega^{n,m+1}_c(G)$ is given by the usual de
Rham differential.

Let $(\Omega^*_c(G), d=d^\prime+d^{\prime\prime})$ be the
differential graded algebra associated with the bicomplex
$(\Omega^{*,*}_c(G), d^\prime,d^{\prime\prime})$.

For an arbitrary cochain $c\in C^{n,m}$ in the bicomplex
$(C^{*,*},d_1,d_2)$ its push-forward by the map
\[
(\gamma_1,\ldots,\gamma_n)\in G^{(n)} \mapsto
((\gamma_1\ldots\gamma_n)^{-1}, \gamma_1,\ldots,\gamma_n)\in
G^{(n+1)},
\]
if $n>0$ and by the natural embedding $G^{(0)}\to G$ if $n=0$
defines a de Rham current of degree $-m$ on $G^{(n+1)}$. Denote by
$\tilde{c}$ the corresponding linear functional on
$\Omega^{n,m}_c(G)$.

The morphism $\Phi$ from the bicomplex $(C^{*,*},d_1,d_2)$ to the
$(b,B)$-bicomplex of the algebra $\cA=C^\infty_c(G)$ is defined as
follows. For any $c\in C^{n,m}$ the corresponding element
$\Phi(c)\in C^{n-m}(\cA,\cA^*)$ is the $(n-m+1)$-linear functional
on $\cA$ given by the formula (with $\ell=n-m+1$)
\begin{multline*}
\Phi(c)(a^0,\ldots,a^\ell)\\ =\lambda_{n,m}\sum_{j=0}^\ell
(-1)^{j(\ell-j)}\tilde{c}(da^{j+1}\ldots da^\ell a^0 da^1\ldots
da^j),\quad a^0,\ldots,a^\ell\in\cA,
\end{multline*}
where $\lambda_{n,m}=\frac{n!}{(\ell+1)!}$. Here we regard the
algebra $\cA$ as a subalgebra of the algebra $\Omega^{0,0}_c(G)$,
and the product and the differential $d$ are taken in the algebra
$\Omega^*_c(G)$.

\begin{ex}
If $G$ is the trivial groupoid associated with a manifold $M$ (see
Example~\ref{ex:trivial}), then any closed current $C$ of degree $k$
on $M$ defines a cocycle $c$ in the complex $(C^*,d)$. This cocycle
has the sole non-zero component $c_{0,-k}=C\in C^{0,-k}$. The
corresponding cyclic cocycle $\Phi(c)$ on the algebra
$C^\infty_c(M)$ coincides with the cocycle given by the
formula~(\ref{e:C}).
\end{ex}

\begin{ex}
More generally, if $(M,\cF)$ is a compact foliated manifold and
$G=G^T_T$ is its reduced holonomy groupoid associated with a
complete transversal $T$ defined by a good cover $\cU=\{U_i\}$, then
any closed holonomy invariant transversal current $C$ of degree $k$
on $M$ defines a cocycle $c$ in the complex $(C^*,d)$. This cocycle
has the sole non-zero component $c_{0,-k}=C\in C^{0,-k}$. The class
in $HC^*(C^\infty_c(G^T_T))$ of the cyclic cocycle $\Phi(c)$ on the
algebra $C^\infty_c(G^T_T)$ corresponds under the isomorphism
defined by Morita equivalence to the class in $HC^*(C^\infty_c(G))$
of the cyclic cocycle on $C^\infty_c(G)$ given by the current $C$
(see Subsection~\ref{s:tfclass}).
\end{ex}

\begin{ex}
This example is a generalization of Examples \ref{ex:C} and
\ref{ex:phi}. It also served as a motivation for the construction of
$\Phi$. Let $M$ be an $n$-dimensional oriented manifold, and
$\Gamma$ a discrete group acting on $M$ by orientation preserving
diffeomorphisms. Consider a $k$-cocycle $\omega$ on $\Gamma$ with
coefficients in the $\Gamma$-module $\Omega^n(M)$ of smooth
differential $n$-forms on $M$, $\omega\in Z^k(\Gamma, \Omega^n(M))$,
such that $\omega(g_1,\ldots,g_k)=0$ in the case when either $g_i =
e$ for some $i$ or $g_1\ldots g_k = e$. As shown in \cite[Lemma
7.1]{Co86}, the following equality defines a cyclic $k$-cocycle
$\tau$ on the algebra $C^\infty_c(M\times \Gamma)\subset
C(M)\rtimes_r \Gamma$ (see Example~\ref{ex:crossed}):
\begin{multline*}
\tau_\omega(f_0,\ldots,f_k)\\ =\sum_{\substack{g_0,\ldots,g_k\in\Gamma\\
g_0\cdot \ldots\cdot g_k=e}} \int_M f_0(x,g_0)f_1(xg_0,g_1)\ldots
f_k(xg_0g_1\ldots g_{k-1},g_k) \omega(g_1,\ldots,g_k)(x).
\end{multline*}
The crossed product groupoid $G=M\times \Gamma$ is an \'etale
groupoid. One can check that the $k$-cocycle $\omega$ defines an
element $\omega\in C^{k,0}$ and that $\Phi(\omega)=\tau_\omega$. For
more details we refer the reader to \cite[Chapter III, Section 2
$\delta$]{Co}.

The construction of the cyclic cocycle, corresponding to the
Godbillon-Vey class of a foliation (see the formula
(\ref{e:GV-cyclic})) is actually a particular case of this
construction (the only point is that the role of the discrete group
$\Gamma$ is played by the reduced holonomy groupoid of the
foliation).
\end{ex}

Let $G$ be a Hausdorff smooth \'etale groupoid. Consider the space
of loops in $G$:
\[
B^{(0)}=\{\gamma\in G : r(\gamma)=s(\gamma)\}.
\]
We say that a subset $\mathcal O\subset B^{(0)}$ is invariant under
the action of $G$ if for any $\gamma\in {\mathcal O}$ and $g\in G$
such that $s(g)=r(\gamma)$ one has an inclusion $g\gamma g^{-1}\in
{\mathcal O}$.

For any invariant subset $\mathcal O\subset B^{(0)}$ one can define
the localized cyclic cohomology $HC^n(C^\infty_c(G))_{{\mathcal
O}}$. Moreover, if $B^{(0)}=\bigsqcup_\alpha {\mathcal O}_\alpha$ is
the representation of $B^{(0)}$ as the disjoint union of open
invariant subsets, then one has the direct sum decomposition
\cite{Brylinsky:N}
\[
HC^n(C^\infty_c(G))=\bigoplus_\alpha HC^n(C^\infty_c(G))_{{\mathcal
O}_\alpha}.
\]
A subset ${\mathcal O}\subset B^{(0)}$ is said to be elliptic if it
is invariant and the order of each element $\gamma\in {\mathcal O}$
is finite. A subset ${\mathcal O}\subset B^{(0)}$ is said to be
hyperbolic if it is invariant and the order of each element
$\gamma\in {\mathcal O}$ is infinite.

One can take the set of units ${\mathcal O}= G^{(0)}$ for an open
and closed subset ${\mathcal O}\subset B^{(0)}$. The localizations
on this set are usually denoted by the subscript $[1]$ instead of
${\mathcal O}$. The morphism $\Phi$ described above provides a
description of the corresponding component
$HP^*(C^\infty_c(G))_{[1]}$ in the periodic cyclic cohomology
$HP^n(C^\infty_c(G))$. Namely, one has the isomorphism
\[
HP^{\rm ev/odd}(C^\infty_c(G))_{[1]}\cong \bigoplus_{k\ {\rm
even/odd}} H^{k+q}_\tau(BG).
\]
A similar description in terms of homologies of some double
complexes is obtained in \cite{Brylinsky:N} for an arbitrary
elliptic component ${\mathcal O}\subset B^{(0)}$. The computation of
the localized cyclic cohomology for hyperbolic components is more
complicated. It uses in a greater extent the combinatorics of the
groupoid.

\section{The index theory of transversally elliptic operators}

\subsection{Equivariant transversally elliptic operators and their distributional index}
Transversally elliptic operators appeared for the first time in the
papers \cite{At:teo,Singer:recent} in the following situation. Let
$M$ be a smooth compact manifold and $G$ a compact Lie group acting
on $M$. Denote by $\mathfrak g$ the Lie algebra of $G$. For any
$X\in \mathfrak g$ denote by $X_M$ the corresponding fundamental
vector field on $M$. Vectors of the form $X_M(x)$ with
$X\in\mathfrak g$ span the tangent space $T_x(G.x)$ to the orbit
passing through $x$. Consider the space
\[
(T^*_GM)_x  =\{\xi\in T^*_xM :\langle \xi, X_M(x)\rangle =0 \
\text{for any}\ X\in \mathfrak g\}.
\]

\begin{defn}
A classical pseudodifferential operator $P$ from $C^\infty(M,\cE^+)$
to $C^\infty(M,\cE^-)$ of order $m$ acting on sections of
$G$-equivariant vector bundles $\cE^\pm$ on $M$ is said to be
$G$-transversally elliptic if it commutes with the action of $G$ and
its principal symbol $\sigma(P)(x,\xi): \pi^*\cE^+\to \pi^*\cE^-$ is
invertible for any $(x,\xi)\in T^*_GM\setminus 0$.
\end{defn}

The choice of a $G$-invariant Riemannian metric $g$ on $M$ and
$G$-invariant Hermitian structures on $\cE^{\pm}$ defines Hilbert
structures in the spaces $L^2(M,\cE^{\pm})$. Let us regard a
$G$-transversally elliptic operator $P\in \Psi^m(M,\cE^+,\cE^-)$ as
a closed unbounded operator from $L^2(M,\cE^{+})$ to
$L^2(M,\cE^{-})$ obtained as its closure from the initial domain
$C^\infty(M,\cE^+)$. Then the kernels $\Ker P$ and $\Ker P^*$ of the
operators $P$ and $P^*$ in $L^2(M,\cE^{+})$ and $L^2(M,\cE^{-})$,
respectively, are $G$-invariant closed subspaces.

For any $G$-equivariant vector bundle $\cE$ on $M$ denote by $T(g)$
the induced action of $g\in G$ in $L^2(M,\cE)$. For any function
$\phi\in C^\infty_c(G)$ define the operator $T(\phi)$ in
$L^2(M,\cE)$ by
\[
T(\phi)=\int_G\phi(g)T(g)dg.
\]

The $G$-equivariant index $\Ind^G (P)$ of $P$ is the $G$-invariant
distribution on $G$ defined by
\[
\langle\Ind^G (P),\phi\rangle= \tr T(\phi)\Pi_{\Ker P} - \tr
T(\phi)\Pi_{\Ker P^*}, \quad \phi\in C^\infty_c(G),
\]
where $\Pi_{\Ker P}$ and $\Pi_{\Ker P^*}$ are the orthogonal
projections to $\Ker P$ and $\Ker P^*$, respectively.

The principal symbol $\sigma(P)$ of a $G$-transversally elliptic
operator $P$ defines an element $[\sigma(P)]$ of the equivariant
$K$-theory $K_G(T^*_GM)$ with compact supports of the space
$T^*_GM$. The $G$-equivariant index $ \Ind^G (P) \in
{\cD}^\prime(G)^G$ of $P$ depends only on the class $[\sigma(P)] \in
K_G(T^*_GM)$. Thus, the $G$-equivariant index induces a homomorphism
of $R(G)$-modules (the analytic index)
\[
\Ind^{G}_a : K_G(T^*_GM)\to {\cD}^\prime(G)^G.
\]

In  \cite{At:teo} an algorithm is given for computing the index of a
$G$-transversally elliptic operator. Using this algorithm, Berline
and Vergne \cite{BerlineVergne96a,BerlineVergne96b} obtained an
explicit cohomological index formula (see also
\cite{Paradan-Vergne08} and the reference therein for recent
advances in this direction).

Let $M$ be a smooth compact manifold and $G$ a compact Lie group
acting on $M$. For a fixed $s\in G$ consider the submanifold
$M(s)=\{x\in M : sx=x\}$ of fixed points of the action of $s$ on
$M$. Denote by $\cN=\cN(M,M(s))$ the normal bundle to $M(s)$ in $M$.
Let $G(s)$ denote the centralizer of $s$, $G(s)=\{t\in G : st=ts\}$,
and ${\mathfrak g}(s)$ its Lie algebra. The action of $G$ induces an
action of $G(s)$ on $M(s)$, and the bundle $\cN=\cN(M,M(s))$ is a
$G(s)$-equivariant vector bundle on $M(s)$. One has the formula
\begin{equation}\label{e:decomp00}
T_xM\cong T_xM(s)\oplus \cN_x,\quad x\in M(s),
\end{equation}
where $T_xM(s)=\{v\in T_xM : sv=v\}$, $\cN_x=(1-s)T_xM$.

The Levi-Civita connection $\nabla^M$ preserves the
decomposition~(\ref{e:decomp00}), which defines a representation
$\nabla^M=\nabla^0\oplus \nabla^1$, where $\nabla^0$ is a
$G(s)$-invariant connection on $TM(s)$, which coincides with the
Levi-Civita connection on $M(s)$, $\nabla^1$ is a $G(s)$-invariant
connection on $\cN$ compatible with the induced metric on $\cN$. A
similar representation holds for the Riemannian curvature form $R^M$
on $TM$:
\[
\left.R^M\right|_{M(s)}\cong R^0\oplus R^1,
\]
where $R^0$ and $R^1$ are the curvatures of the connections
$\nabla^0$ and $\nabla^1$, respectively.

Denote by $\Omega^\infty_G({\mathfrak g},M)=C^\infty({\mathfrak g},
\Omega(M))^G$ the space of $G$-invariant smooth maps from $\mathfrak
g$ to the space $\Omega(M)$ of smooth differential forms on $M$.
Elements of the space $\Omega^\infty_G({\mathfrak g},M)$ will be
called $G$-equivariant forms on $M$. The equivariant differential
$d_\fg: \Omega^\infty_G({\mathfrak g},M) \to
\Omega^\infty_G({\mathfrak g},M)$ is defined by the formula
\[
(d_\fg\alpha)(X)=d(\alpha(X))-\iota(X_M)(\alpha(X)), \quad X\in\fg,
\]
where $\iota(X_M)$ denotes the inner product by $X_M$. We say that a
$G$-equivariant form $\alpha$ on $M$ is equivariantly closed if
$d_\fg\alpha=0$.

Let $\cE$ be a $G$-equivariant vector bundle endowed with a
$G$-invariant connection $\nabla$. The covariant derivative $\nabla$
defines a direct sum decomposition $T\cE=V\cE\oplus H\cE$, where
$H\cE$ is the horizontal space of the connection, and the vertical
space $V\cE$ is isomorphic to $\pi^*\cE$. We define the moment map
$\mu : \fg \to \End(\cE)$ as follows. For any $m\in M$ and $v\in
\cE_m$, the vector $(\mu(Y)v)_m\in \cE_m$ is the projection of the
vector $(Y_\cE)_{v}\in T_{v}\cE$ to $V\cE_{v}\cong \cE_m$. In the
case when $\nabla$ is the Levi-Civita connection on $TM$, the
corresponding moment map  $\mu^M(Y)\in C^\infty(M,so(TM)), Y\in\fg$,
called the Riemannian moment of the manifold $M$, is given by the
formula
\[
\mu^M(Y)Z=-\nabla_ZY_M, \quad Z\in TM.
\]

The equivariant curvature of the connection $\nabla$ is defined by
\[
R(Y)=\mu(Y)+R,
\]
where $R=\nabla^2\in \Omega^2(M,\End \cE)$ is the curvature of the
connection $\nabla$ and $\mu(Y):\cE\to\cE$ is the associated moment
map.

Let $R^0(Y)$ and $R^1(Y)$ be the equivariant curvatures of the
respective connections $\nabla^0$ and $\nabla^1$, and define the
$G(s)$-equivariant forms $J(M(s))$ and $D_s(\cN(M,M(s)))$ on $M(s)$
by
\begin{gather*}
J(M(s))(Y) =\det \frac{e^{R_0(Y)/2}-e^{-R_0(Y)/2}}{R_0(Y)}, \\
D_s(\cN(M,M(s)))(Y) =\det(1-s\, e^{R_1(Y)}).
\end{gather*}
for any $Y\in \fg(s)$. The forms $J(M(s))$ and $D_s(\cN(M,M(s)))$
are equivariantly closed. Moreover, for any $Y$ in a sufficiently
small neighborhood of $0$ in $\fg(s)$ the form $J(M(s))(Y)$ is
invertible in $\Omega(M)$.

Suppose that $P: C^\infty(M,\cE^+) \to C^\infty(M,\cE^-)$ is a
$G$-equivariant transversally elliptic operator acting on sections
of $G$-equivariant vector bundles $\cE^\pm$ on $M$, and let $\sigma
: \pi^*\cE^+\to \pi^*\cE^- $ be its principal symbol. Choose
arbitrary $G$-invariant connections $\nabla^{\cE^\pm}$ on $\cE^\pm$
and define an odd endomorphism $U(\sigma)$ of the $\ZZ_2$-graded
bundle $\pi^*\cE=\pi^*\cE^+ \oplus \pi^*\cE^-$ by the formula
\[
U(\sigma)(x,\xi)=\begin{pmatrix} 0 & \sigma(x,\xi)^* \\
\sigma(x,\xi) & 0 \end{pmatrix}.
\]
Let $\A(\sigma)$ be the superconnection on $\pi^*\cE$ given by
\[
\A(\sigma)=i U(\sigma)+\pi^*\nabla,
\]
or, equivalently,
\[
\A(\sigma)=\begin{pmatrix} \pi^*\nabla^{\cE^+} & i\sigma^* \\
i\sigma & \pi^*\nabla^{\cE^-}
\end{pmatrix}.
\]

For any $Y\in \fg$ the action of the equivariant curvature
$F^\A(Y)\in \Omega(M,\End\ \cE)$ of the superconnection $\A(\sigma)$
on $\Omega(M,\cE)$ has the form
\[
F^\A(Y)=(\A-i(Y_M))^2+\cL^\cE(Y),
\]
where $\cL^\cE(Y)$ is the first order differential operator defined
by the Lie derivative of the action of $G$ on $\Omega(M,\cE)$.

We define the equivariant Chern character
$\operatorname{ch}(\A(\sigma))\in \Omega^\infty_G({\mathfrak g},M)$
by
\[
\operatorname{ch}(\A(\sigma))(Y)=\operatorname{tr}_s \exp F^\A(Y),
\quad Y\in \fg,
\]
and, for any $s\in G$, the equivariant Chern character
$\operatorname{ch}_s(\A(\sigma))\in \Omega^\infty_{G(s)}({\mathfrak
g}(s),M(s))$ by
\[
\operatorname{ch}_s(\A(\sigma))(Y)=\operatorname{tr}_s s^\cE \exp
\left.F^\A(Y)\right|_{M(s)}, \quad Y\in \fg(s),
\]
where $s^\cE$ denotes the action in the fibres of the bundle
$\left.\cE\right|_{M(s)}$ defined by the action of $s$.

Define an equivariant extension of the canonical symplectic form
$\Omega$ on $T^*M$ by
\[
\Omega(Y)=-d_Y\omega^M,\quad Y\in \fg,
\]
where $d_Y$ denotes the operator $d-\iota(Y_M)$ on $\Omega(M)$, or
\[
\Omega(Y)(x,\xi)=\Omega(x,\xi)+\langle\xi, Y_M(x)\rangle,\quad
(x,\xi)\in T^*M.
\]

Finally, we will need the notion of decent of distributions. If $N$
is a submanifold of the group $G$ transverse to $G$-orbits, the
restriction of an arbitrary distribution $\Theta\in \cD'(G)^G$ to
$N$ is well defined as a distribution on $N$. If $U_s(0)$ is a
sufficiently small neighborhood of $0$ in $\fg(s)$, then the
submanifold $s\exp U_s(0)$ is transverse to the orbits of the
adjoint action of $G$ on $G$, and use of the above result allows us
to speak about the restriction of the $G$-invariant distribution
$\Ind^G(P)\in \cD^\prime(G)$ to $s\exp U_s(0)$.

\begin{thm}
Let $s\in G$. For any $Y\in U_s(0)$, one has the equality
\[
\left.\Ind^G(P)\right|_{s\exp U_s(0)}(Y) =
\int_{T^*M(s)}\frac{1}{(2i\pi)^{n_s}} \frac{e^{-id_Y\omega^M}
\operatorname{ch}_s(\A(\sigma))(Y)}{D_s(\cN(M,M(s)))(Y) J(M(s))(Y)}.
\]
\end{thm}

In this formula $n_s=\dim M(s)$, and its right-hand side
$\theta^s\in \cD'(s\exp U_s(0))$ should be understood in the
following way: if $\phi\in C^\infty_c(s\exp U_s(0))$, then
\[
I_\phi=\int_{\fg(s)} \frac{e^{-id_Y\omega^M}
\operatorname{ch}_s(\A(\sigma))(Y)}{D_s(\cN(M,M(s)))(Y) J(M(s))(Y)}
\phi(Y)dY
\]
is rapidly decreasing along the fibres of the bundle $T^*M(s)$ and
\[
\langle \theta^s, \phi \rangle\stackrel{\rm def}{=} \int_{T^*M(s)}
\frac{1}{(2i\pi)^{n_s}} I_\phi.
\]

The index theorem for transversally elliptic operators in the
$K$-theoretic form was proved very recently by Kasparov
\cite{Kas-teo}. In \cite{Yu} the index theorem for transversallly
elliptic operators was proved in the case when the isotropy groups
are finite. In \cite{Kawasaki81} the results of \cite{At:teo} were
applied to prove the index theorem on orbifolds. See the papers
\cite{Fox-Haskell94a,Fox-Haskell94b,Paradan01,Paradan03} for
applications of the index theory in the representation theory of Lie
groups. In \cite{MMS08} the authors connected the equivariant index
of transversally elliptic operators with the fractional analytic
index of projectively elliptic operators associated with an Azumai
bundle and with the twisted $K$-theory (see \cite{MMS06}).

The existence of the index of a transversally elliptic operator in
the case when the group $G$ is non-compact was proved by Hoermander
\cite{At:teo} (see also \cite{Nestke-Z,trans}). In this case the
simplest examples show that the $G$-equivariant index $\Ind^G (P)
\in {\cD}^\prime(G)^G$ of a $G$-transversally elliptic operator $P$
is not a homotopy invariant of its principal symbol $\sigma(P)$.
Therefore, the analytic index map is not well defined. In fact, this
was one of the main motivations for Connes to introduce cyclic
cohomology. The results on computing the index of a transversally
elliptic operator for non-compact Lie groups will be given in the
next subsection.

\subsection{Transversally elliptic operators on foliations}
Let $(M,\cF)$ be a compact foliated manifold and $E$ a Hermitian
vector bundle on $M$. We recall that the principal symbol $p_m$ of a
classical pseudodifferential operator $P\in \Psi ^{m}(M,E)$ is a
smooth section of the bundle $\Hom(\pi^*E)$ on $T^*M\setminus
\{0\}$. The transversal principal symbol $\sigma_{P}$ of $P$ is the
restriction of its principal symbol $p_m$ to
$\widetilde{N}^{\ast}{\mathcal F}=N^{\ast}{\mathcal F}\setminus
\{0\}$. Thus, $\sigma _{P}$ is a smooth section of the bundle
$\Hom(\pi^*E)$ on $\widetilde{N}^{\ast}{\mathcal F}$. An operator
$P\in \Psi^m(M,E)$ is transversally elliptic if its transverse
principal symbol $\sigma _{P}(\nu)$ is invertible for any $\nu\in
\widetilde{N}^*{\mathcal F}$.

Suppose now that $E$ is a holonomy equivariant Hermitian vector
bundle on $M$. Thus, there is given a representation $T$ of the
holonomy groupoid $G$ of the foliation $\cF$ on the fibres of $E$,
that is, for any $\gamma\in G, \gamma:x\rightarrow y$, a linear
operator $T(\gamma):E_x\rightarrow E_y$ is defined. The transversal
principal symbol $\sigma _{P}$ of $P\in \Psi ^{m}(M,E)$ is holonomy
invariant if, for any $\gamma\in G$, $\gamma : x \to y$, one has the
equality
\[
T(\gamma)\circ [\sigma_{P}(dh_{\gamma}^{\ast}(\nu ))] = \sigma
_{P}(\nu )\circ T(\gamma),\quad \nu  \in N^{ \ast}_{y}{\cF},
\]
where we have used the isomorphisms
$(\pi^*E)_{dh_{\gamma}^{\ast}(\nu )}\cong E_x$ and
$(\pi^*E)_\nu\cong E_y$.

Examples of transversally elliptic operators are transverse Dirac
operators (see \cite{matrix-egorov, vanishing} and the references
therein).

Let $M$ be a compact manifold endowed with a Riemannian foliation
$\cF$ of even codimension $q$ and $g_M$ a bundle-like metric on $M$,
and denote by $T^HM$ the orthogonal complement of $T{\cF}$. Let
$\nabla$ be the transverse Levi-Civita connection in $T^HM$.

For any $x\in M$ denote by ${\rm Cl}(Q_x)$ the Clifford algebra of
the Euclidean space $Q_x$. We define a $\ZZ_2$-graded vector bundle
${\rm Cl}(Q)$ on $M$ whose fibre at $x\in M$ is ${\rm Cl}(Q_x)$.
This bundle is associated with the principal $SO(q)$-bundle $O(Q)$
of oriented orthonormal frames in $Q$, ${\rm
Cl}(Q)=O(Q)\times_{O(q)}{\rm Cl}(\RR^q)$. Therefore, the transverse
Levi-Civita connection $\nabla$ induces a natural leafwise flat
connection $\nabla^{{\rm Cl}(Q)}$ in ${\rm Cl}(Q)$, which is
compatible with the Clifford multiplication and preserves the
$\ZZ_2$-grading on ${\rm Cl}(Q)$.

A complex vector bundle $\cE$ on $M$ endowed with an action of the
bundle ${\rm Cl}(Q)$ is called a transverse Clifford module. The
action of an element $a\in C^\infty(M,{\rm Cl}(Q))$ on an $s\in
C^\infty(M,\cE)$ will be denoted by $c(a)s\in C^\infty(M,\cE)$. A
transverse Clifford module $\cE$ is said to be self-adjoint if it is
endowed with a Hermitian metric such that the operator $c(f) :
\cE_x\to \cE_x$ is skew-symmetric for any $x\in M$ and $f\in Q_x$. A
transverse Clifford module $\cE$ has a natural $\ZZ_2$-grading
$\cE=\cE_+\oplus\cE_-$. A connection $\nabla^\cE$ on a transverse
Clifford module $\cE$ is called a Clifford connection if for any
$f\in C^\infty(M,T^HM)$ and $a\in C^\infty(M,{\rm Cl}(Q))$,
\[
[\nabla^{\cE}_f, c(a)]=c(\nabla^{{\rm Cl}(Q)}_fa).
\]
A self-adjoint transverse Clifford module equipped with a Hermitian
Clifford connection is called a transverse Clifford bundle.

Denote by $\tau \in C^\infty(M,T^HM)$ the mean curvature vector of
$\cF$. If $e_1,e_2,\ldots,e_p$ is a local orthonormal base in
$T\cF$, then
\[
\tau=\sum_{i=1}^pP_H(\nabla^L_{e_i}e_i).
\]

Let $\cE$ be a transverse Clifford bundle on $M$ equipped with a
Hermitian Clifford connection $\nabla^\cE$. Let $f_1,\ldots,f_q$ be
a local orthonormal base in $T^HM$. The transverse Dirac operator
$D_\cE$ is defined by
\[
D_\cE= \sum_{\alpha=1}^q
c(f_\alpha)\left(\nabla^{\cE}_{f_\alpha}-\frac12 g_M(\tau, f_\alpha)
\right).
\]
The operator $D_\cE$ is formally self-adjoint in $L^2(M, \cE)$. We
also observe that the operator $D_\cE$ has a holonomy invariant
principal symbol, as follows from the fact that the metric $g_M$ is
bundle-like.

Another example of a transversally elliptic operator is the
transversal de Rham operator acting in the space
$C^{\infty}(M,\Lambda T^HM^{*})$ by the formula
\[
D_H=d_H + d^*_H,
\]
where $d_H$ is the transverse de Rham differential (see
(\ref{e:d})). If the foliation admits a transversal spin structure,
then the operator $D_H$ is connected with the transverse Dirac
operator $D_{F(Q)\otimes F(Q)^*}$, where $F(Q)$ is the associated
transversal spinor bundle, as follows:
\[
D_{F(Q)\otimes F(Q)^*}= D_H-\frac12(\varepsilon_{\tau^*}+i_{\tau}).
\]
Thus, these operators coincide if and only if $\tau=0$, that is, all
the leaves are minimal submanifolds \cite{matrix-egorov}.

In \cite{Lazarov00} formulae for the distributional index of
invariant transversally elliptic operators are obtained in the case
when the group $G=\RR$ acts locally free and isometrically on a
compact Riemannian manifold $M$. In this case the action defines a
non-singular isometric flow $\phi_s: M\to M$, and its orbits define
a foliation $\cF$.

Let $E$ and $F$ be Hermitian vector bundles on $M$. Suppose that
there exists an action of $\RR$ in the fibres of the bundles $E$ and
$F$ that preserves the Hermitian structure and covers the action of
$\RR$ on $M$. Let $D : C^\infty(M,E)\to C^\infty(M,F)$ be an
arbitrary $\RR$-invariant, first-order, transversally elliptic
operator.

We recall that a periodic orbit of the flow $\phi$ is said to be
non-degenerate if $1$ is not an eigenvalue of the associated
Poincar\'e map. Denote by $X$ the generator of the flow. Then there
is a natural action of the flow in the normal bundle $Q$:
\[
d\phi_s(x): Q_x=T_xM/\RR X(x)\to Q_{\phi_s(x)}=T_{\phi_s(x)}M/\RR
X(\phi_s(x)).
\]
A periodic orbit $c$ with a (not necessarily minimal) period $l$ of
the flow $\phi$ is said to be non-degenerate (or simple) if
\[
\det(\id - d\phi_{l}(x) : Q_x\to Q_x)\neq 0,
\]
where $x\in c$ is an arbitrary point on $c$, $\phi_{l}(x)=x$. In
this case, put
\[
b_{l}(c)=\frac{\Tr(\phi_{l}:E_x\to E_x)-\Tr(\phi_{l}:F_x\to F_x)
}{|\det(\id - d\phi_{l}(x) : Q_x\to Q_x)|}.
\]

\begin{thm}[\cite{Lazarov00}]\label{t:dR}
Suppose that all closed orbits of the flow $\phi_s$ are simple. The
restriction of the index $\Ind^\RR (D) \in {\cD}^\prime(\RR)^\RR$ of
the operator $D$ to $\RR\setminus\{0\}$ is given by
\[
\Ind^\RR (D) =\sum_c l(c)\sum_{k\neq 0}b_{kl(c)}(c)
\cdot\delta_{kl(c)},
\]
where $c$ runs over the set of all primitive closed orbits of the
flow $\phi$, and $l(c)$ denotes the length of $c$ (its minimal
positive period).
\end{thm}

Since the set of periods of periodic orbits of the flow is bounded
away from zero, the restriction of $\Ind^\RR (D)$ to some
neighborhood of zero is a distribution supported at $\{0\}$ and,
therefore, is a linear combination of the delta-function at zero
$\delta_0$ and its derivatives.

There is an important particular case when the derivatives of the
delta-function $\delta_0$ do not contribute to the formula for
$\Ind^\RR (D)$, namely, the case when $D$ is the transverse de Rham
operator. Note first of all that in this case the contribution of a
non-degenerate periodic orbit $c$ of $\phi$ is given by
\[
\varepsilon_{l}(c)={\rm sign} \det(\id - d\phi_{l}(x) : Q_x\to Q_x).
\]

Let $T^HM$ denote the orthogonal complement of $T\cF$. The curvature
of the distribution $T^HM$ is a transversal $2$-form $R_H\in
C^\infty(M,\Lambda^2T_HM^*)$ (see (\ref{e:d})), and therefore the
form $\Pf(R_H/2\pi)\in C^\infty(M,\Lambda^qT_HM^*)$ is well defined.
Denote by $X^*$ the $1$-form on $M$ dual to the vector field $X$
with respect to the Riemannian metric. The product
$\Pf(R_H/2\pi)\wedge X^*$ is a top degree differential form on $M$.

\begin{thm}[\cite{Lazarov00}]
In some neighborhood of $0$ in $\RR$,
\[
\Ind^\RR (D)=\int_M \Pf\left(\frac{R_{\cH}}{2\pi}\right)\wedge
X^*\cdot \delta_0.
\]
\end{thm}

Let us now turn to the case when the distribution $T^HM$ is
integrable. Denote by $\cH$ the foliation generated by the
distribution $T^HM$. The action of the flow $\phi$ preserves the
foliation $\cH$, that is, takes any leaf of $\cH$ to a (possibly
different) leaf. There is a holonomy invariant transverse volume
form $\Lambda$ on $\cH$ corresponding to the form $dt$ on $\RR$. Any
transverse Dirac operator $D$ is a tangentially  elliptic operator
with respect to $\cH$, and therefore the $\Lambda$-index
$\Ind_\Lambda(D)$ of $D$ is well defined (see
Subsection~\ref{s:meas}).

\begin{thm}[\cite{Lazarov00}]\label{t:Laz0}
In some neighborhood of $0$ in $\RR$,
\[
\Ind^\RR (D)=\Ind_\Lambda(D)\cdot\delta_0.
\]
\end{thm}

In \cite{AlvKordy2002} the case when the flow $\phi$ is not
necessarily isometric and the operator $D$ is not invariant was
studied. More precisely, the following situation was considered. Let
$M$ be a compact manifold, let $\phi$ be a non-singular flow on $M$,
and let $X$ denote the generator of the flow $\phi$. Suppose that
there exists an integrable distribution $H\subset TM$, which is
invariant under the action of the flow and transversal to orbits of
the flow. If $\cH$ is the foliation on $M$ defined by the
distribution $H$, then $\phi$ preserves the foliation $\cH$.
Consider the leafwise de Rham complex $(\Omega(\cH),d_\cH)$, where
$\Omega(\cH)=C^\infty(M,\Lambda T\cH^*)$ is the space of smooth
leafwise forms on $M$, and $d_\cH$ is the leafwise de Rham
differential. We choose an arbitrary Riemannian metric on the leaves
of $\cH$ and extend it to a Riemannian metric on $M$ by setting
$|X(x)|=1$ for any $x\in M$ and saying that $X$ is orthogonal to
$H$, and we consider the leafwise de Rham operator
\[
D_\cH=d_\cH+d_\cH^*.
\]
It coincides with the transverse de Rham operator for the foliation
$\cF$ defined by the orbits of the flow $\phi$.

\begin{thm}[\cite{AlvKordy2002}]\label{t:lefsch}
In some neighborhood of $0$,
\begin{equation}\label{e:AlvK1}
\Ind^\RR (D_\cH)=\int_M \Pf\left(\frac{R_{\cH}}{2\pi}\right)\wedge
X^*\cdot \delta_0.
\end{equation}
If all closed orbits of the flow $\phi_s$ are simple, then outside
of zero
\begin{equation}\label{e:AlvK2}
\Ind^\RR (D_\cH) =\sum_c l(c)\sum_{k\neq 0}\varepsilon_{kl(c)}(c)
\cdot\delta_{kl(c)},
\end{equation}
where $c$ runs over the set of all primitive closed orbits of the
flow $\phi_t$, $l(c)$ denotes the length of $c$, $x$ is an arbitrary
point on $c$, and
\[
\varepsilon_{l(c)}(c)={\rm sign} \det(\id - d\phi_{l(c)}(x) :
T_x\cH\to T_x\cH).
\]
\end{thm}

In this case the coefficient of $\delta_0$ coincides with the Euler
$\Lambda$-characteristic $\chi_{\Lambda}(\cF)$ of $\cF$ introduced
in \cite{Co79} as the $\Lambda$-index of the leafwise de Rham
operator $D_\cH$. The equality
\[
\chi_{\Lambda}(\cF)=\int_M
\Pf\left(\frac{R_{\cH}}{2\pi}\right)\wedge X^*
\]
is a consequence of the Gauss-Bonnet theorem for measurable
foliations proved in \cite{Co79} as a particular case of the index
theorem for measurable foliations (see Theorem~\ref{t:Co79}).

We`remark also that for an isometric flow $\phi$ the formula
(\ref{e:AlvK2}) was independently proved in~\cite{DenSinghof1}.

We refer the reader to the paper \cite{AlvKordy2002} for connections
of Theorem~\ref{t:lefsch} with the reduced cohomology of $\cH$ and
the Hodge theory for foliations. Taking into account these
relationships, we can understand the statement of
Theorem~\ref{t:lefsch} as a dynamical Lefschetz formula for flows,
that is, a formula, which connects invariants of a foliation with
closed orbits (for a discussion of dynamical Lefschetz formulas for
flows see, for instance, in \cite{Fried87}). We should also mention
a dynamical Lefschetz formula for flows which was proposed as a
conjecture first by Guillemin \cite{Guill77} and later independently
by Patterson \cite{Pat90}. This formula has a form similar to
(\ref{e:AlvK2}), but it is written in the case when the transversal
foliation $\cH$ has codimension greater than $1$ and therefore one
cannot apply analytical results of index theory for transversally
elliptic operators (see a survey of results about the
Guillemin-Patterson conjecture in the book \cite{Juhl}).

The recent interest in the index theory of transversally elliptic
operators, in particular, in the situation when there is a flow on a
compact manifold which preserves a foliation, and in dynamical
Lefschetz formulae for flows, is closely connected with an approach
proposed by Deninger to the study of arithmetic zeta-functions based
on analogies between arithmetic geometry and the theory of dynamical
systems on foliated manifolds (see, for instance, the papers
\cite{Deninger07,Leichtnam05} and the references therein).

In \cite{AlvKordy2008} the results of \cite{AlvKordy2002} were
extended to the case of an arbitrary Lie group action. More
precisely, in this paper a more general situation when $\cH$ is a
Lie foliation of a compact manifold $M$ is considered. In this case
one can define an action of the structural Lie group $G$ on $M$ ``up
to leafwise homotopies'', which enables one to define the index
$\Ind^G (D_\cH)$ of the leafwise de Rham operator $D_\cH$ as a
distribution on $G$. In \cite{AlvKordy2008} a Lefschetz formula is
proved, which gives an expression for $\Ind^G (D_\cH)$ in terms of
the fixed points of the action. It can be regarded as a
generalization of the Selberg formula.

\subsection{Spectral triples associated with transversally
elliptic operators} \label{s:steo} Let $(M,\cF)$ be a compact
foliated manifold, and $E$ a holonomy equivariant Hermitian vector
bundle on $M$. As shown in \cite{Co:nc}, the pair $(H,F)$, where the
Hilbert space $H=L^2(M,E)$ is equipped with an action of the algebra
$C^*(M,\cF)$ by means of the $\ast$-representation $R_E$ and $F\in
\Psi^0(M,E)$ is a transversally elliptic operator with holonomy
invariant transverse principal symbol $\sigma_F$ such that
$\sigma_F^2=1$, and $\sigma^*_F=\sigma_F$, is a Fredholm module over
the $C^*$-algebra $C^*(M,\cF)$. Thus, the transversally elliptic
operator $F$ defines a class $[F]\in K_0(C^*(M,\cF))$. If one
regards the pair $(H,F)$ as a Fredholm module over the algebra
$C^\infty_c(G)$, then this module is summable, and therefore its
Chern character $\tau_n=\operatorname{ch}_\ast (H,F) \in
HC^n(C^\infty_c(G))$ is well defined.

Similarly, one can consider a compact $G$-manifold $M$ ($G$ an
arbitrary Lie group) and a $G$-equivariant Hermitian vector bundle
$E$ on $M$. Then the pair $(H,F)$, where $H=L^2(M,E)$ and $F\in
\Psi^0(M,E)$ is a (not necessarily $G$-invariant) transversally
elliptic operator with invariant transverse principal symbol
$\sigma_F$ such that $\sigma_F^2=1$, and $\sigma^*_F=\sigma_F$, is a
Fredholm module over the $C^*$-algebra $C^*(G)$. The corresponding
Fredholm module over the algebra $C^\infty_c(G)$ is summable, and
its Chern character $\tau_n=\operatorname{ch}_\ast (H,F) \in
HC^n(C^\infty_c(G))$ is defined.

As pointed out in \cite{Co:nc}, if the operator $F$ is invariant,
then its $G$-equivariant index $\chi=\Ind^G(F) \in
{\cD}^\prime(G)^G$ defines a trace on the group algebra
$C^\infty_c(G)$, that is, an element of $HC^0(C^\infty_c(G))$, and
$\tau_{2N}=S^N\chi$ for sufficiently large $N$. In other words, the
Chern character $\tau_{2N}$ and the $G$-equivariant index $\chi$
define the same element in $HP^\ev (C^\infty_c(G))$. In the case
when the operator $F$ is not invariant, the Chern character
$\operatorname{ch}_\ast (H,F) \in HP^*(C^\infty_c(G))$ of the
Fredholm module $(H,\cF)$ is a well-defined homotopy invariant,
unlike its distributional index. This observation served as one of
the main motivations for Connes to introduce the notion of cyclic
cohomology.

In \cite{noncom} spectral triples associated with transversally
elliptic operators on a compact foliated manifold $(M,{\mathcal F})$
were constructed. More precisely, it was proved in \cite{noncom}
that the triple $({\mathcal A},{\mathcal H},D)$ with
\begin{enumerate}
\item ${\mathcal A}$ the algebra $C^{\infty}_c(G)$;
\item ${\mathcal H}$ the Hilbert space $L^2(M,E)$ of square
integrable sections of a holonomy equivariant Hermitian vector
bundle $E$ on which an element $k\in {\mathcal A}$ acts by the
$\ast$-representation $R_E$;
\item $D$ a self-adjoint first order transversally elliptic operator
with holonomy invariant transverse principal symbol such that $D^2$
is self-adjoint and has scalar principal symbol;
\end{enumerate}
is a finite-dimensional spectral triple of dimension $q=\codim \cF$.

As an example, one can assume that the foliation $\cF$ is Riemannian
and $g_M$ is a bundle-like metric on $M$. Let $H$ be the orthogonal
complement of $F=T\cF$. Let ${\mathcal H}=L^2(M,\Lambda^{*}H^{*})$
be the space of transverse differential forms equipped with the
natural action $R_{\Lambda^{*}H^{*}}$ of $\cA$, and $D$ the
transverse signature operator $d_H+d^*_H$. The triple $({\mathcal
A},{\mathcal H},D)$ is a finite-dimensional spectral triple of
dimension $q=\codim \cF$.

It is proved in \cite{noncom} that the spectral triples associated
with transversally elliptic operators are smooth and have simple
discrete dimension spectrum ${\rm Sd}$, which is contained in the
set $\{v\in {\mathbb N}:v\leq q\}$.

One should note that in this case the algebra ${\mathcal A}$ is
non-unital, which can be understood as reflection of the fact that
the space $M/{\mathcal F}$ is non-compact. Therefore, one should
modify the definitions of various geometric and analytic objects
given in the unital case, taking into account behaviour of these
objects at ``infinity''. In \cite{egorgeo}, the algebra
$\Psi^*_0({\mathcal A})$ associated with an arbitrary smooth
spectral triple $({\mathcal A},{\mathcal H},D)$ was constructed.
This algebra can be regarded as an analogue of the algebra of
pseudodifferential operators on a non-compact manifold whose symbols
vanish at infinity along with derivatives of any order. In the same
paper there is a description of the non-commutative
pseudodifferential calculus for spectral triples associated with
transversally elliptic operators in terms of the transversal
pseudodifferential calculus on foliated manifolds, together with a
description of the non-commutative geodesic flow (see also
\cite{matrix-egorov}).

\begin{ex}\label{ex:circle}
Let $\cA=C^\infty(S^1)\rtimes \Gamma$ be the algebraic product of
the algebra $C^\infty(S^1)$ by the group $\Gamma$ of orientation
preserving diffeomorphisms of the circle $S^1$. An arbitrary element
of the algebra $\cA$ is represented as a finite sum
\[
a=\sum_{\phi\in\Gamma}a_\phi U^*_\phi, \quad a_\phi\in
C^\infty(S^1),
\]
the product is given by
\[ (a_\phi U^*_\phi)(b_\psi
U^*_\psi)=a_\phi(b_\psi\circ \phi) U^*_{\phi\psi}
\]
and the involution by
\[
(a_\phi U^*_\phi)^*=U^*_{\phi^{-1}}\bar{a}
\]
(see Example~\ref{ex:crossed}). Define an involutive representation
of the algebra $\cA$ on the Hilbert space $\cH=L^2(S^1)$ by
\[
(\rho(a U^*_\phi)\xi)(x)=a(x)\phi^\prime(x)^{1/2}\xi(\phi(x)), \quad
\xi\in \cH,\quad x\in S^1.
\]
Let the operator $D$ be $D=\frac{1}{i}\frac{\partial}{\partial x}$.

The triple $({\mathcal A},{\mathcal H},D)$ is a spectral triple if
and only if the action of $\Gamma$ is isometric. In the general case
(see \cite{CoM:twisted}) the triple $({\mathcal A},{\mathcal H},D)$
is a $\sigma$-spectral triple, where the automorphism $\sigma$ of
$\cA$ is defined by
\[
\sigma(a U^*_\phi)=\frac{d\phi}{dx} a U^*_\phi.
\]
\end{ex}

\begin{ex}
We give examples of spectral triples $({\mathcal A},{\mathcal H},D)$
associated with the non-commutative torus $T^2_\theta$ (see
\cite{Connes-gravity} and Example~\ref{ex:Atheta}). These triples
are parameterized by a complex number $\tau$ with ${\rm Im}\,\tau
>0$. Put
\[
\cA=\cA_\theta=\left\{a = \sum_{(n,m)\in\ZZ^2}a_{nm}U^nV^m :
a_{nm}\in \cS(\ZZ^2)\right\}.
\]
Define a canonical normalized trace $\tau_0$ on $\cA_\theta$ by
\[
\tau_0(a)=a_{00},\quad a\in \cA_\theta.
\]
Let $L^2(\cA_\theta,\tau_0)$ be the Hilbert space which is the
completion of the space $\cA_\theta$ in the inner product
$(a,b)=\tau_0(b^*a), a,b\in \cA_\theta$. The Hilbert space $\cH$ is
defined as the sum of two copies of the space
$L^2(\cA_\theta,\tau_0)$ equipped with the grading given by the
operator
\[
\gamma=\begin{pmatrix} 1 & 0\\ 0 & -1
\end{pmatrix}.
\]

The representation $\rho$ of the algebra $\cA_\theta$ in $\cH$ is
given by the left multiplication, that is, for any $a\in
\cA_\theta$,
\[
\rho(a)=\begin{pmatrix} \lambda(a) & 0\\ 0 & \lambda (a)
\end{pmatrix},
\]
where the operator $\lambda(a)$ is defined on $\cA_\theta\subset
L^2(\cA_\theta,\tau_0)$ by the formula
\[
\lambda(a)b=ab, \quad b\in \cA_\theta.
\]
We introduce differentiations $\delta_1$ and $\delta_2$ on
$\cA_\theta$ by
\[
\delta_1(U)=2\pi iU,\quad \delta_1(V)=0;\quad \delta_2(U)=0,\quad
\delta_2(V)=2\pi iV.
\]
The operator $D$ explicitly depends on $\tau$ and has the form
\[
D=\begin{pmatrix} 0 & \delta_1+\tau\delta_2\\
-\delta_1-\bar{\tau}\delta_2 & 0
\end{pmatrix}.
\]
The triples constructed above are two-dimensional smooth spectral
triples.

We refer the reader to the book \cite{Naz-Savin-Sternin} and its
references for the index theory of elliptic operators with shifts
and connected non-commutative geometry. The theory of elliptic
operators on the non-commutative torus $T^2_\theta$ (see
\cite{CoCRAS,Co}) is a particular example of this general theory.
\end{ex}

\subsection{Para-Riemannian foliations and transversally
hypoelliptic operators}\label{para}

A foliation $\cF$ on a manifold $M$ is said to be para-Riemannian if
there exists an integrable distribution $V$ on $M$ which contains
the tangent bundle $T\cF$ of $\cF$, and is such that the bundles
$TM/V$ and $V/T\cF$ are holonomy equivariant Riemannian bundles. If
$V$ is an integrable distribution on $M$ which defines the
para-Riemannian structure, and $\cV$ is the corresponding foliation
on $M$, $\cF\subset \cV$, then $\cV$ is Riemannian, and the
restriction of the foliation $\cF$ to each leaf $L$ of $\cV$ is a
Riemannian foliation on $L$.

As mentioned already in Subsection~\ref{s:tfclass}, the interest to
para-Riemannian foliations consists in the fact that in some
problems of index theory the study of arbitrary foliations can be
reduced to the study of para-Riemannian foliations. Let us describe
the corresponding construction by Connes and Moscovici \cite{Co-M}.
In fact, they dealt with a closely related situation (a strongly
Morita equivalent one): they considered an oriented smooth manifold
$W$ endowed with an action of a discrete group $\Gamma $ instead of
a foliated manifold. We consider the fibration $\pi:P(W)\to W$ whose
fibre $P_x(W)$ at $x\in W$ is the space of all Euclidean metrics on
the vector space $T_xW$. Thus, a point $p\in P(W)$ is given by a
point $x\in W$ and a non-degenerate quadratic form on $T_xW$. Let
$F_+(W)$ be the bundle of positive frame in $W$ whose fibre $F_x(W)$
at $x\in W$ is the space of orientation preserving linear
isomorphisms $\RR^n\to T_xW$. Equivalently, the bundle $P(W)$ can be
described as the orbit space of the bundle $F_+(W)$ with respect to
the fibrewise action of the subgroup $SO(n)\subset GL(n,\RR)$,
$P(W)=F_+(W)/SO(n)$. We will use a natural invariant Riemannian
metric on the symmetric space $GL_+(n,\RR)/O(n)$ given by the matrix
Hilbert-Schmidt norm on the tangent space of $GL(n,\RR)/SO(n)$,
which is identified with the space of symmetric $n\times n$
matrices. If we transfer this metric to the fibres $P_x$ of the
bundle $P(W)=P$, then we obtain a Euclidean structure on the
vertical distribution $V\subset TP$. The normal space $N_p=T_pP/V_p$
is naturally identified with the space $T_xW, x=\pi(p).$ Thus, the
quadratic form on $T_xW$ corresponding to $p$ defines a natural
Euclidean structure on $N_p$.

There are natural actions of the group $\Gamma$ on $F_+(W)$ and $P$.
This action takes fibres of the bundle $\pi:P\to W$ to fibres.
Moreover, the Euclidean structures introduced above on the
distributions $V$ and $N$ are invariant under the action of the
group $\Gamma$. In this case one says that there is a triangular
structure on $P$ invariant under the group action. It is an analogue
of a para-Riemannian foliation in this situation. A very essential
aspect of this construction is the fact that the fibres of the
bundle $\pi:P\to W$ are complete Riemannian manifolds of
non-positive curvature.

Another important property of an arbitrary para-Riemannian foliation
$(M,\cF)$ is the existence of a Fredholm module over its
$C^*$-algebra $C^*(M,\cF)$. A construction of this module is given
in \cite{Hil-Skan}. It makes use of transversally hypoelliptic
operators and pseudodifferential operators of type $(\rho, \delta)$.

In \cite{Co-M} Connes and Moscovici described a spectral triple
associated with an invariant triangular structure. Thus, let $P$ be
a smooth manifold equipped with an action of a discrete group
$\Gamma$. Suppose that $P$ is the total space of a bundle $\pi:P\to
W$ over an oriented smooth manifold $W$. The action of $\Gamma$
leaves $W$ invariant and takes fibres of the bundle $\pi:P\to W$ to
fibres. Finally, $\Gamma$ preserves Euclidean structures on the
vertical distribution $V\subset TP$ and on the horizontal
distribution $N=TP/V$.

Consider the Hermitian vector bundle $E=\Lambda^*(V^*\otimes
\CC)\otimes \Lambda^*(N^*\otimes\CC)$ over $P$. The Hermitian
structure in the fibres of $E$ is determined by the Euclidean
structures on $V$ and $N$. The bundles $\Lambda^*(V^*\otimes \CC)$
and $\Lambda^*(N^*\otimes\CC)$ have the grading operators $\gamma_V$
and $\gamma_N$ given by the Hodge operators of the Euclidean
structures and the orientations of $V$ and $N$. The Euclidean
structures on $V$ and $N$ also define a natural volume form $v\in
\Lambda^*V^*\otimes \Lambda^*N^*=\Lambda T^*P$.

Let $\cA$ be the crossed product $C^\infty_c(P)\rtimes \Gamma$. We
recall that this algebra is generated as a linear space by
expressions of the form $fU_g$, where $f\in C^\infty_c(P)$ and $g\in
\Gamma$ (see Example~\ref{ex:actions}).

Let $\cH$ be the space $L^2(P,E)$ equipped by the Hilbert structure
given by the volume form $v$ and the Hermitian structure  on $E$.
The action of $\cA$ in $\cH$ is given in a following way. A function
$f\in C^\infty_c(P)$ acts as the corresponding multiplication
operator in $\cH$. For any $g\in\Gamma$ the unitary operator $U_g$
is given on $\cH$ by the natural actions of $\Gamma$ on sections of
the bundles $V$ and $N$.

Consider the foliation $\cV$ given by the fibres of the fibration
$P$. Then $V=T\cV$, and $N$ is the normal bundle to $\cV$. Denote by
$d_L:C^\infty(P,E)\to C^\infty(P,E)$ the tangential de Rham
differential associated with the foliation $\cV$ (see (\ref{e:d})).
Let $Q_L$ be the second order tangential differential operator
acting in $C^\infty(P,E)$ given by
\[
Q_L=d_Ld^*_L-d^*_Ld_L.
\]
As shown in \cite{Co-M}, the principal symbol of $Q_L$ is homotopic
to the principal symbol of the signature operator $d_L+d^*_L$.

We choose an arbitrary distribution $H$ on $P$ transverse to $V$,
and consider the corresponding transverse de Rham differential $d_H$
(see (\ref{e:d})) and the transverse signature operator
\[
Q_H=d_H+d^*_H.
\]
This operator depends on the choice of $H$, but its transverse
principal symbol is independent of $H$. Let us define a mixed
signature operator $Q$ on $C^\infty(P,E)$ as
\[
Q=Q_L(-1)^{\partial_N}+Q_N,
\]
where $(-1)^{\partial_N}$ denotes the parity operator in the
transverse direction, that is, it coincides with $1$ on
$\Lambda^{\rm ev}N^*$ and with $-1$ on $\Lambda^{\rm odd}N^*$.
Assume that $Q$ is essentially self-adjoint in $\cH$. Using
functional calculus, we define the operator $D$ as $Q=D |D|$.

One should note that, although the above construction can be applied
to any manifold endowed with an invariant triangular structure, the
question of essential self-adjointness of the operator $Q$ is a
difficult analytic question, since the manifold $P$ is non-compact.
In \cite{Co-M} this question was answered just for the example
described above of a triangular structure on the Euclidean metrics
bundle $P(W)$ associated with an arbitrary group $\Gamma$ of
diffeomorphisms of a smooth manifold $W$.

\begin{thm}\cite{Co-M}\label{t:para}
For the triangular structure associated with an arbitrary group
$\Gamma$ of diffeomorphisms of a smooth manifold $W$, the operator
$Q$ is self-adjoint and the triple $(\cA,\cH,D)$ constructed above
is a spectral triple of dimension $\dim V+2\dim N$.
\end{thm}

The proof of this theorem makes essential use of the
pseudodifferential calculus constructed by Beals and Greiner
\cite{BG} on Heisenberg manifolds. It is also shown in this paper
that the non-commutative integral $\bint$ defined by such a spectral
triple coincides with the Wodzicki-Guillemin type trace defined on
the Beals-Greiner algebra of pseudodifferential operators.

In~\cite{CoM:Hopf} Connes and Moscovici computed the Chern character
of the spectral triple constructed in Theorem~\ref{t:para} by using
the non-commutative local index theorem, Theorem~\ref{t:cmodd}. In
fact, a direct computation of the Chern character of a spectral
triple associated with a triangular structure on a smooth manifold
using the formulae given in Theorem~\ref{t:cmodd} is quite
cumbersome even in the one-dimensional case. One gets formulae
involving thousands terms, most of which give zero contribution. To
simplify the computations a priori, Connes and Moscovici introduced
a Hopf algebra $\cH_n$ of transverse vector fields in $\RR^n$ which
plays the role of a quantum symmetry group. They constructed the
cyclic cohomology $HC^*(\cH)$ for an arbitrary Hopf algebra $\cH$
and a map
\[
HC^*(\cH_n)\to HC^*(C^\infty_c(P)\rtimes\Gamma).
\]
Moreover, they showed that the cyclic cohomology $HC^*(\cH_n)$ is
canonically isomorphic to the Gel'fand-Fuchs cohomology
$H^*(W_n,SO(n))$ (see Subsection~\ref{s:high-ind}). Therefore, there
is defined a characteristic homomorphism
\[
\chi_{SO(n)}^* : H^*(W_n,SO(n))\to HP^*(C^\infty(P)\rtimes \Gamma).
\]
It is the composition of the map (\ref{e:second}) and the
homomorphism $\Phi_*$ (see (\ref{e:Phi*})).

The following theorem is the main result of \cite{CoM:Hopf} (cf.
also \cite{CoMEssays} and the survey \cite{skandal-hopf}).

\begin{thm}
Let $(\cA,\cH,D)$ be the spectral triple introduced in
Theorem~\ref{t:para}. The Chern character
$\operatorname{ch}_*(\cA,\cH,D)\in HP^*(C^\infty(PW)\rtimes \Gamma)$
is the image of a universal class $\cL_n\in H^*(W_n,SO(n))$ under
the characteristic homomorphism $\chi_{SO(n)}^*$:
\[
\operatorname{ch}_*(\cA,\cH,D)=\chi_{SO(n)}^*(\cL_n).
\]
\end{thm}

There is one more computation, given in \cite{Perrot}, illustrating
the Connes-Moscovici local index theorem. Let $\Sigma$ be a closed
Riemann surface and $\Gamma$ a discrete pseudogroup of local
conformal maps of $\Sigma$ without fixed points. Using methods of
\cite{CoM:Hopf} (the bundle of Kaehler metrics $P$ on $\Sigma$,
hypoelliptic operators, Hopf algebras), in \cite{Perrot} a spectral
triple, which is a generalization of the classical Dolbeaut complex
to this setting, is constructed. The Chern character of this
spectral triple as a cyclic cocycle on the crossed product
$C^\infty_c(\Sigma)\rtimes \Gamma$ is computed in terms of the
fundamental class $[\Sigma]$ and a cyclic 2-cocycle which is a
generalization of the class Poincar\'e dual to the Euler class. This
formula can be regarded as a non-commutative version of the
Riemann-Roch theorem.

A Hopf algebra of the same type as the algebra $\cH_1$ was
constructed by Kreimer \cite{Kreimer98} for the study of the
algebraic structure of the perturbative quantum field theory. This
connection was further elaborated in \cite{CoKr1,CoKr2,CoKr3}. One
should also mention the papers \cite{CoM:mmj1} and \cite{CoM:mmj2}
by Connes and Moscovici on modular Hecke algebras, where it is
shown, in particular, that such an important algebraic structure on
the modular forms as the Rankin-Cohen brackets has a natural
interpretation in the language of non-commutative geometry in terms
of the Hopf algebra $\cH_1$.

\section{The index theory of tangentially elliptic operators}
Let $(M,{\mathcal F})$ be a compact foliated manifold, and $E$ a
smooth vector bundle on $M$. A linear differential operator $D$ of
order $\mu$ acting in $C^{\infty}(M,E)$ is called a tangential
differential operator if in any foliated chart $\phi : U\subset M\to
I^p\times I^q$ and any trivialization of $E$ over it the operator
$D$ has the form
\begin{equation}
D =\sum _{|\alpha|\leq \mu}a_{\alpha}(x,y)
\frac{\partial^{|\alpha|}}{\partial x_1^{\alpha_1}\ldots \partial
x_p^{\alpha_p}}, \quad (x,y)\in I^{p} \times I^{q},
\label{rev:(1.14)}
\end{equation}
where $\alpha=(\alpha_1,\ldots,\alpha_p)\in \ZZ_+^p$ is a
multiindex, $|\alpha|=\alpha_1+\ldots+\alpha_p$, and the
$a_{\alpha}$ are smooth matrix-valued functions on $I^p\times I^q$.

For a tangential differential operator $D$ given by
(\ref{rev:(1.14)}) in some foliated chart $\phi : U\subset M\to
I^p\times I^q$ and a trivialization of $E$ over it, we define its
tangential (complete) symbol $$ \sigma (x,y,\xi) = \sum_{| \alpha |
\leq \mu} a_{\alpha}(x,y) (i\xi)^{\alpha},\quad (x,y)\in I^p\times
I^q,\quad \xi\in\RR^p, $$ and its tangential principal symbol
$$ \sigma _{\mu}(x,y,\xi) = \sum _{ | \alpha | =\mu}
a_{\alpha}(x,y)(i\xi)^{\alpha},\quad (x,y)\in I^p\times I^q,\quad
\xi\in\RR^p. $$ The tangential principal symbol is invariantly
defined as a section of the bundle $\Hom(\pi_F^*E)$ on $T^*\cF$
(where $\pi_F : T^*\cF\to M$ is the natural projection).

A tangential differential operator $D$ is said to be tangentially
elliptic if its tangential principal symbol $\sigma_{\mu}$ is
invertible for $\xi \neq 0$.

Let $D : C^\infty(M,E)\to C^\infty(M,E)$ be a tangential
differential operator on a compact foliated manifold $(M,\cF)$. The
restrictions of $D$ to the leaves of $\cF$ define a family
$(D_L)_{L\in M/\cF}$, where for any leaf $L$ of $\cF$ the operator
$D_L: C^\infty(L,E\left|_L\right.)\to C^\infty(L,E\left|_L\right.)$
is a differential operator on $L$. For any $x\in M$ the lift of the
operator $D_{L_x}$ by the holonomy covering  $s : G^x\to L_x$
defines a differential operator $D_x : C^\infty(G^x,s^*E)\to
C^\infty(G^x,s^*E)$. The operator family $\{D_x, x\in M\}$ is a
$G$-operator (see (\ref{e:G})). Families of this form will be called
differential $G$-operators. If $D$ is a tangentially elliptic
operator, then the corresponding $G$-operator will be said to be
elliptic. In \cite{Co79} the corresponding algebra of
pseudodifferential $G$-operators on a foliated manifold is
constructed.

\subsection{The index theory for measurable foliations}\label{s:meas}
Let $D$ be a tangentially elliptic operator on a compact foliated
manifold $(M,\cF)$, and $\{D_x: x\in M\}$ the corresponding elliptic
$G$-operator. Suppose that the foliation $\cF$ has a holonomy
invariant transverse measure $\Lambda$. The families $(P_{\Ker
D_x})_{x\in M}$ and $(P_{\Ker D^*_x})_{x\in M}$ consisting of the
orthogonal projections onto $\Ker D_x$ and $\Ker D^*_x$ in the space
$L^2(G^x)$, respectively, define elements $P_{\Ker D}$ and $P_{\Ker
D^*}$ of the foliation von Neumann algebra $W^*_\Lambda(M,\cF)$. The
holonomy invariant measure $\Lambda$ defines a faithful normal
semifinite trace  $\trLambda$ on $W^*_\Lambda(M,\cF)$. It is proved
in \cite{Co79} that the dimensions
\[
\dim_\Lambda \Ker D = \trLambda P_{\Ker D},\quad \dim_\Lambda \Ker
D^* = \trLambda P_{\Ker D^*},
\]
are finite, and therefore the index of the tangentially elliptic
operator $D$ is well defined by
\[
\operatorname{Ind}_\Lambda (D)=\dim_\Lambda \Ker D - \dim_\Lambda
\Ker D^*.
\]

Suppose that the bundle $T\cF$ is oriented. As in the index theory
for families of elliptic operators, the tangential principal symbol
$\sigma_D$ of $D$ defines an element of $K(T^*\cF)$. Denote by
$\pi_{F!}: H^{*}(T\cF)\to H^{*}(M)$ the map given by integration
along the fibres of the bundle $\pi_F:T\cF\to M$.

\begin{thm}\cite{Co79}\label{t:Co79}
One has the formula
\[
\operatorname{Ind}_\Lambda (D)=(-1)^{p(p+1)/2}\langle C,
\pi_{F!}\operatorname{ch}(\sigma_D) \operatorname{Td}(T^*\cF\otimes
\CC) \rangle,
\]
where $C$ is the Ruelle-Sullivan current, corresponding to the
transverse measure $\Lambda$.
\end{thm}

This theorem is completely analogous to the Atiyah-Singer index
theorem in the cohomological form, Theorem~\ref{t:AScoh}, with the
only difference being  that here one uses the pairing with the
Ruelle-Sullivan current $C$ instead of integration over the compact
manifold on the right-hand side of the Atiyah-Singer formula.

An odd version of Theorem~\ref{t:Co79} is proved in \cite{DHK2} (see
also \cite{DHK}). Suppose that $D: C^\infty(M,E)\to C^\infty(M,E)$
is a first order tangentially elliptic operator on a compact
foliated manifold $(M,\cF)$, and let $D_x : C^\infty(G^x,s^*E)\to
C^\infty(G^x,s^*E), x\in M,$ be the corresponding differential
$G$-operator. Assume that the foliation $\cF$ has a holonomy
invariant transverse measure $\Lambda$. Assume also that the
operator $D_x$ is self-adjoint in the space $L^2(G^x,s^*E)$. Denote
by $P_x$ its spectral projection corresponding to the positive
semi-axis. For any $\phi\in C(M,U(N))$ denote by $M_\phi$ the
corresponding $G$-operator of multiplication by $\phi$. The leafwise
Toeplitz operator associated with $\phi$ is a bounded $G$-operator
\[
T_{\phi}=\{T_{\phi,x} : L^2(G^x,s^*E\otimes \CC^N)\to
L^2(G^x,s^*E\otimes \CC^N), x\in M\}
\]
given by
\[
T_{\phi,x}=P_xM_\phi P_x.
\]
It is proved in \cite{DHK2} that if $\phi\in C(M,U(N))$ is
invertible, then the operator $T_\phi$ is a Breuer-Fredholm
operator, and the Breuer-Fredholm index of $T_\phi$ is defined by
\[
\operatorname{Ind}_\Lambda (T_\phi)=\dim_\Lambda (\Ker T_\phi) -
\dim_\Lambda (\Ker T^*_\phi).
\]

Denote by $E_+$ the subbundle of the bundle $\pi_F^*E$ on $ST^*\cF$
spanned by the eigenvectors of the principal symbol $\sigma_1(D)$ of
the operator $D$, corresponding to positive eigenvalues. As above,
denote by $\pi_{F!}: H^{*}(ST^*\cF)\to H^{*}(M)$ the map given by
integration along the fibres of the bundle $\pi_F:ST^*\cF\to M$.

\begin{thm}[\cite{DHK,DHK2}]
\[
\operatorname{Ind}_\Lambda (T_\phi)=(-1)^p
 \langle C, \ch([\phi]) \pi_{F!}\ch(E_+) \Td(T\cF\otimes \CC)\rangle,
\]
where $C$ is the Ruelle-Sullivan current corresponding to the
transverse measure $\Lambda$.
\end{thm}

One should also note the papers
\cite{CDSS,Fedosov-Shubin1,Fedosov-Shubin2,Shubin:UMN79,CoCRAS,CMX,Lesch91,Phillips-Raeburn,CPS2},
where analogous problems of index theory (in both even and odd
settings) were studied in closely related situations --- for
almost-periodic and random operators in $\RR^n$.

In \cite{Ramachandran93} an analogue of the Atiyah-Patodi-Singer
theorem for measurable foliations was proved. In particular, the
eta-invariant of tangentially elliptic operators was introduced
there (see also \cite{Peric}). In \cite{H-L} (see also \cite{H-L91})
an analogue of the Atiyah-Bott-Lefschetz fixed point formula
\cite{Atiyah-Bott} was proved for maps of a compact manifold
equipped with a measurable foliation, which take each leaf to
itself.

\subsection{The $K$-theoretic index theory}\label{s:Kth}
In \cite{Co-skandal} (see also \cite{Co:survey}) a $K$-theoretic
version of the index theorem is proved for tangentially elliptic
operators on an arbitrary compact foliated manifold $(M,\cF)$. Let
$D$ be a tangentially elliptic operator on a compact manifold $M$
acting on sections of a vector bundle $E$ on $M$. Using operator
techniques, one constructs the analytic index ${\rm Ind}_a(D)\in
K_0(C^*_r(M,\cF))$ of the corresponding right-invariant elliptic
$G$-operator $\{D_x : C^\infty_c(G_x, r^*E) \to C^\infty_c(G_x,
r^*E), x\in M\}$ \cite{Co79,Co:survey}, and starting from the class
$[\sigma_D]\in K(T^*\cF)$ defined by the tangential principal symbol
$\sigma_D$ of $D$ one constructs its topological index ${\rm
Ind}_t(D)\in K_0(C^*_r(M,\cF))$.

The analytic index is constructed in the following way. Recall that
the bundle $E$ gives rise to the $C(M)$-$C^*_r(M,\cF)$-bimodule
$\cE_{M,E}$ defined as the completion of the pre-Hilbert
$C^\infty_c(G)$-module $\mathcal E_\infty=C^\infty_c(G, r^*E)$ (see
Subsection~\ref{s:vect}). Any right-invariant pseudodifferential
$G$-operator naturally defines an endomorphism of the Hilbert
$C^*_r(M,\cF)$-module $\cE_{M,E}$. Therefore, the operator $D$
defines an element $[D]\in KK(C(M),C^*_r(M,\cF))$ given by the pair
$(\cE_{M,E},F)$, where $F=D(I+D^2)^{1/2}$. The image of this element
under the map
\[
KK(C(M),C^*_r(M,\cF))\to KK(\CC, C^*_r(M,\cF))=K_0(C^*_r(M,\cF))
\]
is the analytic index of $D$, ${\rm Ind}_a(D)\in K_0(C^*_r(M,\cF))$.
In \cite{DHK2} another construction of the element $[D]\in
KK^1(C(M), C^*_r(M,\cF))$ corresponding to a self-adjoint
tangentially elliptic operator $D$ is given. This construction uses
an extension of the $C^*$-algebra $C^*_r(M,\cF)$ generated by
smoothed leafwise Toeplitz operators.

We now describe the construction of the topological index of a
tangentially elliptic operator $D$. Let $i$ be an embedding of the
manifold $M$ into $\RR^{2n}$. Denote by $N$ the total space of the
normal bundle to the leaves: $N_x=(i_*(F_x))^\bot\subset \RR^{2n}$.
We consider the foliation $\tilde{F}$ on the manifold
$\tilde{M}=M\times \RR^{2n}$ with fibres $\tilde{L}=L\times\{t\}$,
where $L$ is a leaf of $\cF$ and $t\in\RR^{2n}$. The map $N \ni
(x,\xi) \mapsto (x,i(x)+\xi)$ maps some open neighborhood of the
zero section in $N$ to an open transversal $T$ of the foliation
$(\tilde{M},\tilde{\cF})$. For a suitable open neighborhood $\Omega$
of the transversal $T$ in $\tilde{M}$, the $C^*$-algebra
$C^*_r(\Omega,\tilde{\cF}\left|_{\Omega}\right.)$ of the restriction
of the foliation $\tilde{\cF}$ to $\Omega$ is Morita equivalent to
the algebra $C_0(T)$. Hence, the embedding
$C^*_r(\Omega,\tilde{\cF}\left|_{\Omega}\right.) \subset
C^*_r(\tilde{M},\tilde{\cF})$ defines a map $K^0(N)\to
K^0(C^*_r(\tilde{M},\tilde{\cF}))$. Since
$C^*_r(\tilde{M},\tilde{\cF})=C^*_r(M,\cF)\otimes C_0(\RR^{2n})$,
the Bott periodicity implies that
$K^0(C^*_r(\tilde{M},\tilde{\cF}))=K^0(C^*_r(M,\cF))$. With use of
the Thom isomorphism, $K^0(T^*\cF)$ is identified with $K^0(N)$.
Thus, one gets the topological index map
\[
{\rm Ind}_t : K^0(T^*\cF)\to K_0(C^*_r(M,\cF)).
\]

\begin{thm}\cite{Co-skandal,Co:survey}\label{t:CS}
For any tangentially elliptic operator $D$ on a compact foliated
manifold $(M,\cF)$,
\[
{\rm Ind}_a(D)={\rm Ind}_t(D).
\]
\end{thm}

If a foliation $\cF$ is given by the fibres of a fibration $M\to B$,
then $K_0(C^*_r(M,\cF))\cong K^0(B)$ and Theorem~\ref{t:CS} reduces
to the Atiyah-Singer index theorem for families of elliptic
operators, Theorem~\ref{t:ASfamily}. If the foliation has a holonomy
invariant measure $\Lambda$, then the trace ${\rm tr}_\Lambda$ on
the $C^*$-algebra $C^*_r(M,\cF)$ is well defined. In turn, it
defines a map ${\rm Tr}_\Lambda$ from $K^0(C^*_r(M,\cF))$ to $\RR$.
It is not difficult to show that ${\rm Tr}_\Lambda ({\rm
Ind}_a(D))={\rm Ind}_\Lambda (D)$. The composition ${\rm Tr}_\Lambda
\circ {\rm Ind}_t$ can be computed by topological methods, and, as a
consequence of Theorem~\ref{t:CS}, one gets the index theorem for
measurable foliations (see Subsection~\ref{s:meas}). We note also
the paper \cite{Wang07}, which gives a generalization of
Theorem~\ref{t:CS} in the twisted $K$-theory.

Following the ideas of the paper \cite{Baum-D82}, one can give an
equivalent formulation of Theorem~\ref{t:CS} in terms of the map
$\mu$ \cite{Co}. First of all, observe that the principal symbol
$\sigma_D$ of a tangentially elliptic operator $D$ determines in the
geometric $K$-homology group $K^*_{\rm top}(M,\cF)$ a class
$[\sigma_D]$ given by the $K$-cycle $(T\cF,[\sigma_D],p\circ \pi)$,
where $[\sigma_D]\in K^0(T^*\cF)$ is the class given by $\sigma_D$,
and the map $p\circ \pi : T^*\cF\to M/\cF$ is obtained as the
composition of the natural projections $\pi:T^*\cF\to M$ and $p:M\to
M/\cF$. Then
\begin{equation}\label{e:mu}
\mu_r([\sigma_D])= {\rm Ind}_a(D)\in K_0(C^*_r(M,\cF)).
\end{equation}

In \cite{Benameur97} there is an equivariant generalization of
Theorem~\ref{t:CS} to the case of an action of a compact Lie group
$H$ taking each leaf of $\cF$ to itself. As a consequence, the
author extended the Lefschetz theorem proved in~\cite{H-L} to the
case of arbitrary tangentially elliptic complexes under the
assumption that the diffeomorphism $f\colon M\to M$ is included into
an action of a compact Lie group $H$ taking each leaf of $\cF$ to
itself.

Finally, we note that in \cite{Vassout} (see also
\cite{Benameur-Fack06,Vassout07}) semifinite spectral triples
associated with a tangentially elliptic operators are constructed.
Let $(M,\cF)$ be a smooth compact foliated manifold whose leaves are
even-dimensional spin manifolds. Denote by $S$ the associated spinor
bundle, and suppose that there is a holonomy invariant measure
$\Lambda$. We regard the involutive algebra ${\mathcal A}=C(M)$ as a
subalgebra of the semifinite von Neumann algebra
$\cN=W^{\ast}_\Lambda(M,\mathcal F)$ acting in the Hilbert space
\[
\mathcal H=L^2(G,s^*S,m)=\int^\oplus_M L^2(G^x,s^*S,\nu^x)\,d\mu(x),
\]
For any function $a\in C(M)$ its action in $\cH$ is given by the
operator of multiplication by the function $a\circ s$ in each space
$L^2(G^x,s^*S,\nu^x)$. Then the set $({\mathcal A}, {\mathcal H},
D)$, where $D$ is the leafwise Dirac spin operator, is a semifinite
spectral triple whose dimension is finite and equals the dimension
of the foliation.

In \cite{Ben-C-Ph-R-S} (see also \cite{Carey-S}) semifinite spectral
triples associated with differential operators with almost-periodic
coefficients are constructed. Using the semifinite non-commutative
local index theorem \cite{CPRS2,CPRS3} (see
Subsection~\ref{s:semifinite}), the authors computed the spectral
flow along a family of almost-periodic Dirac operators.

\subsection{The higher index theory}\label{s:high-ind}
The approach proposed by Connes and Moscovici in \cite{Co-M90} to
the higher index theory for tangentially elliptic operators is based
on the idea of representing higher indices as pairings of the
$K$-theoretic index of a tangentially elliptic operator with some
cyclic cocycle on the algebra $C^\infty_c(G)$. For this, one should
first modify the definition of the index and give a definition of
the index as an element of $K_0(C^\infty_c(G))$.

The construction of the pseudodifferential calculus on the holonomy
groupoid $G$ \cite{Co79} enables one to construct the short exact
sequence
\[ 0\longrightarrow C^*_r(M,\cF)\longrightarrow
\bar{\Psi}^0(\cF)\stackrel{\sigma_0}{\longrightarrow} C
(S^*\cF)\longrightarrow 0,
\]
where $\Psi^0(\cF)$ is the algebra of zero order pseudodifferential
$G$-operators, $\bar{\Psi}^0(\cF)$ is the uniform completion of the
algebra $\Psi^0(\cF)$ in the space $L^2(G,m)=\int^\oplus_M
L^2(G^x,\nu^x)\,d\mu(x)$, and $\sigma_0$ is the principal symbol
map. The corresponding exact sequence in $K$-theory consists of six
terms:
\[
  \begin{CD}
K_1(C^*_r(M,\cF)) @>>> K_1(\bar{\Psi}^0(\cF)) @>>> K_1(C(S^*\cF)) \\
@AAA @. @VVV
\\ K_0(C(S^*\cF)) @<<<
K_0(\bar{\Psi}^0(\cF)) @<<< K_0(C^*_r(M,\cF))
  \end{CD}
\]
The connecting homomorphism $\partial : K_1(C
(S^*\cF))\longrightarrow K_0(C^*_r(M,\cF))$ can be shown
\cite{Co-skandal} to coincide with the analytical index
homomorphism. More precisely, if $A$ is a zero order elliptic
$G$-pseudodifferential operator, then its principal symbol defines a
class $[\sigma_A]\in K_1(C^\infty (S^*\cF))$, and the element
$\partial [\sigma_A]\in K_0(C^*_r(M,\cF))$ coincides with the
analytical index of $A$ defined in Subsection~\ref{s:Kth}. Explicit
algebraic construction of the connecting homomorphism in $K$-theory
leads to the following rule for constructing the analytic index of
an elliptic symbol $a\in C^\infty(S^*\cF, \End(\pi^*E,\pi^*F))$.
After the bundle $E\oplus F$ is embedded in a trivial bundle, the
function $\tilde{a}=
\begin{pmatrix}
  0 & -a^{-1} \\
  a & 0 \\
\end{pmatrix}
\in C^\infty(S^*\cF, \Hom(\pi^*(E\oplus F)))$ defines an element of
$GL_N(C^\infty(S^*\cF))$ for sufficiently large $N$. This element
can lifted to $GL_N(\Psi^0(\cF))$. For example, choose any operator
$A\in \Psi^0(\cF;E,F)$ such that $\sigma_0(A)=a$ and $B\in
\Psi^0(\cF;F,E)$ such that $\sigma_0(B)=a^{-1}$. Then $S_0=I-BA\in
\Psi^{-1}(\cF;E)$ and $S_1=I-AB\in \Psi^{-1}(\cF;F)$. The operator
\[
L=\begin{pmatrix}
  S_0 & -B-S_0B \\
  A & S_1 \\
\end{pmatrix}\in \Psi^{0}(\cF;E\oplus F)
\]
provides the desired lift. Thus, $\sigma_0(L)=\tilde{a}$ and
\[
L^{-1}=\begin{pmatrix}
  S_0 & B+BS_0 \\
  -A & S_1 \\
\end{pmatrix}\in \Psi^{0}(\cF;E\oplus F)
\]
By definition, put
\begin{equation}\label{e:ind1}
\partial [a]=[P]-[e],
\end{equation}
where $P$ and $e$ are the idempotents defined as follows:
\begin{equation}\label{e:ind2}
P=L\begin{pmatrix}
  I_E & 0 \\
  0 & 0 \\
\end{pmatrix}
L^{-1}= \begin{pmatrix}
  S_0^2 & S_0(I+S_0)B \\
  S_1A & I_F-S^2_1 \\
\end{pmatrix}, \quad e=\begin{pmatrix}
  0 & 0 \\
  0 & I_F \\
\end{pmatrix}.
\end{equation}

The construction of the $G$-pseudodifferential calculus and the
definition of the analytic index have been extended to the case of
an arbitrary Lie groupoid in \cite{MP97,NWX99}.

In fact, the calculus of pseudodifferential operators associated
with the holonomy groupoid enables one to construct a parametrix
$B\in \Psi^0(\cF;F,E)$ such that $S_0=I-BA\in \Psi^{-\infty}(\cF;E)$
and $S_1=I-AB\in \Psi^{-\infty}(\cF;F)$. Then the formulae
(\ref{e:ind1}) and (\ref{e:ind2}) define the analytic index of $A$
as an element ${\rm Ind}\,A \in K_0(C^\infty_c(G))$ (see
\cite{Co-M90,Co}).

The natural embedding $j: C^\infty_c(G)\to C^*_r(M,\cF)$ induces a
map $j_*: K_0(C^\infty_c(G))\to K_0(C^*_r(M,\cF))$. One can show
(see \cite{Co-M90}) that $j_*({\rm Ind}\,A)={\rm Ind}_a\,A$. The map
$j_*$, in general, is not an isomorphism, and therefore we lose some
information in passing from ${\rm Ind}\,A \in K_0(C^\infty_c(G))$ to
${\rm Ind}_a\,A \in K_0(C^*_r(M,\cF))$. Namely, for the analytic
index with values in $K_0(C^*_r(M,\cF))$ there are results like
vanishing or homotopy invariance. For instance, as shown in
\cite{Baum-Connes83} (see also \cite{Hil-Skan92,Hurder93}), the
analytic index ${\rm Ind}_a\,A \in K_0(C^*_r(M,\cF))$ of the
leafwise signature operator on a compact manifold equipped with an
even-dimensional oriented foliation is invariant under leafwise
oriented homotopy equivalences. On the other hand, the analytic
index ${\rm Ind}\,A \in K_0(C^\infty_c(G))$ of an elliptic
$G$-pseudodifferential operator $A$ depends, in general, not just on
the class $[\sigma_A]\in K_0(C^\infty (T^*\cF))$ defined by its
principal symbol. In particular, the analytic index does not, in
general, define a map $K_1(C (S^*\cF))\longrightarrow
K_0(C^\infty_c(G))$. A corresponding example for a certain Lie
groupoid is given by Connes in \cite[Chapter II, Section
10.$\gamma$, Proposition 10]{Co} (see also \cite{Carrillo08c}).

The significance of the analytic index ${\rm Ind}\,A$ of an elliptic
$G$-pseudo\-dif\-fe\-ren\-ti\-al operator $A$ with values in
$K_0(C^\infty_c(G))$ consists in the fact that it gives rise to
numerical invariants (higher indices) upon taking the pairing
$\langle {\rm Ind}\,A, \tau \rangle\in \CC$ with an arbitrary
(periodic) cyclic cocycle $\tau$ on $C^\infty_c(G)$.

In \cite[Chapter III, Section 7.$\gamma$, Corollary 13]{Co} a higher
index theorem for tangentially elliptic operators is formulated. It
is an analogue of the higher index $\Gamma$-index theorem proved in
\cite{Co-M90} for $\Gamma$-invariant elliptic operators on a
$\Gamma$-covering of a compact manifold.

\begin{thm}\label{t:l-ind-Connes}
Let $A$ be a tangentially elliptic operator on a compact manifold
equipped with a transversally oriented foliation $(M,\cF)$. Then for
any $\omega\in H^*(BG)$
\begin{equation}\label{e:Phi}
\langle {\rm Ind}\,A, \Phi_*(\omega)\rangle=(2\pi
i)^{-q}\langle\omega,\ch_\tau(\sigma_A)\rangle.
\end{equation}
\end{thm}

Here ${\rm Ind}\,A\in K_0(C^\infty_c(G))$ is the analytic index of
$A$, $[\sigma_A]$ is the class defined in geometric $K$-homology
group $K^*_{\rm top}(M,\cF)$ by the principal symbol of $A$ (see
Subsection~\ref{s:Kth}), $\Phi_* : H^{*}(BG) \to
HP^*(C^\infty_c(G))$ is the map introduced in
Subsection~\ref{s:second}, and $\ch_\tau(\sigma_A)$ is the twisted
Chern character
\[
\ch_\tau(\sigma_A)=\Td(\tau)^{-1}\Phi\ch(\sigma_A).
\]
It follows from (\ref{e:Fch}) that for any $y=[X,E,f]\in
K_{\ast,\tau}(BG)$ we have
\[
\ch_\tau(y)=\bar{f}_*(\ch(E)\cup\Td(X)\cap [X]).
\]
Applying this formula to the class $[\sigma_A]\in K^*_{\rm
top}(M,\cF)$, one can rewrite the formula (\ref{e:Phi}) as
\[
\langle {\rm Ind}\,A, \Phi_*(\omega)\rangle=\langle
\pi_{F!}\ch(\sigma_A) \Td(T\cF\otimes\CC) \bar{p}^*\omega,
[M]\rangle,
\]
where $\bar{p}: M\to BG$ corresponds to the map $p:M\to M/\cF$, and
$\pi_F:T^*\cF\to M$ is the natural projection.

In \cite{Nistor96} a proof of Theorem~\ref{t:l-ind-Connes} is given
by methods of algebraic topology in the case when the foliation is
the horizontal foliation of a flat foliated bundle (see
Example~\ref{ex:suspension}). This proof is based on a formalism of
cyclic type cohomology developed in papers of Cuntz and Quillen
together with results of the papers \cite{Brylinsky:N} and
\cite{Nistor94}.

Various particular cases of Theorem~\ref{t:l-ind-Connes} have been
studied in several papers. In the case when the foliation is given
by the fibres of a fibration $M\to B$, Theorem~\ref{t:l-ind-Connes}
reduces to the index theorem for families of elliptic operators,
Theorem~\ref{t:ASfamily2}. If the foliation $\cF$ has a holonomy
invariant measure $\Lambda$ and $\omega\in H^*(BG)$ is the
corresponding class, then Theorem~\ref{t:l-ind-Connes} reduces to
the measurable index theorem, Theorem~\ref{t:Co79} (see also the
remarks after Theorem~\ref{t:CS}). In the case when $\omega$
corresponds to an invariant form, the given statement was proved in
\cite{Heitsch95}. For the horizontal foliation of a flat foliated
bundle, particular cases were considered in
\cite{DHK,Jiang,Mor-Natsume,Douglas-Kaminker}. We mention the paper
\cite{DHK}, which establishes a connection between the higher index
of the tangential Dirac operator on the horizontal foliation of a
flat foliated bundle $M\to B$ associated with some transverse
cocycle and the relative eta-invariant of the Dirac operator on the
base $B$.

We should note once more that Theorem~\ref{t:l-ind-Connes} has no
topological and geometrical consequences such as a vanishing theorem
or homotopy invariance. For this, one needs to establish a statement
similar to Theorem~\ref{t:l-ind-Connes} for the analytic index with
values in $K_0(C^*_r(M,\cF))$. Topological invariance of the
cocycles defined by the pairing with elements in the image of the
mao $\Phi_*$ was investigated in \cite{Co86} (see also
Subsection~\ref{s:tfclass}).

Before we formulate the basic result of the paper \cite{Co86}, let
recall some information about secondary characteristic classes of
foliations. These classes are given by the characteristic
homomorphism (see, for instance, \cite{Bernshtein-R73})
\[
\chi_{\cF}: H^*(W_q; O(q)) \to H^*(M,\RR)
\]
which is defined for any codimension $q$ foliation $\cF$ on a smooth
manifold $M$, where $H^*(W_q; O(q))$ denotes the relative cohomology
of the Lie algebra $W_q$ of formal vector fields in $\RR^q$. A basic
property of the secondary characteristic classes is their
functoriality: if a smooth map $f: N\to M$ is transverse to a
transversely oriented foliation $\cF$ and $f^*\cF$ is the foliation
on $N$ induced by $f$ (by definition, the leaves of $f^*\cF$ are the
connected components of the pre-images of the leaves of $\cF$ under
the map $f$), then
\[
f^*(\chi_{\cF}(\alpha))=\chi_{f^*\cF}(\alpha), \quad \alpha\in
H^*(W_q; O(q)).
\]

The classifying space $B\Gamma_q$ of the groupoid $\Gamma_q$
classifies codimension $q$ foliations on a given manifold $M$ in the
sense that any foliation $\cF$ on $M$ defines a map $M\to
B\Gamma_q$, and moreover, in the case when $M$ is compact a homotopy
class of maps $M\to B\Gamma_q$ corresponds to a concordance class of
foliations on $M$ \cite{Haefliger72}. For any codimension $q$
foliation $\cF$ on a manifold $M$, the classifying map $M\to
B\Gamma_q$ is obtained as the composition of the map $\bar{p} : M
\to BG$ associated with the projection $p:M\to M/\cF$ and the
universal map $BG\to B\Gamma_q$ (see Subsection~\ref{s:BC}).

By the functoriality of the characteristic homomorphism, it suffices
to know it for the universal foliation on $B\Gamma_q$. This gives
rise to the universal characteristic homomorphism
\[
\chi: H^*(W_q; O(q)) \to H^*(B\Gamma_q,\RR).
\]
For any complete transversal $T$, $\chi$ is represented as the
composition
\[
\chi: H^*(W_q; O(q)) \to H^*(BG^T_T,\RR) \to H^*(B\Gamma_q,\RR).
\]
Since the groupoids $G^T_T$ and $G$ are equivalent,
$H^*(BG^T_T,\RR)\cong H^*(BG,\RR)$, which defines a map
\begin{equation}\label{e:second}
 H^*(W_q; O(q)) \to H^*(BG,\RR).
\end{equation}
Elements of the image of $H^*(W_q; O(q))$ by the map
(\ref{e:second}) will be called secondary characteristic classes.
For computations of $H^*(W_q; O(q))$ see, for instance, \cite{Fuks}.

In \cite{Gor:topology} Gorokhovsky generalized the construction in
Subsection~\ref{s:GV} (which uses the group of modular
automorphisms) of the cyclic cocycle associated with the
Godbillon-Vey class to the case of arbitrary secondary
characteristic classes. This construction makes essential use of the
cyclic cohomology theory developed in the paper \cite{CoM:Hopf} for
Hopf algebras (see Subsection~\ref{para}).

The main result of \cite{Co86} is formulated as follows.

\begin{thm}\label{t:Co86}
Let $(M,\cF)$ be a (not necessarily compact) foliated manifold,
which is transversally oriented. Let $G$ be its holonomy groupoid
and $\pi : BG\to B\Gamma_q$ the classifying map for the
$\Gamma_q$-structure defined by the foliation. Let $\tau$ be the
bundle on $BG$ defined by the normal bundle $\tau$ of $\cF$. Denote
by $\cR$ the subring in $H^*(BG,\CC)$ generated by the Pontryagin
classes of $\tau$, the Chern classes of holonomy equivariant bundles
on $M$ and the secondary characteristic classes.

For any $P\in \cR$ there is an additive map $\varphi_P$ from
$K_*(C^*_r(M,\cF))$ to $\CC$ such that
\begin{equation}\label{e:varphi}
\varphi_P(\mu_r(x))=\langle\Phi\ch(x), P \rangle, \quad x\in
K^*_{\rm top}(M,\cF).
\end{equation}
\end{thm}

We recall that the Chern character $K^*_{\rm top}(M,\cF)\to
H_{*,\tau}(BG)=H_*(B\tau,S\tau)$ is denoted by $\ch$, and $\Phi:
H_{*,\tau}(BG)\to H_*(BG)$ is the Thom isomorphism.
Using~(\ref{e:Fch}), we can write the formula (\ref{e:varphi}) as
follows for any $x=[X,E,f]\in K_{*,\tau}(BG)$:
\begin{equation}\label{e:mur}
\varphi_P(\mu_r(x))  =\langle \ch(E) \Td(TX\oplus
\bar{f}^*\tau)\bar{f}^*P, [X]\rangle.
\end{equation}
In particular, note that the formulae (\ref{e:muP}) and
(\ref{e:muGV}) are particular cases of the formulae
(\ref{e:varphi}).

As a consequence, Theorem~\ref{t:Co86} leads to information about
injectivity of the map $\mu_r$. It also implies that the map
$\varphi_P$ takes integer values on the image of $\mu$. (One should
note that $\varphi_P(K_*(C^*_r(M,\cF)))$ is not, in general,
contained in $\ZZ$.)

Let $A$ be a tangentially elliptic operator on a compact manifold
equipped with a transversally oriented foliation $(M,\cF)$, and let
${\rm Ind}_a(A)\in K_0(C^*_r(M,\cF))$ be its analytic index.
Applying (\ref{e:mur}) to the class $[\sigma_A]$ in the geometric
$K$-homology group $K^*_{\rm top}(M,\cF)$ defined by the principal
symbol of $A$ (see Subsection~\ref{s:Kth}) and taking into account
(\ref{e:mu}), we obtain the higher index theorem for $A$:
\[
\varphi_P({\rm Ind}_a(A))= \langle \pi_{F!}\ch(\sigma_A)
\Td(T\cF\otimes\CC)\Td (\bar\tau\otimes\CC)\bar{p}^*P, [M]\rangle,
\]
where $P\in \cR$, $\bar{p}: M\to BG$ is the map corresponding to the
map $p:M\to M/\cF$, and $\pi_F:T^*\cF\to M$ is the natural
projection.

In \cite{Hil-Skan} a $K$-theoretic analogue of Theorem~\ref{t:Co86}
is obtained in a particular case. Denote by $p$ the natural
projection $M/\cF\to pt$. Suppose that it is a $K$-oriented map
(which is equivalent to the normal bundle $\tau$ being a
$K$-oriented bundle). Let $p!\in KK^*(C^*(M,\cF),\CC)$ be the
corresponding class, which defines the Gysin homomorphism (see
Subsection~\ref{s:BC}). A holonomy equivariant complex vector bundle
$L$ on $(M,\cF)$ defines an element $[L]\in
KK^*(C^*(M,\cF),C^*(M,\cF))$ (see Subsection~\ref{s:vect}). Define
an element $p_L\in KK^*(C^*(M,\cF),\CC)$ by $p_L=[L]\otimes p!$.

\begin{thm}[\cite{Hil-Skan}]\label{t:HS88}
Let $L$ be a holonomy equivariant complex vector bundle on
$(M,\cF)$. For any $y=[X,\cF_X,x,f]\in K^*_{\rm top}(M,\cF)$,
\[
\mu(y)\otimes p_L= x \otimes  [f^*L]\otimes  (p\circ f)!\in \ZZ.
\]
In particular, for any $y=[X,E,f]\in K_{*,\tau}(BG)$,
\begin{equation}\label{e:pl}
\mu(y)\otimes p_L=\langle \ch(E) \ch(\bar{f}^*L) \Td(TX), [X]
\rangle.
\end{equation}
\end{thm}

Theorem~\ref{t:HS88} implies Theorem~\ref{t:Co86} in the case when
$p!\in KK^*(C_r^*(M,\cF),\CC)$, for instance, when the foliation is
amenable (and, therefore, $C_r^*(M,\cF)= C^*(M,\cF)$). One should
note that the homomorphisms $\varphi\circ \lambda_*$ (where $\lambda
: K_*(C^*(M,\cF)) \to K_*(C^*_r(M,\cF))$ is the natural projection)
and $\otimes p_L$ do not, in general, coincide as homomorphisms from
$K_*(C^*_r(M,\cF))$ to $\CC$ and coincide only on the image of
$\mu$.

As shown in \cite{DGKY}, if a foliation $\cF$ is Riemannian and the
normal bundle $\tau$ has a holonomy invariant complex spin
structure, then the element $p_L$ coincides with the $K$-homology
class $[D_L]\in K_*(C^*(M,\cF))$ defined by the transverse $\rm
Spin^c$ Dirac operator $D_{L}$ with coefficients in $L$. Therefore,
the equality (\ref{e:pl}) can be rewritten as follows. Denote by
$\bar{L}$ the bundle on $BG$, corresponding to $L$. For any
$y=[X,E,f]\in K_{*,\tau}(BG)$,
\[
\mu(y)\otimes [D_{L}]= \langle \ch(E) \ch (\bar{f}^*L) \Td(TX), [X]
\rangle.
\]
The authors of \cite{DGKY} propose viewing this formula as an index
formula for the transverse Dirac operator $D_{L}$ and give examples
where this formula could be useful for computation of the
distributional index of transversally elliptic operators.

In \cite{Gor-Lott03,Gor-Lott04} there is a proof of
Theorem~\ref{t:l-ind-Connes} for a tangential Dirac operator which
generalizes Bisnut's proof in \cite{BismutInv85} of the local index
theorem for families of Dirac operators (see
Subsection~\ref{s:fam}).

Let $(M,\cF)$ be a compact foliated manifold. Suppose that the
leaves of $\cF$ are even-dimensional spin manifolds. Choose a
Riemannian metric in the fibres of the bundle $T\cF$, and denote by
$F(T\cF)$ the associated spinor bundle. Let $V$ be a Hermitian
vector bundle on $M$ equipped with a Hermitian connection
$\nabla^V$. Consider the Clifford bundle $\cE=F(T\cF)\otimes V$ over
the Clifford algebra of $T\cF$ and the associated tangential Dirac
operator $D_{\cE}$.

In \cite{Gor-Lott03} (see also a more precise formulation in
\cite{Gor-Lott04}) the authors proved the local index theorem for
families of Dirac operators invariant under a free, proper,
cocompact action of an \'etale groupoid. The case of a tangential
Dirac operator is reduced to this case using the following
construction. Let $T$ be a complete transversal, and $G^T_T$ the
associated reduced holonomy groupoid. The set $P=G_T$ is a smooth
manifold which is equipped with the natural free proper action of
$G^T_T$ given by right multiplication (see
Subsection~\ref{s:morita1}). We remark that the orbit space
$P/G^T_T$ of this action coincides with $M$. The map $s$ defines a
submersion $\pi : G_T\to (G^T_T)^{(0)}=T$. The submersion $\pi :
G_T\to T$ is a $G^T_T$-equivariant map if we consider the action of
$G^T_T$ on $T$ given by right multiplication. Denote by
$Z_x=\pi^{-1}(x)=G_x$ the fibre of the fibration $\pi : G_T\to T$ at
$x\in T$, and by $TZ$ the vertical tangent space. The projection
$P=G_T\to M=P/G^T_T$ takes fibres of the fibration $\pi$ to leaves
of the foliation $\cF$ on $M$ and $TZ$ to the tangent bundle $F$ of
$\cF$.

The leafwise spin structure in $T\cF$ is lifted to a
$G^T_T$-invariant spin structure in $TZ$. Denote by $F(TZ)$ the
associated spinor bundle and by $\widehat{V}$ the lift of $V$ to
$P$. Consider the leafwise Clifford bundle
$\widehat{\cE}=F(TZ)\otimes \widehat{V}$ and the associated
tangential Dirac operator $D_{\widehat{\cE}}$ on $P$. The operator
$D_{\widehat{\cE}}$ determines a family of Dirac operators acting
along the fibres of the fibration $\pi : P\to T$ and invariant under
the action of $G^T_T$.

The index $\Ind D_{\widehat{\cE}}$ of the operator
$D_{\widehat{\cE}}$ is well defined as an element of the group
$K_0(C^\infty_c(G^T_T)\otimes \cR)$, where $\cR$ is the algebra of
rapidly decreasing infinite real matrices. Consider the differential
graded algebra $\Omega^*_c(G^T_T)$ introduced in
Subsection~\ref{s:second}. For an arbitrary closed graded trace
$\eta$ on $\Omega^*_c(G^T_T)$ the authors of \cite{Gor-Lott03}
define a pairing of the Chern character $\ch(\Ind D_{\cE})\in
HP(C^\infty_c(G^T_T))$ with $\eta$ and an element $\Phi_\eta\in
H^*_\tau(BG)$ associated with $\eta$. As above, denote by $\bar{p} :
M\to BG$ the map corresponding to the natural projection $p : M\to
M/\cF$. Then the following theorem holds.

\begin{thm}[\cite{Gor-Lott03}]\label{t:gorlott}
\[
\langle \ch(\Ind D_{\widehat{\cE}}), \eta\rangle =
\int_M\hat{A}(T\cF)\ch(V)\bar{p}^*\Phi_\eta
\]
\end{thm}

The proof of this Theorem is based on a non-commutative equivariant
version of the Bismut constructions (see Subsection~\ref{s:fam})
applied to the $G^T_T$-equivariant submersion $\pi : P \to T$. Here
an important role is played by a construction of a certain
differentiation
\[
\nabla^{0,1} : C^\infty_c(P,\widehat{\cE})\to
\Omega^{0,1}_c(G^T_T)\otimes_{C^\infty_c(G^T_T)}C^\infty_c(P,\widehat{\cE}).
\]
in ``non-commutative'' directions. It is the use of this
non-commutative connection that enables one to include into
considerations cohomology classes in dimension greater than the
codimension of the foliation (such as the Godbillon-Vey class whivh
is a three-dimensional cohomology class for a codimension one
foliation).

In \cite{Gor-Lott04} a direct proof of Theorem~\ref{t:Co86} for the
tangential Dirac operator is given in the particular case when the
class $\omega$ is defined by a holonomy invariant transverse
current, without using an auxiliary choice of a complete
transversal. The authors make use of a differential graded algebra
similar to the algebra constructed in Subsection~\ref{s:tfclass}. In
the papers \cite{Leichtnam-Piazza05,Leichtnam-Piazza05a}
Theorem~\ref{t:gorlott} was extended to the case when $M$ is a
manifold with boundary and the foliation is transversal to the
boundary.

As mentioned above, a tangentially elliptic operator $D$ on a
compact foliated manifold $(M,\cF)$ defines a class $[D]\in KK(C(M),
C^*_r(M,\cF))$. This class can be represented by an explicit
$p$-summable quasi-homomorphism $\psi_D$ from $C^\infty(M)$ to
$C^\infty_c(G)$ in the sense of the papers \cite{Nistor91,Nistor93}.
In \cite{Gor06} the bivariant Chern character of this
quasi-homomorphism introduced in \cite{Nistor91,Nistor93} is
computed.

There is another approach to higher index theorems for tangentially
elliptic operators based on the use of the Haefliger cohomology (see
Subsection~\ref{s:holonomy}).

Let $(M,\cF)$ be a compact foliated manifold. Suppose that the
dimension of $\cF$ is even, it is oriented, and it has a spin
structure. Let $\cE$ be a Hermitian vector bundle on $M$, $D_\cE$
the corresponding leafwise Dirac operator on $M$ with coefficients
in the bundle $\cE$, and ${\rm Ind}(D_\cE)\in K_0(C^\infty_c(G))$
its analytic index. In \cite{Heitsch95} the Chern character
$\overline{\ch}({\rm Ind}(D_\cE))$ is defined as an element of the
Haefliger cohomology $H^*_c(M/\cF)$ of $\cF$. The construction of
the Chern character is a direct modification of the Bismut
construction (see Subsection~\ref{s:fam}). It makes use of an
analogue of the Bismut superconnection associated with the operator
$D_\cE$ and of the heat operator determined by the curvature of this
superconnection.

Under some additional restrictions on the foliation it is proved in
\cite{Heitsch-L99} that the following equality holds in the
Haefliger cohomology $H^*_c(M/\cF)$:
\[
\overline{\ch}({\rm Ind}(D_\cE))=\frac{1}{(2\pi i)^{p/2}}\int_\cF
\hat{A}(T\cF,\nabla^{T\cF})\ch(\cE,\nabla^\cE).
\]

The Chern character $\ch_a: K_0(C^\infty_c(G))\to H^*_c(M/\cF)$ with
values in the Haefliger cohomology was constructed by the authors of
\cite{Benameur-Heitsch04}, and this enables them to translate the
Connes-Skandalis index theorem for tangentially elliptic operators
to the language of the Haefliger cohomology, using the results of
the paper \cite{Co86}.

Denote by $\hat{\cF}$ the dimension $p$ foliation on the manifold
$M\times \RR^{2k}$ induced by $\cF$. The holonomy groupoid
$G^{\RR^{2k}}$ of this foliation coincides with $G\times \RR^{2k}$.
Taking the composition $\ch_a: K_0(C^\infty_c(G^{\RR^{2k}}))\to
H^*_c(M\times \RR^{2k}/\hat{\cF})$ with integration along
$\RR^{2k}$, one obtains a map
\[
\ch^{\RR^{2k}}_a: K_0(C^\infty_c(G^{\RR^{2k}}))\to H^*_c(M/\cF).
\]
We`note that in general there is no Bott isomorphism between
$K_0(C^\infty_c(G^{\RR^{2k}}))$ and $K_0(C^\infty_c(G))$ in this
case. For sufficiently large $k$ the map $\pi_! : K^0_c(T\cF)\to
K_0(C^\infty_c(G^{\RR^{2k}}))$ is well defined.

\begin{thm}
For any $u\in K^0_c(T\cF)$
\[
\ch^{\RR^{2k}}_a\circ
\pi_!(u)=(-1)^p\int_\cF\pi_{\cF!}(\ch(u))\Td(T\cF\otimes \CC)\in
H^*_c(M/\cF),
\]
where $\pi_{\cF!} : H^*_c(T\cF,\RR)\to H^*(M,\RR)$ is integration
along the fibres of the bundle $\pi_\cF : T\cF\to M$.
\end{thm}

Suppose that the dimension of the foliation $\cF$ is even, it is
oriented and has a spin structure. Let $\cE$ be a Hermitian vector
bundle on $M$, and $D_\cE$ the corresponding leafwise Dirac operator
on $M$ with coefficients in the bundle $\cE$. In
\cite{Benameur-Heitsch2}, assuming that the foliation $\cF$ is
Riemannian and the bundle ${\rm Ind}(D_\cE)$ is transversally
smooth, the authors proved the coincidence of the two Chern
characters defined in previous papers:
\[
\overline{\ch}({\rm Ind}(D_\cE)) = \ch_a({\rm Ind}(D_\cE)) \in
H^*_c(M/\cF).
\]

Finally, in  \cite{Benameur-Heitsch3}, under the assumption that the
family of projections to leafwise harmonic forms in the middle
dimension is transversally smooth in a certain sense, the authors
defined a higher harmonic signature of an even-dimensional oriented
Riemannian foliation on a compact Riemannian manifold and proved its
invariance under leafwise homotopies.

We also mention the papers \cite{Benameur02,Benameur03}, which
concern cyclic versions of the Lefschetz formula for diffeomorphisms
which take each leaf of the foliation to itself.

In all the theorems mentioned above, the higher indices ${\rm
Ind}_\tau (A)$ of an elliptic $G$-operator $A$ depend only on the
class $[\sigma_A]\in K_0(T^*\cF)$ defined by its principal symbol,
and this, as has already been noted above, does not hold for just
any cyclic cocycle $\tau$. Following the paper \cite{Carrillo08c},
we say that a (periodic) cyclic cocycle $\tau$ on the algebra
$C^\infty_c(G)$ can be localized if the value of the functional
$A\mapsto \langle {\rm Ind}\,A, \tau \rangle\in \CC$, defined on the
set of elliptic $G$-operators, depends only on $[\sigma_A]$. In this
case there is a well defined map ${\rm Ind}_\tau : K_0(T^*\cF)\to
\CC$ satisfying the condition
\[
\langle {\rm Ind}\,A, \tau \rangle={\rm Ind}_\tau([\sigma_A]).
\]
The map ${\rm Ind}_\tau$ is called a higher localized index
associated with $\tau$.

We say that a $(k+1)$-linear functional $\tau$ on the space
$C^\infty_c(G)$ is bounded if it extends to a continuous
$(k+1)$-linear functional $\tau_m$ on the space $C^m_c(G)$ for some
$m\in \NN$. Many geometric cocycles (such as the group cocycles, the
transversal fundamental class, and the cocycles defined by the
Godbillon-Vey class and by secondary characteristic classes of the
foliation) are bounded cocycles. It is proved in \cite{Carrillo08a}
that any bounded cyclic cocycle on the algebra $C^\infty_c(G)$ can
be localized. The proof of this fact is based on the construction of
the tangent groupoid $G^T$ associated with the holonomy groupoid $G$
and of a certain algebra $\cS_c(G^T)$ of functions on $G^T$. This
algebra was constructed in \cite{Carrillo08b}. It is a strict
deformation quantization of the Schwartz algebra $\cS(T^*\cF)$. One
can construct an analytic index map
\[
{\rm Ind}_a : K(T^*\cF)\longrightarrow K_0(C^*_r(M,\cF))
\]
by means of the tangent groupoid $G^T$ and its $C^*$-algebra
$C^*_r(G^T)$ (see \cite{Hil-Skan,MP97,NWX99}). Using this
construction, the author of \cite{Carrillo08c} derived a formula for
the higher localized index ${\rm Ind}_\tau$ associated with a
bounded cocycle $\tau$ in terms of an asymptotic limit of cocycles
on the algebra $\cS_c(G^T)$.

All the facts about higher localized indices mentioned above hold
for an arbitrary Lie groupoid $G$. In this case the role of the
cotangent space $T^*\cF$ is played by the space $A^*G$, where $AG$
is the Lie algebroid of the groupoid $G$. It is proved in
\cite{Carrillo08a} that if the groupoid $G$ is \'etale, then any
cyclic cocycle on the algebra $C^\infty_c(G)$ is bounded. This
statement can be applied, for instance, to the reduced holonomy
groupoid $G^T_T$ associated with some complete transversal $T$.


\begin{thebibliography}{100}
\bibitem{AlvKordy2002}
J.~A. {\'A}lvarez~L{\'o}pez, Yu.~A. Kordyukov,
 ``Distributional {B}etti numbers of transitive foliations of
  codimension one'',
{\em Foliations: geometry and dynamics (Warsaw, 2000)}, World Sci.
Publ., River Edge, NJ, 2002, 159--183

\bibitem{AlvKordy2008}
J.~A. {\'A}lvarez~L{\'o}pez, Yu.~A. Kordyukov,
 ``Lefschetz distribution for {Lie} foliations'',
{\em $C^*$-algebras and elliptic theory. II}, Trends in
  Mathematics, Birkh\"auser, Basel, 2008, 1--40

\bibitem{Atiyah-global}
M.~F. Atiyah,
 ``Global theory of elliptic operators'',
{\em Proceedings of the International Symposium on Functional
Analysis, Tokyo}, University Tokyo Press, Tokyo, 1969, 21--30

\bibitem{At:teo}
M.~F. Atiyah,
 {\em Elliptic operators and compact groups},
 Lect. Notes in Math. {\bf 401}, Springer,
  Berlin, 1974

\bibitem{At:discrete}
M. F. Atiyah,  ``Elliptic operators, discrete groups and von Neumann
algebras'',  {\em Colloque ``Analyse et Topologie'' en l'Honneur de
Henri Cartan (Orsay, 1974),}   Asterisque, {\bf 32-33}, Soc. Math.
France, Paris, 1976, 43--72

\bibitem{Atiyah-Bott}
M.~F. Atiyah, R.~Bott,
 ``A {L}efschetz fixed point formula for elliptic complexes. {I}'',
 {\em Ann. of Math. (2)}, {\bf 86} (1967), 374--407

\bibitem{Atiyah-BottII}
M.~F. Atiyah, R.~Bott,
 ``A {L}efschetz fixed point formula for elliptic complexes. {II}.
  {A}pplications'',
 {\em Ann. of Math. (2)}, {\bf 88} (1968), 451--491

\bibitem{AtiyahBottPatodi}
M. F.~Atiyah, R.~Bott, V.~K. Patodi,
 ``On the heat equation and the index theorem'',
 {\em Invent. Math.}, {\bf 19} (1973), 279--330

\bibitem{APS:SARG1}
M.~F. Atiyah, V.~K. Patodi, I.~M. Singer,
 ``Spectral asymmetry and {R}iemannian geometry. {I}'',
 {\em Math. Proc. Cambridge Philos. Soc.}, {\bf 77} (1975), 43--69

\bibitem{APS:SARG3}
M.~F. Atiyah, V.~K. Patodi, I.~M. Singer,
 ``Spectral asymmetry and {R}iemannian geometry. {III}'',
 {\em Math. Proc. Cambridge Philos. Soc.}, {\bf 79}:1 (1976), 71--99

\bibitem{ASII}
M.~F. Atiyah, G.~B. Segal,
 ``The index of elliptic operators. {II}'',
 {\em Ann. of Math. (2)}, {\bf 87} (1968), 531--545

\bibitem{ASI}
M.~F. Atiyah, I.~Singer,
 ``The index of elliptic operators. {I}'',
 {\em Ann. of Math.}, {\bf 87} (1968), 484--530

\bibitem{ASIII}
M.~F. Atiyah, I.~Singer,
 ``The index of elliptic operators. {III}'',
 {\em Ann. of Math.}, {\bf 87}(1968), 546--604

\bibitem{ASIV}
M.~F. Atiyah, I.~Singer,
 ``The index of elliptic operators. {IV}'',
 {\em Ann. of Math.}, {\bf 93} (1971), 119--138

\bibitem{Baaj-Julg}
S.~Baaj, P.~Julg,
 ``Th\'eorie bivariante de {K}asparov et op\'erateurs non born\'es dans
  les {$C\sp{\ast} $}-modules hilbertiens'',
 {\em C. R. Acad. Sci. Paris S\'er. I Math.} {\bf 296}:21 (1983), 875--878

\bibitem{Baum-Connes83}
P. Baum, A. Connes,
 ``Leafwise homotopy equivalence and rational {P}ontrjagin
 classes'', {\em Foliations (Tokyo, 1983)}, Adv. Stud. Pure
  Math., {\bf 5},  North-Holland, Amsterdam, 1985, 1--14

\bibitem{Baum-Connes}
P.~Baum, A.~Connes,
 ``Geometric {$K$}-theory for {L}ie groups and foliations'',
 {\em Enseign. Math. (2)}, {\bf 46}:1-2 (2000), 3--42

\bibitem{Baum-D82}
P.~Baum, R.~G. Douglas,
 ``{$K$} homology and index theory'',
{\em Operator algebras and applications, Part I (Kingston, Ont.,
  1980)}, Proc. Sympos. Pure Math., {\bf 38}, Amer.
  Math. Soc., Providence, R.I., 1982, 117--173

\bibitem{baum-higson-schick}
P. Baum, N. Higson, T. Schick, ``On the equivalence of geometric and
analytic $K$-homology'', {\em Pure Appl. Math. Q.}  {\bf 3}:1
(2007), 1--24.


\bibitem{BG}
R.~Beals, P.~Greiner, {\em Calculus on {H}eisenberg manifolds},
Annals of Math. Studies, {\bf 119}, Princeton Univ. Press,
Princeton, 1988


\bibitem{Benameur97}
M.-T. Benameur,
 ``A longitudinal {L}efschetz theorem in {$K$}-theory'',
 {\em $K$-Theory}, {\bf 12}:3 (1997), 227--257

\bibitem{Benameur02}
M.-T. Benameur,
 ``Cyclic cohomology and the family {L}efschetz theorem'',
 {\em Math. Ann.}, {\bf 323}:1 (2002), 97--121

\bibitem{Benameur03}
M.-T. Benameur,
 ``A higher {L}efschetz formula for flat bundles'',
 {\em Trans. Amer. Math. Soc.}, {\bf 355}:1 (2003), 119--142

\bibitem{Ben-C-Ph-R-S}
M.-T. Benameur, A.~L. Carey, J. Phillips, A. Rennie, F.~A.
  Sukochev, K.~P. Wojciechowski,
 ``An analytic approach to spectral flow in von {N}eumann
 algebras'', {\em Analysis, geometry and topology of elliptic operators},
 World Sci. Publ., Hackensack, NJ, 2006, 297--352

\bibitem{Benameur-Fack06}
M.-T. Benameur, T. Fack,
 ``Type {II} non-commutative geometry. {I}. {D}ixmier trace in von
  {N}eumann algebras'',
 {\em Adv. Math.}, {\bf 199}:1 (2006), 29--87

\bibitem{Benameur-Heitsch04}
M.-T. Benameur, J.~L. Heitsch,
 ``Index theory and non-commutative geometry. {I}. {H}igher families
  index theory'', {\em $K$-Theory}, {\bf 33}:2 (2004), 151--183

\bibitem{Benameur-Heitsch2}
M.-T. Benameur, J.~L. Heitsch,
 ``Index theory and {Non-Commutative Geometry}. {II.} {D}irac operators
  and index bundles'', {\em J. K-Theory} {\bf 1}:2 (2008), 305--356

\bibitem{Benameur-Heitsch3}
M.-T. Benameur, J.~L. Heitsch,
 ``The higher harmonic signature for foliations {I}. {The} untwisted
  case'', preprint arXiv:0711.0352, 2007

\bibitem{BGV}
N.~Berline, E.~Getzler, M.~Vergne,
 {\em Heat Kernels and the Dirac Operator}, Grundl. der Math. Wiss., {\bf 298},
 Springer, Berlin, 1992

\bibitem{BerlineVergne87}
N. Berline, M. Vergne,
 ``A proof of {B}ismut local index theorem for a family of {D}irac
  operators'',
 {\em Topology}, {\bf 26}:4 (1987), 435--463

\bibitem{BerlineVergne96a}
N. Berline, M. Vergne,
 ``The {C}hern character of a transversally elliptic symbol and the
  equivariant index'',
 {\em Invent. Math.}, {\bf 124}:1-3 (1996), 11--49

\bibitem{BerlineVergne96b}
N. Berline, M. Vergne,
 ``L'indice \'equivariant des op\'erateurs transversalement
 elliptiques'',
 {\em Invent. Math.}, {\bf 124}:1-3 (1996), 51--101

\bibitem{Bernshtein-R73}
J. H. Bernstein, B. I. Rozenfel'd, ``Homogeneous spaces of
infinite-dimensional Lie algebras and characharacteristic classes of
foliations'', {\em UMN}, {\bf 28}:4 (1973), 103--138; translation
in: {\em Russian Mathematical Surveys}, {\bf 28}:4 (1973), 107--142

\bibitem{BismutInv85}
J.-M. Bismut,
 ``The {Atiyah}-{Singer} index theorem for families of {Dirac}
  operators: {Two} heat equation proofs'',
 {\em Invent. Math.}, {\bf 83} (1985), 91--151

\bibitem{Bost}
J.~Bost,
 ``Principe {d'Oka}, {$K$}-th\'eorie et syst\`emes dynamiques non
  commutatifs'',
 {\em Invent. Math.}, {\bf 101} (1990), 261--333

\bibitem{Bott71}
R.~Bott,
 ``On topological obstructions to integrability'',
{\em Actes du Congr\`es Intern. Math. Nice, 1970, tome 1},
  Gauthiers-Villars, Paris, 1971, 27--36

\bibitem{Breuer1}
M. Breuer,
 ``Fredholm theories in von {N}eumann algebras. {I}'',
 {\em Math. Ann.}, {\bf 178} (1968), 243--254

\bibitem{Breuer2}
M. Breuer,
 ``Fredholm theories in von {N}eumann algebras. {II}'',
 {\em Math. Ann.}, {\bf 180} (1969), 313--325

\bibitem{BDF}
L. G. Brown, R. G. Douglas, P. A. Fillmore,  ``Extensions of
$C\sp*$-algebras and $K$-homology'', {\em Ann. of Math. (2)}, {\bf
105} (1977), 265--324

\bibitem{Brylinsky:N}
J.-L. Brylinski, V.~Nistor,
 ``Cyclic cohomology of \'etale groupoids'',
 {\em $K$-Theory}, {\bf 8}:4 (1994), 341--365

\bibitem{CPS2}
A.~L. Carey, J.~Phillips, F.~A. Sukochev,
 ``Spectral flow and {D}ixmier traces'',
 {\em Adv. Math.}, {\bf 173}:1 (2003), 68--113

\bibitem{CPRS2}
A.~L. Carey, J.~Phillips, A. Rennie, F.~A. Sukochev,
 ``The local index formula in semifinite von {N}eumann algebras. {I}.
  {S}pectral flow'',
 {\em Adv. Math.}, {\bf 202}:2 (2006), 451--516

\bibitem{CPRS3}
A.~L. Carey, J.~Phillips, A. Rennie, F.~A. Sukochev,
 ``The local index formula in semifinite von {N}eumann algebras. {II}.
  {T}he even case'',
 {\em Adv. Math.}, {\bf 202}:2 (2006), 517--554

\bibitem{Carey-S}
A. L. Carey, F. A. Sukochev, ``Dixmier traces and some applications
in non-commutative geometry'', {\em UMN}, {\bf 61}:6 (2006),
45--110; translation in: {\em Russian Mathematical Surveys}, {\bf
61}:6 (2006), 1039--1099

\bibitem{Carriere84}
Y. Carri{\`e}re,
 ``Flots riemanniens'', {\em Transversal structure of foliations} (Toulouse,
 1982), Ast\'erisque, {\bf 116}, 1984, 31--52


\bibitem{Carrillo08a}
P. Carrillo Rouse, ``A Schwartz type algebra for the Tangent
Groupoid'', {\em K-Theory and Non-commutative Geometry}, European
Mathematical Society (EMS), Z\"urich, 2008, 181--200

\bibitem{Carrillo08b}
P. Carrillo Rouse, ``Compactly supported analytic indices for Lie
groupoids'', preprint arXiv:0803.2060, 2008

\bibitem{Carrillo08c}
P. Carrillo Rouse, ``Higher localized analytic indices and strict
deformation quantization'', preprint arXiv:0810.4480, 2008

\bibitem{CDSS}
L. A. Coburn, R. G. Douglas, D. G. Schaeffer, I. M. Singer,
``$C\sp{\ast}$-algebras of operators on a half-space. II. Index
theory'' {\em Publ. Math.} {\bf 40} (1971), 69--79

\bibitem{CoLNP80}
A.~Connes,
 ``The von {N}eumann algebra of a foliation'',
{\em Mathematical problems in theoretical physics (Proc. Internat.
  Conf., Univ. Rome, Rome, 1977)}, Lecture Notes in Phys.,
  {\bf 80}, Springer, Berlin, 1978, 145--151

\bibitem{Co79}
A.~Connes,
 ``Sur la th\'eorie non commutative de l'int\'egration'',
{\em Alg\`ebres d'op\'erateurs (S\'em., Les Plans-sur-Bex, 1978)},
Lecture Notes in Math., {\bf 725}, Springer, Berlin, 1979, 19--143

\bibitem{CoCRAS}
A.~Connes,
 ``{$C\sp{\ast} $} alg\`ebres et g\'eom\'etrie diff\'erentielle'',
 {\em C. R. Acad. Sci. Paris S\'er. A-B}, {\bf 290}:13 (1980), A599--A604

\bibitem{Co81}
A.~Connes,
 ``An analogue of the {Thom} isomorphism for crossed products of a
  {$C^*$} algebra by an action of {$\mathbb R$}'',
 {\em Adv. in Math.}, {\bf 39} (1981), 31--55

\bibitem{Co:survey}
A.~Connes,
 ``A survey of foliations and operator algebras'',
{\em Operator algebras and applications, Part I (Kingston, Ont.,
  1980)}, Proc. Sympos. Pure Math., {\bf 38}, Amer.
  Math. Soc., Providence, R.I., 1982, 521--628

\bibitem{Co86}
A.~Connes,
 ``Cyclic cohomology and the transverse fundamental class of a
  foliation'',
{\em Geometric methods in operator algebras (Kyoto, 1983)}, Pitman
Res. Notes in Math., {\bf 123}, Longman, Harlow,
  1986, 52--144

\bibitem{Co:nc}
A.~Connes,
 ``Noncommutative differential geometry'',
 {\em Publ. Math.}, {\bf 62} (1986), 41--144

\bibitem{Co-action}
A.~Connes,
 ``The action functional in noncommutative geometry'',
 {\em Commun. Math. Phys.}, {\bf 117} (1988), 673--683

\bibitem{Co}
A.~Connes,
 {\em Noncommutative geometry}.
 Academic Press Inc., San Diego, CA, 1994.

\bibitem{Sp-view}
A.~Connes,
 ``Geometry from the spectral point of view'',
 {\em Lett. Math. Phys.}, {\bf 34} (1995), 203--238

\bibitem{Connes-gravity}
A.~Connes,
 ``Gravity coupled with matter and the foundation of non commutative
  geometry'',
 {\em Commun. Math. Phys.}, {\bf 182} (1996), 155--176

\bibitem{Connes2000}
A.~Connes,
 ``Noncommutative geometry---year 2000'',
 {\em Geom. Funct. Anal.}, (Special Volume, Part II) (2000), 481--559

\bibitem{ConnesLNM1831}
A.~Connes,
 ``Cyclic cohomology, noncommutative geometry and quantum group
  symmetries'',
{\em Non-commutative geometry}, Lecture Notes
  in Math., {\bf 1831},  Springer, Berlin, 2004, 1--71

\bibitem{Connes:reconstruction}
A.~Connes, ``On the spectral characterization of manifolds'',
preprint arXiv:0810.2088, 2008

\bibitem{CoKr1}
A. Connes, D. Kreimer,
 ``Hopf algebras, renormalization and noncommutative geometry'',
 {\em Comm. Math. Phys.}, {\bf 199}:1 (1998), 203--242

\bibitem{CoKr2}
A. Connes, D. Kreimer,
 ``Renormalization in quantum field theory and the {R}iemann-{H}ilbert
  problem. {I}. {T}he {H}opf algebra structure of graphs and the main
  theorem'',
 {\em Comm. Math. Phys.}, {\bf 210}:1 (2000), 249--273

\bibitem{CoKr3}
A. Connes, D. Kreimer,
 ``Renormalization in quantum field theory and the {R}iemann-{H}ilbert
  problem. {II}. {T}he {$\beta$}-function, diffeomorphisms and the
  renormalization group'',
 {\em Comm. Math. Phys.}, {\bf 216}:1 (2001), 215--241

\bibitem{Co-M90}
A.~Connes, H.~Moscovici,
 ``Cyclic cohomology, the {N}ovikov conjecture and hyperbolic
 groups'',
 {\em Topology}, {\bf 29}:3 (1990), 345--388

\bibitem{Co-M}
A.~Connes, H.~Moscovici,
 ``The local index formula in noncommutative geometry'',
 {\em Geom. and Funct. Anal.}, {\bf 5} (1995), 174--243

\bibitem{CoM:Hopf}
A.~Connes, H.~Moscovici,
 ``Hopf algebras, cyclic cohomology and the transverse index
 theorem'',
 {\em Commun. Math. Phys.}, {\bf 198} (1998), 199--246

\bibitem{CoMEssays}
A.~Connes, H.~Moscovici,
 ``Differentiable cyclic cohomology and {H}opf algebraic structures in
  transverse geometry'',
{\em Essays on geometry and related topics, Vol. 1, 2}, Monogr.
Enseign. Math., {\bf 38}, Enseignement Math., Geneva, 2001, 217--255

\bibitem{CoM:mmj1}
A.~Connes, H.~Moscovici,
 ``Modular {H}ecke algebras and their {H}opf symmetry'',
 {\em Mosc. Math. J.}, {\bf 4}:1 (2004), 67--109, 310

\bibitem{CoM:mmj2}
A.~Connes, H.~Moscovici,
 ``Rankin-{C}ohen brackets and the {H}opf algebra of transverse
  geometry'',
 {\em Mosc. Math. J.}, {\bf 4}:1 (2004) 111--130, 311

\bibitem{CoM:twisted}
A.~Connes, H.~Moscovici,
 ``Type {III} and spectral triples'',
{\em Traces in Geometry, Number Theory and Quantum Fields},
  Aspects of Mathematics E 38, Vieweg Verlag, 2008, 57--71

\bibitem{Co-skandal}
A.~Connes, G.~Skandalis,
 ``The longitudinal index theorem for foliations'',
 {\em Publ. Res. Inst. Math. Sci.}, {\bf 20}:6 (1984), 1139--1183

\bibitem{Connes-Tak}
A.~Connes, M.~Takesaki,
 ``The flow of weights on factors of type {III}'',
 {\em T\^{o}hoku Math. J.}, {\bf 29} (1977), 473--575


\bibitem{CrainicM01}
M.~Crainic, I.~Moerdijk,
 ``Foliation groupoids and their cyclic homology'',
 {\em Adv. Math.}, {\bf 157}:2 (2001), 177--197

\bibitem{CMX}
R.E. Curto, P.S. Muhly, J, Xia, ``Toeplitz operators on flows'',
{\em J. Funct. Anal.} {\bf 93}:2  (1990), 391--450

\bibitem{Deninger07}
Ch. Deninger, ``Analogies between analysis on foliated spaces and
arithmetic geometry'', preprint arXiv:0709.2801, 2007

\bibitem{DenSinghof1}
Ch. Deninger, W. Singhof, ``A note on dynamical trace formulas'',
{\em Dynamical, spectral, and arithmetic zeta functions (San
Antonio, TX, 1999).} Contemp. Math., {\bf 290}, Amer. Math. Soc.,
Providence, 2001, 41--55


\bibitem{Dixmier66}
J.~Dixmier,
 ``Existence de traces non normales'',
 {\em C.R. Acad. Sci. Paris Ser A-B}, {\bf 262} (1966), A1107--A1108

\bibitem{DHK}
R.~G. Douglas, S.~Hurder, J.~Kaminker,
 ``Cyclic cocycles, renormalization and eta invariants'',
 {\em Invent Math.}, {\bf 103} (1991), 101--181

\bibitem{DHK2}
R.~G. Douglas, S.~Hurder, J.~Kaminker,
 ``The longitudinal cocycle and the index of the {T}oeplitz
 operators'', {\em J. Funct. Anal.}, {\bf 101} (1991), 120--144

\bibitem{DGKY}
R.~G. Douglas, J.~F. Glazebrook, F.~W. Kamber, G.~L. Yu,
 ``Index formulas for geometric {D}irac operators in {R}iemannian
  foliations'',
 {\em $K$-Theory}, {\bf 9}:5 (1995), 407--441

\bibitem{Douglas-Kaminker}
R.~G. Douglas, J. Kaminker,
 ``Induced cyclic cocycles and higher eta invariants'',
 {\em J. Funct. Anal.}, {\bf 147}:2 (1997), 301--326

\bibitem{F-Sk}
T.~Fack, G.~Skandalis,
 ``Sur les repr\'{e}sentations et ideaux de la {$C^*$}-alg\`{e}bre d'un
  feuilletage'',
 {\em J.Operator Theory}, {\bf 8} (1982), 95--129

\bibitem{Fedosov-Shubin1}
B. V. Fedosov, M. A. Shubin, ``The index of random elliptic
operators. I'', {\em Mat. Sb. (N.S.)}, {\bf 106}:1 (1978), 108--140;
translation in: {\em Mathematics of the USSR-Sbornik}, {\bf 34}:5
(1978), 671--699

\bibitem{Fedosov-Shubin2}
B. V. Fedosov, M. A. Shubin, ``The index of random elliptic
operators. II'', {\em Mat. Sb. (N.S.)}, {\bf 106}:3 (1978),
455--483; translation in: {\em Mathematics of the USSR-Sbornik},
{\bf 35}:1 (1979), 131--156

\bibitem{Ferry:RR} S. Ferry, A. Ranicki, J. Rosenberg, ``A history
and survey of the Novikov conjecture'', {\em Novikov conjectures,
index theorems and rigidity, Vol. 1 (Oberwolfach, 1993),} London
Math. Soc. Lecture Note Ser., {\bf 226}, Cambridge Univ. Press,
Cambridge, 1995, 7--66

\bibitem{Fox-Haskell94a}
J. Fox, P. Haskell,
 ``The index of transversally elliptic operators for locally free
  actions'',
 {\em Pacific J. Math.}, {\bf 164}:1 (1994), 41--85

\bibitem{Fox-Haskell94b}
J. Fox, P. Haskell,
 ``The index of transversally elliptic operators on locally homogeneous
  spaces of finite volume'',
 {\em Michigan Math. J.}, {\bf 41}:2 (1994), 323--336

\bibitem{Fried87}
D. Fried, ``Lefschetz formulas for flows'', {\em The Lefschetz
centennial conference, Part III (Mexico City, 1984)}, Contemp.
Math., {\bf 58}, III, Amer. Math. Soc., Providence, RI, 1987, 19--69

\bibitem{Fuks}
D.~B. Fuchs. {\em Cohomology of infinite-dimensional {L}ie
algebras}. ``Nauka'', Moscow, 1984; translation in:  Contemporary
Soviet Mathematics. Consultants Bureau, New York, 1986.

\bibitem{Getzler86}
E. Getzler,
 ``A short proof of the local {A}tiyah-{S}inger index theorem'',
 {\em Topology}, {\bf 25}:1 (1986), 111--117

\bibitem{Getzler93}
E. Getzler,
 ``The odd {C}hern character in cyclic homology and spectral flow'',
 {\em Topology}, {\bf 32}:3 (1993), 489--507

\bibitem{Gilkey73}
P.~B. Gilkey,
 ``Curvature and the eigenvalues of the {L}aplacian for elliptic
  complexes'',
 {\em Advances in Math.}, {\bf 10} (1973), 344--382

\bibitem{Gilkey95}
P.~B. Gilkey, {\em Invariance theory, the heat equation, and the
Atiyah-Singer index theorem.} Second edition. Studies in Advanced
Mathematics. CRC Press, Boca Raton, FL, 1995


\bibitem{Gor-CMP}
A.~Gorokhovsky,
 ``Characters of cycles, equivariant characteristic classes and
  {Fredholm} modules'',
 {\em Commun. Math. Phys.}, {\bf 208} (1999), 1--23

\bibitem{Gor:topology}
A.~Gorokhovsky,
 ``Secondary characteristic classes and cyclic cohomology of {Hopf}
  algebras'',
 {\em Topology}, {\bf 41} (2002), 993--1016

\bibitem{Gor06}
A. Gorokhovsky,
 ``Bivariant {C}hern character and longitudinal index'',
 {\em J. Funct. Anal.}, {\bf 237}:1 (2006), 105--134

\bibitem{Gor-Lott03}
A.~Gorokhovsky, J.~Lott,
 ``Local index theory over \'etale groupoids'',
 {\em J. Reine und Angew. Math.}, {\bf 560} (2003), 151--198


\bibitem{Gor-Lott04}
A. Gorokhovsky, J. Lott,
 ``Local index theory over foliation groupoids'',
 {\em Adv. Math.}, {\bf 204}:2 (2006), 413--447

\bibitem{Gracia:book}
J.~M. Gracia-Bond{\'{\i}}a, J.~C. V{\'a}rilly, H. Figueroa,
 {\em Elements of non-commutative geometry}.
 Birkh\"auser Advanced Texts: Basler Lehrb\"ucher. Birkh\"auser, Boston, MA, 2001

\bibitem{Gromov95} M. Gromov, Positive curvature, macroscopic dimension,
spectral gaps and higher signatures.  {\em Functional analysis on
the eve of the 21st century, Vol. II (New Brunswick, NJ, 1993),}
 Progr. Math., 132, Birkha"user Boston, Boston, MA, 1996, 1--213

\bibitem{Guill77}
V.~Guillemin,
 V. Guillemin. ``Lectures on spectral theory of elliptic
 operators'',  {\em Duke Math. J.}  {\bf 44} (1977), 485--517

\bibitem{Gu85}
V.~Guillemin,
 ``A new proof of {W}eyl's formula on the asymptotic distribution of
  eigenvalues'',
 {\em Adv. Math.}, {\bf 55} (1985), 131--160


\bibitem{Haefliger72}
A.~Haefliger,
 ``Homotopy and integrability'',
{\em Manifolds--Amsterdam 1970 (Proc. Nuffic Summer School)},
  Lect. Notes in Math., {\bf 197}, Springer, Berlin,
  1971, 133--163

\bibitem{Haefliger:CIME}
A.~Haefliger,
 ``Differential cohomology'',
{\em Differential topology (Varenna, 1976)},
  Liguori, Naples, 1979, 19--70

\bibitem{Haefliger80}
A.~Haefliger,
 ``Some remarks on foliations with minimal leaves'',
 {\em J. Diff. Geom.}, {\bf 15} (1980), 269--384

\bibitem{Haefliger84}
A.~Haefliger,
 ``Groupo\"\i des d'holonomie et classifiants'',
 {\em Transversal structure of foliations (Toulouse, 1982)},
 Ast\'erisque, {\bf 116}, 1984, 70--97

\bibitem{Heitsch95}
J.~L. Heitsch,
 ``Bismut superconnections and the {C}hern character for {D}irac
  operators on foliated manifolds'',
 {\em $K$-Theory}, {\bf 9}:6 (1995), 507--528

\bibitem{H-L}
J.~L. Heitsch, C.~Lazarov,
 ``A {L}efschetz theorem for foliated manifolds'',
 {\em Topology}, {\bf 29} (1990), 127--162

\bibitem{H-L91}
J.~L. Heitsch, C.~Lazarov,
 ``Rigidity theorems for foliations by surfaces and spin
 manifolds'',
 {\em Michigan Math. J.}, {\bf 38} (1991), 285--297

\bibitem{Heitsch-L99}
J.~L. Heitsch, C.~Lazarov,
 ``A general families index theorem'',
 {\em $K$-Theory}, {\bf 18}:2 (1999), 181--202

\bibitem{Hil-Skan83}
M.~Hilsum, G.~Skandalis,
 ``Stabilit\'e des {$C\sp{\ast} $}-alg\`ebres de feuilletages'',
 {\em Ann. Inst. Fourier (Grenoble)}, {\bf 33}:3 (1983), 201--208

\bibitem{Hil-Skan}
M.~Hilsum, G.~Skandalis,
 ``Morphismes {$K$}-orient\'es d'espaces de feuilles et fonctorialite en
  th\'{e}orie de {Kasparov}'',
 {\em Ann. scient. Ec. Norm. Sup.}, {\bf 20} (1987), 325--390

\bibitem{Hil-Skan92}
M.~Hilsum, G.~Skandalis,
 ``Invariance par homotopie de la signature \`a coefficients dans un
  fibr\'e presque plat'',
 {\em J. Reine Angew. Math.}, {\bf 423} (1992), 73--99

\bibitem{Hurder93}
S. Hurder,
 ``Coarse geometry of foliations'',
 {\em Geometric study of foliations (Tokyo, 1993)},
  World Sci. Publ., River Edge, NJ, 1994, 35--96

\bibitem{Hurder08}
S. Hurder,
 ``Classifying foliations'', preprint arXiv:0804.1240, 2008

\bibitem{Hurder-KatokBAMS}
S.~Hurder, A.~Katok,
 ``Secondary classes and transverse measure theory of a foliation'',
 {\em Bull. Amer. Math. Soc. (N.S.)}, {\bf 11} (1984), 347--350

\bibitem{Ji}
R. Ji,
 ``Smooth dense subalgebras of reduced group {$C\sp *$}-algebras,
  {S}chwartz cohomology of groups, and cyclic cohomology'',
 {\em J. Funct. Anal.}, {\bf 107}:1 (1992), 1--33

\bibitem{Jiang}
X. Jiang,
 ``An index theorem on foliated flat bundles'',
 {\em $K$-Theory}, {\bf 12}:4 (1997), 319--359

\bibitem{Juhl}
A. Juhl, {\em Cohomological theory of dynamical zeta functions.}
Progr. in Math., {\bf 194}, Birkha"user, Basel, 2001

\bibitem{KamberTondLNM}
F.~Kamber, Ph. Tondeur,
 {\em Foliated bundles and characteristic classes},
  Lect. Notes in Math., {\bf 493}, Springer, Berlin, 1975.



\bibitem{Karoubi87}
M. Karoubi,
 {\em Homologie cyclique et {$K$}-th\'eorie},
Ast\'erisque, {\bf 149}, 1987.

\bibitem{Kasparov75}
G. G. Kasparov, ``Topological invariants of elliptic operators. I.
$K$-homology'', {\em Izv. Akad. Nauk SSSR Ser. Mat.}, {\bf 39}:4
(1975), 796--838; translation in: {\em Mathematics of the
USSR-Izvestiya}, {\bf 9}:4 (1975), 751--792


\bibitem{Kasparov80}
G. G. Kasparov, ``The operator $K$-functor and extensions of
$C^*$-algebras'', {\em Izv. Akad. Nauk SSSR Ser. Mat.}, {\bf 44}:3
(1980), 571--636; translation in: {\em Mathematics of the
USSR-Izvestiya}, {\bf 16}:3 (1981), 513--572

\bibitem{Kasparov-conspectus}
G. G. Kasparov, ``$K$-theory, group $C\sp *$-algebras, and higher
signatures (conspectus)'', {\em Novikov conjectures, index theorems
and rigidity, Vol. 1 (Oberwolfach, 1993)},  London Math. Soc.
Lecture Note Ser., {\bf 226}, Cambridge Univ. Press, Cambridge,
1995, 101--146

\bibitem{Kas93}
G. G. Kasparov, ``Novikov's conjecture on higher signatures: the
operator $K$-theory approach'', {\em Representation theory of groups
and algebras,}  Contemp. Math., {\bf 145}, Amer. Math. Soc.,
Providence, RI, 1993, 79--99

\bibitem{Kas-teo}
G. G. Kasparov, ``$K$-theoretic index theorems for elliptic and
transversally elliptic operators'', preprint, 2008.

\bibitem{Kawasaki81}
T. Kawasaki,
 ``The index of elliptic operators over {$V$}-manifolds'',
 {\em Nagoya Math. J.}, {\bf 84} (1981), 135--157

\bibitem{invitation}
{\em An invitation to non-commutative geometry}. M. Khalkhali, M.
Marcolli, editors. World Sci. Publ., River Edge, NJ, 2008.

\bibitem{konechny:schw}
A.~Konechny, A.~Schwarz,
 ``Introduction to {M}(atrix) theory and noncommutative geometry'',
 {\em Phys. Rep.}, {\bf 360}:5-6 (2002), 353--465

\bibitem{trans}
Yu.~A. Kordyukov,
 ``Transversally elliptic operators on {$G$}-manifolds of bounded
  geometry'',
 {\em Russ. J. Math. Ph.}, {\bf 2} (1994), 175--198; {\bf 3} (1995), 41--64

\bibitem{noncom}
Yu.~A. Kordyukov,
 ``Noncommutative spectral geometry of {Riemannian} foliations'',
 {\em Manuscripta Math.}, {\bf 94} (1997), 45--73

\bibitem{egorgeo}
Yu.~A. Kordyukov,
 ``Egorov's theorem for transversally elliptic operators on foliated
  manifolds and noncommutative geodesic flow'',
 {\em Math. Phys. Anal. Geom.}, {\bf 8}:2 (2005), 97--119

\bibitem{matrix-egorov}
Yu.~A. Kordyukov,
 ``The {E}gorov theorem for transverse {D}irac-type operators on
  foliated manifolds'',
 {\em J. Geom. Phys.}, {\bf 57}:11 (2007), 2345--2364

\bibitem{vanishing}
Yu.~A. Kordyukov, ``Vanishing theorems for transverse Dirac
operators on Riemannian foliations'', {\em Ann. Global Anal. Geom.}
{\bf 34} (2008), 195--211

\bibitem{survey}
Yu.~A. Kordyukov,
 ``Noncommutative geometry of foliations'',
 {\em J. K-Theory}, {\bf 2} (2008), 219-327

\bibitem{Kreimer98}
D. Kreimer,
 ``On the {H}opf algebra structure of perturbative quantum field
  theories''
 {\em Adv. Theor. Math. Phys.}, {\bf 2}:2(1998), 303--334



\bibitem{Landi:book}
G.~Landi,
 {\em An introduction to noncommutative spaces and their geometries},
  Lecture Notes in Physics. New Series m:
  Monographs, {\bf 51}, Springer, Berlin, 1997

\bibitem{Lawson-M}
H.B. Lawson, M.-L. Michelsohn, {\em Spin geometry.} Princeton
Mathematical Series, 38. Princeton Univ. Press, Princeton, NJ, 1989

\bibitem{Lazarov00}
C.~Lazarov,
 ``Transverse index and periodic orbits'',
 {\em Geom. and Funct. Anal.}, {\bf 10} (2000), 124--159

\bibitem{Leichtnam05}
E. Leichtnam, ``An invitation to Deninger's work on arithmetic zeta
functions'',  {\em Geometry, spectral theory, groups, and dynamics},
Contemp. Math., {\bf 387}, Amer. Math. Soc., Providence, RI, 2005,
201--236

\bibitem{Leichtnam-Piazza04}
E. Leichtnam, P. Piazza, ``Elliptic operators and higher
signatures'', {\em Ann. Inst. Fourier (Grenoble)}  {\bf 54}:5
(2004), 1197--1277

\bibitem{Leichtnam-Piazza05a}
E. Leichtnam, P. Piazza,
 ``Cut-and-paste on foliated bundles'',
 {\em Spectral geometry of manifolds with boundary and
  decomposition of manifolds}, Contemp. Math., {\bf 366},
  Amer. Math. Soc., Providence, RI, 2005, 151--192

\bibitem{Leichtnam-Piazza05}
E. Leichtnam, P. Piazza,
 \'{E}tale groupoids, eta invariants and index theory.
 {\em J. Reine Angew. Math.}, 587:169--233, 2005.

\bibitem{Lesch91}
M. Lesch, ``On the index of the infinitesimal generator of a flow'',
{\em J. Operator Theory} {\bf 26}:1 (1991), 73--92

\bibitem{Loday}
J.-L. Loday,
 {\em Cyclic homology}, Grundl. der Math. Wiss., {\bf 301}, Springer, Berlin, 1998.

\bibitem{Lott:gafa92}
J. Lott, ``Superconnections and higher index theory'', {\em Geom.
Funct. Anal.} {\bf 2} (1992), 421--454

\bibitem{Lott:eta}
J. Lott,  ``Higher eta-invariants'', {\em $K$-Theory} {\bf 6}
(1992), 191--233

\bibitem{MMS06}
V.~Mathai, R.~B. Melrose, I.~M. Singer.
 ``Fractional analytic index''
 {\em J. Differential Geom.}, {\bf 74}:2 (2006), 265--292

\bibitem{MMS08}
V.~Mathai, R.~B. Melrose, I.~M. Singer,
 ``Equivariant and fractional index of projective elliptic
 operators''
 {\em J. Differential Geom.}, {\bf 78}:3 (2008), 465--473

\bibitem{McKeanSinger}
H.~P. McKean, Jr., I.~M. Singer,
 ``Curvature and the eigenvalues of the {L}aplacian'',
 {\em J. Differential Geometry}, {\bf 1}:1 (1967), 43--69

\bibitem{Melrose93}
R. B. Melrose, {\em The Atiyah-Patodi-Singer index theorem.}
Research Notes in Mathematics, {\bf 4}, A K Peters, Ltd., Wellesley,
MA, 1993

\bibitem{Mishch-Fomenko}
A. S. Mishchenko, A. T. Fomenko, ``The index of elliptic operators
over $C^*$-algebras'', {\em Izv. Akad. Nauk SSSR Ser. Mat.}, {\bf
43}:4 (1979), 831--859; translation in: {\em Mathematics of the
USSR-Izvestiya}, {\bf 15}:1 (1980), 87--112

\bibitem{MP97} B. Monthubert, F. Pierrot,
``Indice analytique and Groupoides de Lie'', {\em C.R. Acad. Sci.
Paris, Ser. I}, {\bf 325}:2 (1997), 193--198

\bibitem{M-S}
C.~C. Moore, C.~Schochet,
 {\em Global analysis on foliated spaces}, With appendices
 by S. Hurder, Moore, Schochet and Robert J. Zimmer,
  Mathematical Sciences Research Institute Publications, {\bf 9}
 Springer, New York, 1988


\bibitem{Mor}
H.~Moriyoshi,
 ``On cyclic cocycles associated with the {Godbillon}-{Vey}
 classes'',
{\em Geometric study of foliations (Tokyo, 1993)},
  World Sci. Publishing, River Edge, NJ, 1994, 411--423

\bibitem{Mor-Natsume}
H.~Moriyoshi, T.~Natsume,
 ``The {Godbillon}-{Vey} cyclic cocycle and longitudinal {Dirac}
  operators'',
 {\em Pacific J. Math.}, {\bf 172}:2 (1996), 483--539

\bibitem{Mosc09} H.~Moscovici, ``Local index formula and twisted
spectral tripes'', preprint arXiv:0902.0835, 2009


\bibitem{MRW87}
P.~S. Muhly, J.~N. Renault, D.~P. Williams,
 ``Equivalence and isomorphism for groupoid {$C\sp
 \ast$}-algebras'',
 {\em J. Operator Theory}, {\bf 17}:1 (1987), 3--22

\bibitem{Naz-Savin-Sternin}
V.~E. Nazaikinskii, A.~Yu. Savin, B.~Yu. Sternin,
 {\em Elliptic theory and non-commutative geometry. Nonlocal elliptic operators},
  Operator Theory: Advances and Applications, Advances in Partial Differential
  Equations. {\bf 183},  Birkh\"auser, Basel, 2008

\bibitem{Nestke-Z}
A.~Nestke, P.~Zickermann,
 ``The index of transversally elliptic complexes'',
 {\em Rend. Circ. Mat. Palermo}, {\bf 34}:Suppl.9 (1985), 165--175

\bibitem{Nistor91}
V. Nistor, ``A bivariant Chern character for $p$-summable
quasihomomorphisms'', {\em $K$-Theory}  {\bf 5}:3  (1991), 193--211

\bibitem{Nistor93}
V. Nistor, ``A bivariant Chern-Connes character'', {\em Ann. of
Math. (2)} {\bf 138}:3 (1993), 555--590


\bibitem{Nistor94}
V. Nistor, ``On the Cuntz-Quillen boundary map'', {\em C. R. Math.
Rep. Acad. Sci. Canada},  {\bf 16}:5  (1994), 203--208

\bibitem{Nistor96}
V.~Nistor,
 ``The index of operators on foliated bundles'',
 {\em J. Funct. Anal.}, {\bf 141}:2 (1996), 421--434

\bibitem{NWX99}
V. Nistor, A. Weinstein, P. Xu, ``Pseudodifferential operators on
differential groupoids'', {\em Pacific J. Math.} {\bf 189}:1 (1999),
117--152

\bibitem{Palais}
R. S. Palais, {\em Seminar on the Atiyah-Singer index theorem.}
Annals of Mathematics Studies, {\bf 57}, Princeton Univ. Press,
Princeton, N.J. 1965

\bibitem{Paradan01}
P.-{\'E}. Paradan,
 ``Localization of the {R}iemann-{R}och character'',
 {\em J. Funct. Anal.}, {\bf 187}:2 (2001), 442--509

\bibitem{Paradan03}
P.-{\'E}. Paradan,
 ``{${\rm Spin}\sp c$}-quantization and the {$K$}-multiplicities of the
  discrete series''.
 {\em Ann. Sci. \'Ecole Norm. Sup. (4)}, {\bf 36}:5 (2003), 805--845

\bibitem{Paradan-Vergne08}
P.-{\'E}. Paradan, M. Vergne,
 ``Index of transversally elliptic operators'', preprint arXiv:0804.1225, 2008

\bibitem{Pat90}
S. J. Patterson, ``On Ruelle's zeta-function'', {\em Festschrift in
honor of I. I. Piatetski-Shapiro on the occasion of his sixtieth
birthday, Part II (Ramat Aviv, 1989)},  Israel Math. Conf. Proc., 3,
Weizmann, Jerusalem, 1990, 163--184

\bibitem{Patodi71}
V.~K. Patodi,
 ``Curvature and the eigenforms of the {L}aplace operator'',
 {\em J. Differential Geometry}, {\bf 5} (1971), 233--249

\bibitem{Peric}
G. Peri{\'c},
 ``Eta invariants of {D}irac operators on foliated manifolds'',
 {\em Trans. Amer. Math. Soc.}, {\bf 334}:2 (1992), 761--782

\bibitem{Perrot}
D.~Perrot,
 ``A {Riemann}-{Roch} theorem for one-dimensional complex
 groupoids'', {\em Commun. Math. Phys.}, {\bf 218} (2001), 373--391

\bibitem{phillips87}
J.~Phillips,
 ``The holonomic imperative and the homotopy groupoid of a foliated
  manifold'', {\em Rocky Mountain J. Math.}, {\bf 17}:1 (1987), 151--165

\bibitem{Phillips-Raeburn}
J. Phillips, I. Raeburn,
 ``An index theorem for {T}oeplitz operators with non-commutative symbol
  space'', {\em J. Funct. Anal.}, {\bf 120}:2 (1994), 239--263

\bibitem{Ramachandran93}
M. Ramachandran,
 ``von {N}eumann index theorems for manifolds with boundary'',
 {\em J. Differential Geom.}, {\bf 38}:2 (1993), 315--349


\bibitem{RieffelRot}
M.~A. Rieffel,
 ``{$C^*$}-algebras associated with irrational rotations'',
 {\em Pac. J. Math.}, {\bf 93} (1981), 415--429

\bibitem{Rieffel82}
M.~A. Rieffel,
 ``Morita equivalence for operator algebras'',
 {\em Operator algebras and applications, Part I (Kingston, Ont.,
  1980)}, Proc. Sympos. Pure Math., {\bf 38},  Amer.
  Math. Soc., Providence, R.I., 1982, 285--298

\bibitem{Roe98}
J. Roe, {\em Elliptic operators, topology and asymptotic methods.}
Second edition. Pitman Research Notes in Mathematics Series, {\bf
395}. Longman, Harlow, 1998

\bibitem{Sau}
J.-L. Sauvageot,
 ``Semi-groupe de la chaleur transverse sur la ${C}^*$-algebre d'un
  feuilletage riemannien''
 {\em J. Funct. Anal.}, {\bf 142} (1996), 511--538

\bibitem{Schw}
L.~Schweitzer,
 ``A short proof that {$M_n(A)$} is local if {$A$} is local and
  {Fr\'echet}'',
 {\em Internat. J. Math.}, {\bf 3} (1992), 581--589

\bibitem{Segal68}
G. Segal, ``Classifying spaces and spectral sequences'', {\em Inst.
Hautes \'{E}tudes Sci. Publ. Math.} {\bf 34} (1968), 105--112

\bibitem{Shubin:UMN79}
M. A. Shubin, ``The spectral theory and the index of elliptic
operators with almost periodic coefficients'', {\em UMN}, {\bf 34}:2
(1979), 95--135; translation in: {\em Russian Mathematical Surveys},
{\bf 34}:2 (1979), 109--157

\bibitem{Singer:recent}
I.~M. Singer,
 ``Recent applications of index theory for elliptic operators'',
{\em Partial differential equations (Proc. Sympos. Pure Math., Vol.
XXIII, Univ. California, Berkeley, Calif., 1971)}, Amer. Math.
  Soc., Providence, R. I., 1973, 11--31

\bibitem{Singha-Goswami}
K. B. Sinha, D. Goswami, {\em Quantum stochastic processes and
noncommutative geometry}. Cambridge Tracts in Mathematics, {\bf
169}. Cambridge Univ. Press, Cambridge, 2007.

\bibitem{skandal-hopf}
G. Skandalis,
 ``Noncommutative geometry, the transverse signature operator, and
  {H}opf algebras [after {A}. {C}onnes and {H}. {M}oscovici]'',
{\em Cyclic homology in non-commutative geometry}, Encyclopaedia
Math. Sci., {\bf 121}, Springer, Berlin, 2004, 115--134

\bibitem{Soloviev-Troitsky}
Yu. P. Solov{$'$}{\"e}v,  E. V. Troitski{\u\i},  $C\sp
*$-algebras and elliptic operators in differential topology,
Izdatel{$'$}stvo ``Faktorial'', Moscow, 1996; translation in:
Translations of Mathematical Monographs, 192. American Mathematical
Society, Providence, RI, 2001.

\bibitem{Varilly}
J. C. V\'arilly, {\em An introduction to noncommutative geometry},
EMS Series of Lectures in Mathematics. European Mathematical Society
(EMS), Z\"urich, 2006

\bibitem{Vassout}
S.~Vassout,
 {\em Feuilletages et r\'esidu non commutatif longitudinal}.
 PhD thesis, Institut Jussieu, Paris, 2001.

\bibitem{Vassout07}
S. Vassout, ``Unbounded pseudodifferential calculus on Lie
groupoids'', {\em J. Funct. Anal.}  {\bf 236}:1  (2006), 161--200

\bibitem{Wang07}
B.-L. Wang, ``Geometric cycles, index theory and twisted
K-homology'', {\em J. Noncommutative Geometry}, {\bf 2} (2008),
497--552

\bibitem{Winkeln}
H.~E. Winkelnkemper,
 ``The graph of a foliation'',
 {\em Ann. Glob. Anal. Geom.}, {\bf 1} (1983), 53--75

\bibitem{Wo}
M.~Wodzicki, ``Noncommutative residue. {P}art {I}. {F}undamentals'',
{\em K-theory, arithmetic and geometry (Moscow, 1984-86)}, Lect.
  Notes in Math. {\bf 1289},  Springer, Berlin,
  1987, pages 320--399


\bibitem{Yu}
G.~L. Yu,
 ``Cyclic cohomology and the index theory of transversely elliptic
  operators'',
 {\em Selfadjoint and nonselfadjoint operator algebras and operator
  theory (Fort Worth, TX, 1990)}, Contemp. Math., {\bf 120},
   Amer. Math. Soc., Providence, RI, 1991, 189--192
\end{thebibliography}
\end{document}